\documentstyle[12pt]{article}
\newtheorem{theorem}{Theorem}[section]
\newtheorem{lemma}[theorem]{Lemma}
\newtheorem{question}[theorem]{Question}
\newtheorem{corollary}[theorem]{Corollary}
\newtheorem{definition}[theorem]{Definition}

\newtheorem{claim}[theorem]{Claim}
\newtheorem{remark}[theorem]{Remark}

\begin{document}

\title{ {\bf   On Gorenstein log del Pezzo surfaces }}
\author{Qiang Ye} 
 \date{}
\input prepictex
\input pictex
\input postpictex
\maketitle
\begin{abstract} 
In this paper, we first present the complete list of the singularity types of the Picard number one Gorenstein log del Pezzo surface and the number of the isomorphism classes with the given singularity type. Then we give out a method to find out all singularity types of Gorenstein log del Pezzo surface. As an application, we present the complete list of the Dynkin type of the Picard number two Gorenstein log del Pezzo surfaces. Finally we present the complete list of the singularity type of the relatively minimal Gorenstein log del Pezzo surface and the number of the isomorphism classes with the given singularity type.
\end{abstract}

\section{Introduction}

Let $V$ be a normal projective surface defined over the complex numbers field {\bf C}. We say that $V$ is a {\it log del Pezzo surface} if $V$ has at worst quotient singularities and if the anticanonical divisor $-K_V$ of $V$ is an ample divisor. Then $V$ is {\it Gorenstein}, i.e., $K_V$ is a Cartier divisor, if and only if $V$ has at worst rational double singularities.  We note that a log del Pezzo surface is a rational surface (cf. [GZ]). The log del Pezzo surface has been studied from the various points of view. A general theory on the structure of such singular surfaces is developed in  [M1] , [Z1] and [Z2].
The topological properties of the smooth part of the log del Pezzo surfaces and log Fano varieties are studied in [Fujiki], [GZ], [GZ1], [T1] and [Z3-Z6].

Let $V$ be a Gorenstein log del Pezzo surface. Denote by $V^0:=V-({\rm Sing}V)$ the smooth part of $V$. Let $g: U \longrightarrow V$ be a minimal resolution of singularities, we denote by $D:=g^{-1}({\rm Sing}V)$ the exceptional 
divisor and $\sharp (D)$ the number of the irreducible components of $D$. 

In [Fu], the possible singularities of a Picard number one Gorenstein log del Pezzo surface are classified (see also [Brenton] ).

\begin{question} Let $V$ be a Picard number one Gorenstein log del Pezzo surface,
 is $V$ then determined uniquely by its singularities?
\end{question}

In general,  $V$ is not determined by its singularities. However, when $V^0$ is simply connected, this is almost the case. In [GPZ], [KM] and [MZ1],  they have shown that
for  a Picard number one Gorenstein log del Pezzo surface $V$  with simply connected smooth part,  the singularity type of $V$ is one of the following:
$$
A_1, \ A_1+A_2, \ A_4, \ D_5, \ E_6, \ E_7, \ E_8.
$$
Furthermore, $V$ is uniquely determined by its singularities except in the case of $V(E_8)$ where there are exactly two isomorphism classes.
\vskip 0.5cm

One of our main results is the following:
\vskip 0.5cm

\begin{theorem} Let $V$ be a Picard number one Gorenstein
log del Pezzo surface $(\not= {\bf P}^2)$, then the singularity types of $V$ and the number $m$ of the isomorphism classes of the Picard number one Gorenstein log del Pezzo surfaces with the given singularity type are listed in Table $1.1$.
\vskip 0.5cm
\begin{center}
Table $1.1$
\end{center}
\begin{center}
\begin{tabular}{| c|c|c|c|c|c|}\hline
Dynkin type of $V$ & $m$ & Dynkin type of $V$  & $m$ & Dynkin type of $V$ & $m$
\\ \hline
$E_8$  & $2$ & $2A_4$& $1$ & $A_1+2A_3$ & $1$ 
\\ \hline
$A_1+E_7$  &$2$ &$A_1+A_2+A_5$   & $ 1$ & $E_6$& $1$  
\\ \hline
  $D_8$ & $1$& $2A_1+2A_3$   & $1$ & $A_1+A_5$& $1$       \\ \hline
$A_8$  &$1$ & $4A_2$ &  
$1$  & $3A_2$&$1$ 
\\ \hline
$A_2+E_6$  & $2$ & $E_7$   & $1$  & $D_5$& $1$ 
\\ \hline
$2A_1+D_6$  & $1$ & $A_1+D_6$  & $ 1$ &  $2A_1+A_3$ & $1$
\\ \hline
$2D_4$  & $\infty$& $A_7$  & $1$  & $A_4$ &$1$ 
\\ \hline
$A_3+D_5$  & $1$ &$A_2+A_5$  & $1$  & $A_1+A_2$ &$1$ 
\\ \hline

$A_1+A_7$  & $1$ &$3A_1+D_4$  & $1$  & $A_1$ &$1$ 
\\ \hline

\end{tabular} 
\end{center} 
All isomorphism classes of any given singularity type are realizable (see the proof).
\end{theorem}
\begin{corollary}  Let $V$ be a Picard number one Gorenstein
log del Pezzo surface. Then $V$ is  determined uniquely up to isomorphisms by the Dynkin type of $V$,  except for $V(E_8)$, $V(A_1+E_7)$, $V(A_2+E_6)$ and $V(2D_4)$.
\end{corollary}
\begin{remark} There are several results which are related to our Theorem $1.2$. 

$(1)$ \ In {\rm [Fu] }, the possible singularities of $V$ are classified (see also {\rm  [Brenton]} ).

$(2)$ \ In {\rm  [MZ1]}, a classification of Picard number one Gorenstein log del Pezzo surfaces is given.

In the proof,  we need the result about the extremal rational elliptic surfaces which is proved in {\rm [MP]} but we do not use the classification mentioned in {\rm [Brenton]},  {\rm [Fu]} and {\rm [MZ1]} above. 
\end{remark}
We also find out the singularity types of the Picard number two Gorenstein log del Pezzo surfaces.

\begin{theorem}  Let $V$ be a Picard number two Gorenstein log del Pezzo surface $(\not= {\bf P}^1 \times {\bf P}^1 \ and  \   \Sigma_1 )$. Then there exist exactly  $45$ singularity types which  are listed in  
Table $1.2$.
\begin{center}
Table $1.2$
\end{center}
\begin{center}
\begin{tabular}{| c|c|c|c|c|c|}\hline
$E_7$ &  $A_7$ &   $2A_3+A_1$&  $A_4+A_2$ &  $A_4+A_1$ & $A_3$    \\ \hline
$E_6+A_1$ & $A_6+A_1$ & $3A_2+A_1$  & $2A_3$ & $A_3+2A_1$  & $A_2+A_1$\\ \hline                                                      
$D_7$  & $A_5+A_2$ & $E_6$ & $A_3+A_2+A_1$ & $2A_2+A_1$ &  $A_2$\\ \hline
$D_6+A_1$ & $A_5+2A_1$ &  $D_6$ & $A_3+3A_1$ & $D_4$ &  $2A_1$\\ \hline
$D_5+A_2$ & $A_4+A_3$ & $D_5+A_1$& $3A_2$ & $A_4$ & $A_1$\\ \hline
$D_5+2A_1$ & $A_4+A_2+A_1$ & $D_4+2A_1$ & $6A_1$ & $A_3+A_1$ &  \\ \hline
 $D_4+A_3$ & $A_3+ A_2+2A_1$& $A_6$ &$D_5$ & $A_2+2A_1$ &   \\ \hline
$D_4+3A_1$ & $A_3+4A_1$ & $A_5+A_1$ & $A_5$ & $4A_1$ &   \\ \hline
\end{tabular} 
\end{center} 
All the above singularity types are realizable (see the proof).
\end{theorem}
\vskip 0.5cm
The relatively minimal Gorenstein log del Pezzo surface (see Definition $2.1$) is introduced in [MZ2] where the Dynkin type of the Picard number two relatively minimal Gorenstein log del Pezzo surface is given but the number of the isomorphism classes with the given singularity is not provided. 
In this paper, we prove the following theorem. 
\vskip 0.5cm
\begin{theorem}  Let $V$ be a Picard number two relatively minimal Gorenstein
log del Pezzo surface $(\not= {\bf P}^1 \times {\bf P}^1 \  and \ \Sigma_1 )$. Then $V$ has one of the $10$ singularity types in Table $1.3$.
The number $m$ of the isomorphism classes of the Picard number two relatively minimal Gorenstein log del Pezzo surfaces with the given singularity type are also listed in Table $1.3$.
\begin{center}
Table $1.3$
\end{center}
\begin{center}
\begin{tabular}{| c|c|c|c|c|c|}\hline
Dynkin type of $V$ & $m$ & Dynkin type of $V$  & $m$ & Dynkin type of $V$ & $m$
\\ \hline
$D_7$  & $\infty$ & $D_6$& $\infty$  & $D_4 $ & $1$ 
\\ \hline
$D_5+2A_1$  & $\infty$ &$ D_4+2A_1$   & $1$ &   $4A_1$ & $1$  
\\ \hline
$A_3+4A_1$  & $1$  & $ 2A_3$ & $\infty$  &    &   
\\ \hline
$A_3+D_4$  & $\infty$  & $ 6A_1$ & $1$  &    &   
\\ \hline

\end{tabular} 
\end{center} 
All isomorphism classes of the above singularity types are realized  in the course of  the proof.
\end{theorem}
\begin{corollary}  Let $W$ be a Gorenstein log del Pezzo surface.

$(1)$  Then there are disjoint curves $C_1$, $\cdots$, $C_k$ and a birational morphism 
$\sigma : W \longrightarrow W_{\rm min}$ such that $\sigma(C_i)$ are smooth points , $\sigma$ is an isomorphism outside $\cup C_i$ and $W_{\rm min}$ is isomorphism to one of the surfaces $V$ in Theorem $1.2$ or $1.6$.

$(2)$  The singularity type of $W$ is that of $V$ plus a few type $A_{n_i}$ $(1\leq i \leq k)$ singularities lying on $C_i$.

\end{corollary} 
 
\begin{remark} From {\rm [MZ2, Lemma 2]} or Lemma $2.3$, we know that Theorems $1.2$ and $1.6$ present the complete list of the singularity types of the relatively minimal Gorenstein log del Pezzo surfaces and the number of the isomorphism classes with the given singularity type.
\end{remark}

This paper is organized as follows. In section $2$, we introduce some terminology and preliminaries which will be used throughout the paper. In section $3$, we establish the connection between Gorenstein log del Pezzo surfaces and rational elliptic surfaces. We will also give out a method to classify all singularity types of the Gorenstein log del Pezzo surfaces. In sections $4$ and $5$, we prove Theorems $1.2$, $1.5$ and $1.6$.
\vskip 0.5cm
\hskip -0.6cm {\bf Acknowledgement.} I would like to thank Professor D.-Q. Zhang for many enlightening advice during the preparation of this paper.

\section{Preliminaries}

{\bf A. Relatively minimal Gorenstein log del Pezzo surfaces}
\vskip 0.2cm
\hskip -0.6cm 

Let $V$ be a Gorenstein log del Pezzo surface. Let $f : U \longrightarrow V$ be a minimal resolution of singularities and let $D:=f^{-1}({\rm Sing} V)$. Then there are no $(-n)$-curves on $U$ with $n \geq 2$ except for those $(-2)$-curves contained in $D$. The morphism $f$ is the contraction of all $(-2)$-curves on $U$. 

Let $\Sigma = \{ F_1, \cdots, F_t \}$ be a maximum set of disjoint union of $(-1)$-curves on $U$ such that the connected component of $\sum F_j +D$ containing $F_i$ is a linear chain $F_i+\Delta_i$ of $F_i$ and $n_i$ ($n_i \geq 0$) $(-2)$-curves with the dual graph:
$$
(-1)-(-2)-(-2)-\cdots -(-2).
$$
Let $U \longrightarrow U_{\rm min}$ be the smooth blow-down of all connected components of $\sum_{i=1}^{t} (F_i +\Delta_i )$. Let $D'$ be the image on $U_{\rm min}$ of $D$ and let $U_{\rm min} \longrightarrow V_{\rm min}$ be the contraction of $D'$. Then we have a commutative diagram:
\begin{center}
\begin{tabular}{rrr}

$ U$ & $ \longrightarrow $ & $U_{\rm min}$ \\
 $\downarrow$ & & $\downarrow$ \\
$ V$ & $\longrightarrow $& $  V_{\rm min}$
\end{tabular}
\end{center}

Here the map $V \longrightarrow V_{\rm min}$ contracts ${\bar F}_i :=f(F_i)$ to $t$ smooth points
and is isomorphism outside $\cup {\bar F}_i$. Note that $\rho (V_{\rm min})=\rho (V) - t$.

\begin{definition} The surface $V_{\rm min}$ is  called a 
relatively minimal model of $V$. $V$ is relatively minimal if $V=V_{\rm min}$. Note that $V_{\rm min}$ is not uniquely determined by $V$.
\end{definition}

From the above discussion, we have

\begin{lemma} Let $V$ be a Gorenstein log del Pezzo surface, $U$ the minimal resolution of $V$ and $D$ the exceptional divisor.
  
$(1)$ $V$ is not relatively minimal if and only if there is a $(-1)$-curve $E$ on $U$ such that the connected component of $E+D$ has dual graph:
$$
(-1) \ \ or \ \ (-1)-(-2)-(-2)-\cdots -(-2).
$$
$(2)$ If Sing$V_{\rm min}$ has Dynkin type $\sum a_iA_i+\sum d_iD_i + \sum e_iE_i$, then SingV has Dynkin type  $\sum a_iA_i+\sum d_iD_i + \sum e_iE_i + \sum_{i=1}^t A_{n_i}$.

\end{lemma}

Thanks to the virtue of Miyanishi-Zhang (cf. [25, Lemma 2]), we have the following Lemma.

\begin{lemma}  Let $V$ be a Gorenstein log del Pezzo surface. Then there exists a birational morphism $\eta : V \longrightarrow W$ such that the following assertions hold.

$(1)$ $W$ is a relatively minimal Gorenstein log del Pezzo surfaces of Picard number one or two.

$(2)$  Let $E_1$,...,$E_t$ be all exceptional curves of $\eta$. Then $w_i:=\eta (E_i) (1\leq i\leq t)$ are $t$ distinct smooth points of $W$ for each $i$. Moreover $V$ consists of the same singular points as on $W$ and additional points $x_i$ of the Dynkin type $A_{n_i}$ with $n_i \geq 1$, where  $1\leq i\leq t$.

\end{lemma} 
\vskip 0.5cm
{\bf B. Mordell-Weil groups of  elliptic surfaces}
\vskip 0.5cm
\hskip -0.6cm $(a)$ \ \  Definitions
\vskip 0.5cm
\hskip -0.6cm
Let $C$ be a smooth projective curve over an algebraically
closed field of characteristic 0, and $Y$ a (Jacobian)  {\em elliptic
surface} over $C$. By this we mean the following: $Y$ is a 
smooth projective surface with a relatively minimal elliptic
fibration
$$
  f : Y \longrightarrow C 
$$ 
such that

$(i)$\ $f$ has a global section $\cal O$, and

$(ii)$\ $f$ is not smooth, i.e., there is at 
least one singular fibre.
\vskip 0.5cm
\hskip -0.6cm For every (Jacobian) elliptic fibration $Y$ there 
is a group structure of sections $MW (Y)$ with the distinguished section 
$\cal O$  as zero which we call {\it Mordell-Weil group}. 
Up to a finite group, $MW (Y)$ is identified 
with the relative automorphism group of the fibration. The basic results on the Mordell-Weil groups are proved in [Sh].
\vskip 0.5cm
\hskip -0.6cm Due to a formula of Shioda-Tate we have  the basic 
inequalities
$$
0\leq {\rm rank}MW (Y) \leq \rho (Y)-2 
$$
where $\rho (Y)$ is the Picard number and the discrepancy
in the upper bound is related to the degree of reducibility
of the fibre.

\begin{definition}
An elliptic fibration $Y$ is called {\rm extremal} if $\rho (Y)=h^{1,1} (Y)$  
{\rm (maximal\  Picard \  number) } and
rank $MW (Y)=0$.
\end{definition}

The extremal elliptic surfaces have  been studied from various points of view. In [MP],  the complete classification of extremal rational elliptic surfaces is given. For extremal elliptic $K3$ surfaces, we refer to [ATZ], [MP1],
[MP2],  [SZ] and [Ye].
\begin{definition}\hskip 0.5cm 
$Y$ is called  a  {\rm Mordell-Weil rank $r$} rational elliptic surface if $Y$ is a rational surface and rank $MW (Y)=r$.
\end{definition}

In [M] and [P],  the complete classification of the  rational elliptic surfaces is given. In [OS], they describe the complete structure theorem of the Mordell-Weil groups of the rational elliptic surfaces. 
\vskip 0.5cm
\hskip -0.6cm The following result is proved in [OS, Corollary 2.3 and Theorem 2.5].

\begin{theorem} The Mordell-Weil rank one rational elliptic surface is generated by the section $P$ which is disjoint from the zero-section ${\cal O}$ with  $<P,P> \leq 2$,  where $<\cdot ,\cdot> $ is the height pairing on Mordel-Weil groups which is defined below.
\end{theorem}

\vskip 0.5cm
\hskip -0.6cm $(b)$ \ \  Height Pairing

\vskip 0.5cm
\hskip -0.6cm Let $f$: $Y \longrightarrow {\bf P}^1$ be a 
(relatively minimal) elliptic surface over ${\bf P}^1$ with
a distinguished section ${\cal O}$. The complete list of 
possible fibres has been given by Kodaira [K1]. It 
encompasses two infinite families $(I_n, I_n^*, n\geq 0)$
and six exceptional cases $(II, III, IV, II^*, III^*, IV^*)$.
And they can be considered as sublattices of the N$\acute{e}$ron-Severi group of $Y$.
\vskip 0.5cm
\hskip -0.6cm  Given  an ellipric surface $f$: $Y\longrightarrow {\bf P}^1$,
let $F_{\nu} = f^{-1} ({\nu})$ denote the fibre over $\nu \in
{\bf P}^1$, and let

$  Sing(f) =\{ \nu \in {\bf P}^1 | F_{\nu} \ {\rm is\  singular} \}$.

$  {\bf {\cal R}}= {\rm Red}(f) = \{ \nu \in {\bf P}^1 | F_{\nu}\  {\rm is \ 
reducible}\}. $

\hskip -0.6cm For each $\nu \in {\bf {\cal R}}$, let
$$
  F_{\nu} = f^{-1} (\nu ) = \Theta_{\nu ,0} + \sum_{i=1}^{m_{\nu}
-1} \mu_{\nu ,i} \Theta_{\nu ,i}  \ \ (\mu_{\nu ,0} =1 )
$$
where $ \Theta_{\nu ,i}$ ($ 0\leq i \leq m_{\nu} -1$) are the 
irreducible components of $F_{\nu}$, $m_{\nu}$ being their
number, such that $ \Theta_{\nu ,0}$ is the unique component
of $F_{\nu}$ meeting the zero section. 

\begin{lemma} {\rm (cf. [Sh] )}
For each $\nu \in {\bf P}^1$, the intersection matrix
$$
A_{\nu} = ( (\Theta_{\nu ,i} \cdot \Theta_{\nu ,j}))_{1\leq i,j\leq m_{\nu}-1}
$$
is negative-definite.
\end{lemma} 

\hskip -0.6cm Here we denote by $K=K({\bf P}^1)$ the function field,
$$
  E(K) = {\rm the \ group\ of \ sections \ of\ }\  f,
$$
and
$$
  NS(Y) = {\rm the \ group \ of \ divisors \ on}\  Y \ \ {\rm
 modulo \ algebraic \ equivalence.}
$$
In general, for a rational elliptic surface $Y$ with a section and with at least one singular fibre, there is a unique homomorphism $\varphi$ of $E(K)$ into the 
N$\acute{e}$ron-Severi group $NS(Y)\otimes {\bf Q}$ such that 
$\varphi (P) \equiv (P)$ {\rm mod} $T_{\bf Q}$ and $\varphi (P) \perp T$ where $T$ is the subgroup of $NS(Y)$ generrated by the zero section ${\cal O}$ and all the irreducible components of fibres (cf. [Sh, Lemma 8.1]).

\begin{definition} 
The pairing 
$$
  <P,Q> = - (\varphi (P) \cdot \varphi (Q))
$$
on the Mordell-Weil group $E(K)$ is called
the {\rm height pairing}, and the lattice
$$
   (E(K)/E(K)_{\it tor} , <,>)
$$
is called the {\rm Mordell-Weil Lattice} of the elliptic
curve $E/K$ or of the elliptic surface $f$: $Y \longrightarrow
{\bf P}^1$.
\end{definition}
\begin{theorem} {\rm (Explicit formula for the height pairing) 
[30, Theorem 8.6]} For any $P$,$Q \in E(K)$, 
we have
$$
<P ,Q> = \chi + (P{\cal O}) + (Q{\cal O}) - (PQ) - \sum_{\nu \in R} contr_{\nu}
(P,Q),
$$
$$
<P,P>=2\chi + 2(P{\cal O}) - \sum_{\nu \in R} contr_{\nu}
(P).
$$
\end{theorem}
\begin{remark} {\rm Here $\chi =\chi ({\cal O}_S )$, and $(P{\cal O})$ is the intersection number of the sections $(P)$
and $({\cal O})$, and similarly for $(Q{\cal O})$,$(PQ)$. The term  $contr_{\nu}
(P,Q)$ stands for the local contribution at $\nu \in R$, which is
 defined as follows: suppose that $(P)$ interests $\Theta_{\nu ,i}$
and $(Q)$ intersects $\Theta_{\nu ,j}$. Then we let
$$
contr_{\nu} (P,Q) = \left \{\begin{array}{ll} 
                (-A_{\nu}^{-1})_{i,j} ,& if \ i\geq 1,j\geq 1, \\
                               0                 ,& otherwise.
 \end{array}
\right.
$$
where the first one means the $(i,j)$-entry of the matrix $(-A_{\nu}^{-1})$.
Further we set
$$
 contr_{\nu} (P) =  contr_{\nu} (P,P).
$$

\hskip -0.5cm Arrange $\Theta_i = \Theta_{\nu ,i}$ $(i=0,1,\cdots ,m_{\nu} -1)$ so 
that the simple components are numbered as in the figure below.
\vskip 0.5cm
\centerline{\font\thinlinefont=cmr5
\begingroup\makeatletter\ifx\SetFigFont\undefined
% extract first six characters in \fmtname
\def\x#1#2#3#4#5#6#7\relax{\def\x{#1#2#3#4#5#6}}%
\expandafter\x\fmtname xxxxxx\relax \def\y{splain}%
\ifx\x\y   % LaTeX or SliTeX?
\gdef\SetFigFont#1#2#3{%
  \ifnum #1<17\tiny\else \ifnum #1<20\small\else
  \ifnum #1<24\normalsize\else \ifnum #1<29\large\else
  \ifnum #1<34\Large\else \ifnum #1<41\LARGE\else
     \huge\fi\fi\fi\fi\fi\fi
  \csname #3\endcsname}%
\else
\gdef\SetFigFont#1#2#3{\begingroup
  \count@#1\relax \ifnum 25<\count@\count@25\fi
  \def\x{\endgroup\@setsize\SetFigFont{#2pt}}%
  \expandafter\x
    \csname \romannumeral\the\count@ pt\expandafter\endcsname
    \csname @\romannumeral\the\count@ pt\endcsname
  \csname #3\endcsname}%
\fi
\fi\endgroup
\mbox{\beginpicture
\setcoordinatesystem units <0.60000cm,0.60000cm>
\unitlength=0.60000cm
\linethickness=1pt
\setplotsymbol ({\makebox(0,0)[l]{\tencirc\symbol{'160}}})
\setshadesymbol ({\thinlinefont .})
\setlinear
%
% Fig POLYLINE object
%
\linethickness= 0.500pt
\setplotsymbol ({\thinlinefont .})
\putrectangle corners at  1.746 18.415 and  6.350 18.415
%
% Fig POLYLINE object
%
\linethickness= 0.500pt
\setplotsymbol ({\thinlinefont .})
\putrectangle corners at  1.905 16.351 and  6.350 16.351
%
% Fig POLYLINE object
%
\linethickness= 0.500pt
\setplotsymbol ({\thinlinefont .})
\putrectangle corners at  1.746 16.351 and  1.905 16.351
%
% Fig POLYLINE object
%
\linethickness= 0.500pt
\setplotsymbol ({\thinlinefont .})
\putrectangle corners at  2.381 19.209 and  2.381 15.558
%
% Fig POLYLINE object
%
\linethickness= 0.500pt
\setplotsymbol ({\thinlinefont .})
\putrectangle corners at  5.239 19.209 and  5.239 15.716
%
% Fig POLYLINE object
%
\linethickness= 0.500pt
\setplotsymbol ({\thinlinefont .})
\putrectangle corners at 10.001 16.351 and 14.287 16.351
%
% Fig POLYLINE object
%
\linethickness= 0.500pt
\setplotsymbol ({\thinlinefont .})
\putrectangle corners at 13.335 18.098 and 16.669 18.098
%
% Fig POLYLINE object
%
\linethickness= 0.500pt
\setplotsymbol ({\thinlinefont .})
\putrectangle corners at 15.716 16.351 and 20.479 16.351
%
% Fig POLYLINE object
%
\linethickness= 0.500pt
\setplotsymbol ({\thinlinefont .})
\putrectangle corners at 13.811 18.733 and 13.811 15.558
%
% Fig POLYLINE object
%
\linethickness= 0.500pt
\setplotsymbol ({\thinlinefont .})
\putrectangle corners at 16.034 18.733 and 16.034 15.716
%
% Fig POLYLINE object
%
\linethickness= 0.500pt
\setplotsymbol ({\thinlinefont .})
\putrectangle corners at 10.478 17.145 and 10.478 15.558
%
% Fig POLYLINE object
%
\linethickness= 0.500pt
\setplotsymbol ({\thinlinefont .})
\putrectangle corners at 11.906 17.145 and 11.906 15.558
%
% Fig POLYLINE object
%
\linethickness= 0.500pt
\setplotsymbol ({\thinlinefont .})
\putrectangle corners at 17.304 17.145 and 17.304 15.716
%
% Fig POLYLINE object
%
\linethickness= 0.500pt
\setplotsymbol ({\thinlinefont .})
\putrectangle corners at 18.891 17.145 and 18.891 15.716
%
% Fig POLYLINE object
%
\linethickness= 0.500pt
\setplotsymbol ({\thinlinefont .})
\putrectangle corners at 10.478 17.145 and 10.478 16.828
%
% Fig POLYLINE object
%
\linethickness= 0.500pt
\setplotsymbol ({\thinlinefont .})
\putrectangle corners at 10.478 16.034 and 10.478 15.558
%
% Fig TEXT object
%
\put{\SetFigFont{7}{8.4}{rm}$I_b$} [lB] at  1.429 20.479
%
% Fig TEXT object
%
\put{\SetFigFont{7}{8.4}{rm}  $\Theta_1$ } [lB] at  1.270 17.462
%
% Fig TEXT object
%
\put{\SetFigFont{7}{8.4}{rm}$\Theta_0$} [lB] at  3.334 15.716
%
% Fig TEXT object
%
\put{\SetFigFont{7}{8.4}{rm}$\Theta_{b-1}$} [lB] at  5.397 17.304
%
% Fig TEXT object
%
\put{\SetFigFont{7}{8.4}{rm}$\Theta_0$} [lB] at  9.684 16.510
%
% Fig TEXT object
%
\put{\SetFigFont{7}{8.4}{rm}$\Theta_1$} [lB] at 11.589 16.510
%
% Fig TEXT object
%
\put{\SetFigFont{7}{8.4}{rm}$\Theta_2$} [lB] at 16.669 16.669
%
% Fig TEXT object
%
\put{\SetFigFont{7}{8.4}{rm}$\Theta_3$} [lB] at 19.050 16.669
%
% Fig TEXT object
%
\put{\SetFigFont{7}{8.4}{rm}$I_b^*$} [lB] at  9.366 20.161
\linethickness=0pt
\putrectangle corners at  1.270 20.892 and 20.504 15.532
\endpicture}
}  
\vskip 0.5cm
\begin{center}
Figure $2.1$
\end{center}
 For the other types of reducible fibres, the numbering is irrelevant.
Assume that $(P)$ intersects $\Theta_{\nu ,i}$ and $(Q)$ intersect
$\Theta_{\nu ,j}$ with $i>1$,$j>1$. Then we have the following 
Table 2.1: the fourth row is for the case $i<j$ (interchange $P$, $Q$
if necessary).
\begin{center}
Table $2.1$
\end{center}
\vskip 0.5cm
\begin{center}
\begin{tabular}{|c|c|c|c|c|c|c|}\hline
$T_{\nu}^-$ & $A_1$ & $E_7$ & $A_2$ & $E_6$ & $A_{b-1}$ &  $D_{b+4}$ 
\\ \hline 
type of $F_{\nu}$ & $III$ & $III^*$ & $IV$ & $IV^*$ & $I_b (b\geq 2 )$
& $I_b^* (b\geq 0)$ \\  \hline
$contr_{\nu} (P)$ &$\frac1 2$ & $\frac 3 2$ & $\frac2 3$ & $\frac4 3$&
$\frac{i(b-i)} b$ & $\left \{\begin{array}{ll} 
                 1 , & i=1 \\  
                  1+\frac b 4, & i>1
           \end{array}
\right.$ \\  \hline
$contr_{\nu} (P,Q) (i<j)$& $-$ & $-$ & $\frac1 3$ &$\frac2 3$ &
$\frac{i(b-j)} b$ &  $\left \{\begin{array}{ll} 
                         \frac1 2 , & i=1 \\  
                         \frac{(2+b)}4 ,& i>1
                  \end{array}
\right.  $ \\ \hline
\end{tabular}
\end{center}
}
\end{remark}

\vskip 0.5cm
\hskip -0.6cm {\bf C. Weierstrass model of elliptic surfaces}
\vskip 0.5cm
\hskip -0.6cm Let $f: Y \longrightarrow C$ be a relatively minimal elliptic surface with the global section ${\cal O}$, and $K(Y)$ and $K(C)$ the function fields of $Y$ and of $C$ respectively. $f$ induces a homomorphism $f^* : K(C) \longrightarrow K(Y)$, and $K(Y)$ is a transcendential extension of $K(C)$ of transcendence degree one. $Y$ is birationaly equivalent to the subscheme $Y^*$ in $Proj({\cal O}\oplus {\cal O}(2L)\oplus {\cal O}(3L))$, which is given by the equation
$$
V^2W=4U^3-g_2UW^2-g_3W^3 ,
$$
where ${\cal O}$ is the structure sheaf of $C$, $L$ a line bundle on $C$ with $deg L= p_g -q+1$,  $g_2 \in H^0(C, {\cal O}(4L))$ and $g_3 \in H^0(C, {\cal O}(6L))$ are the sections with $\Delta =g_2^3-27g_3^2 
\not\equiv 0$.
\begin{theorem} {\rm (cf. [Kas] )} \  \  $Y^*$ is an algebraic surface with rational double points as the only singularities. $Y$ is the minimal resolution of $Y^*$. $Y^*$ is determined by $g_2$, $g_3$ up to ${\bf C}^*$-operation
$$
(g_2,g_3) \longrightarrow (\lambda^4 g_2 , \lambda^6 g_3),  \lambda \in {\bf C}^*={\bf C}-\{0\}.
$$
$g_2$ , $g_3$ satisfy

$(i)$ \ $\Delta =g_2^3 -g_3^2 \not\equiv 0$,

$(ii)$ {\rm min}$(3\nu_p (g_2), 2\nu_p (g_3)) <12$ for all $p \in C$,

where $v_p(g_2)$, $v_p(g_3)$ and $v_p(\Delta)$ are the orders of the zeros of $g_2$, $g_3$ and $\Delta$ at $p$. The singular fiber in $Y^*$ over $p$ consists of the minimal resolution of the rational double point and the rational curve, which is defined by the section ${\cal O}$. The type of rational double and thereby the type of the singular fibre determines $v_p(g_2)$, $v_p(g_3)$ and $v_p(\Delta)$. $Y^*$ is call the {\rm Weierstrass model} of the elliptic surface.
\vskip 0.2cm
\hskip -0.6cm The ${\cal J}$-invariant of the model is ${\cal J} = \frac{g_2^3}{\Delta}$.
\end{theorem}
For the value of $v_p(g_2)$, $v_p(g_3)$ and $v_p(\Delta)$ of the singular fibre, we refer to [Her].
\vskip 0.5cm
\hskip -0.6cm {\bf D. Classification of elliptic surfaces with at most $4$ singular fibers}
\vskip 0.5cm
In [Hir], U.Schmickler-Hirzebruch classifies all relatively minimal elliptic fibrations $Y$ over the complex projective line ${\bf P}^1$ which have at most three singular fibers. There are $36$ possible combinations of singular fibers (each of the $36$ potential representations yields different fiber types). Furthermore, U.Schmickler-Hirzebruch exhibits each of the $36$ combinations of singular fibers on an elliptic surface of the desired type and shows that except for the trivial case $(I_0,I_0,I_0)$ and the case $(I_0^*,I_0^*)$, the surface is determined by the singular fiber types. On the other hand, except for the trivial case, those surfaces corresponding to rational or $K3$ surfaces. There are exactly {\em eleven} elliptic $K3$ surfaces with $3$ singular fibers (cf. [SZ], [Ye]). Thus we get the following theorem.

\begin{theorem} Let $Y$ be a rational elliptic surface with at least one and at most three singular fibers. Then except for singular fiber type $(I_0^*,I_0^*)$, the surface is unique up to isomorphisms.
\end{theorem}

In [Her], S.Herfurtner studies the elliptic surfaces with four singular fibers and nonconstant ${\cal J}$-invariant over ${\bf P}^1$. He also distinguish two sets:
$$
T^+=\{I_n (n\geq 0), II,III, IV\} \ \ {\rm and} \ \ T^-=\{I_n^*(n\geq 0), IV^*,III^*,II^*\}.
$$ 
\begin{theorem} {\rm (cf. [Her] )} \ \ Let $Y\longrightarrow {\bf P}^1$ be a minimal elliptic surface with section,  nonconstant ${\cal J}$-invariant and four singular fibers, of which at most one is in $T^-$.

$(i)$ \ If one singular fiber is in $T^-$, three in $T^+$, the Weierstrass model depend on a parameter. Given four different base points, in the case $(I_1^*,I_1,I_1,II)$, there exists four, in the cases $(I_1^*,I_1,I_2,II)$ and $(I_1,I_1,II,IV^*)$, there exist two elliptic surfaces, each depending on the ${\cal J}$-invariant, and for all other fiber combinations there exists precisely one elliptic surface.

$(ii)$\ \ If all four singular fibers are in $T^+$, then the Weierstrass models are determined uniquely up to isomorphisms, except for one combination: $(I_1,I_6,II,III)$, where there are two nonisomorphic models.
\end{theorem}

\vskip 0.5cm
\hskip -0.6cm {\bf E. Rational surfaces}
\vskip 0.5cm
\begin{lemma} Let $Y$ be a smooth rational surface. Then we have

$(1)$ Suppose that $L$ is a smooth rational curve with $L^2 \geq 0$. Then $|L|$ is base point free and $h^0(Y,L)=2+L^2$.

$(2)$ Suppose that $L$ is an irreducible curve with $p_a(L)=1$ and $L^2 \geq 1$. Then a general member of $|L|$ is smooth and $h^0(Y,L)=L^2+1$.
If $L^2 \geq 2$, then $Bs|L|=\emptyset$. If $L^2=1$, then $|L|$ has exactly one base point.

\end{lemma}

{\em Proof.} (1) is clear when $L^2=0$. So we assume $L^2 \geq 1$. 
Consider the exact sequence below
$$
0 \longrightarrow {\cal O}_Y \longrightarrow {\cal O}_Y (L) \longrightarrow {\cal O}_{{\bf P}^1}(L^2) \longrightarrow 0. 
$$
Then,  we have $h^0(Y,L) = 1 + h^0({\bf P}^1, {\cal O} (L^2)) = 1 + (L^2 +1)
=L^2 +2$. 

To see $|L|$ is base point free, we let $\sigma : L_1 \longrightarrow L$ be the blow-up at $q\in L$ with $E$ the exceptional curve and let 
$L'$ be the proper transform of $L$. Then  $L'^2=L^2-1$ and $h^0(L')=L'^2 +2=L^2 + 1= h^0(L)-1$. On the other hand, $|L'|=|\sigma^* (L) - E| = |L-q|$; thus $|L|$ is base point free for
$$
h^0(L-q)= h^0(L)-1.
$$

(2) Considering the similar exact sequence as in (1) and the equality
$$
\chi ({\cal O}_L (L)) = h^0(L, {\cal O}_L (L))- h^1(L, {\cal O}_L (L)) =
\chi({\cal O}_L) + deg {\cal O}_L (L),
$$
we get
\begin{eqnarray*}
h^0 (Y,L)&=&1+ h^0(L, {\cal O}_L (L)) = 1+ h^0(L, \omega_L (-L))+ 1-p_a(L) + L^2  \\
&=& 1+L^2 + h^0 (L, \omega_L (-L)) = 1+L^2.
\end{eqnarray*} 

Here $\omega_L$ is the dualizing sheaf and we have applied Serre  duality for $L$. From this precise formula for $h^0(Y, L)$ and the similar discussion as in (1), we see that $Bs|L|=\emptyset$ if $L^2 \geq 2$. For the last assertion, see [DZ].

\section{Reduction to the Rational Elliptic Surfaces }
\vskip 0.5cm
Before starting the results, we explain our terminology.
\vskip 0.5cm
$(a)$ \ \  Let $V$ be a projective, normal surface with only rational double points as singularities. As usual, ratinal double points are indicated by their Dynkin types $A_n$, $D_n$, $E_6$, $E_7$ and $E_8$. When we say $V$ a surface of type $2A_1+D_6$ for example, it means that $V$ has $3$ singular points, one of which is of type $D_6$ and the other two are of type $A_1$. We indicate this by writing $V(2A_1+D_6)$.
\vskip 0.5cm
$(b)$ \ For  a log del Pezzo surface $V(2A_1+D_6)$ for example,  we denote by $U(2A_1 +D_6)$ the minimal resolution of  $V(2A_1 +D_6)$.
\vskip 0.3cm

For readers' convenience, we give a short proof of a result in [5] using Kawamata-Viehweg Vanishing Theorem which was not available at that time.

\begin{lemma}  Let $U$ be the minimal resolution of a Gorenstein log del Pezzo surface $V$. Then $|-K_U|$ has a reduced irreducible member. 
\end{lemma}

{\em Proof.}\ Note that $K_U$ is the pull back of $K_V$ and $K_U^2=K_V^2$ (cf. [1, 2]). Since $-K_U$, being the pull back of $-K_V$, is nef and big,  by the Kawamata-Viehweg Vanishing theorem,
we have $H^1(U, -K_U)= (0)$ (cf. [KMM]). Hence by the Riemann-Roch theorem, we get
$\dim |-K_U| = K_U^2$. Suppose that a member $A$ of $|-K_U|$ contains an
arithmetic genus $\geq 1$ irreducible component $A_0$. The Riemann-Roch theorem implies that $ |A_0 + K_U| \not= \emptyset$. From this and $0=A+K_U=(A_0 +K_U) +(A-A_0)$, we deduce that $A_0=A$ with $p_a (A_0)=1$ 
and Lemma 3.1 is true.

  Thus we may assume that every member of $|-K_U|$ is a union of smooth rational curves. The Stein factorization and the fact that $q(U)=0$ imply that a general member of $|-K_U|$ is of the form $M_1 + \cdots + M_k +F$, where $F$ is the fixed part of the linear system, $M_i \cong {\bf P}^1$, and $M_i \sim M_j$.  

  Now suppose that $K_U^2 =1$. If $F=0$, then $k=1$ and $M_1^2 =1$ and Lemma 3.1 is true. Since $-K_U$ is nef and big, it is $1$-connected by a result of C.P.Ramanujam (cf. [Reid]). Hence $1=K_U^2 \geq (kM_1 + F).kM_1 \geq 1+ k^2 M_1^2$. Thus $M_1^2 =0$, $k=1$, $M_1.F=1$ and $K_U.F=0$. Now intersecting the relation $-K_U \sim M_1 +F$ with the smooth rational curve $M_1$ of self-intersection $0$, one gets a contraction. So Lemma 3.1 is true when $K_U^2=1$.

For $K_U^2 \geq 2$, let $U_1 \longrightarrow U$ be the blow-up of a point on $M_1-(M_2 + \cdots + M_k +F)$. then $-K_{U_1}$ is linearly equivalent to the proper transform of $M_1+\cdots +M_k+F$ and hence $nef$ and $big$. If $U_1 \longrightarrow V_1$ is the contraction of all $(-2)$-curves then $V_1$ is a Gorenstein log del Pezzo surface. For $U_1$, we may argue as in the case of $U$ and have reduced the proof of Lemma 3.1 to the case $K_U^2=1$, which has been dealt with in the previous paragraph. This completes the proof of Lemma 3.1.

\begin{lemma} Let $U$ be the minimal resolution of a Picard number one Gorenstein log del Pezzo surface $V$ with $K_V^2=d$ where $2\leq d \leq 7$. Then we can construct another surface $U_{d-1}$ which is the minimal resolution of a Picard number one Gorenstein log del Pezzo surface $V'$ with $K^2_{V'}=d-1$.
\end{lemma}

{\em Proof.}  \ From Lemma 3.1,  we may choose an irreducible member of $|-K_U|$, say $A$ and let $E_d$ be a $(-1)$-curve on $U$. Such a $(-1)$-curve exists because $K_V^2 \leq 7$. Then $A$ meets $E_d$ at a point $q$. We may choose $A$ general such that $E_d$ is the only $(-1)$-curve through $q$.  Let $U_{d-1} \longrightarrow U$ be the blow-up of $q$ with $E_{d-1}$ the except curve. Then $-K_{U_{d-1}}$ is linearly equivalent to the proper transform $A_{d-1}$ of $A$. Since $K_{U_{d-1}}^2 >0$, the divisor $-K_{U_{d-1}}$ is nef and big. And the curves having $0$ intersection with $A_{d-1}$ are precisely the inverse image 
of the $(-2)$-curves on $U$ (contractible to singular points on $V$)
and the proper transform $E'_d$ of $E_d$. It follows that the contraction of all the $(-2)$-curves on $U_{d-1}$ gives a Picard number one Gorenstein log del Pezzo surface $V'$.

\begin{lemma}  Let $V$ be a Picard number one Gorenstein log del Pezzo surface with $K_V^2=1$ and $U$ the minimal resolution of $V$ with $D$ the exceptional divisor. Then there exists an extremal rational elliptic surface $Y$ such that $g : Y \longrightarrow U$ is the
blow-down of one  $(-1)$-curve.
\end{lemma}

{\em Proof.}\hskip 0.5cm Since $K_U $ is the pull back of $K_V$, we have $K_U^2=K_V^2=1$, $K_U^2 =10-\rho (U)$, $\rho (U)=\rho (V)+ \sharp (D)$ and $\sharp (D) =8$. By Lemma 2.14, we know $Bs|-K_U|= \{p\}$.  

Let $f: Y \longrightarrow U$ be the blow-up at $p$, then with Lemma 2.14, 
we have

(1)   $Y$ is a rational surface.

(2)   $K_Y^2 =0$ and $|-K_Y|$ induces an elliptic fibration of $Y$.

(3)   $\rho (Y) = 10 - K_Y^2 =10$.

(4)  By the formula of Shioda (cf. [S])
$$
  10=\rho (Y)= \sharp (D) + {\rm rank \ MW(Y)} +2,
$$
we get  ${\rm rank \ MW(Y)}=0$.
\vskip 0.3cm
Thus $Y$ is an extremal rational elliptic surface. This proves the lemma.
\vskip 0.5cm
Combining Lemmas 3.2 and 3.3, we get the following Theorem 3.4.

\begin{theorem}  Let $V$ be a Picard number one Gorenstein log del Pezzo surface with $K_V^2 \leq 7$ and $U \longrightarrow V$  the minimal resolution . Then there exists an extremal rational elliptic surface $Y$ such that $f: Y \longrightarrow U$ is a succession of blow downs of some $(-1)$-curves. 
\end{theorem}

\vskip 0.3cm
\hskip -0.6cm  We may generalize the above Theorem $3.4$ to the Picard number $r$ ($r\geq 2$) Gorenstein log del Pezzo surfaces.

 \begin{theorem}
 Let $V$ be a Picard number $r$ ($r\geq 2$) singular Gorenstein log del Pezzo surface and $U \longrightarrow V$ the minimal resolution. Then there exists a Mordell-Weil rank $(r-1)$ rational elliptic surface $Y$ such that $f: Y \longrightarrow U$ is a succession of blow downs of some $(-1)$-curves.
\end{theorem}
\vskip 0.3cm
The definition of the Mordell-Weil rank $r$ rational elliptic surface is given in the Secition $2$. 
\vskip 0.3cm

{\em Proof.} \ \ For simplicity, we only prove Theorem 3.5 for $r=2$. The general case can be
proved in the similar way.
\vskip 0.3cm
We let $U \longrightarrow V$  be the minimal resolution of a Picard number two Gorenstein log del Pezzo surface $V$ with $D$ the exceptional divisor and $d:=K_V^2\geq 2$. From Lemma 3.1, we know that $|-K_U|$ has a reduced irreducible member. 

Since  $10-K_V^2=\rho (V)+ \sharp (D)$,  $\rho (V)=2$ and  $\sharp (D) \geq 1$ beacuse $V$ is a singular surface, 
we get $ 10- K_V^2 \geq 3$ and $K_V^2 \leq 7$. On the other hand, since $K_U$ is the pull back of $K_V$, $K_U^2 =K_V^2 \leq 7$.  

Thus we may choose an irreducible member of $|-K_U|$ say $A$ and let $E_d$ be a $(-1)$-curve on $U$. Such a $(-1)$-curve exists because $K_V^2 \leq 7$. Then $A$ meets $E_d$ at a point $q$. 

Since $K_U^2 \geq 2$, we know $|-K_U|$ is base point free by Lemma 2.14. Furthermore since the smooth reduced irreducible member of $|-K_U|$ is dense in $|-K_U|$ and since there exists only finitely many $(-1)$-curves on $U$ (cf. the proof of [6, Lemma $13$]), we get the following Claim.

\begin{claim} We may choose $A \in |-K_U|$ such that there doesn't exist another $(-1)$ curve passing through the point $q$.
\end{claim}
Let $U_{d-1}$ be the blow-up of $q$ with $E_{d-1}$ the exceptional curve and $E_d'$ the proper transform of $E_d$. Then by the similar discussion as Lemma 3.2, we know that $U_{d-1}$ is the minimal resolution of a Gorenstein log del Pezzo surface $V'$ with $K^2_{V'}=d-1$. 

\begin{claim} 
 $\rho (V')=2$.
\end{claim}

{\em Proof of the Claim.}\ It is clearly that $\rho (V') \leq 2$. Thus to prove the Claim, we need to show that $\rho (V')\not= 1$. Suppose $\rho (V')=1$. Then from the above construction, we know that there must exist a curve $C$ on $U$ such that $C^2 \geq 0$ and $mult_p C \geq 2$ and the proper transform $C'$ of $C$ in $U_{d-1}$ is an $(-2)$-curve. Thus we have $C'.E_{d-1} = mult_p C=m \geq 2$, $E_d'.E_{d-1}=1$ where $C'$ and $E_d'$ are $(-2)$-curves.
  
  Now we pay attention to the set $I= \{E_d', E_{d-1} ,C'\}$ in $U_{d-1}$.  
Clearly we may do the same operation (blow-up) on $U_{d-1}$ as before and finally get an rational elliptic surface $Y$.  Since $C'$ and $E_d'$ are already $(-2)$-curves  we will keep the same symbols in $Y$ and let $E'$ be the proper transform of $E_{d-1}$. To prove the Claim, we need to consider two cases.

1.  $E'$ is an $(-1)$-curve.

In this case and from the elliptic fibration of $Y$, we see that $C'$ and
$E_d'$ will be in two different elliptic fibers and $E'$ is a section of $Y$ which is a contradiction since $C'.E'=C'.E_{d-1} \geq 2$ .

2. $E'$ is an $(-2)$-curve.

In this case we see that $E_d'$, $E'$ and $C'$ will be in one elliptic fiber. We also get a contracdiction since  $C'.E' \geq 2$ and we can't find an elliptic fiber contains $C'$ and $E'$.

Thus we prove the Claim. We note that the above proof  can easily be generalized to prove the general case.

\vskip 0.3cm
\hskip -0.6cm   Continue this process. Since $K_{U_{d-1}}^2=d-1 < K_U^2$, this process must stop at the $(d-1)$-th step. Thus,  we obtain a surface $U_1$ which is the minimal resolution of the Picard number two Gorenstein log del Pezzo surface $V_1$ with $K^2_{V_1}=1$. On $U_1$ we can not do the above process again. But since $K_{U_1}^2=1$, $|-K_{U_1}|$ has exactly one base point $q$ by Lemma 2.14. Blowing up $q$ will give rise to a Mordell-Weil rank one rational elliptic surface $Y$ as in Lemma 3.3. Thus Theorem 3.5 is proved.  

\begin{remark}
 From Theorems $3.4$ and   $3.5$, together with {\rm [OS]} and {\rm [P]}, we may classify all singularity types of the Gorenstein log del Pezzo surfaces.

\end{remark}

\begin{question} From the proofs of Theorems $3.4$ and $3.5$, it is natural for us to ask the converse question: If $U \longrightarrow V$ is the minimal resolution of a Picard number  one (resp. Picard number $r \geq 2$) Gorenstein log del Pezzo surface $V$ . Let $E$ be a $(-1)$-curve on $U$ and $\phi : U \longrightarrow U'$  the blow-down of $E$, then could $U'$ be a minimal resoultion of a Picard number one (resp. Picard number $r \geq 2$) Gorenstein log del Pezzo surface $V'$?  
\end{question}

\begin{lemma} Let the notation be the same as in Question $3.9$. Let $f: U \longrightarrow U'$ be the blow-down of a $(-1)$-curve $E$ and $D$ the exceptional divisor of the minimal resolution $U \longrightarrow V$. Let $U' \longrightarrow V'$ be the contraction of all $(-2)$-curves and let $D'$ be the exceptional divisor.  Then

$(1)$ \ \ $V'$ is a Gorenstein log del Pezzo surface.

$(2)$ \ \ If $D.E=1$, $\rho (V')= \rho (V)$. The converse is also true.
\end{lemma}

{\em Proof.} \ We fix a nonsingular member  $A \in |-K_U|$, the existence
of which is proven in [D] . We also let $A'=f_*(A)$.
Then $A' \in |-K_{U'}|$ and $(A'^2)\geq (A^2) \geq 1$. Since $-K_{U}$ is nef and big, so is $-K_{U'}$. To vertify $(1)$, we have only to show that $-K_{V'}$ is ample. In other words, for every curve $C$ on $U'$ with $(-K_{U'},C)=0$, we must show that $C \leq D'$. Since $-K_{U'}$ is nef and big, by Hodge Index Theorem, we get $C^2 <0$. On the other hand,
$$
2p_a(C)-2 =C^2+C.K_{U'}=C^2 <0.
$$
Thus we have $p_a(C)=0$ and so $C={\bf P}^1$, $C^2=-2$ and $C \leq D'$. 
$(1)$ is proved. 
\vskip 0.5cm
\hskip -0.6cm
 $(2)$ Since $E.D=1$, we have $K_{U'}^2=K_U^2+1$ and
$\sharp (D)= \sharp (D') +1$. Thus the first part in the assertion  $(2)$ follows from the relations $K_U^2=K_V^2=10-\rho (U)$, $K_{U'}^2=K_{V'}^2=10-\rho (U')$, $\rho (U) = \rho (V) + \sharp (D)$and $\rho (U')=\rho (V') + \sharp (D')$. The converse part is clear from the above discussion.
\vskip 0.5cm
\hskip -0.6cm
We will call a $(-1)$-curve  which satisfies the part $(2)$ of Lemma $3.10$ a {\it Nice Exceptional Curve} or simply an {\it NEC}.

\section{Complete Classification of the Picard number one Gorenstein log del Pezzo Surfaces}

We shall prove Theorem $1.2$ in the present section.
\vskip 0.5cm
\hskip -0.6cm Let $V$ be a Picard number one Gorenstein log del Pezzo surfaces  and $g$ : $U \longrightarrow V$  a minimal resolution of singularities. Then $K_U^2=K_V^2 >0$, $|-nK_U| \not= \phi$ for $n>0$ and $U$ is a rational surface. Hence $K_U^2 + \rho (U)=10$. This implies that $1\leq K_U^2 \leq 9$.
We observe:
\vskip 0.5cm
(i) If $K_U^2=9$, then $V$ is smooth and isomorphic to ${\bf P}^2$.
\vskip 0.5cm
(ii) If $K_U^2=8$, then $U$ is the Hirzebruch surface $\Sigma_2$ and $V$ is obtained by contracting the $(-2)$-curve to an $A_1$-singularity. In this case $V$ is isomorphic to the quadric cone $Q := \{ X^2+Y^2+Z^2=0 \}$ in ${\bf P}^3$.
\begin{lemma}  $K_U^2\not= 7$. 
\end{lemma}

{\em Proof.} Assume the contrary, i.e., $K_U^2=7$. In this case
$U$ is obtained from a Hirzebruch surface $\Sigma_n$ by one blow-up with $E$ the exceptional curve. Let $S$ be a curve on $\Sigma_n$ with self-intersection $-n$  and $L$ a fibre of the ${\bf P}^1$-fibration on $\Sigma_n$. We may assume that the blow-up point $p\in L$. We have $K_{\Sigma_n} \sim -2S-(2+n)L$. 

If $p\not\in S$, then $K_U\sim -2S-(2+n)L'-(1+n)E$, where $L'$ is the proper transform of $L$. Clearly, $U$ contains a $(-2)$-curve different from $S$, say $C$. Since $K_U .C=0$, the curve $C$ is disjoint from $S$, $L'$ and $E$. This is impossible. 

If $p\in S$, then $K_U \sim -2S'-(2+n)L'-(3+n)E$ where $S'$ is the proper transform of $S$. Again there is a $(-2)$-curve different from $S$. Thus we get a contradiction as above. This completes the proof of the lemma.

\vskip 0.5cm
\hskip -0.6cm In the rest of this section, we assume $1\leq K_V^2 \leq 6$.
\vskip 0.3cm
\hskip -0.6cm 
We start with $K_V^2 =1$.
\vskip 0.3cm
\hskip -0.6cm
From Lemma 3.3, we know that if $U$ is the minimal resolution of a Picard number one Gorenstein log del Pezzo surface $V$,  then every NEC on $U$ comes from either the $(-1)$-curve on an extremal rational elliptic  surface $Y$ or the 
$(-2)$-curve which intersects with the exceptional curve of the blow-down $Y \longrightarrow U$. The $(-1)$-curves on $Y$ are the sections of the elliptic fibration and form the Mordell-Weil group of  the elliptic fibration $Y \longrightarrow {\bf P}^1$.   

 From [MP], we get a complete list of the extremal rational elliptic surfaces $Y$ which are listed in Table $4.1$. For convenience, we label the surfaces. By contracting the zero section $\cal O$ in $MW(Y)$, we get a surface $U$ which is the minimal resolution of a Picard number one Gorenstein log del Pezzo surface $V$ with $K_V^2 = 1$.   By a surface $U_j^i$, we mean that it is the blow-down of some $(-1)$-curves from the No.$i$ extremal rational elliptic  surface and $K^2=j$. We let $n$ be the number of the NEC on $U_1$.
\vskip 0.5cm
\begin{center}
Table $4.1$
\end{center}
\begin{center}
\begin{tabular}{|c|c| c|c|c|c|c|c|}\hline
$\sharp$ & $Y$  &  $U_1$  & $n$& $\sharp$ & $Y$ & $U_1$ & $n$
\\ \hline
$1$& $II, II^*$  & $U^1_1(E_8)$& $1$& $9$&$I_1^*,I_4,I_1$  & $U^9_1(D_5+A_3)$ &$1$
\\ \hline
$2$& $III, III^*$  &$U^2_1(A_1+E_7)$  &$1$ & $10$&$ I_2^*,I_2,I_2$ & $U^{10}_1(D_6+2A_1)$ &$1$ 
\\ \hline
 $3$& $IV, IV^*$ & $U^3_1(A_2+E_6)$ &$1$  & $11$&$ I_9,I_1,I_1,I_1  $ & $U^{11}_1(A_8)$ &$2$      \\ \hline
$4$&$I_0^*,I_0^*$   & $U^4_1(2D_4)$ & $2$&$12$& $I_8,I_2,I_1,I_1$  &$U^{12}_1(A_7+A_1)$ &$1$
\\ \hline
$5$&$II^*,I_1,I_1$   & $U^5_1(E_8)$  &$1$ &$13$& $ I_5,I_5,I_1,I_1$  &$U^{13}_1(2A_4)$ &$0$
\\ \hline
$6$&$ III^*,I_2,I_1$   & $U^6_1(E_7+A_1)$ &$1$ &$14$& $ I_6,I_3,I_2,I_1$  & $U^{14}_1(A_5+A_2+A_1)$ &$0$
\\ \hline
$7$&$IV^*,I_3,I_1$  & $U^7_1(E_6+A_2)$  &$1$& $15$&$I_4,I_4,I_2,I_2$  &$U^{15}_1(2A_1+2A_3)$ &$0$
\\ \hline
$8$&$I_4^*,I_1,I_1$   &$U^8_1(D_8)$ & $2$ & $16$&$I_3,I_3,I_3,I_3$   &$U^{16}_1(4A_2)$ &$0$
\\ \hline
\end{tabular} 
\end{center} 

From the proof of Lemma 3.3, we know there is a one to one correspondence between the No.$i$ extremal rational elliptic surfaces $Y$ and $U^i_1$.  On the other hand, for every $i\not=4$, there is only one isomorphism class of the No.$i$ extremal rational elliptic surface (cf. [MP, Theorem $5.4$]). The extremal rational elliptic surfaces of No.$4$, i.e., with singular fibers $(I_0^*, I_0^*)$, are parameterized by $j\in {\bf C}$ [ibid.].  Thus we get the following lemma:

\begin{lemma} Let  $V$ be the Picard number one  Gorenstein log del Pezzo surface with $K^2_V=1$. Then the singularity type of $V$, the number $m$ of the isomorphism class with the given singularity type are listed in Table $4.2$.

\vskip 0.5cm
\begin{center}
Table $4.2$
\end{center}
\begin{center}
\begin{tabular}{| c|c|c|c|}\hline
 Dynkin type of $V$ & $m$ & Dynkin type of $V$  & $m$\\ \hline
$E_8$ & $2$  &  $A_8$ & $1$\\ \hline
$A_1+E_7$ & $2$  &  $A_1+A_7$& $1$  \\ \hline
$A_2+E_6$ & $2$ & $2A_4$ & $1$ \\ \hline
$2D_4$ & $\infty$ &$A_1+A_2+A_5$ & $1$  \\ \hline
$D_8$ & $1$ & $2A_1+2A_3$ & $1$ \\ \hline
$A_3+D_5$ & $1$ & $4A_2$ & $1$  \\ \hline
$2A_1+D_6$ & $1$  & & \\ \hline
\end{tabular} 
\end{center} 
\end{lemma}

Next we consider $K_V^2=2$.
\vskip 0.5cm

By blowing down an NEC on $U_1$ which is listed in Table $4.3$, we get a surface $U_2$ which is the minimal resolution of a Picard number one Gorenstein log del Pezzo surface $V_2$. They are listed in the following table where we let $n$ be number of the NECs on $U_2$.
\vskip 0.5cm
\begin{center}
Table $4.3$
\end{center}
\begin{center}
\begin{tabular}{| c|c|c| c|c|c|c|} \hline
  $U_1$ & $U_2$ &$n$ & $U_1$ & $U_2$ &$n$ \\ \hline
  $U^1_1 (E_8)$ & $U^1_2(E_7)$ &$1$ &$U^8_1(D_8)$ & $U^8_2(A_1+D_6)$& $1$ \\ \hline
$U^2_1(A_1+E_7)$ & $U^2_2(A_1+D_6)$ &$1$ &$U^8_1(D_8)$ & $U^8_2(A_7)$ &$2$ \\ \hline
  $U^3_1(A_2+E_6)$ & $U^3_2(A_2+E_5)$ &$1$ &$U^9_1(A_3+D_5)$ & $U^9_2(A_1+2A_3)$& $0$\\ \hline
 $U^4_1(2D_4)$ & $U^4_2(3A_1+D_4)$& $0$ &$U^{10}_1(2A_1+D_6)$ & $U^{10}_2(3A_1+D_4)$ &$0$  \\ \hline
 $U^5_1(E_8)$ & $U^5_2(E_7)$ &$1$ &$U^{11}_1(A_8)$& $U^{11}_2(A_2+A_5)$& $1$ \\ \hline
 $U^6_1(A_1+E_7)$ & $U^6_2(A_1+D_6)$& $1$ &$U^{12}_1(A_1+A_7)$ & $U^{12}_2(A_1+2A_3)$ &$0$ \\ \hline
 $U^7_1(A_2+E_6)$ & $U^7_2(A_2+A_5)$& $1$  & & & \\ \hline

\end{tabular} 
\end{center} 

\begin{remark} \ $(i)$ \ In the following discussion, the Figure $i$ $(i=1,2,3,4)$ are listed in the configurations of  Appendix.

(ii) To see $U_2^{11}(A_2+A_5)$ also doesn't depend on the choice of the blow-down curve $E_i$($i=1,2$) on $U_1^{11} (A_8)$ (cf. Figure $2$),  we look at the No.$11$ extremal rational elliptic surface $Y$ (cf. Table $2$). A simple calculation shows that the Mordell-Weil group of $Y$ is ${\bf Z}/3{\bf Z}$ and there are three sections $E_0$,$E_1$ and $E_2$ (cf. Figure $2$). The section $E_1$ defines an automorphism $g$ of $Y$ such that $g(E_0)=E_1$, $g(E_1)=E_2$ and $g(E_2)=E_0$. For simplicity, we denote by $bd_i$ the blow-down of the $E_i$ $(i=0,1,2)$. Then
$$
bd_2 (U_1^{11}(A_8)) = bd_2 \circ bd_0 (Y) \cong bd_0 \circ bd_1 (Y) \cong bd_1 \circ bd_0 (Y) = bd_1 (U_1^{11} (A_8)).
$$
\end{remark}

\begin{claim} $U^1_2(E_7) \cong U^5_2(E_7)$.
\end{claim}
{\em Proof.} \ \ \ For simplicity,  we let 
$$
\begin{array}{llll}
U_1 =U_1^1(E_8), &   V_1 = V_1^1(E_8), &  U_2 =U_2^1(E_7), & V_2 = V_2^1(E_7),  \\
U'_1 = U_1^5(E_8), &  V'_1 =V_1^5(E_8), &  U'_2 = U_2^5(E_7), & V'_2 =V_2^5(E_7).
\end{array}
$$
By the construction of $U_2$, $U_2'$, there is a unique $(-1)$-curve $E_2$ on $U_2$. To prove $U_2 \cong U'_2$, we only need to show that $|-K_{U_2}|$ has a cuspidal member $A$ which is the image of the type $II$ fiber and also a nodal member $A'$; in fact if $U_1 \rightarrow U_2$ (resp. $U'_1 \rightarrow U_2$) is the blow-up of $q=A \cap E_2$ (resp. $q'=A' \cap E_2$) 
then $|-K_{U_2}|$ (resp.$|-K_{U'_2}|$) has a cuspidal (resp. nodal) member, whence 
$$U_2 \cong {\rm \ (the \ blow-down\ of \ the \ unique \ (-1)-curve  \ 
on \ U_1' \  )} = U'_2.
$$
 We let $Y$ (resp. $Y'$) be the extremal rational elliptic  surface with singular fiber type $II^*$, $II$ (resp. $II^*$,$I_1$,$I_1$). There is a   composition $Y \longrightarrow {\bf P}^2$ (resp. $Y' \longrightarrow {\bf P}^2$) of blow-down of the section $E$ and all components in the type $II^*$ fiber except for a multiplicity-$3$ component $C'_3$.

  Clearly the above birational morphism $Y \longrightarrow {\bf P}^2$ factors as $Y \longrightarrow U_2 \longrightarrow {\bf P}^2$, where $Y\longrightarrow U_2$ is the blow-down of $E$ and the multiplicity-1 component $C_1$ of $II^*$. 

  Thus the existence of $A$, $A'$ on $U_2$ is equivalent to the existence of a cuspidal cubic $B$ (the image of the type $II$ fibre) and a nodal cubic $B'$ on ${\bf P}^2$ with a common inflex (the image of the multiplicity-3 component $C'_3$ in the type $II^*$ fibre) defined by $Z=0$, with local intersection $I(B,B';[0,1,0])=7$ at the inflex point $[0,1,0]$. We may assume that $B=\{ Y^2Z=X^3\}$ and just let
$B'=\{ Y^2Z=X^3+XZ^2+\sqrt{-4/27} \}$. This completes the proof of the Claim.

\begin{claim}
$(1)$ $U^2_2(A_1+D_6)\cong U^6_2(A_1+D_6) \cong U^8_2(A_1+D_6)$  \ \ {\rm (cf. Figure $1$)}.

$(2)$ $U^3_2(A_2+A_5)\cong U^7_2(A_2+A_5) \cong U^{11}_2(A_2+A_5)$ \ \ {\rm (cf. Figure $2$)}.

$(3)$  $U^4_2(3A_1+D_4) \cong U^{10}_2(3A_1+D_4)$ \ \ {\rm (cf. Figure $3$)}.

$(4)$  $U_2^9(A_1+2A_3) \cong U_2^{12}(A_1+2A_3)$\ \ {\rm (cf. Figure $4$)}.
\end{claim}

{\em Proof.} 
(1)\ \ From the construction of $U^2_2(A_1+D_6)$,  there is an NEC $E_2$ 
on $U^2_2(A_1+D_6)$ and there is a smooth rational curve $L$ and $L^2=0$ such that $L\cap E_2 = \{p \}$. On the other hand, by Lemma 2.14, we know that $|L|$ is base point free and defines a ${\bf P}^1$-fibration. We choose a general member $L_1$ of $|L|$, and let  $\{q \}=L_1 \cap E_2$. It is easy to see that $U_1^2(A_1+E_7)$ is the blow-up of $U_2^2(A_1+D_6)$ at point $p$ and $U_1^6(A_1+E_7)$ is the blow-up of $U_2^2(A_1+D_6)$ at point $q$. Thus
$$
U_2^2(A_1+D_6) \cong ({\rm \ the \ blow-down \  of \ the \ uniqe \ NEC \ on\ } U_1^6(A_1+E_7)) = U^6_2(A_1+D_6)
$$
On the other hand, blowing up at point $r$ in $U_2^6(A_1+D_6)$ will produce
$U_1^8(D_8)$, and there is unique way from $U_1^8(D_8)$ to $U_2^8(A_1+D_6)$.
Thus
$$
U^2_2(A_1+D_6) \cong U^8_2(A_1+D_6).
$$
By the similar argument as (1), we may prove (2), (4).

\vskip 0.5cm
(3) \  Since $U_1^{10}(2A_1+D_6)$ is unique and there is only one $NEC$ on $U_1^{10}(2A_1+D_6)$,  we know $U_2^4(3A_1+D_4)$ is unique and
$$
U^4_2(3A_1+D_4) \cong U^{10}_2(3A_1+D_4).
$$
This proves the Claim.
\vskip 0.5cm
Since $U_i^j \cong U_{i'}^{j'}$ if and only if $V_i^j \cong V_{i'}^{j'}$, we get the following lemma:

\begin{lemma} Let  $V$ be the Picard number one  Gorenstein log del Pezzo surface with $K^2_V=2$. Then the singularity type of $V$ is one of the  following:
$$
E_7,  A_7,  A_1+D_6, A_1+2A_3, A_2+A_5,  3A_1+D_4.
$$
Furthermore, $V$ is uniquely determined by its singularity type.
\end{lemma}

Next we deal with $K_V^2=3$.
\vskip 0.5cm
By blowing down an NEC on $U_2$ which is listed in Table $4.4$, we get a surface $U_3$ which is the minimal resolution of a Picard number one Gorenstein log del Pezzo surface $V_3$. They are listed in the following table where we let $n$ be the number of the NEC on $U_3$.
\vskip 0.5cm
\begin{center}
Table $4.4$
\end{center}
\begin{center}
\begin{tabular}{| c|c|c|c|c|c|} \hline
  $U_2$ & $U_3$ &$n$   &$U_2$ & $U_3$ &$n$ \\ \hline
  $U^1_2 (E_7)$ & $U^1_3(E_6)$ &$1$ & $U^3_2(A_2+A_5)$ & $U^3_3(3A_2)$ &$0$ \\ \hline
$U^6_2(A_1+D_6)$ & $U^6_3(A_1+A_5)$ &$1$   & $U^8_2(A_7)$ & $U^8_3(A_1+A_5)$& $1$ \\ \hline
\end{tabular} 
\end{center} 

\begin{remark} $(i)$ There are two NECs, $E_0'$ and $E_1'$ on $U_2^8(A_7)$ (cf. Figure $1$).  To see $U_3^8(A_1+A_5)$ doesn't depend on the choice of the blow-down of the
$E_i'$ on $U_2^8(A_7)$ where $i=1,2$, we look at the No.8 extremal rational elliptic surface $Y$ (cf. Table $2$). A direct calculation shows that the Mordell-Weil group of $Y$ is ${\bf Z}/2{\bf Z}$ and there are two sections $E_0$, $E_1$ (cf. Figure $1$). The section $E_1$ define an automorphism $g$ of $Y$ which interchanges $E_0$ and $E_1$. And this $g$ will induce an automorphism of $U_2^8(A_7)$ which interchanges $E_i'$ ($i=0,1$). Thus blowing down $E_0$ or $E_1$ will produce the same surface. 

$(ii)$ \  By the similar argument as in Claim $4.5$, we have
 $$
U^8_3(A_1+A_5) \cong U^6_3(A_1+A_5) \ \ {\rm (cf. Figure \ 1)}.
$$ 
\end{remark}
Thus we have the following lemma:
\begin{lemma} Let  $V$ be the Picard number one  Gorenstein log del Pezzo surface with $K^2_V=3$. Then the singularity type of $V$ is one of the  following:
$$
E_6,  3A_2, A_1+A_5.
$$
$V$ is uniquely determined by its singularity type.

\end{lemma}

Finally we deal with \ $ 4\leq K_V^2 \leq 6$.
\vskip 0.5cm 
By blowing down an NEC on $U_3$ which is listed in Table $4.5$, we get a surface $U_4$ which is the minimal resolution of a Picard number one Gorenstein log del Pezzo surface $V_4$. 
\vskip 0.5cm
\begin{center}
Table $4.5$
\end{center}
\begin{center}
\begin{tabular}{| c|c|c|} \hline
  $U_3$ & $U_4$ &$n$   \\ \hline
  $U^1_3 (E_6)$ & $U^1_4(D_5)$ &$1$   \\ \hline
$U^8_3(A_1+A_5)$ & $U^8_4(2A_1+A_3)$ &$0$  \\ \hline
\end{tabular} 
\end{center}

\begin{lemma} Let  $V$ be the Picard number one  Gorenstein log del Pezzo surface with $K^2_V=4$. Then the singularity type of $V$ is one of the following:
$$
D_5, 2A_1+A_3
$$
Furthermore,  $V$ is uniquely determined by its singularity type.
\end{lemma}

By blowing down the unique NEC on $U^1_4(D_5)$, we get $U_5^1(A_4)$ and so the Picard number one Gorenstein log del Pezzo surface $V(A_4)$.  Then blowing down the unique NEC on $U_5^1(A_4)$, we get $U(A_1+A_2)$ with $\sharp \{{\it NEC} \ {\rm on}\  U(A_1+A_2) \}=0$, and so the Picard number one Gorenstein log del Pezzo surface $V(A_1+A_2)$.
\vskip 0.5cm
Combining Lemmas $4.2$, $4.6$, $4.8$, $4.9$ and the above discussion, we prove  Theorem $1.2$.

\section{Complete classification of the Picard number two relatively minimal Gorenstein log del Pezzo Surfaces}

We shall prove Theorems $1.5$ and $1.6$ in the present section. 
\vskip 0.3cm
\hskip -0.6cm
Let $V$ be a Picard number two Gorenstein log del Pezzo surface and $g$ : $U \longrightarrow V$  a minimal resolution of singularities. We also denote by $D:=g^{-1}({\rm Sing}V)$ the exceptional divisor.
\vskip 0.3cm
\hskip -0.6cm
We start with $K_V^2=1$.
\vskip 0.3cm
\hskip -0.6cm
 From Theorem $3.5$,we know that any rational surface $U$ which is the minimal resolution of a Picard number two Gorenstein log del Pezzo surface $V$ with $K_V^2=1$  comes from the Mordell-Weil {\em rank one} rational elliptic surface $Y$. Meanwhile, every $NEC$ on $U$ comes from either the $(-1)$-curve on  
$Y$ or the $(-2)$-curve which intersects the NEC  of the blow-down $Y \longrightarrow U$. The $(-1)$-curves on $Y$ are the sections of the elliptic fibration and form the Mordell-Weil group of  the elliptic fibration $Y \longrightarrow {\bf P}^1$.   
\vskip 0.3cm
\hskip -0.6cm
   From [M], [OS] and [P], we get a complete list of the Mordell-Weil {\em rank one} rational elliptic surfaces $Y$ and its Mordell-Weil group.  By contracting the zero section ${\cal O}$ in $MW(Y)$, we get a surface $U_1$ which is the minimal resolution of a Picard number two Gorenstein log del Pezzo surfaces $V$ with $K_V^2 = 1$.  For convenience, we label the surfaces which are listed in Table $1$ of  Appendix. By a surface $U_j^i$, we mean that it is the smooth blow-down of some $(-1)$-curves from the No.$i$ Mordell-Weil rank one rational elliptic  surface and $K^2=j$. Thus we get the following lemma.

\begin{lemma} Let  $V$ be the Picard number two  Gorenstein log del Pezzo surface with $K^2_V=1$. Then there exist $18$ singularity types which are listed in Table $5.1$.

\vskip 0.3cm
\begin{center}
Table $5.1$
\end{center}
\begin{center}
\begin{tabular}{| c|c|c|c|c|}\hline
$E_7$ & $D_5+A_2$  &  $A_7$ & $A_4+A_3$  & $2A_3+A_1$ \\ \hline
$E_6+A_1$ & $D_5+2A_1$  &  $A_1+A_6$& $A_4+A_2+A_1$  & $3A_2+A_1$\\ \hline
$D_7$ &  $D_4+A_3$ & $A_5+A_2$ & $A_3+A_2+2A_1$ &\\ \hline
$D_6+A_1$ & $ D_4+3A_1$ &$A_5+2A_1$ & $A_3+4A_1$ & \\ \hline
\end{tabular} 
\end{center} 
\end{lemma}
In order to find out the list of all Picard number two {\em relatively minimal}  Gorenstein log del Pezzo surfaces, we need to exclude the surfaces on which there is a $(-1)$-curve $E$ and a connected component $R$ of $D$ ($R$ may be zero) such that $E+R$ is a linear chain and $E+R$ is a connected component of $E+D$ (cf. Table $1$ of  Appendix).
In fact we have the following lemma:

\begin{lemma} None of $U_1^i$ ($7\leq i\leq 38$) is  the minimal resolution of a Picard number two relatively minimal Gorenstein log del Pezzo surface.
\end{lemma}

{\em Proof.} We verify the above assertions by checking each case separately. To exhibit our arguments, we treat only  No.$12$ and No.$24$,  i.e., 
$U_1^{12}(2A_3+A_1)$ and $U_1^{24}(3A_2+A_1)$ respectively. The other cases can be treated  similarly. We also write out the linear chain exists in the corresponging  surface in  Table $1$ of the Appedenix.
\vskip 0.3cm
\hskip -0.6cm
(1) $U_1^{12}(2A_3+A_1)$. 
\vskip 0.3cm
\hskip -0.6cm
The Mordell-Weil group of the No.$12$ Mordell-Weil rank one rational elliptic surface $Y$ (with fibre type $2I_4,I_2,2I_1$)  is $A_1^*\oplus {\bf Z}/4{\bf Z}$ where $A_1^*$ is the rank one lattice generated by an element $e$ such that $<e,e>=1/2$.  We let $P$ be the generator of the lattice $A_1^*$ and $Q$ the generator of the torsion part of the Mordell-Weil group,i.e., ${\bf Z}/4{\bf Z}$.  A direct calculation shows that there is one element in the set $\{ P_2, P_2+Q, P_2+2Q, P_2+3Q \}$ where $P_2:=2P$, say $P_2$, passing through the $0$-th components of the singular fibres $2I_4, I_2$
and $(P_2{\cal O} )=0$. Since $U_1^{12}(2A_3+A_1)$ is the blow-down of the zero section ${\cal O}$ of $Y$, the image of $P_2$ will be the $(-1)$-curve disjoint from $D=2A_3+A_1$. Thus $U_1^{12}(2A_3+A_1)$ is not the minimal resolution of a Picard number two relatively minimal Gorenstein log del Pezzo surface.
\vskip 0.3cm
\hskip -0.6cm      
$(2)$ $U_1^{24}(3A_2+A_1)$.
\vskip 0.3cm
\hskip -0.6cm
To start with, we let $P$ be the generator of the rank $1$ lattice $<1/6>$. Firstly we need to find out all the $(-1)$-curves on $U_1^{24}(3A_2+A_1)$. Since we blow down the zero section ${\cal O}$ on the No.24 Mordell-Weil rank one rational elliptic surface with the fibre type $(3I_3, I_2, I_1)$, the section $P_m :=mP$ ($m\in {\bf Z}$) which intersects
 ${\cal O}$ will not be a  $(-1)$-curves on $U_1^{24}(3A_2+A_1)$ again. Thus we get the following claim:

\begin{claim} A $(-1)$-curve on $U_1$ is either the image of the zero component of the singular fibre or the image of the section disjoint from the ${\cal O}$. 
\end{claim}
Suppose $(P{\cal O})=0$. Then from $m^2/6=<P_m,P_m>=2\chi + 2(P_m{\cal O}) - \sum_{\nu \in R} contr_{\nu} 
(P_m) = 2 -\sum_{\nu \in R} contr_{\nu}(P_m) \leq 2$,
we get $m^2 \leq 12 $ and $|m|\leq 3$. Conversely,  a simple calculation shows that $(P_3{\cal O})=0$ and $P_3$ only passes through the $1$-st component of $I_2$. Thus $P_3+A_1$ will form a linear chain required. 
\vskip 0.3cm
 
By a similar discussion as in Lemma 4.2, we know there is a one to one correspondence between No.$i$ Mordell-Weil rank one rational elliptic surfaces and $U^i_1$, where $1\leq i \leq 6$. Thus finding the number of the isomorphism classes of the relatively minimal Dynkin type on $K_V^2=1$ is equivalent to finding the number of the isomorphism classes of the rational elliptic surfaces with the corresponding singular fibre type.

\begin{remark} \ \ $(1)$ \ \ From {\rm [Hir]} or Theorem $2.12$ , we know 
the Nos.$1$ and $4$ Mordell-Weil rank one rational elliptic surfaces are unique up to isomorphism,i.e., are determined by the singular fibre types. 

$(2)$  \ \ From {\rm [Her]} or Theorem $2.13$, we know that the number of the isomorphism classes of each of the Nos. $2$, $5$ and $6$ Mordell-Weil rank one rational elliptic surfaces is  infinity. In fact, the isomorphism classes depend on a single parameter.
 
\end{remark}

\begin{lemma} The No.$3$ Mordell-Weil rank one rational elliptic surface $S$ is uniquely determined by its singular fibre 
type $(I_4, I_2,I_2,I_2, I_2)$, up to isomorphisms.
\end{lemma}

{\em Proof.} \ \ To see the uniqueness of $S$, we need to calculate the polynomial $g_2$, $g_3$ and $\Delta$ in the homogeneous coordinates $(X,Y)$ of ${\bf P}^1$ where $g_2$, $g_3$ are the sections which determine the {\it Weierstrass model} of $S$ (cf. Theorem $2.11$). For simplicity,
we let $G_2=g_2$ and $G_3=3\sqrt{3}g_3$. Then $\Delta=G_2^3-G_3^2$.

The orders $\nu_p$ of zeros for the fibre type $(I_4, I_2,I_2,I_2, I_2)$ (cf. [Her]) have to be:
\begin{center}
Table $5.2$
\end{center}     
\begin{center}
\begin{tabular}{ c||c | c||c}\hline
$p$   & $\nu_p(G_2)$  &  $\nu_p(G_3)$ &  $\nu_p (\Delta)$ \\ \hline
$p_1$ & $0$  &  $0$& $4$   \\
$p_2$ &  $0$ & $0$ & $2$ \\ 
$p_3$ & $ 0$ &$0$ & $2$ \\ 
$p_4$ &  $0$ & $0$ & $2$ \\ 
$p_5$ &  $0$ & $0$ & $2$ \\ \hline
{\rm sum} & $0$ & $0$ & $12$ \\ \hline 
\end{tabular} 
\end{center} 

Therefore 
\begin{equation}  
-{\bar \Delta}^2 = \Delta =G_2^3-G_3^2    
\end{equation}

with ${\bar \Delta}$, $G_3 \in H^0({\bf P}^1, {\cal O}(6L))$ where $deg L=p_g-q+1=1$ since $S$ is a rational surface. Thus $(1)$ is equivalent to
$$
G_2^3 = G_3^2 -  {\bar \Delta}^2 = (G_3-{\bar \Delta} )(G_3 + {\bar \Delta}).
$$
It follows from Table $5.2$ that $G_2$, $G_3$ and ${\bar \Delta}$ are relatively 
prime. We let
$$
H_1^3=G_3-{\bar \Delta} \ \ {\rm and} \ \ H_2^3=G_3+ {\bar \Delta}
$$
It is easy to see that  $H_1$, $H_2 \in H^0({\bf P}^1, {\cal O}(2L))$ and they are relatively prime. Thus
$$
-{\bar \Delta} = \frac{1}{2} (H_1^3-H_2^3)=\frac{1}{2}(H_1-H_2)(\eta^2 H_2-\eta H_1)(\eta H_2-\eta^2 H_1)
$$
where $\eta = e^{\frac{2\pi i}{3}}$. By $(1)$, ${\bar \Delta}$ has one double zero at $p_1$ and four single zeros at $p_i$ where $i=2$,$3$,$4$,$5$. Let
\begin{eqnarray*}
J_1 &:=&  H_1-H_2 \\
J_2 &:=&   \eta^2 H_2-\eta H_1 \\
J_3 &:=& \eta H_2-\eta^2 H_1 
\end{eqnarray*}
 with $J_i \in H^0({\bf P}_1, {\cal O}(2L))$, $i=1$, $2$, $3$. Similarly, $J_1$, $J_2$ and $J_3$ are relatively prime.

It can be assumed that double point of ${\bar \Delta}$ is that of  $J_2$ at $p_1=0$, two of the four single zeros of ${\bar \Delta}$ are that of  $J_3$ at $p_2=\infty$ and $p_3=1$. Therefore we may write $J_2=X^2$ and $J_3=Y(X-Y)$. It follows
\begin{eqnarray*}
J_1 &=& X^2+XY-Y^2 \\
H_1 &=&   \frac{1}{1-\eta}(X^2- \eta XY +\eta Y^2) \\
H_2 &=& \frac{1}{1-\eta} (\eta X^2-XY+Y^2) \\
{\bar \Delta} &=& -\frac{1}{2}(Y^2-XY-X^2)\cdot X^2\cdot Y(X-Y)
\end{eqnarray*}
and
\begin{eqnarray*}
G_2 &=& H_1\cdot H_2 \\
G_3 &=&  H_1^3+  {\bar \Delta}  \\
 \Delta &=& -\frac{1}{4}(Y^2-XY-X^2)^2\cdot X^4\cdot Y^2(X-Y)^2
\end{eqnarray*}
This proves the lemma.
\vskip 0.3cm
Thus we get the following lemma (cf. the proof of Lemma $3.3$):

\begin{lemma} Let  $V$ be a Picard number two relatively minimal  Gorenstein log del Pezzo surface with $K^2_V=1$. Then the singularity type of $V$ and the number $m$ of the isomorphism classes of the given singularity type are listed in Table $5.3$.
\vskip 0.5cm
\begin{center}
Table $5.3$
\end{center}
\begin{center}
\begin{tabular}{| c|c|c|c|}\hline
 Dynkin type of $V$ & $m$ & Dynkin type of $V$  & $m$\\ \hline
$D_7$ & $\infty $  &  $D_5+2A_1$ & $\infty $\\ \hline
$A_3+4A_1$ & $1$  &  $D_4+A_3$& $\infty $  \\ \hline
\end{tabular} 
\end{center} 
\end{lemma}
\vskip 0.3cm
Next we deal with $K_V^2=2$.
\vskip 0.3cm
By blowing down an NEC on $U_1$, we get a surface $U_2$ which is the minimal resolution of the Picard number two Gorenstein log del Pezzo surface $V$ with $K_V^2=2$.

\begin{theorem} Let  $V$ be the Picard number two  Gorenstein log del Pezzo surface with $K^2_V=2$. 

\hskip -0.6cm $(1)$ There exist $12$ singularity types which are listed in the following:
$E_6$, $D_4+2A_1$, $A_4+A_2$, $A_3+3A_1$, $D_6$, $A_6$, $2A_3$, $3A_2$, $D_5+A_1$, $A_5+A_1$, $A_3+A_2+A_1$, $6A_1$.

\hskip -0.6cm 
$(2)$ If $V$ is a relatively minimal surface, then the singularity type of $V$  and the number $m$ of the isomorphism classes with the given singularity type are listed in Table $5.4$.
\vskip 0.5cm
\begin{center}
Table $5.4$
\end{center}
\begin{center}
\begin{tabular}{| c|c|c|c|}\hline
 Dynkin type of $V$ & $m$ & Dynkin type of $V$  & $m$\\ \hline
$D_6$ & $\infty$  &  $D_4+2A_1$ & $1 $\\ \hline
$6A_1$ & $1$  &  $2A_3$& $\infty $  \\ \hline
\end{tabular} 
\end{center} 
\end{theorem}
{\em Proof.} \ The singularity types and the relatively minimal singularity types can be found by a direct calculation together with the definition of the relatively minimal surfaces (cf. Definition $2.1$). For the details, we refer to Table $3$ of  Appendix. 
\vskip 0.5cm
\hskip -0.6cm
Now we shall deal with the number of isomorphism classes with the given relatively minimal singularity type separately. Firstly, we prove the following lemma.

\begin{lemma} The number of the isomorphism classes of the Picard number two Gorenstein log del Pezzo surface with the singularity type $D_6$ is infinity.
\end{lemma}
{\em Proof.} \ We let $U$ be the minimal resolution of $V$. By Table $3$ of  Appendix, 
there are four possibilities for $U$: start from No.$1$ ($2$, $7$ or $8$) Mordell-Weil rank one rational elliptic surface.  Since each of these surfaces has at least one type $I_1$ fibre, we have
\begin{claim} There exists a reduced irreducible nodal member $A \in |-K_U |$.
\end{claim}
From the configuration of the negative curves in $U$ which is listed below, we know there are three $(-1)$-curves: $l_1$, $l_2$ and $l_3$.
\vskip 0.5cm
\centerline{\font\thinlinefont=cmr5
\begingroup\makeatletter\ifx\SetFigFont\undefined%
\gdef\SetFigFont#1#2#3#4#5{%
  \reset@font\fontsize{#1}{#2pt}%
  \fontfamily{#3}\fontseries{#4}\fontshape{#5}%
  \selectfont}%
\fi\endgroup%
\mbox{\beginpicture
\setcoordinatesystem units <0.93000cm,0.93000cm>
\unitlength=0.93000cm
\linethickness=1pt
\setplotsymbol ({\makebox(0,0)[l]{\tencirc\symbol{'160}}})
\setshadesymbol ({\thinlinefont .})
\setlinear
%
% Fig POLYLINE object
%
\linethickness= 0.500pt
\setplotsymbol ({\thinlinefont .})
\plot  3.175 21.590  5.874 22.701 /
%
% Fig POLYLINE object
%
\linethickness= 0.500pt
\setplotsymbol ({\thinlinefont .})
\plot  5.080 22.701  7.461 21.590 /
%
% Fig POLYLINE object
%
\linethickness= 0.500pt
\setplotsymbol ({\thinlinefont .})
\plot  6.668 21.590  9.049 22.543 /
%
% Fig POLYLINE object
%
\linethickness= 0.500pt
\setplotsymbol ({\thinlinefont .})
\plot  8.572 22.543 10.954 21.590 /
%
% Fig POLYLINE object
%
\linethickness= 0.500pt
\setplotsymbol ({\thinlinefont .})
\plot  9.207 21.907 10.954 23.336 /
%
% Fig POLYLINE object
%
\linethickness= 0.500pt
\setplotsymbol ({\thinlinefont .})
\plot 10.001 21.590 11.589 23.019 /
%
% Fig POLYLINE object
%
\linethickness= 0.500pt
\setplotsymbol ({\thinlinefont .})
\setdashes < 0.1270cm>
\plot  2.540 22.860  4.286 21.431 /
%
% Fig POLYLINE object
%
\linethickness= 0.500pt
\setplotsymbol ({\thinlinefont .})
\plot  3.493 23.495  5.239 21.907 /
%
% Fig POLYLINE object
%
\linethickness= 0.500pt
\setplotsymbol ({\thinlinefont .})
\plot  9.684 23.654 10.636 22.701 /
%
% Fig TEXT object
%
\put{\SetFigFont{11}{13.2}{rm}$U$} [lB] at  6.191 20.161
%
% Fig TEXT object
%
\put{\SetFigFont{11}{13.2}{rm}$l_1$} [lB] at  2.540 23.019
%
% Fig TEXT object
%
\put{\SetFigFont{11}{13.2}{rm}$l_2$} [lB] at  3.651 23.336
%
% Fig TEXT object
%
\put{\SetFigFont{11}{13.2}{rm}$l_3$} [lB] at  9.842 23.654
%
% Fig TEXT object
%
\put{\SetFigFont{11}{13.2}{rm}$G_5$} [lB] at  6.350 21.273
%
% Fig TEXT object
%
\put{\SetFigFont{11}{13.2}{rm}$G_4$} [lB] at  7.620 21.273
%
% Fig TEXT object
%
\put{\SetFigFont{11}{13.2}{rm}$G_1$} [lB] at  5.715 22.860
%
% Fig TEXT object
%
\put{\SetFigFont{11}{13.2}{rm}$G_2$} [lB] at 10.795 23.495
%
% Fig TEXT object
%
\put{\SetFigFont{11}{13.2}{rm}$G_3$} [lB] at 11.748 22.701
%
% Fig TEXT object
%
\put{\SetFigFont{11}{13.2}{rm}$G_6$} [lB] at 10.636 21.273
\linethickness=0pt
\putrectangle corners at  2.515 24.067 and 11.748 20.161
\endpicture}
}
\vskip 0.5cm
\begin{center}
Figure $5.1$
\end{center}
\begin{claim} There exists an automorphism on  $U$ which 
interchanges $l_1$ and $l_2$.
\end{claim}
{\em Proof of the claim.} \ \ Let $W \longrightarrow U$ be the blow up of 
$p=A \cap l_3$. It is easy to see that $W$ is a surface which is the minimal resolution of the Picard number two Gorenstein log del Pezzo surface with the singularity type $E_7$. Thus $W$ is either $U_1^7(E_7)$ or $U_1^8(E_7)$ and the proper transform of $l_i$ (i=1,2) intersects the $1$st component of $E_7$. From Table $2$ of  Appendix, we know that $l_i$ (i=1,2) comes from the Mordell-Weil section $P_1$ or $P_{-1}$ of the No.$7$ or $8$ Mordell-Weil rank one rational elliptic surface where $P_1$ is the generator of the rank one Mordell-Weil lattice $A_1^*$. In any case, we may define an automorphism $\Phi$ on No.7 (No.8) Mordell-Weil rank one rational elliptic surface which is the involution on a general fibre. Then
$\Phi (P_1)= P_{-1}$ and  $\Phi$ will induce an automorphism on $U$ which interchanges $l_1$ and $l_2$. This proves the claim.
\vskip 0.3cm
\hskip -0.6cm Let $U \longrightarrow {\bf P}^2$ be the composition of blowing down of $l_3+G_2+G_6+G_3$, $l_1+G_1+G_4$ to $p$, $q$ respectively where $G_i$ 
$(i=1, \cdots , 6)$ is the component of $D_6$ (cf. Figure $5.1$). 
\vskip 0.5cm
\centerline{\font\thinlinefont=cmr5
\begingroup\makeatletter\ifx\SetFigFont\undefined%
\gdef\SetFigFont#1#2#3#4#5{%
  \reset@font\fontsize{#1}{#2pt}%
  \fontfamily{#3}\fontseries{#4}\fontshape{#5}%
  \selectfont}%
\fi\endgroup%
\mbox{\beginpicture
\setcoordinatesystem units <0.73000cm,0.73000cm>
\unitlength=0.73000cm
\linethickness=1pt
\setplotsymbol ({\makebox(0,0)[l]{\tencirc\symbol{'160}}})
\setshadesymbol ({\thinlinefont .})
\setlinear
%
% Fig POLYLINE object
%
\linethickness= 0.500pt
\setplotsymbol ({\thinlinefont .})
\plot  1.429 23.812  2.540 22.384 /
%
% Fig POLYLINE object
%
\linethickness= 0.500pt
\setplotsymbol ({\thinlinefont .})
\plot  1.746 22.701  2.699 23.495 /
%
% Fig POLYLINE object
%
\linethickness= 0.500pt
\setplotsymbol ({\thinlinefont .})
\plot  2.064 22.384  3.016 23.178 /
%
% Fig POLYLINE object
%
\linethickness= 0.500pt
\setplotsymbol ({\thinlinefont .})
\plot  1.270 23.178  2.381 24.130 /
%
% Fig POLYLINE object
%
\linethickness= 0.500pt
\setplotsymbol ({\thinlinefont .})
\plot  2.540 23.019  4.763 22.860 /
%
% Fig POLYLINE object
%
\linethickness= 0.500pt
\setplotsymbol ({\thinlinefont .})
\plot  4.604 24.130  4.445 21.749 /
%
% Fig POLYLINE object
%
\linethickness= 0.500pt
\setplotsymbol ({\thinlinefont .})
\setdashes < 0.1270cm>
\plot  2.223 24.765  2.223 23.654 /
%
% Fig POLYLINE object
%
\linethickness= 0.500pt
\setplotsymbol ({\thinlinefont .})
\plot  3.969 25.082  4.763 23.495 /
%
% Fig POLYLINE object
%
\linethickness= 0.500pt
\setplotsymbol ({\thinlinefont .})
\setsolid
\putrule from  6.668 22.860 to  7.303 22.860
%
% arrow head
%
\plot  7.048 22.796  7.303 22.860  7.048 22.924 /
%
%
% Fig POLYLINE object
%
\linethickness= 0.500pt
\setplotsymbol ({\thinlinefont .})
\setdashes < 0.1270cm>
\plot  4.763 22.384  3.016 20.796 /
%
% Fig POLYLINE object
%
\linethickness= 0.500pt
\setplotsymbol ({\thinlinefont .})
\setsolid
\plot 12.224 25.082 11.906 20.637 /
%
% Fig POLYLINE object
%
\linethickness= 0.500pt
\setplotsymbol ({\thinlinefont .})
\plot  9.684 23.495 13.494 22.384 /
%
% Fig POLYLINE object
%
\linethickness= 0.500pt
\setplotsymbol ({\thinlinefont .})
\putrule from  1.905 24.289 to  1.909 24.289
\plot  1.909 24.289  1.916 24.291 /
\plot  1.916 24.291  1.928 24.293 /
\plot  1.928 24.293  1.952 24.297 /
\plot  1.952 24.297  1.981 24.304 /
\plot  1.981 24.304  2.024 24.310 /
\plot  2.024 24.310  2.076 24.320 /
\plot  2.076 24.320  2.140 24.331 /
\plot  2.140 24.331  2.216 24.344 /
\plot  2.216 24.344  2.303 24.359 /
\plot  2.303 24.359  2.402 24.378 /
\plot  2.402 24.378  2.510 24.395 /
\plot  2.510 24.395  2.629 24.416 /
\plot  2.629 24.416  2.756 24.435 /
\plot  2.756 24.435  2.889 24.456 /
\plot  2.889 24.456  3.029 24.479 /
\plot  3.029 24.479  3.173 24.500 /
\plot  3.173 24.500  3.317 24.522 /
\plot  3.317 24.522  3.465 24.543 /
\plot  3.465 24.543  3.611 24.562 /
\plot  3.611 24.562  3.755 24.579 /
\plot  3.755 24.579  3.895 24.596 /
\plot  3.895 24.596  4.032 24.608 /
\plot  4.032 24.608  4.163 24.621 /
\plot  4.163 24.621  4.290 24.630 /
\plot  4.290 24.630  4.411 24.636 /
\plot  4.411 24.636  4.523 24.638 /
\putrule from  4.523 24.638 to  4.629 24.638
\plot  4.629 24.638  4.727 24.634 /
\plot  4.727 24.634  4.815 24.625 /
\plot  4.815 24.625  4.898 24.613 /
\plot  4.898 24.613  4.972 24.598 /
\plot  4.972 24.598  5.040 24.577 /
\plot  5.040 24.577  5.099 24.551 /
\plot  5.099 24.551  5.152 24.522 /
\plot  5.152 24.522  5.199 24.488 /
\plot  5.199 24.488  5.239 24.448 /
\plot  5.239 24.448  5.275 24.403 /
\plot  5.275 24.403  5.304 24.354 /
\plot  5.304 24.354  5.330 24.299 /
\plot  5.330 24.299  5.349 24.238 /
\plot  5.349 24.238  5.366 24.172 /
\plot  5.366 24.172  5.378 24.100 /
\plot  5.378 24.100  5.387 24.024 /
\plot  5.387 24.024  5.393 23.942 /
\plot  5.393 23.942  5.395 23.855 /
\putrule from  5.395 23.855 to  5.395 23.762
\plot  5.395 23.762  5.391 23.664 /
\plot  5.391 23.664  5.385 23.563 /
\plot  5.385 23.563  5.376 23.459 /
\plot  5.376 23.459  5.366 23.351 /
\plot  5.366 23.351  5.355 23.241 /
\plot  5.355 23.241  5.340 23.127 /
\plot  5.340 23.127  5.326 23.012 /
\plot  5.326 23.012  5.309 22.896 /
\plot  5.309 22.896  5.292 22.782 /
\plot  5.292 22.782  5.273 22.665 /
\plot  5.273 22.665  5.251 22.549 /
\plot  5.251 22.549  5.230 22.435 /
\plot  5.230 22.435  5.209 22.320 /
\plot  5.209 22.320  5.188 22.210 /
\plot  5.188 22.210  5.165 22.102 /
\plot  5.165 22.102  5.141 21.999 /
\plot  5.141 21.999  5.116 21.897 /
\plot  5.116 21.897  5.091 21.800 /
\plot  5.091 21.800  5.065 21.706 /
\plot  5.065 21.706  5.038 21.620 /
\plot  5.038 21.620  5.008 21.537 /
\plot  5.008 21.537  4.978 21.461 /
\plot  4.978 21.461  4.949 21.389 /
\plot  4.949 21.389  4.915 21.323 /
\plot  4.915 21.323  4.881 21.262 /
\plot  4.881 21.262  4.845 21.207 /
\plot  4.845 21.207  4.805 21.158 /
\plot  4.805 21.158  4.763 21.114 /
\plot  4.763 21.114  4.714 21.071 /
\plot  4.714 21.071  4.663 21.035 /
\plot  4.663 21.035  4.606 21.006 /
\plot  4.606 21.006  4.544 20.980 /
\plot  4.544 20.980  4.479 20.959 /
\plot  4.479 20.959  4.409 20.942 /
\plot  4.409 20.942  4.335 20.932 /
\plot  4.335 20.932  4.257 20.925 /
\plot  4.257 20.925  4.174 20.921 /
\plot  4.174 20.921  4.087 20.923 /
\plot  4.087 20.923  3.998 20.925 /
\plot  3.998 20.925  3.905 20.934 /
\plot  3.905 20.934  3.808 20.942 /
\plot  3.808 20.942  3.711 20.955 /
\plot  3.711 20.955  3.611 20.970 /
\plot  3.611 20.970  3.509 20.987 /
\plot  3.509 20.987  3.406 21.004 /
\plot  3.406 21.004  3.304 21.023 /
\plot  3.304 21.023  3.200 21.042 /
\plot  3.200 21.042  3.099 21.061 /
\plot  3.099 21.061  2.999 21.080 /
\plot  2.999 21.080  2.900 21.099 /
\plot  2.900 21.099  2.802 21.118 /
\plot  2.802 21.118  2.707 21.137 /
\plot  2.707 21.137  2.616 21.154 /
\plot  2.616 21.154  2.529 21.171 /
\plot  2.529 21.171  2.445 21.186 /
\plot  2.445 21.186  2.364 21.201 /
\plot  2.364 21.201  2.290 21.213 /
\plot  2.290 21.213  2.220 21.224 /
\plot  2.220 21.224  2.155 21.234 /
\plot  2.155 21.234  2.095 21.245 /
\plot  2.095 21.245  2.040 21.251 /
\plot  2.040 21.251  1.990 21.260 /
\plot  1.990 21.260  1.945 21.266 /
\plot  1.945 21.266  1.905 21.273 /
\plot  1.905 21.273  1.863 21.281 /
\plot  1.863 21.281  1.829 21.289 /
\plot  1.829 21.289  1.799 21.300 /
\plot  1.799 21.300  1.778 21.311 /
\plot  1.778 21.311  1.761 21.321 /
\plot  1.761 21.321  1.750 21.334 /
\plot  1.750 21.334  1.744 21.349 /
\plot  1.744 21.349  1.742 21.364 /
\plot  1.742 21.364  1.746 21.380 /
\plot  1.746 21.380  1.753 21.397 /
\plot  1.753 21.397  1.763 21.414 /
\plot  1.763 21.414  1.776 21.433 /
\plot  1.776 21.433  1.793 21.452 /
\plot  1.793 21.452  1.810 21.471 /
\plot  1.810 21.471  1.829 21.491 /
\plot  1.829 21.491  1.848 21.510 /
\plot  1.848 21.510  1.869 21.527 /
\plot  1.869 21.527  1.888 21.546 /
\plot  1.888 21.546  1.909 21.560 /
\plot  1.909 21.560  1.928 21.575 /
\plot  1.928 21.575  1.947 21.588 /
\plot  1.947 21.588  1.964 21.598 /
\plot  1.964 21.598  1.981 21.607 /
\plot  1.981 21.607  1.996 21.613 /
\plot  1.996 21.613  2.009 21.618 /
\putrule from  2.009 21.618 to  2.021 21.618
\plot  2.021 21.618  2.032 21.615 /
\plot  2.032 21.615  2.043 21.609 /
\plot  2.043 21.609  2.053 21.601 /
\plot  2.053 21.601  2.064 21.590 /
\plot  2.064 21.590  2.079 21.569 /
\plot  2.079 21.569  2.095 21.541 /
\plot  2.095 21.541  2.115 21.505 /
\plot  2.115 21.505  2.136 21.461 /
\plot  2.136 21.461  2.161 21.408 /
\plot  2.161 21.408  2.187 21.349 /
\plot  2.187 21.349  2.216 21.281 /
\plot  2.216 21.281  2.248 21.207 /
\plot  2.248 21.207  2.282 21.126 /
\plot  2.282 21.126  2.316 21.044 /
\plot  2.316 21.044  2.349 20.957 /
\plot  2.349 20.957  2.383 20.872 /
\plot  2.383 20.872  2.415 20.792 /
\plot  2.415 20.792  2.445 20.718 /
\plot  2.445 20.718  2.472 20.650 /
\plot  2.472 20.650  2.493 20.595 /
\plot  2.493 20.595  2.512 20.551 /
\plot  2.512 20.551  2.525 20.517 /
\plot  2.525 20.517  2.534 20.496 /
\plot  2.534 20.496  2.538 20.485 /
\plot  2.538 20.485  2.540 20.479 /
%
% Fig POLYLINE object
%
\linethickness= 0.500pt
\setplotsymbol ({\thinlinefont .})
\putrule from  9.842 24.765 to  9.842 24.761
\plot  9.842 24.761  9.845 24.750 /
\plot  9.845 24.750  9.849 24.733 /
\plot  9.849 24.733  9.855 24.706 /
\plot  9.855 24.706  9.864 24.668 /
\plot  9.864 24.668  9.874 24.619 /
\plot  9.874 24.619  9.887 24.560 /
\plot  9.887 24.560  9.904 24.492 /
\plot  9.904 24.492  9.921 24.414 /
\plot  9.921 24.414  9.942 24.329 /
\plot  9.942 24.329  9.963 24.238 /
\plot  9.963 24.238  9.986 24.145 /
\plot  9.986 24.145 10.010 24.052 /
\plot 10.010 24.052 10.033 23.961 /
\plot 10.033 23.961 10.058 23.872 /
\plot 10.058 23.872 10.084 23.787 /
\plot 10.084 23.787 10.107 23.709 /
\plot 10.107 23.709 10.132 23.635 /
\plot 10.132 23.635 10.158 23.569 /
\plot 10.158 23.569 10.183 23.512 /
\plot 10.183 23.512 10.209 23.463 /
\plot 10.209 23.463 10.234 23.421 /
\plot 10.234 23.421 10.262 23.385 /
\plot 10.262 23.385 10.289 23.357 /
\plot 10.289 23.357 10.319 23.336 /
\plot 10.319 23.336 10.357 23.319 /
\plot 10.357 23.319 10.397 23.309 /
\putrule from 10.397 23.309 to 10.439 23.309
\plot 10.439 23.309 10.486 23.317 /
\plot 10.486 23.317 10.535 23.334 /
\plot 10.535 23.334 10.585 23.357 /
\plot 10.585 23.357 10.638 23.387 /
\plot 10.638 23.387 10.693 23.423 /
\plot 10.693 23.423 10.751 23.463 /
\plot 10.751 23.463 10.808 23.508 /
\plot 10.808 23.508 10.867 23.552 /
\plot 10.867 23.552 10.926 23.597 /
\plot 10.926 23.597 10.983 23.641 /
\plot 10.983 23.641 11.041 23.683 /
\plot 11.041 23.683 11.096 23.721 /
\plot 11.096 23.721 11.151 23.753 /
\plot 11.151 23.753 11.204 23.781 /
\plot 11.204 23.781 11.252 23.802 /
\plot 11.252 23.802 11.301 23.815 /
\plot 11.301 23.815 11.345 23.821 /
\putrule from 11.345 23.821 to 11.390 23.821
\plot 11.390 23.821 11.430 23.812 /
\plot 11.430 23.812 11.468 23.798 /
\plot 11.468 23.798 11.506 23.779 /
\plot 11.506 23.779 11.544 23.753 /
\plot 11.544 23.753 11.580 23.721 /
\plot 11.580 23.721 11.618 23.686 /
\plot 11.618 23.686 11.656 23.643 /
\plot 11.656 23.643 11.697 23.599 /
\plot 11.697 23.599 11.735 23.550 /
\plot 11.735 23.550 11.775 23.497 /
\plot 11.775 23.497 11.813 23.444 /
\plot 11.813 23.444 11.851 23.389 /
\plot 11.851 23.389 11.889 23.332 /
\plot 11.889 23.332 11.925 23.277 /
\plot 11.925 23.277 11.957 23.222 /
\plot 11.957 23.222 11.989 23.169 /
\plot 11.989 23.169 12.014 23.116 /
\plot 12.014 23.116 12.037 23.067 /
\plot 12.037 23.067 12.054 23.021 /
\plot 12.054 23.021 12.067 22.976 /
\plot 12.067 22.976 12.073 22.936 /
\plot 12.073 22.936 12.071 22.896 /
\plot 12.071 22.896 12.065 22.860 /
\plot 12.065 22.860 12.052 22.828 /
\plot 12.052 22.828 12.031 22.797 /
\plot 12.031 22.797 12.002 22.767 /
\plot 12.002 22.767 11.966 22.739 /
\plot 11.966 22.739 11.919 22.712 /
\plot 11.919 22.712 11.864 22.686 /
\plot 11.864 22.686 11.803 22.663 /
\plot 11.803 22.663 11.735 22.640 /
\plot 11.735 22.640 11.661 22.619 /
\plot 11.661 22.619 11.582 22.598 /
\plot 11.582 22.598 11.502 22.578 /
\plot 11.502 22.578 11.419 22.557 /
\plot 11.419 22.557 11.339 22.536 /
\plot 11.339 22.536 11.261 22.517 /
\plot 11.261 22.517 11.187 22.494 /
\plot 11.187 22.494 11.119 22.473 /
\plot 11.119 22.473 11.060 22.447 /
\plot 11.060 22.447 11.009 22.422 /
\plot 11.009 22.422 10.969 22.394 /
\plot 10.969 22.394 10.939 22.365 /
\plot 10.939 22.365 10.924 22.333 /
\plot 10.924 22.333 10.920 22.299 /
\plot 10.920 22.299 10.930 22.263 /
\plot 10.930 22.263 10.954 22.225 /
\plot 10.954 22.225 10.977 22.197 /
\plot 10.977 22.197 11.005 22.170 /
\plot 11.005 22.170 11.041 22.140 /
\plot 11.041 22.140 11.081 22.109 /
\plot 11.081 22.109 11.127 22.077 /
\plot 11.127 22.077 11.180 22.043 /
\plot 11.180 22.043 11.239 22.007 /
\plot 11.239 22.007 11.305 21.969 /
\plot 11.305 21.969 11.375 21.931 /
\plot 11.375 21.931 11.449 21.888 /
\plot 11.449 21.888 11.529 21.848 /
\plot 11.529 21.848 11.612 21.804 /
\plot 11.612 21.804 11.699 21.759 /
\plot 11.699 21.759 11.790 21.715 /
\plot 11.790 21.715 11.883 21.670 /
\plot 11.883 21.670 11.978 21.624 /
\plot 11.978 21.624 12.073 21.577 /
\plot 12.073 21.577 12.171 21.531 /
\plot 12.171 21.531 12.268 21.484 /
\plot 12.268 21.484 12.366 21.440 /
\plot 12.366 21.440 12.463 21.393 /
\plot 12.463 21.393 12.556 21.351 /
\plot 12.556 21.351 12.649 21.306 /
\plot 12.649 21.306 12.740 21.266 /
\plot 12.740 21.266 12.827 21.226 /
\plot 12.827 21.226 12.910 21.188 /
\plot 12.910 21.188 12.988 21.154 /
\plot 12.988 21.154 13.064 21.120 /
\plot 13.064 21.120 13.134 21.090 /
\plot 13.134 21.090 13.200 21.061 /
\plot 13.200 21.061 13.261 21.038 /
\plot 13.261 21.038 13.318 21.014 /
\plot 13.318 21.014 13.369 20.995 /
\plot 13.369 20.995 13.415 20.978 /
\plot 13.415 20.978 13.456 20.966 /
\plot 13.456 20.966 13.494 20.955 /
\plot 13.494 20.955 13.542 20.944 /
\plot 13.542 20.944 13.583 20.942 /
\plot 13.583 20.942 13.614 20.947 /
\plot 13.614 20.947 13.640 20.959 /
\plot 13.640 20.959 13.659 20.978 /
\plot 13.659 20.978 13.672 21.004 /
\plot 13.672 21.004 13.678 21.033 /
\plot 13.678 21.033 13.680 21.071 /
\putrule from 13.680 21.071 to 13.680 21.112
\plot 13.680 21.112 13.674 21.156 /
\plot 13.674 21.156 13.667 21.203 /
\plot 13.667 21.203 13.657 21.249 /
\plot 13.657 21.249 13.646 21.298 /
\plot 13.646 21.298 13.633 21.344 /
\plot 13.633 21.344 13.621 21.391 /
\plot 13.621 21.391 13.608 21.433 /
\plot 13.608 21.433 13.595 21.471 /
\plot 13.595 21.471 13.583 21.505 /
\plot 13.583 21.505 13.570 21.535 /
\plot 13.570 21.535 13.557 21.558 /
\plot 13.557 21.558 13.542 21.575 /
\plot 13.542 21.575 13.528 21.586 /
\plot 13.528 21.586 13.513 21.590 /
\putrule from 13.513 21.590 to 13.494 21.590
\plot 13.494 21.590 13.473 21.584 /
\plot 13.473 21.584 13.447 21.573 /
\plot 13.447 21.573 13.415 21.556 /
\plot 13.415 21.556 13.382 21.533 /
\plot 13.382 21.533 13.341 21.503 /
\plot 13.341 21.503 13.295 21.467 /
\plot 13.295 21.467 13.244 21.427 /
\plot 13.244 21.427 13.189 21.380 /
\plot 13.189 21.380 13.130 21.330 /
\plot 13.130 21.330 13.068 21.275 /
\plot 13.068 21.275 13.003 21.217 /
\plot 13.003 21.217 12.939 21.158 /
\plot 12.939 21.158 12.874 21.101 /
\plot 12.874 21.101 12.812 21.044 /
\plot 12.812 21.044 12.753 20.991 /
\plot 12.753 20.991 12.700 20.942 /
\plot 12.700 20.942 12.656 20.902 /
\plot 12.656 20.902 12.617 20.866 /
\plot 12.617 20.866 12.588 20.839 /
\plot 12.588 20.839 12.567 20.820 /
\plot 12.567 20.820 12.552 20.807 /
\plot 12.552 20.807 12.545 20.800 /
\plot 12.545 20.800 12.541 20.796 /
%
% Fig TEXT object
%
\put{\SetFigFont{9}{10.8}{rm}$A$} [lB] at  4.445 25.082
%
% Fig TEXT object
%
\put{\SetFigFont{9}{10.8}{rm}$K^2=2$} [lB] at  3.016 20.003
%
% Fig TEXT object
%
\put{\SetFigFont{9}{10.8}{rm}$A$} [lB] at 10.795 23.971
%
% Fig TEXT object
%
\put{\SetFigFont{9}{10.8}{rm}$p$} [lB] at 10.160 22.701
%
% Fig TEXT object
%
\put{\SetFigFont{9}{10.8}{rm}${\bf P}^2$} [lB] at 10.478 19.685
%
% Fig TEXT object
%
\put{\SetFigFont{9}{10.8}{rm}$q$} [lB] at 12.383 22.860
\linethickness=0pt
\putrectangle corners at  1.245 25.432 and 13.705 19.685
\endpicture}
}
\vskip 0.5cm
\begin{center}
Figure $5.2$
\end{center}
From the above construction and Claim $5.10$, we know that $U$ is uniquely determined by a pair $M=(S,C,p)$ where $S={\bf P}^2$, $C$ is a nodal cubic  curve on $S$ and $p$ is a point of $C$ which is neither a singular point of $C$ nor an inflex of $C$. An automorphism of $U$ will induce an isomorphism on $M=(S,C,p)$ by making use of Claim $5.10$ if necessary. In order to prove the lemma, we have only to show that there exist infinitly many isomorphism classes of such pairs. 
\vskip 0.3cm
\hskip -0.6cm
Before stating the next result, we employ the following notation.

$S$ = ${\bf P}^2$ with coordinates $X$,$Y$ and $Z$.

 $C$ = $\{[X:Y:Z] | Y^2 Z=X^2(X+Z) \}$.

$S'$ = ${\bf P}^2$ with coordinates $X'$,$Y'$ and $Z'$. 

$C'$ = $\{[X':Y':Z'] | Y'^2 Z'=X'^2(X'+Z') \}$.
  
\begin{claim} If $\Phi : S \longrightarrow S'$ is an isomorphism such that $\Phi (C)=C'$, then either $\Phi ([X:Y:Z]) =[X':Y':Z']$ or $\Phi ([X:Y:Z])=[X':-Y':Z']$.
\end{claim}
{\em Proof of the claim.}  Clearly, $\Phi$ is linear map. We set: 
$$
\left(\begin{array}{r}
   X' \\
   Y' \\
   Z'
   \end{array}
   \right)  =
\Phi \left(\begin{array}{r}
   X \\
   Y \\
   Z 
   \end{array}
   \right) =
 \left( \begin{array}{rrr}
   a_{11} & a_{12} & a_{13} \\
   a_{21} & a_{22} &a_{23} \\
   a_{31} & a_{32} & a_{33}
   \end{array}
   \right) \cdot
 \left( \begin{array}{r}
   X \\
   Y  \\
   Z
    \end{array}
   \right) = A \cdot
 \left( \begin{array}{r}
   X \\
   Y \\
   Z
   \end{array}
   \right)
$$
Then $det A \not=0$. Since $\Phi (C)=C'$, we have
$$
\Phi ([0:0:1])=[0:0:1]   \ \ {\rm and} \ \ \Phi ([1:0:-1])=[1:0:-1].
$$
From the above equalities, we may set
$$
a_{33}=1, a_{11}=1- a_{31}, a_{13}=a_{23}=a_{21}=0.
$$
For simplicity, we denote $a_{12}=a$, $a_{22}=b$, $a_{32}=c$ and $a_{11}=d$. Thus we have
$$
X'=dX+aY, \ \ Y'=bY, \ \ Z'=(1-d)X+cY+Z.
$$
Substituting above equalities to $Y'^2Z'=X'^2(X'+Z')$, we get
\begin{equation}  
(bY)^2\cdot ((1-d)X+cY+Z)=(dX+aY)^2\cdot (X+aY+cY+Z).  
\end{equation}
(2) is hold for all $[X:Y:Z] \in C$. On the other hand, if $[X:Y:Z]\in C$, then
$[X:-Y:Z] \in C$, thus we also have
\begin{equation}  
(bY)^2\cdot ((1-d)X-cY+Z)=(dX-aY)^2\cdot (X-aY-cY+Z). 
\end{equation}
From $(2)+(3)$, we may get
\begin{eqnarray*}
b^2Y^2[(1-d)X+Z] &=&(X+Z)(d^2X^2+a^2Y^2)+2ad(a+c)Y^2X \\
&=& d^2Y^2Z+a^2Y^2(X+Z)+2ad(a+c)XY^2.
\end{eqnarray*}
So for all $[X:Y:Z]\in C-\{ [0:0:1], [1:0:-1] \}$, we must have
$$
b^2[(1-d)X+Z]=d^2Z+a^2(X+Z)+2ad(a+c)X
$$
Since we can take $Z=1$ and arbitrarily many $X$,  we get
\begin{equation}  
b^2=d^2+a^2  
\end{equation}

\begin{equation}  
b^2(1-d)=a^2+2ad(a+c)
\end{equation}
By $(2)-(3)$, we get
$$
cb^2Y^3=2adXY(X+Z)+(a+c)d^2YX^2+a^2(a+c)Y^3.
$$
Since above equality holds for all $[X:Y:Z] \in C-\{ [0:0:1], [1:0:-1] \} $, we get
$$
cb^2Y^2=2adX(X+Z)+(a+c)d^2X^2+a^2(a+c)Y^2.  
$$
We let $g(X,Y,Z)=cb^2Y^2-[2adX(X+Z)+(a+c)d^2X^2+a^2(a+c)Y^2]$. Then we know $C$ is included in the closure of the zero set of $g$. On the other hand, $C$ is an irreducible cubic. Thus we know
$g = 0$ and so
\begin{equation}  
2ad=0,
\end{equation}
\begin{equation}  
 b^2c=a^2(a+c).
\end{equation}
By $(6)$, we have $a=0$ or $d=0$. But $det A \not=0$, we see that $a$ and $d$ can't both equal to $0$. If $d=0$, then by $(4)$ and $(7)$, we get $a=0$, a contradiction. Thus $a=0$. Then by $(5)$ and $det A \not =0$ which implies $b\not=0$, we have $d=1$ and so by $(4)$, $b=\pm 1$.
This proves the claim. 
\vskip 0.3cm
\hskip -0.6cm
Let $M=(S, C, p_1 )$ and $M_1=(S, C, p_2)$ where $S$ and $C$ are defined above and $p_1=[X_1:Y_1:Z_1]$ and $p_2=[X_2:Y_2:Z_2]$. From the above claim, we see that $M_1$ and
$M_2$ are not isomorphic if $[X_1:Y_1:Z_1]\not=[X_2:\pm Y_2:Z_2]$. Moving $p \in C $, we may get infinitely many non-isomorphism classes of the pair $(S,C,p)$ and so the infinite many non-ismorphism classes of $U$. This proves the lemma.
\vskip 0.5cm
\hskip -0.6cm
Next we consider the singularity type $6A_1$.
\vskip 0.5cm
\hskip -0.6cm
From Table $3$ of the Appedenix, there are at most two possibilities for  
$U(D_6)$: Start with either $U_1^3(A_3+4A_1)$ or $U_1^{17}(D_4+3A_1)$ and blow down an NEC.  
Since the configuration of $U_2^3(6A_1)$ and $U_2^{17}(6A_1)$ are the same, we only discuss one of these, say, $U_2^3(6A_1)$. By Tables $1$ and $2$ in the Appedenix, we find that blowing up a point on $Q_1$ or $Q_2$ will give rise to the surface $U_1$ which is the minimal resolution of a Picard number two log del Pezzo surface with singularity type $A_3+4A_1$. By Lemma 5.6, we know $U_1=U_1^3(A_3+4A_1)$. On the other hand, there is a unique way to blow down an NEC on $U_1^3(A_3+4A_1)$ to $U_2^{3}(6A_1)$, modular an automorphism switching the two NECs (as in Claim $5.10$). Thus the uniqueness of the surface with the singularity type $6A_1$ is due to the uniqueness of the Picard number two log del Pezzo surface with the singularity type $A_3+4A_1$. 
\vskip 0.5cm
\hskip -0.6cm
Now we deal with the singularity type $D_4+2A_1$.
\vskip 0.5cm
\hskip -0.6cm
 We let $V(D_4+2A_1)$ be the  Picard number two Gorenstein log del Pezzo surface with the singularity type $D_4+2A_1$, $U(D_4+2A_1)$ the minimal resolution of $V(D_4+2A_1)$ and $D$ the exceptional divisor. The configuration of $D$ is in the left part of the following figure. Blowing up $p$, $q$ of $D$ will give rise to the No.$9$ extremal rational elliptic surface $Y$ with the singular fibre type $(I_1^*, I_4, I_1)$ which is unique up to isomorphisms (cf. [MP, Theorem 5.4]). On the other hand, a direct calculation shows that the Mordell-Weil group of $Y$ is ${\bf Z}/4{\bf Z}$ with the generator $Q$ and $Y \longrightarrow U(D_4+2A_1)$ is the blow-down of ${\cal O}$ and $2Q$ on $Y$.
\vskip 0.5cm 
 \centerline{\font\thinlinefont=cmr5
\begingroup\makeatletter\ifx\SetFigFont\undefined%
\gdef\SetFigFont#1#2#3#4#5{%
  \reset@font\fontsize{#1}{#2pt}%
  \fontfamily{#3}\fontseries{#4}\fontshape{#5}%
  \selectfont}%
\fi\endgroup%
\mbox{\beginpicture
\setcoordinatesystem units <0.58000cm,0.58000cm>
\unitlength=0.58000cm
\linethickness=1pt
\setplotsymbol ({\makebox(0,0)[l]{\tencirc\symbol{'160}}})
\setshadesymbol ({\thinlinefont .})
\setlinear
%
% Fig POLYLINE object
%
\linethickness= 0.500pt
\setplotsymbol ({\thinlinefont .})
\plot  4.128 21.273  8.414 19.367 /
%
% Fig POLYLINE object
%
\linethickness= 0.500pt
\setplotsymbol ({\thinlinefont .})
\plot  7.144 22.225  6.032 19.685 /
%
% Fig POLYLINE object
%
\linethickness= 0.500pt
\setplotsymbol ({\thinlinefont .})
\plot  7.938 21.431  6.985 19.367 /
%
% Fig POLYLINE object
%
\linethickness= 0.500pt
\setplotsymbol ({\thinlinefont .})
\plot 10.319 22.543 12.859 20.796 /
%
% Fig POLYLINE object
%
\linethickness= 0.500pt
\setplotsymbol ({\thinlinefont .})
\plot  9.684 21.590 12.224 19.844 /
%
% Fig POLYLINE object
%
\linethickness= 0.500pt
\setplotsymbol ({\thinlinefont .})
\setdashes < 0.1270cm>
\plot 12.065 23.654 10.001 20.479 /
%
% Fig POLYLINE object
%
\linethickness= 0.500pt
\setplotsymbol ({\thinlinefont .})
\plot 12.700 21.749 10.160 17.939 /
%
% Fig POLYLINE object
%
\linethickness= 0.500pt
\setplotsymbol ({\thinlinefont .})
\plot  6.509 21.907 10.954 22.543 /
%
% Fig POLYLINE object
%
\linethickness= 0.500pt
\setplotsymbol ({\thinlinefont .})
\plot  7.303 20.637 10.160 21.590 /
%
% Fig POLYLINE object
%
\linethickness= 0.500pt
\setplotsymbol ({\thinlinefont .})
\setsolid
\plot 15.716 19.209 19.050 21.431 /
%
% Fig POLYLINE object
%
\linethickness= 0.500pt
\setplotsymbol ({\thinlinefont .})
\plot 17.780 21.273 20.955 18.891 /
%
% Fig POLYLINE object
%
\linethickness= 0.500pt
\setplotsymbol ({\thinlinefont .})
\plot 20.320 21.907 18.733 19.685 /
%
% Fig POLYLINE object
%
\linethickness= 0.500pt
\setplotsymbol ({\thinlinefont .})
\plot 20.955 21.114 19.685 19.050 /
%
% Fig POLYLINE object
%
\linethickness= 0.500pt
\setplotsymbol ({\thinlinefont .})
\plot 16.034 22.543 16.510 19.209 /
%
% Fig POLYLINE object
%
\linethickness= 0.500pt
\setplotsymbol ({\thinlinefont .})
\plot 16.986 21.431 17.621 17.939 /
%
% Fig POLYLINE object
%
\linethickness= 0.500pt
\setplotsymbol ({\thinlinefont .})
\plot 22.543 22.066 25.718 20.637 /
%
% Fig POLYLINE object
%
\linethickness= 0.500pt
\setplotsymbol ({\thinlinefont .})
\plot 22.225 20.479 25.082 19.209 /
%
% Fig POLYLINE object
%
\linethickness= 0.500pt
\setplotsymbol ({\thinlinefont .})
\plot 23.336 23.178 23.019 19.367 /
%
% Fig POLYLINE object
%
\linethickness= 0.500pt
\setplotsymbol ({\thinlinefont .})
\plot 24.924 22.225 23.971 17.780 /
%
% Fig POLYLINE object
%
\linethickness= 0.500pt
\setplotsymbol ({\thinlinefont .})
\setdashes < 0.1270cm>
\plot 17.145 18.415 24.606 18.256 /
%
% Fig POLYLINE object
%
\linethickness= 0.500pt
\setplotsymbol ({\thinlinefont .})
\plot 15.558 22.066 24.924 22.701 /
%
% Fig POLYLINE object
%
\linethickness= 0.500pt
\setplotsymbol ({\thinlinefont .})
\setsolid
\plot  5.239 23.178  4.921 18.256 /
%
% Fig POLYLINE object
%
\linethickness= 0.500pt
\setplotsymbol ({\thinlinefont .})
\setdashes < 0.1270cm>
\plot  4.604 22.860 12.383 23.336 /
%
% Fig POLYLINE object
%
\linethickness= 0.500pt
\setplotsymbol ({\thinlinefont .})
\plot  4.445 18.574 10.954 18.256 /
%
% Fig POLYLINE object
%
\linethickness= 0.500pt
\setplotsymbol ({\thinlinefont .})
\plot 23.178 22.066 23.178 22.066 /
%
% Fig POLYLINE object
%
\linethickness= 0.500pt
\setplotsymbol ({\thinlinefont .})
\plot 19.685 21.273 23.019 22.066 /
%
% Fig POLYLINE object
%
\linethickness= 0.500pt
\setplotsymbol ({\thinlinefont .})
\plot 20.003 20.161 22.860 20.320 /
%
% Fig POLYLINE object
%
\linethickness= 0.500pt
\setplotsymbol ({\thinlinefont .})
\setsolid
%
% arrow head
%
\plot 13.906 20.384 13.652 20.320 13.906 20.256 /
\putrule from 13.652 20.320 to 15.240 20.320
%
% Fig TEXT object
%
\put{\SetFigFont{7}{8.4}{rm}$p$} [lB] at 11.430 23.654
%
% Fig TEXT object
%
\put{\SetFigFont{7}{8.4}{rm}$q$} [lB] at  9.842 18.574
%
% Fig TEXT object
%
\put{\SetFigFont{7}{8.4}{rm}$(I_1^*, I_4, I_1)$} [lB] at 18.733 23.812
%
% Fig TEXT object
%
\put{\SetFigFont{7}{8.4}{rm}$U(D_4+2A_1)$} [lB] at  7.144 23.654
\linethickness=0pt
\putrectangle corners at  4.102 24.225 and 25.743 17.755
\endpicture}
}
\vskip 0.5cm
\begin{center}
Figure $5.3$
\end{center}
\vskip 0.5cm
\hskip -0.6cm 
Finally, we look at the singularity type $2A_3$. 
\vskip 0.5cm
\hskip -0.6cm
 We let $V(2A_3)$ be the  Picard number two Gorenstein log del Pezzo surface with the singularity type $2A_3$, $U(2A_3)$ the minimal resolution of $V(2A_3)$ and $D$ the exceptional divisor. The configuration of $D$ is in the left part of the following figure.
\vskip 0.5cm
  \centerline{\font\thinlinefont=cmr5
\begingroup\makeatletter\ifx\SetFigFont\undefined%
\gdef\SetFigFont#1#2#3#4#5{%
  \reset@font\fontsize{#1}{#2pt}%
  \fontfamily{#3}\fontseries{#4}\fontshape{#5}%
  \selectfont}%
\fi\endgroup%
\mbox{\beginpicture
\setcoordinatesystem units <0.73000cm,0.73000cm>
\unitlength=0.73000cm
\linethickness=1pt
\setplotsymbol ({\makebox(0,0)[l]{\tencirc\symbol{'160}}})
\setshadesymbol ({\thinlinefont .})
\setlinear
%
% Fig POLYLINE object
%
\linethickness= 0.500pt
\setplotsymbol ({\thinlinefont .})
\plot  4.763 23.812  2.699 21.590 /
%
% Fig POLYLINE object
%
\linethickness= 0.500pt
\setplotsymbol ({\thinlinefont .})
\putrule from  3.334 22.860 to  3.334 19.050
%
% Fig POLYLINE object
%
\linethickness= 0.500pt
\setplotsymbol ({\thinlinefont .})
\plot  3.334 19.050  3.334 19.050 /
%
% Fig POLYLINE object
%
\linethickness= 0.500pt
\setplotsymbol ({\thinlinefont .})
\plot  2.857 20.161  5.239 17.621 /
%
% Fig POLYLINE object
%
\linethickness= 0.500pt
\setplotsymbol ({\thinlinefont .})
\plot  6.826 23.971  8.731 21.907 /
%
% Fig POLYLINE object
%
\linethickness= 0.500pt
\setplotsymbol ({\thinlinefont .})
\putrule from  8.255 23.178 to  8.255 19.209
%
% Fig POLYLINE object
%
\linethickness= 0.500pt
\setplotsymbol ({\thinlinefont .})
\plot  8.731 20.479  6.668 17.621 /
%
% Fig POLYLINE object
%
\linethickness= 0.500pt
\setplotsymbol ({\thinlinefont .})
\setdashes < 0.1270cm>
\plot  3.810 23.495  7.779 23.495 /
%
% Fig POLYLINE object
%
\linethickness= 0.500pt
\setplotsymbol ({\thinlinefont .})
\plot  4.128 18.256  8.096 18.256 /
%
% Fig POLYLINE object
%
\linethickness= 0.500pt
\setplotsymbol ({\thinlinefont .})
\plot  2.857 21.114  6.191 22.860 /
%
% Fig POLYLINE object
%
\linethickness= 0.500pt
\setplotsymbol ({\thinlinefont .})
\plot  5.239 22.860  8.731 21.590 /
%
% Fig POLYLINE object
%
\linethickness= 0.500pt
\setplotsymbol ({\thinlinefont .})
\plot  2.857 20.637  6.509 19.367 /
%
% Fig POLYLINE object
%
\linethickness= 0.500pt
\setplotsymbol ({\thinlinefont .})
\plot  5.556 19.367  8.890 21.114 /
%
% Fig POLYLINE object
%
\linethickness= 0.500pt
\setplotsymbol ({\thinlinefont .})
\setsolid
\putrule from 14.287 21.114 to 17.462 21.114
%
% Fig POLYLINE object
%
\linethickness= 0.500pt
\setplotsymbol ({\thinlinefont .})
\putrule from 14.764 22.701 to 14.764 18.891
%
% Fig POLYLINE object
%
\linethickness= 0.500pt
\setplotsymbol ({\thinlinefont .})
\putrule from 16.034 23.495 to 16.034 20.320
%
% Fig POLYLINE object
%
\linethickness= 0.500pt
\setplotsymbol ({\thinlinefont .})
\putrule from 16.669 22.701 to 16.669 20.320
%
% Fig POLYLINE object
%
\linethickness= 0.500pt
\setplotsymbol ({\thinlinefont .})
\putrule from 15.399 22.701 to 15.399 19.526
%
% Fig POLYLINE object
%
\linethickness= 0.500pt
\setplotsymbol ({\thinlinefont .})
\putrule from 19.050 21.114 to 21.749 21.114
%
% Fig POLYLINE object
%
\linethickness= 0.500pt
\setplotsymbol ({\thinlinefont .})
\putrule from 19.367 22.701 to 19.367 20.479
%
% Fig POLYLINE object
%
\linethickness= 0.500pt
\setplotsymbol ({\thinlinefont .})
\putrule from 19.844 23.654 to 19.844 20.479
%
% Fig POLYLINE object
%
\linethickness= 0.500pt
\setplotsymbol ({\thinlinefont .})
\putrule from 20.479 22.701 to 20.479 19.685
%
% Fig POLYLINE object
%
\linethickness= 0.500pt
\setplotsymbol ({\thinlinefont .})
\putrule from 21.273 22.701 to 21.273 18.891
%
% Fig POLYLINE object
%
\linethickness= 0.500pt
\setplotsymbol ({\thinlinefont .})
\setdashes < 0.1270cm>
\plot 16.351 22.225 19.526 22.225 /
%
% Fig POLYLINE object
%
\linethickness= 0.500pt
\setplotsymbol ({\thinlinefont .})
\plot 15.716 23.178 20.320 23.178 /
%
% Fig POLYLINE object
%
\linethickness= 0.500pt
\setplotsymbol ({\thinlinefont .})
\plot 15.081 20.003 20.796 20.003 /
%
% Fig POLYLINE object
%
\linethickness= 0.500pt
\setplotsymbol ({\thinlinefont .})
\plot 14.446 19.209 21.907 19.209 /
%
% Fig POLYLINE object
%
\linethickness= 0.500pt
\setplotsymbol ({\thinlinefont .})
\setsolid
%
% arrow head
%
\plot 10.573 21.177 10.319 21.114 10.573 21.050 /
\putrule from 10.319 21.114 to 12.859 21.114
%
% Fig TEXT object
%
\put{\SetFigFont{9}{10.8}{rm}$p$} [lB] at  5.397 21.907
%
% Fig TEXT object
%
\put{\SetFigFont{9}{10.8}{rm}$q$} [lB] at  5.715 19.844
%
% Fig TEXT object
%
\put{\SetFigFont{9}{10.8}{rm}$(I_0^* ,I_0^*)$} [lB] at 16.510 23.971
%
% Fig TEXT object
%
\put{\SetFigFont{9}{10.8}{rm}$U(2A_3)$} [lB] at  5.080 24.130
\linethickness=0pt
\putrectangle corners at  2.673 24.543 and 21.933 17.596
\endpicture}
}
\vskip 0.5cm
\begin{center}
Figure $5.4$
\end{center} 
We denote
\begin{center}
\begin{tabular}{ccc}
A  &= & {\rm  \ the \ set \ of \  the \ isomorphism \ classes \ of } $U(2A_3)$. \\
B  &= &  {\rm the \  set \ of \ the \ extremal \ rational \ elliptic \ surfaces } \\
   & & {\rm  with \ the \ singular \ fibre \ type } 
$(I_0^*,I_0^*)$.

\end{tabular}
\end{center}

 We define an surjective map $\Phi : A \longrightarrow B$ where $\Phi$ is the blow-up $p$ and $q$ of $D$. Thus the infinity of the number of the isomorphism classes of $U(2A_3)$ is due to the infinity of the  number of the isomorphism classes of the extremal rational elliptic surface with the singular fibre type $(I_0^*, I_0^*)$ (cf. [MP, Theorem $5.4$]).

We  may do the same operation as above to get a surface  which is the minimal resolution of the Picard number two Gorenstein log del Pezzo surface 
$V$ with $K_V^2\geq 3$. For details, we refer to the Table $3$ of the Appedenix.
\begin{lemma} Let  $V$ be the Picard number two  Gorenstein log del Pezzo surface with $K^2_V\geq 3$. Then the singularity type of $V$ is one of the following  listed in Table $5.5$.
\begin{center}
Table $5.5$
\end{center}
\begin{center}
\begin{tabular}{| c||c|}\hline
 $K_V^2$  & Dynkin type of $V$ \\ \hline

$3$ &   $D_5$, $A_5$, $A_4+A_1$, $A_3+2A_1$, $2A_2+A_1$    \\ \hline
$4$ &   $D_4$, $A_4$, $A_3+A_1$, $A_2+2A_1$, $4A_1$      \\ \hline
$5$ &  $A_3$, $A_2+A_1$   \\ \hline

$6$  &  $A_2$, $2A_1$   \\ \hline

$7$ &  $A_1$  \\ \hline

\end{tabular} 
\end{center} 
\end{lemma}
\begin{remark}  Let $V$ be a Picard number two Gorenstein log del Pezzo surface and $U$ the minimal resolution of $V$.

$(1)$ \ \ From a direct calulation, we find that the configurations of the negative curves with the same singularity type are the same except for the type $A_5+A_1$ which has two different configurations where one has two NECs (cf. $U_2^9(A_5+A_1)$) and the other has three NECs (cf. $U_2^{10}(A_5+A_1)$ or $U_2^{11}(A_5+A_1)$).

$(2)$ The configurations of the negative curves of $U$ with $K_U^2 \geq 3$ are listed in the Figure $6$ of  Appendix. From these configurations, it is easy to see  that  there are only two relatively minimal Dynkin types in the above table which are $4A_1$ and $D_4$.
\end{remark}

\begin{lemma} The Picard number two relatively minimal  Gorenstein log del Pezzo surface with the singularity type $4A_1$ is unique.
\end{lemma}
{\em Proof.} We let $U$ be the minimal resolution of $V(4A_1)$ and $D$ the exceptional divisor. Then the configuation of $D$ is in the left part  of the  following figure. 
\vskip 0.5cm
\centerline{\font\thinlinefont=cmr5
\begingroup\makeatletter\ifx\SetFigFont\undefined%
\gdef\SetFigFont#1#2#3#4#5{%
  \reset@font\fontsize{#1}{#2pt}%
  \fontfamily{#3}\fontseries{#4}\fontshape{#5}%
  \selectfont}%
\fi\endgroup%
\mbox{\beginpicture
\setcoordinatesystem units <0.78000cm,0.78000cm>
\unitlength=0.78000cm
\linethickness=1pt
\setplotsymbol ({\makebox(0,0)[l]{\tencirc\symbol{'160}}})
\setshadesymbol ({\thinlinefont .})
\setlinear
%
% Fig POLYLINE object
%
\linethickness= 0.500pt
\setplotsymbol ({\thinlinefont .})
\setdashes < 0.1270cm>
\plot  1.429 22.701  1.429 20.796 /
%
% Fig POLYLINE object
%
\linethickness= 0.500pt
\setplotsymbol ({\thinlinefont .})
\plot  2.223 20.637  4.286 20.637 /
%
% Fig POLYLINE object
%
\linethickness= 0.500pt
\setplotsymbol ({\thinlinefont .})
\plot  2.064 23.178  4.128 23.178 /
%
% Fig POLYLINE object
%
\linethickness= 0.500pt
\setplotsymbol ({\thinlinefont .})
\plot  4.128 23.178  4.128 23.178 /
%
% Fig POLYLINE object
%
\linethickness= 0.500pt
\setplotsymbol ({\thinlinefont .})
\plot  4.763 23.019  4.763 20.955 /
%
% Fig POLYLINE object
%
\linethickness= 0.500pt
\setplotsymbol ({\thinlinefont .})
\plot  4.763 20.955  4.763 20.955 /
%
% Fig POLYLINE object
%
\linethickness= 0.500pt
\setplotsymbol ({\thinlinefont .})
\setsolid
\plot  2.699 23.654  1.111 22.066 /
%
% Fig POLYLINE object
%
\linethickness= 0.500pt
\setplotsymbol ({\thinlinefont .})
\plot  1.111 21.590  3.016 20.161 /
%
% Fig POLYLINE object
%
\linethickness= 0.500pt
\setplotsymbol ({\thinlinefont .})
\plot  5.239 21.749  3.493 20.320 /
%
% Fig POLYLINE object
%
\linethickness= 0.500pt
\setplotsymbol ({\thinlinefont .})
\plot  3.334 23.654  5.239 22.384 /
%
% Fig POLYLINE object
%
\linethickness= 0.500pt
\setplotsymbol ({\thinlinefont .})
\putrule from  6.509 21.907 to  7.779 21.907
%
% arrow head
%
\plot  7.525 21.844  7.779 21.907  7.525 21.971 /
%
%
% Fig POLYLINE object
%
\linethickness= 0.500pt
\setplotsymbol ({\thinlinefont .})
\setdashes < 0.1270cm>
\plot  8.255 22.543  9.684 20.955 /
%
% Fig POLYLINE object
%
\linethickness= 0.500pt
\setplotsymbol ({\thinlinefont .})
\plot  8.890 23.336 11.113 23.336 /
%
% Fig POLYLINE object
%
\linethickness= 0.500pt
\setplotsymbol ({\thinlinefont .})
\plot  9.525 23.812  8.255 21.907 /
%
% Fig POLYLINE object
%
\linethickness= 0.500pt
\setplotsymbol ({\thinlinefont .})
\plot 10.478 23.812 11.906 22.225 /
%
% Fig POLYLINE object
%
\linethickness= 0.500pt
\setplotsymbol ({\thinlinefont .})
\plot 11.748 22.860 10.636 20.955 /
%
% Fig POLYLINE object
%
\linethickness= 0.500pt
\setplotsymbol ({\thinlinefont .})
\plot 10.636 20.955 10.636 20.955 /
%
% Fig POLYLINE object
%
\linethickness= 0.500pt
\setplotsymbol ({\thinlinefont .})
\plot  8.890 21.273 11.430 21.273 /
%
% Fig POLYLINE object
%
\linethickness= 0.500pt
\setplotsymbol ({\thinlinefont .})
\plot 11.430 21.273 11.430 21.273 /
%
% Fig POLYLINE object
%
\linethickness= 0.500pt
\setplotsymbol ({\thinlinefont .})
\setsolid
\putrule from 13.176 22.066 to 14.287 22.066
%
% arrow head
%
\plot 14.034 22.003 14.287 22.066 14.034 22.130 /
%
%
% Fig POLYLINE object
%
\linethickness= 0.500pt
\setplotsymbol ({\thinlinefont .})
\putrule from 15.399 20.955 to 18.733 20.955
%
% Fig POLYLINE object
%
\linethickness= 0.500pt
\setplotsymbol ({\thinlinefont .})
\putrule from 19.367 23.812 to 19.367 21.431
%
% Fig POLYLINE object
%
\linethickness= 0.500pt
\setplotsymbol ({\thinlinefont .})
\setdots < 0.0953cm>
\plot 15.399 21.590 18.891 21.590 /
%
% Fig POLYLINE object
%
\linethickness= 0.500pt
\setplotsymbol ({\thinlinefont .})
\plot 15.399 23.336 18.733 23.336 /
%
% Fig POLYLINE object
%
\linethickness= 0.500pt
\setplotsymbol ({\thinlinefont .})
\plot 15.716 23.812 15.716 21.273 /
%
% Fig POLYLINE object
%
\linethickness= 0.500pt
\setplotsymbol ({\thinlinefont .})
\plot 18.256 23.812 18.256 21.273 /
%
% Fig TEXT object
%
\put{\SetFigFont{9}{10.8}{rm}$E_1$} [lB] at  2.223 23.495
%
% Fig TEXT object
%
\put{\SetFigFont{9}{10.8}{rm}$f_2$} [lB] at  3.651 23.495
%
% Fig TEXT object
%
\put{\SetFigFont{9}{10.8}{rm}$l_1$} [lB] at  1.587 21.590
%
% Fig TEXT object
%
\put{\SetFigFont{9}{10.8}{rm}$l_2$} [lB] at  4.445 21.749
%
% Fig TEXT object
%
\put{\SetFigFont{9}{10.8}{rm}$f_1$} [lB] at  2.064 20.955
%
% Fig TEXT object
%
\put{\SetFigFont{9}{10.8}{rm}$E_2$} [lB] at  3.810 20.955
%
% Fig TEXT object
%
\put{\SetFigFont{8}{9.6}{rm} $E_1$} [lB] at  8.414 22.701
%
% Fig TEXT object
%
\put{\SetFigFont{8}{9.6}{rm}$f_1$} [lB] at  8.890 21.749
%
% Fig TEXT object
%
\put{\SetFigFont{8}{9.6}{rm}$f_2$} [lB] at 10.954 22.860
%
% Fig TEXT object
%
\put{\SetFigFont{8}{9.6}{rm}$E_2$} [lB] at 10.795 21.749
%
% Fig TEXT object
%
\put{\SetFigFont{8}{9.6}{rm}${\bf P}^1\times {\bf P}^1$} [lB] at 15.399 19.526
%
% Fig TEXT object
%
\put{\SetFigFont{9}{10.8}{rm}$(\infty ,\infty )$} [lB] at 14.764 23.495
%
% Fig TEXT object
%
\put{\SetFigFont{9}{10.8}{rm}$(0,0)$} [lB] at 18.256 21.114
%
% Fig TEXT object
%
\put{\SetFigFont{9}{10.8}{rm}$\infty$} [lB] at 15.399 20.320
%
% Fig TEXT object
%
\put{\SetFigFont{9}{10.8}{rm}$0$} [lB] at 17.939 20.320
%
% Fig TEXT object
%
\put{\SetFigFont{9}{10.8}{rm}$\infty$} [lB] at 19.685 23.178
%
% Fig TEXT object
%
\put{\SetFigFont{9}{10.8}{rm}$0$} [lB] at 19.685 21.431
%
% Fig TEXT object
%
\put{\SetFigFont{9}{10.8}{rm}$f_1$} [lB] at 15.240 22.384
%
% Fig TEXT object
%
\put{\SetFigFont{9}{10.8}{rm}$f_2$} [lB] at 18.415 22.384
%
% Fig TEXT object
%
\put{\SetFigFont{9}{10.8}{rm}$U$} [lB] at  2.699 19.209
%
% Fig TEXT object
%
\put{\SetFigFont{9}{10.8}{rm}$W$} [lB] at  9.842 19.367
%
% Fig TEXT object
%
\put{\SetFigFont{9}{10.8}{rm}$\sigma_1$} [lB] at  6.509 22.066
%
% Fig TEXT object
%
\put{\SetFigFont{9}{10.8}{rm}$\sigma_2$} [lB] at 13.018 22.225
\linethickness=0pt
\putrectangle corners at  1.086 23.908 and 19.685 19.209
\endpicture}
}
 \vskip 0.5cm
\begin{center}
Figure $5.5$
\end{center}

Let $\sigma_1 : U \longrightarrow W$ be the contraction of $l_1$ and $l_2$ and
$\sigma_2 : W \longrightarrow {\bf P}^1 \times {\bf P}^1$ the contraction of $E_1$ and $E_2$. By changing the coordinates of ${\bf P}^1 \times {\bf P}^1$, we may assume $\sigma_2 \circ \sigma_1 (E_1)=p_1=(\infty ,\infty )$ and 
$\sigma_2 \circ \sigma_1 (E_2)=p_2=(0, 0)$. There is unique way from $({\bf P}^1\times {\bf P}^1, p_1,p_2)$ to $U$. Thus we prove the lemma.
\vskip 0.3cm 

\begin{lemma} The Picard number two relatively minimal Gorenstein log del Pezzo surface with the singularity type $D_4$ is unique.
\end{lemma} 
{\em Proof.} \ From Figure $6$ in  Appendix, we know that there is unique way from the surface $U(D_5)$ which is the minimal resolution of the Picard number two Gorenstein log del Pezzo surface with the singularity type $D_5$ to  $U(D_4)$. 
%\vskip 0.5cm
% \centerline{\input ex46.tex}
%\vskip 0.5cm

Let $U(A_1+D_6) \longrightarrow U(D_5)$ be  the blow up at the point $p$ (see Figure $6$ of  Appendix). Then  $U(A_1+D_6)$ is the minimal resolution of a Picard number one Gorenstein log del Pezzo surface with the singularity type $A_1+D_6$.
On the other hand,  $U(A_1+D_6)$ is unique by Lemma $4.6$ and there is a unique way from $U(A_1+D_6)$ to $U(D_5)$. Thus we know $U(D_5)$ is unique up to isomorphisms. This proves the lemma. 
\vskip 0.5cm
\hskip -0.6cm
Combining Lemma $5.1$, Theorem $5.7$ and Lemma $5.12$, we prove Theorem $1.5$.
\vskip 0.3cm
\hskip -0.6cm Combining Lemmas $5.6$, Theorem $5.7$, Lemmas $5.14$ and $5.15$, we prove Theorem $1.6$.

\appendix
\section*{Appendix. \ \ List of Configurations and Tables}
\vskip 0.5cm
{\bf Configurations which is used in the proof of Theorem $1.2$}
\vskip 0.3cm
\hskip -0.6cm In the following configurations, we exclude the fibres of type $I_1$ or $II$ which does not contribute to the calculation of the Mordell-Weil group of the corresponding fibration [Sh]. For simplicity, we only list the $(-1)$-curves which will be used in the proof of Theorem $1.2$. 

  \vskip 0.5cm
\centerline{\font\thinlinefont=cmr5
\begingroup\makeatletter\ifx\SetFigFont\undefined%
\gdef\SetFigFont#1#2#3#4#5{%
  \reset@font\fontsize{#1}{#2pt}%
  \fontfamily{#3}\fontseries{#4}\fontshape{#5}%
  \selectfont}%
\fi\endgroup%
\mbox
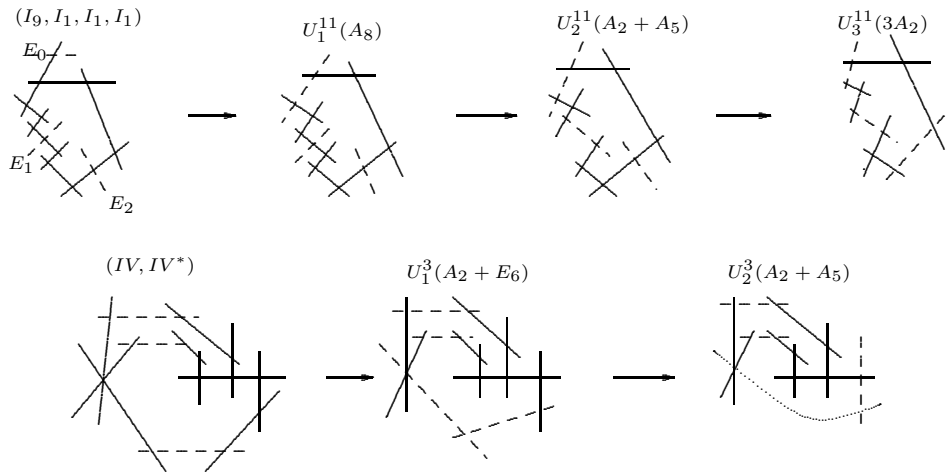
{\beginpicture
\setcoordinatesystem units <0.68000cm,0.68000cm>
\unitlength=0.68000cm
\linethickness=1pt
\setplotsymbol ({\makebox(0,0)[l]{\tencirc\symbol{'160}}})
\setshadesymbol ({\thinlinefont .})
\setlinear
%
% Fig POLYLINE object
%
\linethickness= 0.500pt
\setplotsymbol ({\thinlinefont .})
\putrule from  3.016 24.130 to  7.144 24.130
%
% Fig POLYLINE object
%
\linethickness= 0.500pt
\setplotsymbol ({\thinlinefont .})
\setdashes < 0.1270cm>
\plot  3.016 22.543  7.144 22.543 /
%
% Fig POLYLINE object
%
\linethickness= 0.500pt
\setplotsymbol ({\thinlinefont .})
\setdots < 0.0953cm>
\plot  3.016 20.955  6.985 20.955 /
%
% Fig POLYLINE object
%
\linethickness= 0.500pt
\setplotsymbol ({\thinlinefont .})
\plot  6.985 20.955  7.144 20.955 /
%
% Fig TEXT object
%
\put{\SetFigFont{11}{13.2}{rm}Stands for a curve with self-intersection $0$.} [lB] at  7.620 20.796
%
% Fig TEXT object
%
\put{\SetFigFont{11}{13.2}{rm}Stands for a $(-2)$-curve} [lB] at  7.620 23.971
%
% Fig TEXT object
%
\put{\SetFigFont{11}{13.2}{rm}Stands for a $(-1)$-curve} [lB] at  7.620 22.384
\linethickness=0pt
\putrectangle corners at  2.991 24.543 and  7.620 20.796
\endpicture}
}
   \vskip 0.5cm
   \centerline{\font\thinlinefont=cmr5
\begingroup\makeatletter\ifx\SetFigFont\undefined%
\gdef\SetFigFont#1#2#3#4#5{%
  \reset@font\fontsize{#1}{#2pt}%
  \fontfamily{#3}\fontseries{#4}\fontshape{#5}%
  \selectfont}%
\fi\endgroup%
\mbox{\beginpicture
\setcoordinatesystem units <0.60000cm,0.60000cm>
\unitlength=0.60000cm
\linethickness=1pt
\setplotsymbol ({\makebox(0,0)[l]{\tencirc\symbol{'160}}})
\setshadesymbol ({\thinlinefont .})
\setlinear
%
% Fig CIRCULAR ARC object
%
\linethickness= 0.500pt
\setplotsymbol ({\thinlinefont .})
\circulararc 171.203 degrees from  2.540 21.273 center at  1.581 21.279
%
% Fig CIRCULAR ARC object
%
\linethickness= 0.500pt
\setplotsymbol ({\thinlinefont .})
\circulararc 216.207 degrees from  0.635 22.860 center at  1.374 22.970
%
% Fig CIRCULAR ARC object
%
\linethickness= 0.500pt
\setplotsymbol ({\thinlinefont .})
\circulararc 245.148 degrees from 10.795 20.955 center at 10.522 21.577
%
% Fig POLYLINE object
%
\linethickness= 0.500pt
\setplotsymbol ({\thinlinefont .})
\putrule from  3.969 22.384 to  6.668 22.384
%
% Fig POLYLINE object
%
\linethickness= 0.500pt
\setplotsymbol ({\thinlinefont .})
\putrule from  6.191 23.019 to  6.191 22.066
%
% Fig POLYLINE object
%
\linethickness= 0.500pt
\setplotsymbol ({\thinlinefont .})
\putrule from  5.715 23.019 to  5.715 21.749
%
% Fig POLYLINE object
%
\linethickness= 0.500pt
\setplotsymbol ({\thinlinefont .})
\putrule from  5.239 22.066 to  5.874 22.066
%
% Fig POLYLINE object
%
\linethickness= 0.500pt
\setplotsymbol ({\thinlinefont .})
\putrule from  5.397 22.225 to  5.397 21.114
%
% Fig POLYLINE object
%
\linethickness= 0.500pt
\setplotsymbol ({\thinlinefont .})
\putrule from  4.763 23.019 to  4.763 21.749
%
% Fig POLYLINE object
%
\linethickness= 0.500pt
\setplotsymbol ({\thinlinefont .})
\putrule from  4.128 22.860 to  5.080 22.860
%
% Fig POLYLINE object
%
\linethickness= 0.500pt
\setplotsymbol ({\thinlinefont .})
\putrule from  4.445 23.812 to  4.445 22.701
%
% Fig POLYLINE object
%
\linethickness= 0.500pt
\setplotsymbol ({\thinlinefont .})
\setdashes < 0.1270cm>
\plot  1.587 23.336  4.921 23.336 /
%
% Fig POLYLINE object
%
\linethickness= 0.500pt
\setplotsymbol ({\thinlinefont .})
\plot  2.064 21.590  5.715 21.590 /
%
% Fig POLYLINE object
%
\linethickness= 0.500pt
\setplotsymbol ({\thinlinefont .})
\setsolid
\putrule from  7.303 22.384 to  8.255 22.384
%
% arrow head
%
\plot  8.001 22.320  8.255 22.384  8.001 22.447 /
%
%
% Fig POLYLINE object
%
\linethickness= 0.500pt
\setplotsymbol ({\thinlinefont .})
\putrule from 12.065 22.066 to 14.764 22.066
%
% Fig POLYLINE object
%
\linethickness= 0.500pt
\setplotsymbol ({\thinlinefont .})
\putrule from 14.446 22.701 to 14.446 21.749
%
% Fig POLYLINE object
%
\linethickness= 0.500pt
\setplotsymbol ({\thinlinefont .})
\putrule from 13.970 22.701 to 13.970 21.590
\putrule from 13.970 21.590 to 13.970 21.273
%
% Fig POLYLINE object
%
\linethickness= 0.500pt
\setplotsymbol ({\thinlinefont .})
\putrule from 13.494 21.590 to 14.287 21.590
%
% Fig POLYLINE object
%
\linethickness= 0.500pt
\setplotsymbol ({\thinlinefont .})
\plot 14.129 21.590 14.129 21.590 /
%
% Fig POLYLINE object
%
\linethickness= 0.500pt
\setplotsymbol ({\thinlinefont .})
\putrule from 13.652 21.749 to 13.652 20.796
%
% Fig POLYLINE object
%
\linethickness= 0.500pt
\setplotsymbol ({\thinlinefont .})
\putrule from 13.018 22.701 to 13.018 21.590
%
% Fig POLYLINE object
%
\linethickness= 0.500pt
\setplotsymbol ({\thinlinefont .})
\putrule from 12.383 22.543 to 13.335 22.543
%
% Fig POLYLINE object
%
\linethickness= 0.500pt
\setplotsymbol ({\thinlinefont .})
\setdashes < 0.1270cm>
\plot 12.541 23.495 12.541 22.384 /
%
% Fig POLYLINE object
%
\linethickness= 0.500pt
\setplotsymbol ({\thinlinefont .})
\plot 10.795 21.273 14.287 21.114 /
%
% Fig POLYLINE object
%
\linethickness= 0.500pt
\setplotsymbol ({\thinlinefont .})
\setsolid
\putrule from 15.875 22.225 to 17.304 22.225
%
% arrow head
%
\plot 17.050 22.162 17.304 22.225 17.050 22.288 /
%
%
% Fig POLYLINE object
%
\linethickness= 0.500pt
\setplotsymbol ({\thinlinefont .})
\putrule from 21.114 22.225 to 24.130 22.225
%
% Fig POLYLINE object
%
\linethickness= 0.500pt
\setplotsymbol ({\thinlinefont .})
\putrule from 23.812 22.860 to 23.812 21.907
%
% Fig POLYLINE object
%
\linethickness= 0.500pt
\setplotsymbol ({\thinlinefont .})
\putrule from 23.336 22.860 to 23.336 21.431
%
% Fig POLYLINE object
%
\linethickness= 0.500pt
\setplotsymbol ({\thinlinefont .})
\putrule from 22.860 21.749 to 23.495 21.749
%
% Fig POLYLINE object
%
\linethickness= 0.500pt
\setplotsymbol ({\thinlinefont .})
\putrule from 23.019 21.907 to 23.019 20.955
%
% Fig POLYLINE object
%
\linethickness= 0.500pt
\setplotsymbol ({\thinlinefont .})
\putrule from 22.543 22.860 to 22.543 21.749
%
% Fig POLYLINE object
%
\linethickness= 0.500pt
\setplotsymbol ({\thinlinefont .})
\setdashes < 0.1270cm>
\plot 19.526 21.273 23.336 21.273 /
%
% Fig POLYLINE object
%
\linethickness= 0.500pt
\setplotsymbol ({\thinlinefont .})
\plot 21.590 23.336 22.860 22.384 /
%
% Fig POLYLINE object
%
\linethickness= 0.500pt
\setplotsymbol ({\thinlinefont .})
\setdots < 0.0953cm>
\plot 22.384 23.178 18.574 21.273 /
%
% Fig POLYLINE object
%
\linethickness= 0.500pt
\setplotsymbol ({\thinlinefont .})
\plot 18.415 22.225 19.526 22.384 /
\plot 19.526 22.384 22.066 23.336 /
%
% Fig POLYLINE object
%
\linethickness= 0.500pt
\setplotsymbol ({\thinlinefont .})
\setdashes < 0.1270cm>
\plot  9.525 21.907 13.176 23.495 /
%
% Fig POLYLINE object
%
\linethickness= 0.500pt
\setplotsymbol ({\thinlinefont .})
\setsolid
\plot 20.003 20.955 20.000 20.959 /
\plot 20.000 20.959 19.998 20.970 /
\plot 19.998 20.970 19.992 20.989 /
\plot 19.992 20.989 19.986 21.018 /
\plot 19.986 21.018 19.973 21.059 /
\plot 19.973 21.059 19.958 21.110 /
\plot 19.958 21.110 19.941 21.171 /
\plot 19.941 21.171 19.920 21.243 /
\plot 19.920 21.243 19.897 21.321 /
\plot 19.897 21.321 19.871 21.406 /
\plot 19.871 21.406 19.844 21.493 /
\plot 19.844 21.493 19.816 21.579 /
\plot 19.816 21.579 19.789 21.666 /
\plot 19.789 21.666 19.759 21.751 /
\plot 19.759 21.751 19.732 21.831 /
\plot 19.732 21.831 19.706 21.903 /
\plot 19.706 21.903 19.679 21.971 /
\plot 19.679 21.971 19.653 22.030 /
\plot 19.653 22.030 19.628 22.083 /
\plot 19.628 22.083 19.602 22.130 /
\plot 19.602 22.130 19.577 22.168 /
\plot 19.577 22.168 19.552 22.200 /
\plot 19.552 22.200 19.526 22.225 /
\plot 19.526 22.225 19.494 22.248 /
\plot 19.494 22.248 19.459 22.267 /
\plot 19.459 22.267 19.423 22.280 /
\plot 19.423 22.280 19.382 22.286 /
\plot 19.382 22.286 19.340 22.288 /
\putrule from 19.340 22.288 to 19.293 22.288
\plot 19.293 22.288 19.247 22.284 /
\plot 19.247 22.284 19.198 22.278 /
\plot 19.198 22.278 19.149 22.267 /
\plot 19.149 22.267 19.099 22.257 /
\plot 19.099 22.257 19.050 22.244 /
\plot 19.050 22.244 19.001 22.231 /
\plot 19.001 22.231 18.957 22.217 /
\plot 18.957 22.217 18.912 22.200 /
\plot 18.912 22.200 18.872 22.183 /
\plot 18.872 22.183 18.834 22.166 /
\plot 18.834 22.166 18.802 22.145 /
\plot 18.802 22.145 18.775 22.121 /
\plot 18.775 22.121 18.752 22.096 /
\plot 18.752 22.096 18.733 22.066 /
\plot 18.733 22.066 18.720 22.035 /
\plot 18.720 22.035 18.711 21.999 /
\plot 18.711 21.999 18.707 21.958 /
\putrule from 18.707 21.958 to 18.707 21.910
\plot 18.707 21.910 18.709 21.855 /
\plot 18.709 21.855 18.716 21.793 /
\plot 18.716 21.793 18.726 21.725 /
\plot 18.726 21.725 18.737 21.651 /
\plot 18.737 21.651 18.752 21.573 /
\plot 18.752 21.573 18.768 21.493 /
\plot 18.768 21.493 18.785 21.410 /
\plot 18.785 21.410 18.802 21.330 /
\plot 18.802 21.330 18.821 21.251 /
\plot 18.821 21.251 18.836 21.181 /
\plot 18.836 21.181 18.851 21.118 /
\plot 18.851 21.118 18.864 21.065 /
\plot 18.864 21.065 18.874 21.023 /
\plot 18.874 21.023 18.883 20.991 /
\plot 18.883 20.991 18.887 20.972 /
\plot 18.887 20.972 18.889 20.959 /
\plot 18.889 20.959 18.891 20.955 /
%
% Fig TEXT object
%
\put{\SetFigFont{7}{8.4}{rm}$U_1^2(A_1+E_7)$} [lB] at 10.636 24.448
%
% Fig TEXT object
%
\put{\SetFigFont{7}{8.4}{rm}$U_2^2(A_1+D_6)$} [lB] at 20.003 24.448
%
% Fig TEXT object
%
\put{\SetFigFont{6}{7.2}{rm}$E_2$} [lB] at 22.860 22.543
%
% Fig TEXT object
%
\put{\SetFigFont{6}{7.2}{rm}$p$} [lB] at 21.749 23.336
%
% Fig TEXT object
%
\put{\SetFigFont{6}{7.2}{rm}$q$} [lB] at 21.749 22.543
%
% Fig TEXT object
%
\put{\SetFigFont{6}{7.2}{rm}$L$} [lB] at 18.256 22.860
%
% Fig TEXT object
%
\put{\SetFigFont{6}{7.2}{rm}$L_1$} [lB] at 18.415 20.955
%
% Fig TEXT object
%
\put{\SetFigFont{7}{8.4}{rm}$(III,III^*)$} [lB] at  2.064 24.448
\linethickness=0pt
\putrectangle corners at  0.618 24.860 and 24.155 20.771
\endpicture}
}
   \vskip 0.5cm
   \centerline{\font\thinlinefont=cmr5
\begingroup\makeatletter\ifx\SetFigFont\undefined%
\gdef\SetFigFont#1#2#3#4#5{%
  \reset@font\fontsize{#1}{#2pt}%
  \fontfamily{#3}\fontseries{#4}\fontshape{#5}%
  \selectfont}%
\fi\endgroup%
\mbox{\beginpicture
\setcoordinatesystem units <0.46000cm,0.46000cm>
\unitlength=0.46000cm
\linethickness=1pt
\setplotsymbol ({\makebox(0,0)[l]{\tencirc\symbol{'160}}})
\setshadesymbol ({\thinlinefont .})
\setlinear
%
% Fig POLYLINE object
%
\linethickness= 0.500pt
\setplotsymbol ({\thinlinefont .})
\putrule from  3.969 22.384 to  6.668 22.384
%
% Fig POLYLINE object
%
\linethickness= 0.500pt
\setplotsymbol ({\thinlinefont .})
\putrule from  6.191 23.019 to  6.191 22.066
%
% Fig POLYLINE object
%
\linethickness= 0.500pt
\setplotsymbol ({\thinlinefont .})
\putrule from  5.715 23.019 to  5.715 21.749
%
% Fig POLYLINE object
%
\linethickness= 0.500pt
\setplotsymbol ({\thinlinefont .})
\putrule from  5.239 22.066 to  5.874 22.066
%
% Fig POLYLINE object
%
\linethickness= 0.500pt
\setplotsymbol ({\thinlinefont .})
\putrule from  5.397 22.225 to  5.397 21.114
%
% Fig POLYLINE object
%
\linethickness= 0.500pt
\setplotsymbol ({\thinlinefont .})
\putrule from  4.763 23.019 to  4.763 21.749
%
% Fig POLYLINE object
%
\linethickness= 0.500pt
\setplotsymbol ({\thinlinefont .})
\putrule from  4.128 22.860 to  5.080 22.860
%
% Fig POLYLINE object
%
\linethickness= 0.500pt
\setplotsymbol ({\thinlinefont .})
\putrule from  4.445 23.812 to  4.445 22.701
%
% Fig POLYLINE object
%
\linethickness= 0.500pt
\setplotsymbol ({\thinlinefont .})
\setdashes < 0.1270cm>
\plot  2.064 21.590  5.715 21.590 /
%
% Fig POLYLINE object
%
\linethickness= 0.500pt
\setplotsymbol ({\thinlinefont .})
\setsolid
\putrule from  7.303 22.384 to  8.255 22.384
%
% arrow head
%
\plot  8.001 22.320  8.255 22.384  8.001 22.447 /
%
%
% Fig POLYLINE object
%
\linethickness= 0.500pt
\setplotsymbol ({\thinlinefont .})
\putrule from 12.065 22.066 to 14.764 22.066
%
% Fig POLYLINE object
%
\linethickness= 0.500pt
\setplotsymbol ({\thinlinefont .})
\putrule from 14.446 22.701 to 14.446 21.749
%
% Fig POLYLINE object
%
\linethickness= 0.500pt
\setplotsymbol ({\thinlinefont .})
\putrule from 13.970 22.701 to 13.970 21.590
\putrule from 13.970 21.590 to 13.970 21.273
%
% Fig POLYLINE object
%
\linethickness= 0.500pt
\setplotsymbol ({\thinlinefont .})
\putrule from 13.494 21.590 to 14.287 21.590
%
% Fig POLYLINE object
%
\linethickness= 0.500pt
\setplotsymbol ({\thinlinefont .})
\plot 14.129 21.590 14.129 21.590 /
%
% Fig POLYLINE object
%
\linethickness= 0.500pt
\setplotsymbol ({\thinlinefont .})
\putrule from 13.652 21.749 to 13.652 20.796
%
% Fig POLYLINE object
%
\linethickness= 0.500pt
\setplotsymbol ({\thinlinefont .})
\putrule from 13.018 22.701 to 13.018 21.590
%
% Fig POLYLINE object
%
\linethickness= 0.500pt
\setplotsymbol ({\thinlinefont .})
\putrule from 12.383 22.543 to 13.335 22.543
%
% Fig POLYLINE object
%
\linethickness= 0.500pt
\setplotsymbol ({\thinlinefont .})
\setdashes < 0.1270cm>
\plot 12.541 23.495 12.541 22.384 /
%
% Fig POLYLINE object
%
\linethickness= 0.500pt
\setplotsymbol ({\thinlinefont .})
\plot 10.795 21.273 14.287 21.114 /
%
% Fig POLYLINE object
%
\linethickness= 0.500pt
\setplotsymbol ({\thinlinefont .})
\setsolid
\putrule from 15.875 22.225 to 17.304 22.225
%
% arrow head
%
\plot 17.050 22.162 17.304 22.225 17.050 22.288 /
%
%
% Fig POLYLINE object
%
\linethickness= 0.500pt
\setplotsymbol ({\thinlinefont .})
\putrule from 21.114 22.225 to 24.130 22.225
%
% Fig POLYLINE object
%
\linethickness= 0.500pt
\setplotsymbol ({\thinlinefont .})
\putrule from 23.812 22.860 to 23.812 21.907
%
% Fig POLYLINE object
%
\linethickness= 0.500pt
\setplotsymbol ({\thinlinefont .})
\putrule from 23.336 22.860 to 23.336 21.431
%
% Fig POLYLINE object
%
\linethickness= 0.500pt
\setplotsymbol ({\thinlinefont .})
\putrule from 22.860 21.749 to 23.495 21.749
%
% Fig POLYLINE object
%
\linethickness= 0.500pt
\setplotsymbol ({\thinlinefont .})
\putrule from 23.019 21.907 to 23.019 20.955
%
% Fig POLYLINE object
%
\linethickness= 0.500pt
\setplotsymbol ({\thinlinefont .})
\putrule from 22.543 22.860 to 22.543 21.749
%
% Fig POLYLINE object
%
\linethickness= 0.500pt
\setplotsymbol ({\thinlinefont .})
\setdashes < 0.1270cm>
\plot 19.526 21.273 23.336 21.273 /
%
% Fig POLYLINE object
%
\linethickness= 0.500pt
\setplotsymbol ({\thinlinefont .})
\plot 21.590 23.336 22.860 22.384 /
%
% Fig POLYLINE object
%
\linethickness= 0.500pt
\setplotsymbol ({\thinlinefont .})
\plot  9.842 22.066 13.176 23.336 /
%
% Fig POLYLINE object
%
\linethickness= 0.500pt
\setplotsymbol ({\thinlinefont .})
\setsolid
\plot  3.493 23.654  0.635 22.066 /
%
% Fig POLYLINE object
%
\linethickness= 0.500pt
\setplotsymbol ({\thinlinefont .})
\setdashes < 0.1270cm>
\plot  2.857 23.495  5.080 23.495 /
%
% Fig POLYLINE object
%
\linethickness= 0.500pt
\setplotsymbol ({\thinlinefont .})
\setsolid
\plot 21.273 21.273 21.273 21.273 /
%
% Fig POLYLINE object
%
\linethickness= 0.500pt
\setplotsymbol ({\thinlinefont .})
\plot 21.273 21.273 21.273 21.273 /
%
% Fig POLYLINE object
%
\linethickness= 0.500pt
\setplotsymbol ({\thinlinefont .})
\putrule from 25.082 22.225 to 26.035 22.225
%
% arrow head
%
\plot 25.781 22.162 26.035 22.225 25.781 22.288 /
%
%
% Fig POLYLINE object
%
\linethickness= 0.500pt
\setplotsymbol ({\thinlinefont .})
\putrule from 27.940 22.225 to 30.321 22.225
%
% Fig POLYLINE object
%
\linethickness= 0.500pt
\setplotsymbol ({\thinlinefont .})
\putrule from 29.845 22.860 to 29.845 21.749
%
% Fig POLYLINE object
%
\linethickness= 0.500pt
\setplotsymbol ({\thinlinefont .})
\putrule from 29.369 22.860 to 29.369 21.431
%
% Fig POLYLINE object
%
\linethickness= 0.500pt
\setplotsymbol ({\thinlinefont .})
\putrule from 28.893 21.749 to 29.527 21.749
%
% Fig POLYLINE object
%
\linethickness= 0.500pt
\setplotsymbol ({\thinlinefont .})
\putrule from 29.051 21.907 to 29.051 20.955
%
% Fig POLYLINE object
%
\linethickness= 0.500pt
\setplotsymbol ({\thinlinefont .})
\setdashes < 0.1270cm>
\plot 28.575 22.860 28.575 21.749 /
%
% Fig POLYLINE object
%
\linethickness= 0.500pt
\setplotsymbol ({\thinlinefont .})
\plot 26.988 21.114 29.369 21.114 /
%
% Fig POLYLINE object
%
\linethickness= 0.500pt
\setplotsymbol ({\thinlinefont .})
\setsolid
\putrule from 10.319 22.701 to 10.319 22.697
\putrule from 10.319 22.697 to 10.319 22.689
\plot 10.319 22.689 10.317 22.672 /
\plot 10.317 22.672 10.315 22.646 /
\plot 10.315 22.646 10.312 22.608 /
\plot 10.312 22.608 10.308 22.562 /
\plot 10.308 22.562 10.304 22.502 /
\plot 10.304 22.502 10.300 22.435 /
\plot 10.300 22.435 10.293 22.356 /
\plot 10.293 22.356 10.289 22.267 /
\plot 10.289 22.267 10.283 22.174 /
\plot 10.283 22.174 10.276 22.075 /
\plot 10.276 22.075 10.270 21.973 /
\plot 10.270 21.973 10.266 21.872 /
\plot 10.266 21.872 10.262 21.770 /
\plot 10.262 21.770 10.257 21.668 /
\plot 10.257 21.668 10.255 21.573 /
\plot 10.255 21.573 10.253 21.482 /
\putrule from 10.253 21.482 to 10.253 21.397
\plot 10.253 21.397 10.255 21.319 /
\plot 10.255 21.319 10.257 21.247 /
\plot 10.257 21.247 10.264 21.181 /
\plot 10.264 21.181 10.270 21.124 /
\plot 10.270 21.124 10.279 21.071 /
\plot 10.279 21.071 10.289 21.027 /
\plot 10.289 21.027 10.304 20.989 /
\plot 10.304 20.989 10.319 20.955 /
\plot 10.319 20.955 10.340 20.921 /
\plot 10.340 20.921 10.365 20.894 /
\plot 10.365 20.894 10.395 20.870 /
\plot 10.395 20.870 10.427 20.849 /
\plot 10.427 20.849 10.463 20.834 /
\plot 10.463 20.834 10.503 20.820 /
\plot 10.503 20.820 10.545 20.809 /
\plot 10.545 20.809 10.590 20.800 /
\plot 10.590 20.800 10.636 20.794 /
\plot 10.636 20.794 10.687 20.790 /
\plot 10.687 20.790 10.736 20.786 /
\plot 10.736 20.786 10.787 20.784 /
\putrule from 10.787 20.784 to 10.839 20.784
\plot 10.839 20.784 10.888 20.786 /
\plot 10.888 20.786 10.939 20.788 /
\plot 10.939 20.788 10.986 20.794 /
\plot 10.986 20.794 11.032 20.800 /
\plot 11.032 20.800 11.077 20.811 /
\plot 11.077 20.811 11.117 20.826 /
\plot 11.117 20.826 11.155 20.841 /
\plot 11.155 20.841 11.189 20.862 /
\plot 11.189 20.862 11.218 20.887 /
\plot 11.218 20.887 11.246 20.919 /
\plot 11.246 20.919 11.271 20.955 /
\plot 11.271 20.955 11.288 20.989 /
\plot 11.288 20.989 11.305 21.027 /
\plot 11.305 21.027 11.320 21.074 /
\plot 11.320 21.074 11.333 21.122 /
\plot 11.333 21.122 11.343 21.179 /
\plot 11.343 21.179 11.354 21.243 /
\plot 11.354 21.243 11.364 21.315 /
\plot 11.364 21.315 11.373 21.391 /
\plot 11.373 21.391 11.379 21.474 /
\plot 11.379 21.474 11.388 21.562 /
\plot 11.388 21.562 11.392 21.658 /
\plot 11.392 21.658 11.398 21.755 /
\plot 11.398 21.755 11.402 21.857 /
\plot 11.402 21.857 11.407 21.963 /
\plot 11.407 21.963 11.411 22.066 /
\plot 11.411 22.066 11.415 22.170 /
\plot 11.415 22.170 11.417 22.272 /
\plot 11.417 22.272 11.419 22.367 /
\plot 11.419 22.367 11.422 22.458 /
\plot 11.422 22.458 11.424 22.540 /
\plot 11.424 22.540 11.426 22.614 /
\plot 11.426 22.614 11.428 22.680 /
\putrule from 11.428 22.680 to 11.428 22.733
\putrule from 11.428 22.733 to 11.428 22.777
\plot 11.428 22.777 11.430 22.809 /
\putrule from 11.430 22.809 to 11.430 22.832
\putrule from 11.430 22.832 to 11.430 22.847
\putrule from 11.430 22.847 to 11.430 22.856
\putrule from 11.430 22.856 to 11.430 22.860
%
% Fig POLYLINE object
%
\linethickness= 0.500pt
\setplotsymbol ({\thinlinefont .})
\putrule from 18.733 22.701 to 18.733 22.697
\putrule from 18.733 22.697 to 18.733 22.689
\plot 18.733 22.689 18.730 22.674 /
\plot 18.730 22.674 18.728 22.650 /
\plot 18.728 22.650 18.726 22.619 /
\plot 18.726 22.619 18.722 22.574 /
\plot 18.722 22.574 18.718 22.521 /
\plot 18.718 22.521 18.713 22.456 /
\plot 18.713 22.456 18.707 22.382 /
\plot 18.707 22.382 18.703 22.299 /
\plot 18.703 22.299 18.697 22.208 /
\plot 18.697 22.208 18.690 22.113 /
\plot 18.690 22.113 18.684 22.011 /
\plot 18.684 22.011 18.677 21.907 /
\plot 18.677 21.907 18.671 21.804 /
\plot 18.671 21.804 18.667 21.698 /
\plot 18.667 21.698 18.663 21.596 /
\plot 18.663 21.596 18.658 21.499 /
\plot 18.658 21.499 18.656 21.404 /
\putrule from 18.656 21.404 to 18.656 21.315
\putrule from 18.656 21.315 to 18.656 21.232
\plot 18.656 21.232 18.658 21.156 /
\plot 18.658 21.156 18.663 21.084 /
\plot 18.663 21.084 18.669 21.021 /
\plot 18.669 21.021 18.677 20.963 /
\plot 18.677 20.963 18.688 20.915 /
\plot 18.688 20.915 18.701 20.868 /
\plot 18.701 20.868 18.716 20.830 /
\plot 18.716 20.830 18.733 20.796 /
\plot 18.733 20.796 18.758 20.760 /
\plot 18.758 20.760 18.788 20.729 /
\plot 18.788 20.729 18.821 20.703 /
\plot 18.821 20.703 18.860 20.682 /
\plot 18.860 20.682 18.902 20.665 /
\plot 18.902 20.665 18.948 20.650 /
\plot 18.948 20.650 18.999 20.640 /
\plot 18.999 20.640 19.052 20.631 /
\plot 19.052 20.631 19.109 20.625 /
\plot 19.109 20.625 19.169 20.621 /
\plot 19.169 20.621 19.228 20.618 /
\putrule from 19.228 20.618 to 19.287 20.618
\putrule from 19.287 20.618 to 19.348 20.618
\plot 19.348 20.618 19.408 20.621 /
\plot 19.408 20.621 19.467 20.625 /
\plot 19.467 20.625 19.524 20.631 /
\plot 19.524 20.631 19.577 20.640 /
\plot 19.577 20.640 19.628 20.650 /
\plot 19.628 20.650 19.674 20.665 /
\plot 19.674 20.665 19.717 20.682 /
\plot 19.717 20.682 19.755 20.703 /
\plot 19.755 20.703 19.789 20.729 /
\plot 19.789 20.729 19.818 20.760 /
\plot 19.818 20.760 19.844 20.796 /
\plot 19.844 20.796 19.861 20.830 /
\plot 19.861 20.830 19.875 20.868 /
\plot 19.875 20.868 19.888 20.915 /
\plot 19.888 20.915 19.899 20.963 /
\plot 19.899 20.963 19.907 21.021 /
\plot 19.907 21.021 19.914 21.084 /
\plot 19.914 21.084 19.918 21.156 /
\plot 19.918 21.156 19.920 21.232 /
\putrule from 19.920 21.232 to 19.920 21.315
\putrule from 19.920 21.315 to 19.920 21.404
\plot 19.920 21.404 19.918 21.499 /
\plot 19.918 21.499 19.914 21.596 /
\plot 19.914 21.596 19.909 21.698 /
\plot 19.909 21.698 19.905 21.804 /
\plot 19.905 21.804 19.899 21.907 /
\plot 19.899 21.907 19.892 22.011 /
\plot 19.892 22.011 19.886 22.113 /
\plot 19.886 22.113 19.880 22.208 /
\plot 19.880 22.208 19.873 22.299 /
\plot 19.873 22.299 19.869 22.382 /
\plot 19.869 22.382 19.863 22.456 /
\plot 19.863 22.456 19.859 22.521 /
\plot 19.859 22.521 19.854 22.574 /
\plot 19.854 22.574 19.850 22.619 /
\plot 19.850 22.619 19.848 22.650 /
\plot 19.848 22.650 19.846 22.674 /
\plot 19.846 22.674 19.844 22.689 /
\putrule from 19.844 22.689 to 19.844 22.697
\putrule from 19.844 22.697 to 19.844 22.701
%
% Fig POLYLINE object
%
\linethickness= 0.500pt
\setplotsymbol ({\thinlinefont .})
\putrule from  1.111 22.701 to  1.111 22.697
\plot  1.111 22.697  1.113 22.689 /
\plot  1.113 22.689  1.115 22.672 /
\plot  1.115 22.672  1.120 22.644 /
\plot  1.120 22.644  1.126 22.606 /
\plot  1.126 22.606  1.135 22.557 /
\plot  1.135 22.557  1.143 22.498 /
\plot  1.143 22.498  1.156 22.428 /
\plot  1.156 22.428  1.168 22.348 /
\plot  1.168 22.348  1.183 22.259 /
\plot  1.183 22.259  1.200 22.162 /
\plot  1.200 22.162  1.219 22.062 /
\plot  1.219 22.062  1.238 21.958 /
\plot  1.238 21.958  1.259 21.852 /
\plot  1.259 21.852  1.278 21.749 /
\plot  1.278 21.749  1.302 21.649 /
\plot  1.302 21.649  1.323 21.552 /
\plot  1.323 21.552  1.346 21.461 /
\plot  1.346 21.461  1.369 21.376 /
\plot  1.369 21.376  1.393 21.298 /
\plot  1.393 21.298  1.416 21.226 /
\plot  1.416 21.226  1.441 21.162 /
\plot  1.441 21.162  1.467 21.107 /
\plot  1.467 21.107  1.494 21.059 /
\plot  1.494 21.059  1.524 21.018 /
\plot  1.524 21.018  1.554 20.983 /
\plot  1.554 20.983  1.587 20.955 /
\plot  1.587 20.955  1.623 20.932 /
\plot  1.623 20.932  1.662 20.915 /
\plot  1.662 20.915  1.702 20.900 /
\plot  1.702 20.900  1.748 20.889 /
\plot  1.748 20.889  1.797 20.883 /
\plot  1.797 20.883  1.848 20.879 /
\plot  1.848 20.879  1.905 20.877 /
\plot  1.905 20.877  1.964 20.879 /
\plot  1.964 20.879  2.026 20.881 /
\plot  2.026 20.881  2.089 20.885 /
\plot  2.089 20.885  2.157 20.892 /
\plot  2.157 20.892  2.225 20.900 /
\plot  2.225 20.900  2.292 20.911 /
\plot  2.292 20.911  2.362 20.921 /
\plot  2.362 20.921  2.430 20.932 /
\plot  2.430 20.932  2.498 20.947 /
\plot  2.498 20.947  2.563 20.961 /
\plot  2.563 20.961  2.627 20.976 /
\plot  2.627 20.976  2.686 20.995 /
\plot  2.686 20.995  2.743 21.014 /
\plot  2.743 21.014  2.796 21.035 /
\plot  2.796 21.035  2.845 21.061 /
\plot  2.845 21.061  2.887 21.086 /
\plot  2.887 21.086  2.925 21.116 /
\plot  2.925 21.116  2.957 21.150 /
\plot  2.957 21.150  2.982 21.186 /
\plot  2.982 21.186  3.001 21.226 /
\plot  3.001 21.226  3.016 21.273 /
\plot  3.016 21.273  3.023 21.315 /
\plot  3.023 21.315  3.027 21.359 /
\putrule from  3.027 21.359 to  3.027 21.410
\plot  3.027 21.410  3.020 21.467 /
\plot  3.020 21.467  3.012 21.529 /
\plot  3.012 21.529  2.999 21.594 /
\plot  2.999 21.594  2.980 21.668 /
\plot  2.980 21.668  2.959 21.749 /
\plot  2.959 21.749  2.934 21.833 /
\plot  2.934 21.833  2.904 21.924 /
\plot  2.904 21.924  2.872 22.022 /
\plot  2.872 22.022  2.836 22.121 /
\plot  2.836 22.121  2.798 22.229 /
\plot  2.798 22.229  2.756 22.337 /
\plot  2.756 22.337  2.714 22.449 /
\plot  2.714 22.449  2.669 22.564 /
\plot  2.669 22.564  2.625 22.678 /
\plot  2.625 22.678  2.580 22.790 /
\plot  2.580 22.790  2.536 22.900 /
\plot  2.536 22.900  2.491 23.008 /
\plot  2.491 23.008  2.451 23.108 /
\plot  2.451 23.108  2.411 23.203 /
\plot  2.411 23.203  2.375 23.290 /
\plot  2.375 23.290  2.343 23.368 /
\plot  2.343 23.368  2.314 23.436 /
\plot  2.314 23.436  2.290 23.495 /
\plot  2.290 23.495  2.269 23.544 /
\plot  2.269 23.544  2.252 23.582 /
\plot  2.252 23.582  2.242 23.611 /
\plot  2.242 23.611  2.233 23.630 /
\plot  2.233 23.630  2.227 23.643 /
\plot  2.227 23.643  2.225 23.652 /
\plot  2.225 23.652  2.223 23.654 /
%
% Fig POLYLINE object
%
\linethickness= 0.500pt
\setplotsymbol ({\thinlinefont .})
\putrule from 26.352 22.701 to 26.352 22.697
\putrule from 26.352 22.697 to 26.352 22.689
\putrule from 26.352 22.689 to 26.352 22.672
\plot 26.352 22.672 26.355 22.646 /
\putrule from 26.355 22.646 to 26.355 22.610
\plot 26.355 22.610 26.357 22.564 /
\putrule from 26.357 22.564 to 26.357 22.507
\plot 26.357 22.507 26.359 22.437 /
\plot 26.359 22.437 26.361 22.356 /
\plot 26.361 22.356 26.365 22.265 /
\plot 26.365 22.265 26.367 22.168 /
\plot 26.367 22.168 26.372 22.064 /
\plot 26.372 22.064 26.374 21.954 /
\plot 26.374 21.954 26.380 21.842 /
\plot 26.380 21.842 26.384 21.728 /
\plot 26.384 21.728 26.388 21.615 /
\plot 26.388 21.615 26.395 21.505 /
\plot 26.395 21.505 26.401 21.397 /
\plot 26.401 21.397 26.408 21.296 /
\plot 26.408 21.296 26.414 21.201 /
\plot 26.414 21.201 26.420 21.110 /
\plot 26.420 21.110 26.429 21.027 /
\plot 26.429 21.027 26.437 20.951 /
\plot 26.437 20.951 26.448 20.881 /
\plot 26.448 20.881 26.458 20.820 /
\plot 26.458 20.820 26.469 20.764 /
\plot 26.469 20.764 26.482 20.716 /
\plot 26.482 20.716 26.496 20.673 /
\plot 26.496 20.673 26.511 20.637 /
\plot 26.511 20.637 26.530 20.602 /
\plot 26.530 20.602 26.551 20.570 /
\plot 26.551 20.570 26.575 20.542 /
\plot 26.575 20.542 26.602 20.519 /
\plot 26.602 20.519 26.630 20.500 /
\plot 26.630 20.500 26.659 20.485 /
\plot 26.659 20.485 26.693 20.470 /
\plot 26.693 20.470 26.727 20.460 /
\plot 26.727 20.460 26.763 20.449 /
\plot 26.763 20.449 26.799 20.443 /
\plot 26.799 20.443 26.837 20.436 /
\plot 26.837 20.436 26.877 20.432 /
\plot 26.877 20.432 26.916 20.430 /
\plot 26.916 20.430 26.954 20.428 /
\putrule from 26.954 20.428 to 26.994 20.428
\plot 26.994 20.428 27.032 20.432 /
\plot 27.032 20.432 27.068 20.436 /
\plot 27.068 20.436 27.104 20.443 /
\plot 27.104 20.443 27.136 20.453 /
\plot 27.136 20.453 27.167 20.466 /
\plot 27.167 20.466 27.197 20.483 /
\plot 27.197 20.483 27.225 20.504 /
\plot 27.225 20.504 27.248 20.530 /
\plot 27.248 20.530 27.269 20.559 /
\plot 27.269 20.559 27.288 20.595 /
\plot 27.288 20.595 27.305 20.637 /
\plot 27.305 20.637 27.316 20.678 /
\plot 27.316 20.678 27.326 20.722 /
\plot 27.326 20.722 27.335 20.773 /
\plot 27.335 20.773 27.343 20.832 /
\plot 27.343 20.832 27.347 20.896 /
\plot 27.347 20.896 27.354 20.968 /
\plot 27.354 20.968 27.356 21.046 /
\plot 27.356 21.046 27.358 21.133 /
\plot 27.358 21.133 27.360 21.226 /
\putrule from 27.360 21.226 to 27.360 21.328
\plot 27.360 21.328 27.358 21.433 /
\putrule from 27.358 21.433 to 27.358 21.546
\plot 27.358 21.546 27.356 21.662 /
\plot 27.356 21.662 27.352 21.783 /
\plot 27.352 21.783 27.349 21.903 /
\plot 27.349 21.903 27.345 22.026 /
\plot 27.345 22.026 27.341 22.149 /
\plot 27.341 22.149 27.337 22.267 /
\plot 27.337 22.267 27.333 22.382 /
\plot 27.333 22.382 27.328 22.490 /
\plot 27.328 22.490 27.324 22.589 /
\plot 27.324 22.589 27.320 22.680 /
\plot 27.320 22.680 27.316 22.761 /
\plot 27.316 22.761 27.313 22.828 /
\plot 27.313 22.828 27.311 22.885 /
\plot 27.311 22.885 27.309 22.932 /
\plot 27.309 22.932 27.307 22.966 /
\putrule from 27.307 22.966 to 27.307 22.991
\plot 27.307 22.991 27.305 23.006 /
\putrule from 27.305 23.006 to 27.305 23.015
\putrule from 27.305 23.015 to 27.305 23.019
%
% Fig POLYLINE object
%
\linethickness= 0.500pt
\setplotsymbol ({\thinlinefont .})
\setdots < 0.0953cm>
\plot 18.256 22.066 18.258 22.066 /
\plot 18.258 22.066 18.267 22.066 /
\plot 18.267 22.066 18.280 22.068 /
\plot 18.280 22.068 18.301 22.070 /
\plot 18.301 22.070 18.328 22.073 /
\plot 18.328 22.073 18.366 22.077 /
\plot 18.366 22.077 18.415 22.081 /
\plot 18.415 22.081 18.472 22.087 /
\plot 18.472 22.087 18.540 22.094 /
\plot 18.540 22.094 18.616 22.100 /
\plot 18.616 22.100 18.701 22.109 /
\plot 18.701 22.109 18.792 22.119 /
\plot 18.792 22.119 18.889 22.130 /
\plot 18.889 22.130 18.991 22.140 /
\plot 18.991 22.140 19.094 22.151 /
\plot 19.094 22.151 19.200 22.164 /
\plot 19.200 22.164 19.308 22.174 /
\plot 19.308 22.174 19.412 22.187 /
\plot 19.412 22.187 19.516 22.200 /
\plot 19.516 22.200 19.617 22.212 /
\plot 19.617 22.212 19.715 22.227 /
\plot 19.715 22.227 19.810 22.240 /
\plot 19.810 22.240 19.899 22.255 /
\plot 19.899 22.255 19.983 22.267 /
\plot 19.983 22.267 20.066 22.282 /
\plot 20.066 22.282 20.142 22.297 /
\plot 20.142 22.297 20.216 22.314 /
\plot 20.216 22.314 20.286 22.329 /
\plot 20.286 22.329 20.352 22.346 /
\plot 20.352 22.346 20.415 22.365 /
\plot 20.415 22.365 20.479 22.384 /
\plot 20.479 22.384 20.544 22.405 /
\plot 20.544 22.405 20.610 22.430 /
\plot 20.610 22.430 20.676 22.456 /
\plot 20.676 22.456 20.743 22.483 /
\plot 20.743 22.483 20.811 22.515 /
\plot 20.811 22.515 20.881 22.549 /
\plot 20.881 22.549 20.955 22.585 /
\plot 20.955 22.585 21.029 22.623 /
\plot 21.029 22.623 21.107 22.665 /
\plot 21.107 22.665 21.188 22.710 /
\plot 21.188 22.710 21.270 22.756 /
\plot 21.270 22.756 21.353 22.805 /
\plot 21.353 22.805 21.440 22.854 /
\plot 21.440 22.854 21.524 22.904 /
\plot 21.524 22.904 21.609 22.955 /
\plot 21.609 22.955 21.692 23.006 /
\plot 21.692 23.006 21.772 23.055 /
\plot 21.772 23.055 21.846 23.101 /
\plot 21.846 23.101 21.918 23.144 /
\plot 21.918 23.144 21.982 23.184 /
\plot 21.982 23.184 22.039 23.220 /
\plot 22.039 23.220 22.087 23.249 /
\plot 22.087 23.249 22.130 23.275 /
\plot 22.130 23.275 22.162 23.296 /
\plot 22.162 23.296 22.187 23.313 /
\plot 22.187 23.313 22.204 23.324 /
\plot 22.204 23.324 22.217 23.330 /
\plot 22.217 23.330 22.223 23.334 /
\plot 22.223 23.334 22.225 23.336 /
%
% Fig TEXT object
%
\put{\SetFigFont{6}{7.2}{rm}$(I_1, I_2, III^*)$} [lB] at  2.064 24.765
%
% Fig TEXT object
%
\put{\SetFigFont{6}{7.2}{rm}$U_1^6(A_1+E_7)$} [lB] at 10.954 24.448
%
% Fig TEXT object
%
\put{\SetFigFont{6}{7.2}{rm}$U_2^6(A_1+D_6)$} [lB] at 19.685 24.289
%
% Fig TEXT object
%
\put{\SetFigFont{6}{7.2}{rm}$r$} [lB] at 21.114 20.955
%
% Fig TEXT object
%
\put{\SetFigFont{6}{7.2}{rm}$U_3^6(A_1+A_5)$} [lB] at 26.829 24.130
\linethickness=0pt
\putrectangle corners at  0.610 25.178 and 30.347 20.403
\endpicture}
}
   \vskip 0.5cm
   \centerline{\font\thinlinefont=cmr5
\begingroup\makeatletter\ifx\SetFigFont\undefined%
\gdef\SetFigFont#1#2#3#4#5{%
  \reset@font\fontsize{#1}{#2pt}%
  \fontfamily{#3}\fontseries{#4}\fontshape{#5}%
  \selectfont}%
\fi\endgroup%
\mbox{\beginpicture
\setcoordinatesystem units <0.56000cm,0.56000cm>
\unitlength=0.56000cm
\linethickness=1pt
\setplotsymbol ({\makebox(0,0)[l]{\tencirc\symbol{'160}}})
\setshadesymbol ({\thinlinefont .})
\setlinear
%
% Fig POLYLINE object
%
\linethickness= 0.500pt
\setplotsymbol ({\thinlinefont .})
\putrule from  1.111 23.654 to  2.699 23.654
%
% Fig POLYLINE object
%
\linethickness= 0.500pt
\setplotsymbol ({\thinlinefont .})
\putrule from  3.175 24.448 to  3.175 23.178
%
% Fig POLYLINE object
%
\linethickness= 0.500pt
\setplotsymbol ({\thinlinefont .})
\putrule from  2.381 24.448 to  2.381 23.178
%
% Fig POLYLINE object
%
\linethickness= 0.500pt
\setplotsymbol ({\thinlinefont .})
\putrule from  2.223 24.289 to  3.334 24.289
%
% Fig POLYLINE object
%
\linethickness= 0.500pt
\setplotsymbol ({\thinlinefont .})
\putrule from  3.016 23.654 to  4.604 23.654
%
% Fig POLYLINE object
%
\linethickness= 0.500pt
\setplotsymbol ({\thinlinefont .})
\putrule from  3.810 24.448 to  3.810 23.178
%
% Fig POLYLINE object
%
\linethickness= 0.500pt
\setplotsymbol ({\thinlinefont .})
\putrule from  4.286 25.082 to  4.286 23.178
%
% Fig POLYLINE object
%
\linethickness= 0.500pt
\setplotsymbol ({\thinlinefont .})
\putrule from  1.905 24.289 to  1.905 23.178
%
% Fig POLYLINE object
%
\linethickness= 0.500pt
\setplotsymbol ({\thinlinefont .})
\putrule from  1.429 25.082 to  1.429 23.178
%
% Fig POLYLINE object
%
\linethickness= 0.500pt
\setplotsymbol ({\thinlinefont .})
\setdashes < 0.1270cm>
\plot  3.651 24.924  5.397 24.924 /
%
% Fig POLYLINE object
%
\linethickness= 0.500pt
\setplotsymbol ({\thinlinefont .})
\plot  0.635 24.924  2.064 24.924 /
%
% Fig POLYLINE object
%
\linethickness= 0.500pt
\setplotsymbol ({\thinlinefont .})
\setsolid
\putrule from  6.350 23.654 to  7.303 23.654
%
% arrow head
%
\plot  7.048 23.590  7.303 23.654  7.048 23.717 /
%
%
% Fig POLYLINE object
%
\linethickness= 0.500pt
\setplotsymbol ({\thinlinefont .})
\putrule from  8.572 23.654 to 10.478 23.654
%
% Fig POLYLINE object
%
\linethickness= 0.500pt
\setplotsymbol ({\thinlinefont .})
\putrule from 10.795 23.654 to 12.541 23.654
%
% Fig POLYLINE object
%
\linethickness= 0.500pt
\setplotsymbol ({\thinlinefont .})
\putrule from 10.160 24.448 to 10.160 23.178
%
% Fig POLYLINE object
%
\linethickness= 0.500pt
\setplotsymbol ({\thinlinefont .})
\putrule from 10.954 24.448 to 10.954 23.178
%
% Fig POLYLINE object
%
\linethickness= 0.500pt
\setplotsymbol ({\thinlinefont .})
\putrule from  9.842 24.289 to 11.113 24.289
%
% Fig POLYLINE object
%
\linethickness= 0.500pt
\setplotsymbol ({\thinlinefont .})
\putrule from 11.748 24.448 to 11.748 23.178
%
% Fig POLYLINE object
%
\linethickness= 0.500pt
\setplotsymbol ({\thinlinefont .})
\plot 11.748 23.178 11.748 23.178 /
%
% Fig POLYLINE object
%
\linethickness= 0.500pt
\setplotsymbol ({\thinlinefont .})
\putrule from  9.366 24.289 to  9.366 23.178
%
% Fig POLYLINE object
%
\linethickness= 0.500pt
\setplotsymbol ({\thinlinefont .})
\putrule from 13.494 23.654 to 14.605 23.654
%
% arrow head
%
\plot 14.351 23.590 14.605 23.654 14.351 23.717 /
%
%
% Fig POLYLINE object
%
\linethickness= 0.500pt
\setplotsymbol ({\thinlinefont .})
\putrule from 17.780 24.448 to 17.780 23.178
%
% Fig POLYLINE object
%
\linethickness= 0.500pt
\setplotsymbol ({\thinlinefont .})
\putrule from 18.574 24.448 to 18.574 23.178
%
% Fig POLYLINE object
%
\linethickness= 0.500pt
\setplotsymbol ({\thinlinefont .})
\putrule from 17.462 24.289 to 18.891 24.289
%
% Fig POLYLINE object
%
\linethickness= 0.500pt
\setplotsymbol ({\thinlinefont .})
\putrule from 17.145 24.289 to 17.145 23.178
%
% Fig POLYLINE object
%
\linethickness= 0.500pt
\setplotsymbol ({\thinlinefont .})
\putrule from 10.636 22.384 to 10.636 21.431
%
% arrow head
%
\plot 10.573 21.685 10.636 21.431 10.700 21.685 /
%
%
% Fig POLYLINE object
%
\linethickness= 0.500pt
\setplotsymbol ({\thinlinefont .})
\putrule from 17.939 22.384 to 17.939 21.590
%
% arrow head
%
\plot 17.875 21.844 17.939 21.590 18.002 21.844 /
%
%
% Fig POLYLINE object
%
\linethickness= 0.500pt
\setplotsymbol ({\thinlinefont .})
\putrule from  8.731 20.003 to 10.319 20.003
%
% Fig POLYLINE object
%
\linethickness= 0.500pt
\setplotsymbol ({\thinlinefont .})
\putrule from 10.636 20.003 to 12.541 20.003
%
% Fig POLYLINE object
%
\linethickness= 0.500pt
\setplotsymbol ({\thinlinefont .})
\putrule from 10.954 20.637 to 10.954 19.685
%
% Fig POLYLINE object
%
\linethickness= 0.500pt
\setplotsymbol ({\thinlinefont .})
\putrule from 10.001 20.637 to 10.001 19.685
%
% Fig POLYLINE object
%
\linethickness= 0.500pt
\setplotsymbol ({\thinlinefont .})
\putrule from  9.842 20.479 to 11.113 20.479
%
% Fig POLYLINE object
%
\linethickness= 0.500pt
\setplotsymbol ({\thinlinefont .})
\putrule from  9.366 20.637 to  9.366 19.685
%
% Fig POLYLINE object
%
\linethickness= 0.500pt
\setplotsymbol ({\thinlinefont .})
\plot  9.366 19.685  9.366 19.685 /
%
% Fig POLYLINE object
%
\linethickness= 0.500pt
\setplotsymbol ({\thinlinefont .})
\putrule from 11.589 20.637 to 11.589 19.685
%
% Fig POLYLINE object
%
\linethickness= 0.500pt
\setplotsymbol ({\thinlinefont .})
\setdashes < 0.1270cm>
\plot  8.890 21.114  8.890 19.526 /
%
% Fig POLYLINE object
%
\linethickness= 0.500pt
\setplotsymbol ({\thinlinefont .})
\plot 12.224 21.273 12.224 19.526 /
%
% Fig POLYLINE object
%
\linethickness= 0.500pt
\setplotsymbol ({\thinlinefont .})
\setsolid
\putrule from 17.621 20.637 to 17.621 19.526
%
% Fig POLYLINE object
%
\linethickness= 0.500pt
\setplotsymbol ({\thinlinefont .})
\putrule from 16.986 20.479 to 16.986 19.526
%
% Fig POLYLINE object
%
\linethickness= 0.500pt
\setplotsymbol ({\thinlinefont .})
\putrule from 17.462 20.479 to 18.574 20.479
%
% Fig POLYLINE object
%
\linethickness= 0.500pt
\setplotsymbol ({\thinlinefont .})
\putrule from 19.367 24.448 to 19.367 23.178
%
% Fig POLYLINE object
%
\linethickness= 0.500pt
\setplotsymbol ({\thinlinefont .})
\putrule from 18.256 20.637 to 18.256 19.526
%
% Fig POLYLINE object
%
\linethickness= 0.500pt
\setplotsymbol ({\thinlinefont .})
\putrule from 13.970 20.003 to 15.081 20.003
%
% arrow head
%
\plot 14.827 19.939 15.081 20.003 14.827 20.066 /
%
%
% Fig POLYLINE object
%
\linethickness= 0.500pt
\setplotsymbol ({\thinlinefont .})
\putrule from 20.479 19.844 to 21.114 19.844
%
% arrow head
%
\plot 20.860 19.780 21.114 19.844 20.860 19.907 /
%
%
% Fig POLYLINE object
%
\linethickness= 0.500pt
\setplotsymbol ({\thinlinefont .})
\setdashes < 0.1270cm>
\plot 21.590 19.844 22.701 19.844 /
%
% Fig POLYLINE object
%
\linethickness= 0.500pt
\setplotsymbol ({\thinlinefont .})
\plot 23.019 19.844 24.130 19.844 /
%
% Fig POLYLINE object
%
\linethickness= 0.500pt
\setplotsymbol ({\thinlinefont .})
\setsolid
\putrule from 21.749 20.320 to 21.749 19.367
%
% Fig POLYLINE object
%
\linethickness= 0.500pt
\setplotsymbol ({\thinlinefont .})
\putrule from 22.384 19.367 to 22.384 20.637
%
% Fig POLYLINE object
%
\linethickness= 0.500pt
\setplotsymbol ({\thinlinefont .})
\putrule from 23.812 20.320 to 23.812 19.367
%
% Fig POLYLINE object
%
\linethickness= 0.500pt
\setplotsymbol ({\thinlinefont .})
\putrule from 23.336 20.637 to 23.336 19.367
%
% Fig POLYLINE object
%
\linethickness= 0.500pt
\setplotsymbol ({\thinlinefont .})
\putrule from 22.066 20.479 to 23.495 20.479
%
% Fig POLYLINE object
%
\linethickness= 0.500pt
\setplotsymbol ({\thinlinefont .})
\setdashes < 0.1270cm>
\plot  8.890 24.765  8.890 23.178 /
%
% Fig POLYLINE object
%
\linethickness= 0.500pt
\setplotsymbol ({\thinlinefont .})
\plot 11.906 24.606 12.859 24.606 /
%
% Fig POLYLINE object
%
\linethickness= 0.500pt
\setplotsymbol ({\thinlinefont .})
\setsolid
\putrule from 12.224 24.765 to 12.224 23.178
%
% Fig POLYLINE object
%
\linethickness= 0.500pt
\setplotsymbol ({\thinlinefont .})
\setdashes < 0.1270cm>
\plot 16.669 23.654 18.098 23.654 /
%
% Fig POLYLINE object
%
\linethickness= 0.500pt
\setplotsymbol ({\thinlinefont .})
\plot 19.209 24.289 20.003 24.289 /
%
% Fig POLYLINE object
%
\linethickness= 0.500pt
\setplotsymbol ({\thinlinefont .})
\setsolid
\putrule from 18.415 23.654 to 19.685 23.654
%
% Fig POLYLINE object
%
\linethickness= 0.500pt
\setplotsymbol ({\thinlinefont .})
\putrule from 19.050 24.289 to 19.050 23.178
%
% Fig POLYLINE object
%
\linethickness= 0.500pt
\setplotsymbol ({\thinlinefont .})
\setdashes < 0.1270cm>
\plot 16.669 20.003 17.939 20.003 /
%
% Fig POLYLINE object
%
\linethickness= 0.500pt
\setplotsymbol ({\thinlinefont .})
\plot 19.209 20.637 19.209 19.526 /
%
% Fig POLYLINE object
%
\linethickness= 0.500pt
\setplotsymbol ({\thinlinefont .})
\setsolid
\putrule from 18.891 20.479 to 18.891 19.526
%
% Fig POLYLINE object
%
\linethickness= 0.500pt
\setplotsymbol ({\thinlinefont .})
\putrule from 18.098 20.003 to 19.526 20.003
%
% Fig TEXT object
%
\put{\SetFigFont{7}{8.4}{rm}$(I_1,I_1,I_4^*)$} [lB] at  1.587 25.718
%
% Fig TEXT object
%
\put{\SetFigFont{7}{8.4}{rm}$U_1^8(D_8)$} [lB] at  9.684 25.559
%
% Fig TEXT object
%
\put{\SetFigFont{7}{8.4}{rm}$U_2^8(A_1+D_6)$} [lB] at 16.669 25.559
%
% Fig TEXT object
%
\put{\SetFigFont{7}{8.4}{rm}$U_3^8(A_1+A_5)$} [lB] at 17.304 18.733
%
% Fig TEXT object
%
\put{\SetFigFont{7}{8.4}{rm}$U_4^8(2A_1+A_3)$} [lB] at 22.066 18.733
%
% Fig TEXT object
%
\put{\SetFigFont{7}{8.4}{rm}$U_2^8(A_7)$} [lB] at  9.842 18.891
%
% Fig TEXT object
%
\put{\SetFigFont{7}{8.4}{rm}$E_0$} [lB] at  0.635 24.924
%
% Fig TEXT object
%
\put{\SetFigFont{7}{8.4}{rm}$E_1$} [lB] at  4.604 24.924
%
% Fig TEXT object
%
\put{\SetFigFont{7}{8.4}{rm}$E_1'$} [lB] at 12.383 21.114
%
% Fig TEXT object
%
\put{\SetFigFont{7}{8.4}{rm}$E_0'$} [lB] at  8.414 20.796
\linethickness=0pt
\putrectangle corners at  0.610 26.130 and 24.155 18.733
\endpicture}
}
   \vskip 0.5cm
   \begin{center}
   Figure $1$
   \end{center}
   \vskip 0.5cm
   \centerline{\font\thinlinefont=cmr5
\begingroup\makeatletter\ifx\SetFigFont\undefined%
\gdef\SetFigFont#1#2#3#4#5{%
  \reset@font\fontsize{#1}{#2pt}%
  \fontfamily{#3}\fontseries{#4}\fontshape{#5}%
  \selectfont}%
\fi\endgroup%
\mbox{\beginpicture
\setcoordinatesystem units <0.56000cm,0.56000cm>
\unitlength=0.56000cm
\linethickness=1pt
\setplotsymbol ({\makebox(0,0)[l]{\tencirc\symbol{'160}}})
\setshadesymbol ({\thinlinefont .})
\setlinear
%
% Fig POLYLINE object
%
\linethickness= 0.500pt
\setplotsymbol ({\thinlinefont .})
\putrule from  1.905 23.495 to  3.969 23.495
%
% Fig POLYLINE object
%
\linethickness= 0.500pt
\setplotsymbol ({\thinlinefont .})
\plot  2.699 24.448  1.746 22.701 /
%
% Fig POLYLINE object
%
\linethickness= 0.500pt
\setplotsymbol ({\thinlinefont .})
\plot  1.587 23.178  2.381 22.543 /
%
% Fig POLYLINE object
%
\linethickness= 0.500pt
\setplotsymbol ({\thinlinefont .})
\plot  2.381 22.860  1.905 22.225 /
%
% Fig POLYLINE object
%
\linethickness= 0.500pt
\setplotsymbol ({\thinlinefont .})
\plot  1.905 22.543  2.699 21.749 /
%
% Fig POLYLINE object
%
\linethickness= 0.500pt
\setplotsymbol ({\thinlinefont .})
\plot  2.857 22.066  2.223 21.431 /
%
% Fig POLYLINE object
%
\linethickness= 0.500pt
\setplotsymbol ({\thinlinefont .})
\plot  2.223 21.749  3.175 20.796 /
%
% Fig POLYLINE object
%
\linethickness= 0.500pt
\setplotsymbol ({\thinlinefont .})
\plot  4.286 22.066  2.699 20.796 /
%
% Fig POLYLINE object
%
\linethickness= 0.500pt
\setplotsymbol ({\thinlinefont .})
\plot  3.175 23.812  4.128 21.431 /
%
% Fig POLYLINE object
%
\linethickness= 0.500pt
\setplotsymbol ({\thinlinefont .})
\setdashes < 0.1270cm>
\plot  1.905 24.130  3.175 24.130 /
%
% Fig POLYLINE object
%
\linethickness= 0.500pt
\setplotsymbol ({\thinlinefont .})
\plot  2.699 22.543  1.905 21.749 /
%
% Fig POLYLINE object
%
\linethickness= 0.500pt
\setplotsymbol ({\thinlinefont .})
\plot  3.175 21.907  3.810 20.796 /
%
% Fig POLYLINE object
%
\linethickness= 0.500pt
\setplotsymbol ({\thinlinefont .})
\setsolid
\putrule from  5.715 22.701 to  6.826 22.701
%
% arrow head
%
\plot  6.572 22.638  6.826 22.701  6.572 22.765 /
%
%
% Fig POLYLINE object
%
\linethickness= 0.500pt
\setplotsymbol ({\thinlinefont .})
\putrule from  8.414 23.654 to 10.160 23.654
%
% Fig POLYLINE object
%
\linethickness= 0.500pt
\setplotsymbol ({\thinlinefont .})
\plot  7.938 23.178  8.890 22.543 /
%
% Fig POLYLINE object
%
\linethickness= 0.500pt
\setplotsymbol ({\thinlinefont .})
\plot  8.890 23.019  8.255 22.066 /
%
% Fig POLYLINE object
%
\linethickness= 0.500pt
\setplotsymbol ({\thinlinefont .})
\plot  8.255 22.384  9.207 21.590 /
%
% Fig POLYLINE object
%
\linethickness= 0.500pt
\setplotsymbol ({\thinlinefont .})
\plot  9.207 22.066  8.572 21.273 /
%
% Fig POLYLINE object
%
\linethickness= 0.500pt
\setplotsymbol ({\thinlinefont .})
\plot  8.572 21.590  9.525 20.796 /
%
% Fig POLYLINE object
%
\linethickness= 0.500pt
\setplotsymbol ({\thinlinefont .})
\plot 10.636 22.066  9.049 20.796 /
%
% Fig POLYLINE object
%
\linethickness= 0.500pt
\setplotsymbol ({\thinlinefont .})
\plot  9.525 23.971 10.795 21.273 /
%
% Fig POLYLINE object
%
\linethickness= 0.500pt
\setplotsymbol ({\thinlinefont .})
\setdashes < 0.1270cm>
\plot  9.049 24.130  7.938 22.543 /
%
% Fig POLYLINE object
%
\linethickness= 0.500pt
\setplotsymbol ({\thinlinefont .})
\plot  9.049 22.384  8.255 21.590 /
%
% Fig POLYLINE object
%
\linethickness= 0.500pt
\setplotsymbol ({\thinlinefont .})
\plot  9.684 21.907 10.160 20.796 /
%
% Fig POLYLINE object
%
\linethickness= 0.500pt
\setplotsymbol ({\thinlinefont .})
\setsolid
\putrule from 12.065 22.701 to 13.494 22.701
%
% arrow head
%
\plot 13.240 22.638 13.494 22.701 13.240 22.765 /
%
%
% Fig POLYLINE object
%
\linethickness= 0.500pt
\setplotsymbol ({\thinlinefont .})
\putrule from 14.446 23.812 to 16.192 23.812
%
% Fig POLYLINE object
%
\linethickness= 0.500pt
\setplotsymbol ({\thinlinefont .})
\plot 14.287 23.178 15.240 22.701 /
%
% Fig POLYLINE object
%
\linethickness= 0.500pt
\setplotsymbol ({\thinlinefont .})
\plot 15.081 23.336 14.446 22.225 /
%
% Fig POLYLINE object
%
\linethickness= 0.500pt
\setplotsymbol ({\thinlinefont .})
\plot 14.446 22.225 14.446 22.225 /
%
% Fig POLYLINE object
%
\linethickness= 0.500pt
\setplotsymbol ({\thinlinefont .})
\plot 15.558 22.225 14.922 21.273 /
%
% Fig POLYLINE object
%
\linethickness= 0.500pt
\setplotsymbol ({\thinlinefont .})
\plot 14.922 21.590 15.875 20.796 /
%
% Fig POLYLINE object
%
\linethickness= 0.500pt
\setplotsymbol ({\thinlinefont .})
\plot 16.986 22.225 15.240 20.796 /
%
% Fig POLYLINE object
%
\linethickness= 0.500pt
\setplotsymbol ({\thinlinefont .})
\plot 15.558 24.289 17.145 21.590 /
%
% Fig POLYLINE object
%
\linethickness= 0.500pt
\setplotsymbol ({\thinlinefont .})
\setdashes < 0.1270cm>
\plot 15.081 24.448 14.287 22.701 /
%
% Fig POLYLINE object
%
\linethickness= 0.500pt
\setplotsymbol ({\thinlinefont .})
\plot 14.605 22.701 15.716 21.749 /
%
% Fig POLYLINE object
%
\linethickness= 0.500pt
\setplotsymbol ({\thinlinefont .})
\plot 16.034 22.066 16.828 20.955 /
%
% Fig POLYLINE object
%
\linethickness= 0.500pt
\setplotsymbol ({\thinlinefont .})
\setsolid
\putrule from 18.256 22.701 to 19.526 22.701
%
% arrow head
%
\plot 19.272 22.638 19.526 22.701 19.272 22.765 /
%
%
% Fig POLYLINE object
%
\linethickness= 0.500pt
\setplotsymbol ({\thinlinefont .})
\putrule from 21.273 23.971 to 23.336 23.971
%
% Fig POLYLINE object
%
\linethickness= 0.500pt
\setplotsymbol ({\thinlinefont .})
\plot 22.384 24.606 23.654 21.907 /
%
% Fig POLYLINE object
%
\linethickness= 0.500pt
\setplotsymbol ({\thinlinefont .})
\plot 21.749 21.907 22.701 21.273 /
%
% Fig POLYLINE object
%
\linethickness= 0.500pt
\setplotsymbol ({\thinlinefont .})
\plot 22.384 22.701 21.749 21.273 /
%
% Fig POLYLINE object
%
\linethickness= 0.500pt
\setplotsymbol ({\thinlinefont .})
\plot 21.749 21.431 21.749 21.431 /
%
% Fig POLYLINE object
%
\linethickness= 0.500pt
\setplotsymbol ({\thinlinefont .})
\plot 21.273 23.495 21.907 23.178 /
%
% Fig POLYLINE object
%
\linethickness= 0.500pt
\setplotsymbol ({\thinlinefont .})
\plot 21.749 23.495 21.431 22.543 /
%
% Fig POLYLINE object
%
\linethickness= 0.500pt
\setplotsymbol ({\thinlinefont .})
\setdashes < 0.1270cm>
\plot 21.590 24.448 21.273 23.178 /
%
% Fig POLYLINE object
%
\linethickness= 0.500pt
\setplotsymbol ({\thinlinefont .})
\plot 21.431 22.860 22.543 22.066 /
%
% Fig POLYLINE object
%
\linethickness= 0.500pt
\setplotsymbol ({\thinlinefont .})
\plot 23.654 22.701 22.225 21.114 /
%
% Fig TEXT object
%
\put{\SetFigFont{7}{8.4}{rm}$(I_9,I_1,I_1,I_1)$} [lB] at  1.587 24.924
%
% Fig TEXT object
%
\put{\SetFigFont{7}{8.4}{rm}$U_1^{11}(A_8)$} [lB] at  8.414 24.606
%
% Fig TEXT object
%
\put{\SetFigFont{7}{8.4}{rm}$U_2^{11}(A_2+A_5)$} [lB] at 14.446 24.765
%
% Fig TEXT object
%
\put{\SetFigFont{7}{8.4}{rm}$U_3^{11}(3A_2)$} [lB] at 21.114 24.765
%
% Fig TEXT object
%
\put{\SetFigFont{7}{8.4}{rm}$E_0$} [lB] at  1.746 24.130
%
% Fig TEXT object
%
\put{\SetFigFont{7}{8.4}{rm}$E_1$} [lB] at  1.429 21.431
%
% Fig TEXT object
%
\put{\SetFigFont{7}{8.4}{rm}$E_2$} [lB] at  3.810 20.479
\linethickness=0pt
\putrectangle corners at  1.429 25.337 and 23.679 20.479
\endpicture}
}
   \vskip 0.5cm
   \centerline{\font\thinlinefont=cmr5
\begingroup\makeatletter\ifx\SetFigFont\undefined%
\gdef\SetFigFont#1#2#3#4#5{%
  \reset@font\fontsize{#1}{#2pt}%
  \fontfamily{#3}\fontseries{#4}\fontshape{#5}%
  \selectfont}%
\fi\endgroup%
\mbox{\beginpicture
\setcoordinatesystem units <0.56000cm,0.56000cm>
\unitlength=0.56000cm
\linethickness=1pt
\setplotsymbol ({\makebox(0,0)[l]{\tencirc\symbol{'160}}})
\setshadesymbol ({\thinlinefont .})
\setlinear
%
% Fig POLYLINE object
%
\linethickness= 0.500pt
\setplotsymbol ({\thinlinefont .})
\putrule from  4.128 23.178 to  6.668 23.178
%
% Fig POLYLINE object
%
\linethickness= 0.500pt
\setplotsymbol ({\thinlinefont .})
\putrule from  4.604 23.812 to  4.604 22.543
%
% Fig POLYLINE object
%
\linethickness= 0.500pt
\setplotsymbol ({\thinlinefont .})
\putrule from  5.397 24.448 to  5.397 22.701
%
% Fig POLYLINE object
%
\linethickness= 0.500pt
\setplotsymbol ({\thinlinefont .})
\putrule from  6.032 23.812 to  6.032 21.907
%
% Fig POLYLINE object
%
\linethickness= 0.500pt
\setplotsymbol ({\thinlinefont .})
\plot  3.969 24.289  4.763 23.495 /
%
% Fig POLYLINE object
%
\linethickness= 0.500pt
\setplotsymbol ({\thinlinefont .})
\plot  3.810 24.924  5.556 23.495 /
%
% Fig POLYLINE object
%
\linethickness= 0.500pt
\setplotsymbol ({\thinlinefont .})
\plot  6.509 22.860  4.763 20.955 /
%
% Fig POLYLINE object
%
\linethickness= 0.500pt
\setplotsymbol ({\thinlinefont .})
\plot  3.175 24.130  1.587 22.225 /
%
% Fig POLYLINE object
%
\linethickness= 0.500pt
\setplotsymbol ({\thinlinefont .})
\plot  2.540 25.082  2.223 22.066 /
%
% Fig POLYLINE object
%
\linethickness= 0.500pt
\setplotsymbol ({\thinlinefont .})
\plot  1.746 23.971  3.810 20.955 /
%
% Fig POLYLINE object
%
\linethickness= 0.500pt
\setplotsymbol ({\thinlinefont .})
\setdashes < 0.1270cm>
\plot  2.699 23.971  4.445 23.971 /
%
% Fig POLYLINE object
%
\linethickness= 0.500pt
\setplotsymbol ({\thinlinefont .})
\plot  3.175 21.431  5.715 21.431 /
%
% Fig POLYLINE object
%
\linethickness= 0.500pt
\setplotsymbol ({\thinlinefont .})
\plot  2.223 24.606  4.604 24.606 /
%
% Fig POLYLINE object
%
\linethickness= 0.500pt
\setplotsymbol ({\thinlinefont .})
\setsolid
\putrule from  7.620 23.178 to  8.731 23.178
%
% arrow head
%
\plot  8.477 23.114  8.731 23.178  8.477 23.241 /
%
%
% Fig POLYLINE object
%
\linethickness= 0.500pt
\setplotsymbol ({\thinlinefont .})
\putrule from 10.636 23.178 to 13.176 23.178
%
% Fig POLYLINE object
%
\linethickness= 0.500pt
\setplotsymbol ({\thinlinefont .})
\putrule from 11.271 23.971 to 11.271 22.701
%
% Fig POLYLINE object
%
\linethickness= 0.500pt
\setplotsymbol ({\thinlinefont .})
\putrule from 11.906 24.606 to 11.906 22.701
%
% Fig POLYLINE object
%
\linethickness= 0.500pt
\setplotsymbol ({\thinlinefont .})
\putrule from 12.700 23.812 to 12.700 21.907
%
% Fig POLYLINE object
%
\linethickness= 0.500pt
\setplotsymbol ({\thinlinefont .})
\plot 10.636 24.289 11.430 23.495 /
%
% Fig POLYLINE object
%
\linethickness= 0.500pt
\setplotsymbol ({\thinlinefont .})
\plot 10.636 25.082 12.224 23.654 /
%
% Fig POLYLINE object
%
\linethickness= 0.500pt
\setplotsymbol ({\thinlinefont .})
\plot 10.001 24.289  9.049 22.225 /
%
% Fig POLYLINE object
%
\linethickness= 0.500pt
\setplotsymbol ({\thinlinefont .})
\putrule from  9.525 25.082 to  9.525 22.384
%
% Fig POLYLINE object
%
\linethickness= 0.500pt
\setplotsymbol ({\thinlinefont .})
\setdashes < 0.1270cm>
\plot  9.684 24.130 11.113 24.130 /
%
% Fig POLYLINE object
%
\linethickness= 0.500pt
\setplotsymbol ({\thinlinefont .})
\plot  9.207 24.765 11.430 24.765 /
%
% Fig POLYLINE object
%
\linethickness= 0.500pt
\setplotsymbol ({\thinlinefont .})
\plot  8.890 23.971 11.430 21.273 /
%
% Fig POLYLINE object
%
\linethickness= 0.500pt
\setplotsymbol ({\thinlinefont .})
\plot 10.636 21.749 13.176 22.543 /
%
% Fig POLYLINE object
%
\linethickness= 0.500pt
\setplotsymbol ({\thinlinefont .})
\setsolid
\putrule from 14.446 23.178 to 15.875 23.178
%
% arrow head
%
\plot 15.621 23.114 15.875 23.178 15.621 23.241 /
%
%
% Fig POLYLINE object
%
\linethickness= 0.500pt
\setplotsymbol ({\thinlinefont .})
\putrule from 18.256 23.178 to 20.637 23.178
%
% Fig POLYLINE object
%
\linethickness= 0.500pt
\setplotsymbol ({\thinlinefont .})
\putrule from 19.526 24.448 to 19.526 22.701
%
% Fig POLYLINE object
%
\linethickness= 0.500pt
\setplotsymbol ({\thinlinefont .})
\putrule from 18.733 23.971 to 18.733 22.701
%
% Fig POLYLINE object
%
\linethickness= 0.500pt
\setplotsymbol ({\thinlinefont .})
\plot 18.098 24.289 19.050 23.495 /
%
% Fig POLYLINE object
%
\linethickness= 0.500pt
\setplotsymbol ({\thinlinefont .})
\plot 18.098 25.082 19.844 23.654 /
%
% Fig POLYLINE object
%
\linethickness= 0.500pt
\setplotsymbol ({\thinlinefont .})
\plot 17.780 24.289 16.986 22.701 /
%
% Fig POLYLINE object
%
\linethickness= 0.500pt
\setplotsymbol ({\thinlinefont .})
\putrule from 17.304 25.082 to 17.304 22.543
%
% Fig POLYLINE object
%
\linethickness= 0.500pt
\setplotsymbol ({\thinlinefont .})
\setdashes < 0.1270cm>
\plot 16.986 24.924 18.733 24.924 /
%
% Fig POLYLINE object
%
\linethickness= 0.500pt
\setplotsymbol ({\thinlinefont .})
\plot 17.462 24.130 18.733 24.130 /
%
% Fig POLYLINE object
%
\linethickness= 0.500pt
\setplotsymbol ({\thinlinefont .})
\plot 20.320 24.130 20.320 22.066 /
%
% Fig POLYLINE object
%
\linethickness= 0.500pt
\setplotsymbol ({\thinlinefont .})
\setdots < 0.0953cm>
\plot 16.828 23.812 16.830 23.810 /
\plot 16.830 23.810 16.836 23.806 /
\plot 16.836 23.806 16.847 23.796 /
\plot 16.847 23.796 16.863 23.783 /
\plot 16.863 23.783 16.889 23.762 /
\plot 16.889 23.762 16.921 23.736 /
\plot 16.921 23.736 16.963 23.702 /
\plot 16.963 23.702 17.012 23.660 /
\plot 17.012 23.660 17.071 23.614 /
\plot 17.071 23.614 17.137 23.558 /
\plot 17.137 23.558 17.211 23.499 /
\plot 17.211 23.499 17.293 23.434 /
\plot 17.293 23.434 17.380 23.364 /
\plot 17.380 23.364 17.473 23.290 /
\plot 17.473 23.290 17.568 23.213 /
\plot 17.568 23.213 17.666 23.135 /
\plot 17.666 23.135 17.765 23.057 /
\plot 17.765 23.057 17.867 22.979 /
\plot 17.867 22.979 17.966 22.902 /
\plot 17.966 22.902 18.064 22.828 /
\plot 18.064 22.828 18.159 22.756 /
\plot 18.159 22.756 18.252 22.689 /
\plot 18.252 22.689 18.343 22.623 /
\plot 18.343 22.623 18.428 22.564 /
\plot 18.428 22.564 18.510 22.507 /
\plot 18.510 22.507 18.589 22.456 /
\plot 18.589 22.456 18.665 22.409 /
\plot 18.665 22.409 18.735 22.367 /
\plot 18.735 22.367 18.802 22.331 /
\plot 18.802 22.331 18.868 22.297 /
\plot 18.868 22.297 18.931 22.269 /
\plot 18.931 22.269 18.991 22.246 /
\plot 18.991 22.246 19.050 22.225 /
\plot 19.050 22.225 19.120 22.204 /
\plot 19.120 22.204 19.190 22.189 /
\plot 19.190 22.189 19.262 22.178 /
\plot 19.262 22.178 19.334 22.174 /
\plot 19.334 22.174 19.408 22.172 /
\plot 19.408 22.172 19.484 22.176 /
\plot 19.484 22.176 19.562 22.185 /
\plot 19.562 22.185 19.645 22.197 /
\plot 19.645 22.197 19.729 22.212 /
\plot 19.729 22.212 19.816 22.231 /
\plot 19.816 22.231 19.907 22.253 /
\plot 19.907 22.253 19.998 22.278 /
\plot 19.998 22.278 20.089 22.303 /
\plot 20.089 22.303 20.180 22.331 /
\plot 20.180 22.331 20.269 22.358 /
\plot 20.269 22.358 20.354 22.386 /
\plot 20.354 22.386 20.434 22.413 /
\plot 20.434 22.413 20.508 22.439 /
\plot 20.508 22.439 20.574 22.462 /
\plot 20.574 22.462 20.631 22.483 /
\plot 20.631 22.483 20.680 22.500 /
\plot 20.680 22.500 20.720 22.515 /
\plot 20.720 22.515 20.750 22.526 /
\plot 20.750 22.526 20.771 22.534 /
\plot 20.771 22.534 20.786 22.538 /
\plot 20.786 22.538 20.792 22.540 /
\plot 20.792 22.540 20.796 22.543 /
%
% Fig TEXT object
%
\put{\SetFigFont{7}{8.4}{rm}$(IV,IV^*)$} [lB] at  2.381 25.718
%
% Fig TEXT object
%
\put{\SetFigFont{7}{8.4}{rm}$U_1^3(A_2+E_6)$} [lB] at  9.525 25.559
%
% Fig TEXT object
%
\put{\SetFigFont{7}{8.4}{rm}$U_2^3(A_2+A_5)$} [lB] at 17.145 25.559
\linethickness=0pt
\putrectangle corners at  1.562 26.130 and 20.822 20.930
\endpicture}
}
   \vskip 0.5cm
   \centerline{\font\thinlinefont=cmr5
\begingroup\makeatletter\ifx\SetFigFont\undefined%
\gdef\SetFigFont#1#2#3#4#5{%
  \reset@font\fontsize{#1}{#2pt}%
  \fontfamily{#3}\fontseries{#4}\fontshape{#5}%
  \selectfont}%
\fi\endgroup%
\mbox{\beginpicture
\setcoordinatesystem units <0.56000cm,0.56000cm>
\unitlength=0.56000cm
\linethickness=1pt
\setplotsymbol ({\makebox(0,0)[l]{\tencirc\symbol{'160}}})
\setshadesymbol ({\thinlinefont .})
\setlinear
%
% Fig POLYLINE object
%
\linethickness= 0.500pt
\setplotsymbol ({\thinlinefont .})
\putrule from  4.128 23.178 to  6.668 23.178
%
% Fig POLYLINE object
%
\linethickness= 0.500pt
\setplotsymbol ({\thinlinefont .})
\putrule from  4.604 23.812 to  4.604 22.543
%
% Fig POLYLINE object
%
\linethickness= 0.500pt
\setplotsymbol ({\thinlinefont .})
\putrule from  5.397 24.448 to  5.397 22.701
%
% Fig POLYLINE object
%
\linethickness= 0.500pt
\setplotsymbol ({\thinlinefont .})
\putrule from  6.032 23.812 to  6.032 21.907
%
% Fig POLYLINE object
%
\linethickness= 0.500pt
\setplotsymbol ({\thinlinefont .})
\plot  3.969 24.289  4.763 23.495 /
%
% Fig POLYLINE object
%
\linethickness= 0.500pt
\setplotsymbol ({\thinlinefont .})
\plot  3.810 24.924  5.556 23.495 /
%
% Fig POLYLINE object
%
\linethickness= 0.500pt
\setplotsymbol ({\thinlinefont .})
\plot  6.509 22.860  4.763 20.955 /
%
% Fig POLYLINE object
%
\linethickness= 0.500pt
\setplotsymbol ({\thinlinefont .})
\plot  3.175 24.130  1.587 22.225 /
%
% Fig POLYLINE object
%
\linethickness= 0.500pt
\setplotsymbol ({\thinlinefont .})
\plot  2.540 25.082  2.223 22.066 /
%
% Fig POLYLINE object
%
\linethickness= 0.500pt
\setplotsymbol ({\thinlinefont .})
\setdashes < 0.1270cm>
\plot  2.699 23.971  4.445 23.971 /
%
% Fig POLYLINE object
%
\linethickness= 0.500pt
\setplotsymbol ({\thinlinefont .})
\plot  3.175 21.431  5.715 21.431 /
%
% Fig POLYLINE object
%
\linethickness= 0.500pt
\setplotsymbol ({\thinlinefont .})
\plot  2.223 24.606  4.604 24.606 /
%
% Fig POLYLINE object
%
\linethickness= 0.500pt
\setplotsymbol ({\thinlinefont .})
\setsolid
\putrule from  7.620 23.178 to  8.731 23.178
%
% arrow head
%
\plot  8.477 23.114  8.731 23.178  8.477 23.241 /
%
%
% Fig POLYLINE object
%
\linethickness= 0.500pt
\setplotsymbol ({\thinlinefont .})
\putrule from 10.636 23.178 to 13.176 23.178
%
% Fig POLYLINE object
%
\linethickness= 0.500pt
\setplotsymbol ({\thinlinefont .})
\putrule from 11.271 23.971 to 11.271 22.701
%
% Fig POLYLINE object
%
\linethickness= 0.500pt
\setplotsymbol ({\thinlinefont .})
\putrule from 11.906 24.606 to 11.906 22.701
%
% Fig POLYLINE object
%
\linethickness= 0.500pt
\setplotsymbol ({\thinlinefont .})
\putrule from 12.700 23.812 to 12.700 21.907
%
% Fig POLYLINE object
%
\linethickness= 0.500pt
\setplotsymbol ({\thinlinefont .})
\plot 10.636 24.289 11.430 23.495 /
%
% Fig POLYLINE object
%
\linethickness= 0.500pt
\setplotsymbol ({\thinlinefont .})
\plot 10.636 25.082 12.224 23.654 /
%
% Fig POLYLINE object
%
\linethickness= 0.500pt
\setplotsymbol ({\thinlinefont .})
\plot 10.001 24.289  9.049 22.225 /
%
% Fig POLYLINE object
%
\linethickness= 0.500pt
\setplotsymbol ({\thinlinefont .})
\putrule from  9.525 25.082 to  9.525 22.384
%
% Fig POLYLINE object
%
\linethickness= 0.500pt
\setplotsymbol ({\thinlinefont .})
\setdashes < 0.1270cm>
\plot  9.684 24.130 11.113 24.130 /
%
% Fig POLYLINE object
%
\linethickness= 0.500pt
\setplotsymbol ({\thinlinefont .})
\plot  9.207 24.765 11.430 24.765 /
%
% Fig POLYLINE object
%
\linethickness= 0.500pt
\setplotsymbol ({\thinlinefont .})
\plot 10.636 21.749 13.176 22.543 /
%
% Fig POLYLINE object
%
\linethickness= 0.500pt
\setplotsymbol ({\thinlinefont .})
\setsolid
\putrule from 14.446 23.178 to 15.875 23.178
%
% arrow head
%
\plot 15.621 23.114 15.875 23.178 15.621 23.241 /
%
%
% Fig POLYLINE object
%
\linethickness= 0.500pt
\setplotsymbol ({\thinlinefont .})
\putrule from 18.256 23.178 to 20.637 23.178
%
% Fig POLYLINE object
%
\linethickness= 0.500pt
\setplotsymbol ({\thinlinefont .})
\putrule from 19.526 24.448 to 19.526 22.701
%
% Fig POLYLINE object
%
\linethickness= 0.500pt
\setplotsymbol ({\thinlinefont .})
\putrule from 18.733 23.971 to 18.733 22.701
%
% Fig POLYLINE object
%
\linethickness= 0.500pt
\setplotsymbol ({\thinlinefont .})
\plot 18.098 24.289 19.050 23.495 /
%
% Fig POLYLINE object
%
\linethickness= 0.500pt
\setplotsymbol ({\thinlinefont .})
\plot 18.098 25.082 19.844 23.654 /
%
% Fig POLYLINE object
%
\linethickness= 0.500pt
\setplotsymbol ({\thinlinefont .})
\plot 17.780 24.289 16.986 22.701 /
%
% Fig POLYLINE object
%
\linethickness= 0.500pt
\setplotsymbol ({\thinlinefont .})
\putrule from 17.304 25.082 to 17.304 22.543
%
% Fig POLYLINE object
%
\linethickness= 0.500pt
\setplotsymbol ({\thinlinefont .})
\setdashes < 0.1270cm>
\plot 16.986 24.924 18.733 24.924 /
%
% Fig POLYLINE object
%
\linethickness= 0.500pt
\setplotsymbol ({\thinlinefont .})
\plot 17.462 24.130 18.733 24.130 /
%
% Fig POLYLINE object
%
\linethickness= 0.500pt
\setplotsymbol ({\thinlinefont .})
\plot 20.320 24.130 20.320 22.066 /
%
% Fig POLYLINE object
%
\linethickness= 0.500pt
\setplotsymbol ({\thinlinefont .})
\setsolid
\plot  2.064 24.130  4.128 20.796 /
%
% Fig POLYLINE object
%
\linethickness= 0.500pt
\setplotsymbol ({\thinlinefont .})
\setdashes < 0.1270cm>
\plot  8.890 24.130 12.065 21.749 /
%
% Fig POLYLINE object
%
\linethickness= 0.500pt
\setplotsymbol ({\thinlinefont .})
\setdots < 0.0953cm>
\plot 16.986 24.130 16.988 24.128 /
\plot 16.988 24.128 16.993 24.122 /
\plot 16.993 24.122 17.001 24.111 /
\plot 17.001 24.111 17.014 24.096 /
\plot 17.014 24.096 17.033 24.073 /
\plot 17.033 24.073 17.058 24.043 /
\plot 17.058 24.043 17.092 24.005 /
\plot 17.092 24.005 17.130 23.959 /
\plot 17.130 23.959 17.177 23.904 /
\plot 17.177 23.904 17.230 23.840 /
\plot 17.230 23.840 17.289 23.770 /
\plot 17.289 23.770 17.357 23.694 /
\plot 17.357 23.694 17.427 23.609 /
\plot 17.427 23.609 17.503 23.523 /
\plot 17.503 23.523 17.583 23.429 /
\plot 17.583 23.429 17.666 23.334 /
\plot 17.666 23.334 17.750 23.239 /
\plot 17.750 23.239 17.837 23.142 /
\plot 17.837 23.142 17.924 23.044 /
\plot 17.924 23.044 18.011 22.949 /
\plot 18.011 22.949 18.095 22.856 /
\plot 18.095 22.856 18.180 22.767 /
\plot 18.180 22.767 18.263 22.682 /
\plot 18.263 22.682 18.341 22.600 /
\plot 18.341 22.600 18.419 22.523 /
\plot 18.419 22.523 18.493 22.454 /
\plot 18.493 22.454 18.563 22.388 /
\plot 18.563 22.388 18.633 22.329 /
\plot 18.633 22.329 18.699 22.274 /
\plot 18.699 22.274 18.762 22.227 /
\plot 18.762 22.227 18.821 22.183 /
\plot 18.821 22.183 18.881 22.147 /
\plot 18.881 22.147 18.938 22.115 /
\plot 18.938 22.115 18.995 22.087 /
\plot 18.995 22.087 19.050 22.066 /
\plot 19.050 22.066 19.116 22.045 /
\plot 19.116 22.045 19.183 22.030 /
\plot 19.183 22.030 19.251 22.022 /
\plot 19.251 22.022 19.321 22.020 /
\plot 19.321 22.020 19.393 22.024 /
\plot 19.393 22.024 19.467 22.035 /
\plot 19.467 22.035 19.545 22.049 /
\plot 19.545 22.049 19.626 22.070 /
\plot 19.626 22.070 19.710 22.098 /
\plot 19.710 22.098 19.799 22.128 /
\plot 19.799 22.128 19.888 22.164 /
\plot 19.888 22.164 19.981 22.202 /
\plot 19.981 22.202 20.077 22.244 /
\plot 20.077 22.244 20.170 22.288 /
\plot 20.170 22.288 20.265 22.333 /
\plot 20.265 22.333 20.356 22.380 /
\plot 20.356 22.380 20.445 22.424 /
\plot 20.445 22.424 20.532 22.468 /
\plot 20.532 22.468 20.610 22.511 /
\plot 20.610 22.511 20.682 22.549 /
\plot 20.682 22.549 20.745 22.585 /
\plot 20.745 22.585 20.800 22.614 /
\plot 20.800 22.614 20.847 22.640 /
\plot 20.847 22.640 20.885 22.661 /
\plot 20.885 22.661 20.913 22.678 /
\plot 20.913 22.678 20.932 22.689 /
\plot 20.932 22.689 20.944 22.695 /
\plot 20.944 22.695 20.953 22.699 /
\plot 20.953 22.699 20.955 22.701 /
%
% Fig TEXT object
%
\put{\SetFigFont{7}{8.4}{rm}$(I_1,I_3,IV^*)$} [lB] at  2.381 25.559
%
% Fig TEXT object
%
\put{\SetFigFont{7}{8.4}{rm}$U_1^7(A_2+E_6)$} [lB] at  9.525 25.400
%
% Fig TEXT object
%
\put{\SetFigFont{7}{8.4}{rm}$U_2^7(A_2+A_5)$} [lB] at 17.304 25.559
\linethickness=0pt
\putrectangle corners at  1.562 25.971 and 20.980 20.771
\endpicture}
}
   \vskip 0.5cm
   \begin{center}
    Figure $2$
   \end{center}
   \vskip 0.5cm
   \centerline{\font\thinlinefont=cmr5
\begingroup\makeatletter\ifx\SetFigFont\undefined%
\gdef\SetFigFont#1#2#3#4#5{%
  \reset@font\fontsize{#1}{#2pt}%
  \fontfamily{#3}\fontseries{#4}\fontshape{#5}%
  \selectfont}%
\fi\endgroup%
\mbox{\beginpicture
\setcoordinatesystem units <0.56000cm,0.56000cm>
\unitlength=0.56000cm
\linethickness=1pt
\setplotsymbol ({\makebox(0,0)[l]{\tencirc\symbol{'160}}})
\setshadesymbol ({\thinlinefont .})
\setlinear
%
% Fig POLYLINE object
%
\linethickness= 0.500pt
\setplotsymbol ({\thinlinefont .})
\putrule from  5.080 23.336 to  5.080 22.066
%
% Fig POLYLINE object
%
\linethickness= 0.500pt
\setplotsymbol ({\thinlinefont .})
\putrule from  3.493 22.701 to  5.397 22.701
%
% Fig POLYLINE object
%
\linethickness= 0.500pt
\setplotsymbol ({\thinlinefont .})
\putrule from  4.763 23.178 to  5.874 23.178
%
% Fig POLYLINE object
%
\linethickness= 0.500pt
\setplotsymbol ({\thinlinefont .})
\putrule from  5.397 23.812 to  5.397 23.019
%
% Fig POLYLINE object
%
\linethickness= 0.500pt
\setplotsymbol ({\thinlinefont .})
\putrule from  5.715 23.812 to  5.715 23.019
%
% Fig POLYLINE object
%
\linethickness= 0.500pt
\setplotsymbol ({\thinlinefont .})
\putrule from  4.286 23.336 to  4.286 22.066
%
% Fig POLYLINE object
%
\linethickness= 0.500pt
\setplotsymbol ({\thinlinefont .})
\plot  2.381 25.082  1.111 23.019 /
%
% Fig POLYLINE object
%
\linethickness= 0.500pt
\setplotsymbol ({\thinlinefont .})
\plot  6.350 25.241  8.255 23.654 /
%
% Fig POLYLINE object
%
\linethickness= 0.500pt
\setplotsymbol ({\thinlinefont .})
\putrule from  8.414 22.701 to  9.525 22.701
%
% arrow head
%
\plot  9.271 22.638  9.525 22.701  9.271 22.765 /
%
%
% Fig POLYLINE object
%
\linethickness= 0.500pt
\setplotsymbol ({\thinlinefont .})
\putrule from 13.018 23.495 to 13.018 22.225
%
% Fig POLYLINE object
%
\linethickness= 0.500pt
\setplotsymbol ({\thinlinefont .})
\putrule from 12.700 23.336 to 13.811 23.336
%
% Fig POLYLINE object
%
\linethickness= 0.500pt
\setplotsymbol ({\thinlinefont .})
\putrule from 13.335 23.971 to 13.335 23.019
%
% Fig POLYLINE object
%
\linethickness= 0.500pt
\setplotsymbol ({\thinlinefont .})
\putrule from 13.652 23.971 to 13.652 23.019
%
% Fig POLYLINE object
%
\linethickness= 0.500pt
\setplotsymbol ({\thinlinefont .})
\putrule from 12.224 23.495 to 12.224 22.225
%
% Fig POLYLINE object
%
\linethickness= 0.500pt
\setplotsymbol ({\thinlinefont .})
\setdashes < 0.1270cm>
\plot 11.430 25.400 11.589 22.066 /
%
% Fig POLYLINE object
%
\linethickness= 0.500pt
\setplotsymbol ({\thinlinefont .})
\setsolid
\putrule from 10.954 22.701 to 13.494 22.701
%
% Fig POLYLINE object
%
\linethickness= 0.500pt
\setplotsymbol ({\thinlinefont .})
\putrule from  3.810 22.066 to  3.810 25.241
%
% Fig POLYLINE object
%
\linethickness= 0.500pt
\setplotsymbol ({\thinlinefont .})
\setdashes < 0.1270cm>
\plot  2.064 24.924  7.144 24.924 /
%
% Fig POLYLINE object
%
\linethickness= 0.500pt
\setplotsymbol ({\thinlinefont .})
\plot 12.224 24.765  9.525 23.336 /
%
% Fig POLYLINE object
%
\linethickness= 0.500pt
\setplotsymbol ({\thinlinefont .})
\setsolid
\putrule from 16.510 22.701 to 18.256 22.701
%
% arrow head
%
\plot 18.002 22.638 18.256 22.701 18.002 22.765 /
%
%
% Fig POLYLINE object
%
\linethickness= 0.500pt
\setplotsymbol ({\thinlinefont .})
\setdashes < 0.1270cm>
\plot 19.685 22.701 22.701 22.701 /
%
% Fig POLYLINE object
%
\linethickness= 0.500pt
\setplotsymbol ({\thinlinefont .})
\setsolid
\putrule from 21.273 23.495 to 21.273 22.066
%
% Fig POLYLINE object
%
\linethickness= 0.500pt
\setplotsymbol ({\thinlinefont .})
\putrule from 21.907 23.495 to 21.907 22.066
%
% Fig POLYLINE object
%
\linethickness= 0.500pt
\setplotsymbol ({\thinlinefont .})
\putrule from 21.749 23.336 to 22.701 23.336
%
% Fig POLYLINE object
%
\linethickness= 0.500pt
\setplotsymbol ({\thinlinefont .})
\putrule from 22.225 24.130 to 22.225 23.178
%
% Fig POLYLINE object
%
\linethickness= 0.500pt
\setplotsymbol ({\thinlinefont .})
\putrule from 22.543 24.130 to 22.543 23.178
%
% Fig POLYLINE object
%
\linethickness= 0.500pt
\setplotsymbol ({\thinlinefont .})
\setdots < 0.0953cm>
\plot 19.050 24.606 20.161 22.066 /
%
% Fig POLYLINE object
%
\linethickness= 0.500pt
\setplotsymbol ({\thinlinefont .})
\setdashes < 0.1270cm>
\plot 11.113 24.448 16.192 23.971 /
%
% Fig POLYLINE object
%
\linethickness= 0.500pt
\setplotsymbol ({\thinlinefont .})
\setsolid
\plot  1.111 23.971  1.115 23.967 /
\plot  1.115 23.967  1.124 23.959 /
\plot  1.124 23.959  1.141 23.942 /
\plot  1.141 23.942  1.166 23.918 /
\plot  1.166 23.918  1.200 23.887 /
\plot  1.200 23.887  1.240 23.848 /
\plot  1.240 23.848  1.291 23.802 /
\plot  1.291 23.802  1.344 23.751 /
\plot  1.344 23.751  1.403 23.698 /
\plot  1.403 23.698  1.463 23.645 /
\plot  1.463 23.645  1.522 23.592 /
\plot  1.522 23.592  1.579 23.544 /
\plot  1.579 23.544  1.634 23.499 /
\plot  1.634 23.499  1.687 23.457 /
\plot  1.687 23.457  1.736 23.423 /
\plot  1.736 23.423  1.780 23.393 /
\plot  1.780 23.393  1.825 23.370 /
\plot  1.825 23.370  1.865 23.351 /
\plot  1.865 23.351  1.905 23.336 /
\plot  1.905 23.336  1.945 23.326 /
\plot  1.945 23.326  1.985 23.319 /
\plot  1.985 23.319  2.028 23.315 /
\plot  2.028 23.315  2.074 23.313 /
\putrule from  2.074 23.313 to  2.121 23.313
\putrule from  2.121 23.313 to  2.172 23.313
\plot  2.172 23.313  2.223 23.315 /
\plot  2.223 23.315  2.275 23.317 /
\plot  2.275 23.317  2.330 23.321 /
\plot  2.330 23.321  2.383 23.326 /
\plot  2.383 23.326  2.436 23.332 /
\plot  2.436 23.332  2.487 23.340 /
\plot  2.487 23.340  2.534 23.349 /
\plot  2.534 23.349  2.578 23.360 /
\plot  2.578 23.360  2.616 23.374 /
\plot  2.616 23.374  2.648 23.391 /
\plot  2.648 23.391  2.671 23.410 /
\plot  2.671 23.410  2.690 23.434 /
\plot  2.690 23.434  2.699 23.461 /
\putrule from  2.699 23.461 to  2.699 23.495
\plot  2.699 23.495  2.692 23.525 /
\plot  2.692 23.525  2.682 23.558 /
\plot  2.682 23.558  2.665 23.594 /
\plot  2.665 23.594  2.639 23.639 /
\plot  2.639 23.639  2.610 23.688 /
\plot  2.610 23.688  2.574 23.741 /
\plot  2.574 23.741  2.532 23.800 /
\plot  2.532 23.800  2.483 23.865 /
\plot  2.483 23.865  2.430 23.935 /
\plot  2.430 23.935  2.373 24.009 /
\plot  2.373 24.009  2.311 24.086 /
\plot  2.311 24.086  2.248 24.164 /
\plot  2.248 24.164  2.182 24.244 /
\plot  2.182 24.244  2.119 24.323 /
\plot  2.119 24.323  2.055 24.399 /
\plot  2.055 24.399  1.996 24.471 /
\plot  1.996 24.471  1.941 24.536 /
\plot  1.941 24.536  1.890 24.594 /
\plot  1.890 24.594  1.850 24.644 /
\plot  1.850 24.644  1.814 24.685 /
\plot  1.814 24.685  1.789 24.716 /
\plot  1.789 24.716  1.770 24.740 /
\plot  1.770 24.740  1.757 24.754 /
\plot  1.757 24.754  1.748 24.761 /
\plot  1.748 24.761  1.746 24.765 /
%
% Fig POLYLINE object
%
\linethickness= 0.500pt
\setplotsymbol ({\thinlinefont .})
\plot  7.461 24.765  7.459 24.761 /
\plot  7.459 24.761  7.451 24.754 /
\plot  7.451 24.754  7.440 24.740 /
\plot  7.440 24.740  7.421 24.716 /
\plot  7.421 24.716  7.394 24.687 /
\plot  7.394 24.687  7.360 24.646 /
\plot  7.360 24.646  7.319 24.598 /
\plot  7.319 24.598  7.271 24.539 /
\plot  7.271 24.539  7.218 24.475 /
\plot  7.218 24.475  7.159 24.405 /
\plot  7.159 24.405  7.097 24.329 /
\plot  7.097 24.329  7.034 24.251 /
\plot  7.034 24.251  6.968 24.172 /
\plot  6.968 24.172  6.905 24.094 /
\plot  6.905 24.094  6.845 24.018 /
\plot  6.845 24.018  6.788 23.946 /
\plot  6.788 23.946  6.735 23.876 /
\plot  6.735 23.876  6.687 23.810 /
\plot  6.687 23.810  6.644 23.751 /
\plot  6.644 23.751  6.606 23.696 /
\plot  6.606 23.696  6.576 23.647 /
\plot  6.576 23.647  6.551 23.603 /
\plot  6.551 23.603  6.532 23.563 /
\plot  6.532 23.563  6.517 23.527 /
\plot  6.517 23.527  6.509 23.495 /
\plot  6.509 23.495  6.505 23.461 /
\plot  6.505 23.461  6.507 23.429 /
\plot  6.507 23.429  6.515 23.402 /
\plot  6.515 23.402  6.530 23.374 /
\plot  6.530 23.374  6.549 23.351 /
\plot  6.549 23.351  6.572 23.326 /
\plot  6.572 23.326  6.600 23.305 /
\plot  6.600 23.305  6.629 23.283 /
\plot  6.629 23.283  6.663 23.262 /
\plot  6.663 23.262  6.699 23.243 /
\plot  6.699 23.243  6.737 23.224 /
\plot  6.737 23.224  6.775 23.207 /
\plot  6.775 23.207  6.814 23.192 /
\plot  6.814 23.192  6.854 23.178 /
\plot  6.854 23.178  6.892 23.167 /
\plot  6.892 23.167  6.930 23.156 /
\plot  6.930 23.156  6.966 23.150 /
\plot  6.966 23.150  7.002 23.148 /
\putrule from  7.002 23.148 to  7.038 23.148
\plot  7.038 23.148  7.074 23.152 /
\plot  7.074 23.152  7.108 23.163 /
\plot  7.108 23.163  7.144 23.178 /
\plot  7.144 23.178  7.176 23.194 /
\plot  7.176 23.194  7.207 23.218 /
\plot  7.207 23.218  7.243 23.245 /
\plot  7.243 23.245  7.279 23.279 /
\plot  7.279 23.279  7.319 23.317 /
\plot  7.319 23.317  7.362 23.364 /
\plot  7.362 23.364  7.408 23.415 /
\plot  7.408 23.415  7.457 23.470 /
\plot  7.457 23.470  7.508 23.531 /
\plot  7.508 23.531  7.563 23.597 /
\plot  7.563 23.597  7.618 23.666 /
\plot  7.618 23.666  7.675 23.738 /
\plot  7.675 23.738  7.732 23.810 /
\plot  7.732 23.810  7.787 23.882 /
\plot  7.787 23.882  7.840 23.950 /
\plot  7.840 23.950  7.891 24.016 /
\plot  7.891 24.016  7.938 24.077 /
\plot  7.938 24.077  7.978 24.130 /
\plot  7.978 24.130  8.012 24.177 /
\plot  8.012 24.177  8.041 24.215 /
\plot  8.041 24.215  8.062 24.244 /
\plot  8.062 24.244  8.077 24.263 /
\plot  8.077 24.263  8.088 24.278 /
\plot  8.088 24.278  8.094 24.285 /
\plot  8.094 24.285  8.096 24.289 /
%
% Fig POLYLINE object
%
\linethickness= 0.500pt
\setplotsymbol ({\thinlinefont .})
\plot  9.684 24.130  9.686 24.126 /
\plot  9.686 24.126  9.690 24.115 /
\plot  9.690 24.115  9.696 24.096 /
\plot  9.696 24.096  9.707 24.069 /
\plot  9.707 24.069  9.724 24.031 /
\plot  9.724 24.031  9.743 23.982 /
\plot  9.743 23.982  9.764 23.925 /
\plot  9.764 23.925  9.792 23.859 /
\plot  9.792 23.859  9.819 23.789 /
\plot  9.819 23.789  9.851 23.717 /
\plot  9.851 23.717  9.883 23.643 /
\plot  9.883 23.643  9.912 23.571 /
\plot  9.912 23.571  9.944 23.503 /
\plot  9.944 23.503  9.974 23.440 /
\plot  9.974 23.440 10.003 23.381 /
\plot 10.003 23.381 10.033 23.330 /
\plot 10.033 23.330 10.058 23.285 /
\plot 10.058 23.285 10.086 23.247 /
\plot 10.086 23.247 10.109 23.218 /
\plot 10.109 23.218 10.135 23.194 /
\plot 10.135 23.194 10.160 23.178 /
\plot 10.160 23.178 10.188 23.165 /
\plot 10.188 23.165 10.217 23.158 /
\plot 10.217 23.158 10.249 23.156 /
\plot 10.249 23.156 10.283 23.158 /
\plot 10.283 23.158 10.319 23.165 /
\plot 10.319 23.165 10.357 23.175 /
\plot 10.357 23.175 10.397 23.188 /
\plot 10.397 23.188 10.439 23.203 /
\plot 10.439 23.203 10.480 23.220 /
\plot 10.480 23.220 10.522 23.239 /
\plot 10.522 23.239 10.564 23.260 /
\plot 10.564 23.260 10.602 23.281 /
\plot 10.602 23.281 10.640 23.302 /
\plot 10.640 23.302 10.676 23.326 /
\plot 10.676 23.326 10.708 23.351 /
\plot 10.708 23.351 10.736 23.376 /
\plot 10.736 23.376 10.759 23.402 /
\plot 10.759 23.402 10.776 23.429 /
\plot 10.776 23.429 10.789 23.461 /
\plot 10.789 23.461 10.795 23.495 /
\plot 10.795 23.495 10.797 23.529 /
\plot 10.797 23.529 10.793 23.567 /
\plot 10.793 23.567 10.784 23.611 /
\plot 10.784 23.611 10.770 23.662 /
\plot 10.770 23.662 10.751 23.717 /
\plot 10.751 23.717 10.725 23.779 /
\plot 10.725 23.779 10.698 23.846 /
\plot 10.698 23.846 10.664 23.920 /
\plot 10.664 23.920 10.628 23.999 /
\plot 10.628 23.999 10.590 24.077 /
\plot 10.590 24.077 10.549 24.160 /
\plot 10.549 24.160 10.509 24.238 /
\plot 10.509 24.238 10.471 24.314 /
\plot 10.471 24.314 10.435 24.384 /
\plot 10.435 24.384 10.401 24.448 /
\plot 10.401 24.448 10.376 24.498 /
\plot 10.376 24.498 10.353 24.541 /
\plot 10.353 24.541 10.338 24.570 /
\plot 10.338 24.570 10.327 24.589 /
\plot 10.327 24.589 10.321 24.602 /
\plot 10.321 24.602 10.319 24.606 /
%
% Fig POLYLINE object
%
\linethickness= 0.500pt
\setplotsymbol ({\thinlinefont .})
\plot 15.081 25.082 15.079 25.080 /
\plot 15.079 25.080 15.073 25.072 /
\plot 15.073 25.072 15.060 25.061 /
\plot 15.060 25.061 15.041 25.042 /
\plot 15.041 25.042 15.016 25.015 /
\plot 15.016 25.015 14.982 24.981 /
\plot 14.982 24.981 14.942 24.939 /
\plot 14.942 24.939 14.895 24.890 /
\plot 14.895 24.890 14.842 24.835 /
\plot 14.842 24.835 14.785 24.773 /
\plot 14.785 24.773 14.724 24.710 /
\plot 14.724 24.710 14.660 24.642 /
\plot 14.660 24.642 14.597 24.577 /
\plot 14.597 24.577 14.535 24.509 /
\plot 14.535 24.509 14.476 24.443 /
\plot 14.476 24.443 14.419 24.380 /
\plot 14.419 24.380 14.366 24.318 /
\plot 14.366 24.318 14.317 24.261 /
\plot 14.317 24.261 14.275 24.208 /
\plot 14.275 24.208 14.237 24.160 /
\plot 14.237 24.160 14.205 24.115 /
\plot 14.205 24.115 14.177 24.075 /
\plot 14.177 24.075 14.156 24.037 /
\plot 14.156 24.037 14.139 24.003 /
\plot 14.139 24.003 14.129 23.971 /
\plot 14.129 23.971 14.120 23.937 /
\plot 14.120 23.937 14.118 23.904 /
\plot 14.118 23.904 14.120 23.872 /
\plot 14.120 23.872 14.129 23.840 /
\plot 14.129 23.840 14.139 23.810 /
\plot 14.139 23.810 14.154 23.779 /
\plot 14.154 23.779 14.173 23.749 /
\plot 14.173 23.749 14.194 23.717 /
\plot 14.194 23.717 14.218 23.688 /
\plot 14.218 23.688 14.243 23.658 /
\plot 14.243 23.658 14.271 23.628 /
\plot 14.271 23.628 14.298 23.601 /
\plot 14.298 23.601 14.326 23.575 /
\plot 14.326 23.575 14.355 23.550 /
\plot 14.355 23.550 14.385 23.529 /
\plot 14.385 23.529 14.415 23.510 /
\plot 14.415 23.510 14.444 23.495 /
\plot 14.444 23.495 14.474 23.484 /
\plot 14.474 23.484 14.506 23.480 /
\plot 14.506 23.480 14.537 23.478 /
\plot 14.537 23.478 14.571 23.484 /
\plot 14.571 23.484 14.605 23.495 /
\plot 14.605 23.495 14.637 23.510 /
\plot 14.637 23.510 14.671 23.529 /
\plot 14.671 23.529 14.709 23.556 /
\plot 14.709 23.556 14.749 23.588 /
\plot 14.749 23.588 14.793 23.626 /
\plot 14.793 23.626 14.842 23.671 /
\plot 14.842 23.671 14.895 23.721 /
\plot 14.895 23.721 14.952 23.777 /
\plot 14.952 23.777 15.014 23.838 /
\plot 15.014 23.838 15.077 23.904 /
\plot 15.077 23.904 15.143 23.973 /
\plot 15.143 23.973 15.210 24.045 /
\plot 15.210 24.045 15.278 24.119 /
\plot 15.278 24.119 15.344 24.191 /
\plot 15.344 24.191 15.407 24.263 /
\plot 15.407 24.263 15.469 24.329 /
\plot 15.469 24.329 15.524 24.390 /
\plot 15.524 24.390 15.572 24.445 /
\plot 15.572 24.445 15.615 24.492 /
\plot 15.615 24.492 15.649 24.530 /
\plot 15.649 24.530 15.676 24.560 /
\plot 15.676 24.560 15.695 24.581 /
\plot 15.695 24.581 15.706 24.596 /
\plot 15.706 24.596 15.714 24.602 /
\plot 15.714 24.602 15.716 24.606 /
%
% Fig POLYLINE object
%
\linethickness= 0.500pt
\setplotsymbol ({\thinlinefont .})
\plot 18.415 23.812 18.419 23.810 /
\plot 18.419 23.810 18.428 23.804 /
\plot 18.428 23.804 18.445 23.793 /
\plot 18.445 23.793 18.472 23.779 /
\plot 18.472 23.779 18.506 23.755 /
\plot 18.506 23.755 18.550 23.730 /
\plot 18.550 23.730 18.603 23.696 /
\plot 18.603 23.696 18.663 23.660 /
\plot 18.663 23.660 18.728 23.622 /
\plot 18.728 23.622 18.796 23.582 /
\plot 18.796 23.582 18.866 23.542 /
\plot 18.866 23.542 18.934 23.503 /
\plot 18.934 23.503 18.999 23.467 /
\plot 18.999 23.467 19.063 23.434 /
\plot 19.063 23.434 19.120 23.406 /
\plot 19.120 23.406 19.175 23.381 /
\plot 19.175 23.381 19.221 23.362 /
\plot 19.221 23.362 19.266 23.347 /
\plot 19.266 23.347 19.304 23.338 /
\plot 19.304 23.338 19.338 23.334 /
\plot 19.338 23.334 19.367 23.336 /
\plot 19.367 23.336 19.397 23.343 /
\plot 19.397 23.343 19.427 23.355 /
\plot 19.427 23.355 19.452 23.372 /
\plot 19.452 23.372 19.478 23.396 /
\plot 19.478 23.396 19.503 23.421 /
\plot 19.503 23.421 19.528 23.453 /
\plot 19.528 23.453 19.552 23.487 /
\plot 19.552 23.487 19.577 23.523 /
\plot 19.577 23.523 19.600 23.563 /
\plot 19.600 23.563 19.624 23.601 /
\plot 19.624 23.601 19.645 23.643 /
\plot 19.645 23.643 19.664 23.683 /
\plot 19.664 23.683 19.681 23.724 /
\plot 19.681 23.724 19.696 23.764 /
\plot 19.696 23.764 19.706 23.802 /
\plot 19.706 23.802 19.713 23.838 /
\plot 19.713 23.838 19.715 23.874 /
\plot 19.715 23.874 19.710 23.908 /
\plot 19.710 23.908 19.702 23.939 /
\plot 19.702 23.939 19.685 23.971 /
\plot 19.685 23.971 19.666 23.999 /
\plot 19.666 23.999 19.641 24.024 /
\plot 19.641 24.024 19.609 24.054 /
\plot 19.609 24.054 19.571 24.081 /
\plot 19.571 24.081 19.526 24.113 /
\plot 19.526 24.113 19.475 24.147 /
\plot 19.475 24.147 19.416 24.181 /
\plot 19.416 24.181 19.353 24.217 /
\plot 19.353 24.217 19.283 24.255 /
\plot 19.283 24.255 19.209 24.295 /
\plot 19.209 24.295 19.130 24.333 /
\plot 19.130 24.333 19.052 24.373 /
\plot 19.052 24.373 18.976 24.412 /
\plot 18.976 24.412 18.900 24.450 /
\plot 18.900 24.450 18.830 24.483 /
\plot 18.830 24.483 18.766 24.513 /
\plot 18.766 24.513 18.711 24.541 /
\plot 18.711 24.541 18.667 24.562 /
\plot 18.667 24.562 18.631 24.579 /
\plot 18.631 24.579 18.603 24.591 /
\plot 18.603 24.591 18.586 24.600 /
\plot 18.586 24.600 18.578 24.604 /
\plot 18.578 24.604 18.574 24.606 /
%
% Fig POLYLINE object
%
\linethickness= 0.500pt
\setplotsymbol ({\thinlinefont .})
\plot 23.495 24.606 23.491 24.602 /
\plot 23.491 24.602 23.484 24.591 /
\plot 23.484 24.591 23.472 24.575 /
\plot 23.472 24.575 23.453 24.547 /
\plot 23.453 24.547 23.425 24.513 /
\plot 23.425 24.513 23.393 24.469 /
\plot 23.393 24.469 23.355 24.416 /
\plot 23.355 24.416 23.315 24.361 /
\plot 23.315 24.361 23.273 24.299 /
\plot 23.273 24.299 23.230 24.238 /
\plot 23.230 24.238 23.188 24.177 /
\plot 23.188 24.177 23.150 24.117 /
\plot 23.150 24.117 23.116 24.062 /
\plot 23.116 24.062 23.086 24.009 /
\plot 23.086 24.009 23.061 23.963 /
\plot 23.061 23.963 23.042 23.920 /
\plot 23.042 23.920 23.029 23.880 /
\plot 23.029 23.880 23.021 23.844 /
\plot 23.021 23.844 23.019 23.812 /
\plot 23.019 23.812 23.021 23.779 /
\plot 23.021 23.779 23.029 23.745 /
\plot 23.029 23.745 23.044 23.709 /
\plot 23.044 23.709 23.063 23.675 /
\plot 23.063 23.675 23.084 23.637 /
\plot 23.084 23.637 23.112 23.601 /
\plot 23.112 23.601 23.142 23.563 /
\plot 23.142 23.563 23.171 23.523 /
\plot 23.171 23.523 23.205 23.484 /
\plot 23.205 23.484 23.239 23.448 /
\plot 23.239 23.448 23.273 23.415 /
\plot 23.273 23.415 23.307 23.385 /
\plot 23.307 23.385 23.338 23.357 /
\plot 23.338 23.357 23.372 23.338 /
\plot 23.372 23.338 23.404 23.326 /
\plot 23.404 23.326 23.434 23.319 /
\plot 23.434 23.319 23.465 23.324 /
\plot 23.465 23.324 23.495 23.336 /
\plot 23.495 23.336 23.518 23.353 /
\plot 23.518 23.353 23.544 23.374 /
\plot 23.544 23.374 23.569 23.404 /
\plot 23.569 23.404 23.597 23.440 /
\plot 23.597 23.440 23.626 23.484 /
\plot 23.626 23.484 23.660 23.535 /
\plot 23.660 23.535 23.694 23.594 /
\plot 23.694 23.594 23.730 23.658 /
\plot 23.730 23.658 23.768 23.728 /
\plot 23.768 23.728 23.808 23.804 /
\plot 23.808 23.804 23.848 23.882 /
\plot 23.848 23.882 23.889 23.961 /
\plot 23.889 23.961 23.929 24.039 /
\plot 23.929 24.039 23.967 24.115 /
\plot 23.967 24.115 24.003 24.187 /
\plot 24.003 24.187 24.035 24.251 /
\plot 24.035 24.251 24.062 24.308 /
\plot 24.062 24.308 24.083 24.354 /
\plot 24.083 24.354 24.102 24.390 /
\plot 24.102 24.390 24.115 24.416 /
\plot 24.115 24.416 24.124 24.433 /
\plot 24.124 24.433 24.128 24.443 /
\plot 24.128 24.443 24.130 24.448 /
%
% Fig POLYLINE object
%
\linethickness= 0.500pt
\setplotsymbol ({\thinlinefont .})
\setdots < 0.0953cm>
\plot 19.685 22.066 19.687 22.070 /
\plot 19.687 22.070 19.689 22.077 /
\plot 19.689 22.077 19.696 22.092 /
\plot 19.696 22.092 19.704 22.113 /
\plot 19.704 22.113 19.717 22.142 /
\plot 19.717 22.142 19.734 22.181 /
\plot 19.734 22.181 19.753 22.229 /
\plot 19.753 22.229 19.778 22.286 /
\plot 19.778 22.286 19.806 22.352 /
\plot 19.806 22.352 19.837 22.424 /
\plot 19.837 22.424 19.871 22.502 /
\plot 19.871 22.502 19.907 22.587 /
\plot 19.907 22.587 19.945 22.672 /
\plot 19.945 22.672 19.983 22.761 /
\plot 19.983 22.761 20.024 22.847 /
\plot 20.024 22.847 20.062 22.934 /
\plot 20.062 22.934 20.102 23.017 /
\plot 20.102 23.017 20.140 23.097 /
\plot 20.140 23.097 20.180 23.175 /
\plot 20.180 23.175 20.218 23.249 /
\plot 20.218 23.249 20.254 23.317 /
\plot 20.254 23.317 20.292 23.383 /
\plot 20.292 23.383 20.328 23.444 /
\plot 20.328 23.444 20.364 23.499 /
\plot 20.364 23.499 20.403 23.554 /
\plot 20.403 23.554 20.441 23.605 /
\plot 20.441 23.605 20.479 23.654 /
\plot 20.479 23.654 20.519 23.702 /
\plot 20.519 23.702 20.561 23.749 /
\plot 20.561 23.749 20.604 23.796 /
\plot 20.604 23.796 20.648 23.842 /
\plot 20.648 23.842 20.693 23.889 /
\plot 20.693 23.889 20.737 23.935 /
\plot 20.737 23.935 20.784 23.984 /
\plot 20.784 23.984 20.828 24.033 /
\plot 20.828 24.033 20.875 24.081 /
\plot 20.875 24.081 20.921 24.130 /
\plot 20.921 24.130 20.968 24.181 /
\plot 20.968 24.181 21.016 24.229 /
\plot 21.016 24.229 21.063 24.280 /
\plot 21.063 24.280 21.112 24.329 /
\plot 21.112 24.329 21.160 24.376 /
\plot 21.160 24.376 21.209 24.424 /
\plot 21.209 24.424 21.260 24.469 /
\plot 21.260 24.469 21.311 24.513 /
\plot 21.311 24.513 21.364 24.553 /
\plot 21.364 24.553 21.416 24.591 /
\plot 21.416 24.591 21.471 24.627 /
\plot 21.471 24.627 21.526 24.661 /
\plot 21.526 24.661 21.586 24.689 /
\plot 21.586 24.689 21.645 24.714 /
\plot 21.645 24.714 21.706 24.733 /
\plot 21.706 24.733 21.772 24.748 /
\plot 21.772 24.748 21.838 24.759 /
\plot 21.838 24.759 21.907 24.765 /
\plot 21.907 24.765 21.967 24.765 /
\plot 21.967 24.765 22.028 24.763 /
\plot 22.028 24.763 22.094 24.754 /
\plot 22.094 24.754 22.164 24.744 /
\plot 22.164 24.744 22.238 24.729 /
\plot 22.238 24.729 22.316 24.708 /
\plot 22.316 24.708 22.399 24.685 /
\plot 22.399 24.685 22.485 24.657 /
\plot 22.485 24.657 22.578 24.625 /
\plot 22.578 24.625 22.678 24.589 /
\plot 22.678 24.589 22.780 24.551 /
\plot 22.780 24.551 22.888 24.509 /
\plot 22.888 24.509 22.998 24.462 /
\plot 22.998 24.462 23.112 24.416 /
\plot 23.112 24.416 23.226 24.367 /
\plot 23.226 24.367 23.345 24.316 /
\plot 23.345 24.316 23.461 24.265 /
\plot 23.461 24.265 23.575 24.215 /
\plot 23.575 24.215 23.688 24.164 /
\plot 23.688 24.164 23.796 24.115 /
\plot 23.796 24.115 23.899 24.069 /
\plot 23.899 24.069 23.995 24.024 /
\plot 23.995 24.024 24.081 23.984 /
\plot 24.081 23.984 24.160 23.948 /
\plot 24.160 23.948 24.229 23.914 /
\plot 24.229 23.914 24.289 23.887 /
\plot 24.289 23.887 24.337 23.865 /
\plot 24.337 23.865 24.376 23.846 /
\plot 24.376 23.846 24.403 23.834 /
\plot 24.403 23.834 24.424 23.823 /
\plot 24.424 23.823 24.437 23.817 /
\plot 24.437 23.817 24.445 23.815 /
\plot 24.445 23.815 24.448 23.812 /
%
% Fig TEXT object
%
\put{\SetFigFont{7}{8.4}{rm}$U_1^{10}(2A_1+D_6)$} [lB] at 11.271 25.876
%
% Fig TEXT object
%
\put{\SetFigFont{7}{8.4}{rm}$(I_2, I_2^*,I_2)$} [lB] at  3.016 25.876
%
% Fig TEXT object
%
\put{\SetFigFont{7}{8.4}{rm}$U_2^{10}(3A_1+D_4)$} [lB] at 19.209 25.718
\linethickness=0pt
\putrectangle corners at  1.086 26.289 and 24.473 22.041
\endpicture}
}
   \vskip 0.5cm
   \centerline{\font\thinlinefont=cmr5
\begingroup\makeatletter\ifx\SetFigFont\undefined%
\gdef\SetFigFont#1#2#3#4#5{%
  \reset@font\fontsize{#1}{#2pt}%
  \fontfamily{#3}\fontseries{#4}\fontshape{#5}%
  \selectfont}%
\fi\endgroup%
\mbox{\beginpicture
\setcoordinatesystem units <0.56000cm,0.56000cm>
\unitlength=0.56000cm
\linethickness=1pt
\setplotsymbol ({\makebox(0,0)[l]{\tencirc\symbol{'160}}})
\setshadesymbol ({\thinlinefont .})
\setlinear
%
% Fig POLYLINE object
%
\linethickness= 0.500pt
\setplotsymbol ({\thinlinefont .})
\putrule from  0.794 21.590 to  2.540 21.590
%
% Fig POLYLINE object
%
\linethickness= 0.500pt
\setplotsymbol ({\thinlinefont .})
\putrule from  2.857 21.590 to  4.604 21.590
%
% Fig POLYLINE object
%
\linethickness= 0.500pt
\setplotsymbol ({\thinlinefont .})
\plot  4.445 21.590  4.445 21.590 /
%
% Fig POLYLINE object
%
\linethickness= 0.500pt
\setplotsymbol ({\thinlinefont .})
\putrule from  2.223 22.384 to  2.223 21.273
%
% Fig POLYLINE object
%
\linethickness= 0.500pt
\setplotsymbol ({\thinlinefont .})
\putrule from  3.016 22.384 to  3.016 21.273
%
% Fig POLYLINE object
%
\linethickness= 0.500pt
\setplotsymbol ({\thinlinefont .})
\putrule from  1.905 22.860 to  1.905 21.431
%
% Fig POLYLINE object
%
\linethickness= 0.500pt
\setplotsymbol ({\thinlinefont .})
\putrule from  1.905 21.590 to  1.905 21.273
%
% Fig POLYLINE object
%
\linethickness= 0.500pt
\setplotsymbol ({\thinlinefont .})
\putrule from  3.334 23.019 to  3.334 21.273
%
% Fig POLYLINE object
%
\linethickness= 0.500pt
\setplotsymbol ({\thinlinefont .})
\putrule from  1.587 22.225 to  1.587 20.955
%
% Fig POLYLINE object
%
\linethickness= 0.500pt
\setplotsymbol ({\thinlinefont .})
\putrule from  3.810 22.384 to  3.810 20.955
%
% Fig POLYLINE object
%
\linethickness= 0.500pt
\setplotsymbol ({\thinlinefont .})
\putrule from  1.111 22.225 to  1.111 20.320
%
% Fig POLYLINE object
%
\linethickness= 0.500pt
\setplotsymbol ({\thinlinefont .})
\putrule from  4.286 22.384 to  4.286 20.161
%
% Fig POLYLINE object
%
\linethickness= 0.500pt
\setplotsymbol ({\thinlinefont .})
\setdashes < 0.1270cm>
\plot  2.064 22.066  3.175 22.066 /
%
% Fig POLYLINE object
%
\linethickness= 0.500pt
\setplotsymbol ({\thinlinefont .})
\plot  1.746 22.701  3.651 22.701 /
%
% Fig POLYLINE object
%
\linethickness= 0.500pt
\setplotsymbol ({\thinlinefont .})
\plot  1.429 21.114  4.128 21.114 /
%
% Fig POLYLINE object
%
\linethickness= 0.500pt
\setplotsymbol ({\thinlinefont .})
\plot  0.953 20.637  4.604 20.637 /
%
% Fig POLYLINE object
%
\linethickness= 0.500pt
\setplotsymbol ({\thinlinefont .})
\setsolid
\putrule from  5.397 21.590 to  6.509 21.590
%
% arrow head
%
\plot  6.255 21.526  6.509 21.590  6.255 21.654 /
%
%
% Fig POLYLINE object
%
\linethickness= 0.500pt
\setplotsymbol ({\thinlinefont .})
\putrule from  8.731 22.225 to  8.731 21.273
%
% Fig POLYLINE object
%
\linethickness= 0.500pt
\setplotsymbol ({\thinlinefont .})
\putrule from  9.684 22.225 to  9.684 21.273
%
% Fig POLYLINE object
%
\linethickness= 0.500pt
\setplotsymbol ({\thinlinefont .})
\putrule from  8.414 22.543 to  8.414 21.273
%
% Fig POLYLINE object
%
\linethickness= 0.500pt
\setplotsymbol ({\thinlinefont .})
\putrule from 10.001 22.701 to 10.001 21.273
%
% Fig POLYLINE object
%
\linethickness= 0.500pt
\setplotsymbol ({\thinlinefont .})
\putrule from  7.938 22.225 to  7.938 20.955
%
% Fig POLYLINE object
%
\linethickness= 0.500pt
\setplotsymbol ({\thinlinefont .})
\putrule from 10.319 22.225 to 10.319 20.955
%
% Fig POLYLINE object
%
\linethickness= 0.500pt
\setplotsymbol ({\thinlinefont .})
\setdashes < 0.1270cm>
\plot  8.572 22.066  9.842 22.066 /
%
% Fig POLYLINE object
%
\linethickness= 0.500pt
\setplotsymbol ({\thinlinefont .})
\plot  8.255 22.384 10.160 22.384 /
%
% Fig POLYLINE object
%
\linethickness= 0.500pt
\setplotsymbol ({\thinlinefont .})
\plot  7.779 21.114 10.478 21.114 /
%
% Fig POLYLINE object
%
\linethickness= 0.500pt
\setplotsymbol ({\thinlinefont .})
\setsolid
\putrule from 11.748 21.590 to 13.176 21.590
%
% arrow head
%
\plot 12.922 21.526 13.176 21.590 12.922 21.654 /
%
%
% Fig POLYLINE object
%
\linethickness= 0.500pt
\setplotsymbol ({\thinlinefont .})
\putrule from 16.192 21.590 to 18.415 21.590
%
% Fig POLYLINE object
%
\linethickness= 0.500pt
\setplotsymbol ({\thinlinefont .})
\putrule from 15.558 22.066 to 15.558 21.273
%
% Fig POLYLINE object
%
\linethickness= 0.500pt
\setplotsymbol ({\thinlinefont .})
\putrule from 16.351 22.225 to 16.351 21.273
%
% Fig POLYLINE object
%
\linethickness= 0.500pt
\setplotsymbol ({\thinlinefont .})
\putrule from 15.240 22.860 to 15.240 21.273
%
% Fig POLYLINE object
%
\linethickness= 0.500pt
\setplotsymbol ({\thinlinefont .})
\putrule from 16.669 22.860 to 16.669 21.273
%
% Fig POLYLINE object
%
\linethickness= 0.500pt
\setplotsymbol ({\thinlinefont .})
\putrule from 14.764 22.225 to 14.764 20.955
%
% Fig POLYLINE object
%
\linethickness= 0.500pt
\setplotsymbol ({\thinlinefont .})
\putrule from 17.304 22.384 to 17.304 20.955
%
% Fig POLYLINE object
%
\linethickness= 0.500pt
\setplotsymbol ({\thinlinefont .})
\setdashes < 0.1270cm>
\plot 13.811 21.590 15.875 21.590 /
%
% Fig POLYLINE object
%
\linethickness= 0.500pt
\setplotsymbol ({\thinlinefont .})
\plot 15.399 21.907 16.510 21.907 /
%
% Fig POLYLINE object
%
\linethickness= 0.500pt
\setplotsymbol ({\thinlinefont .})
\plot 15.081 22.701 16.828 22.701 /
%
% Fig POLYLINE object
%
\linethickness= 0.500pt
\setplotsymbol ({\thinlinefont .})
\plot 14.605 21.114 17.621 21.114 /
%
% Fig POLYLINE object
%
\linethickness= 0.500pt
\setplotsymbol ({\thinlinefont .})
\setsolid
\plot  7.144 21.114  8.890 21.907 /
%
% Fig POLYLINE object
%
\linethickness= 0.500pt
\setplotsymbol ({\thinlinefont .})
\plot  9.366 21.907 11.430 20.955 /
%
% Fig POLYLINE object
%
\linethickness= 0.500pt
\setplotsymbol ({\thinlinefont .})
\setdashes < 0.1270cm>
\plot  6.985 21.431  9.207 19.844 /
%
% Fig POLYLINE object
%
\linethickness= 0.500pt
\setplotsymbol ({\thinlinefont .})
\plot 11.271 21.431  8.572 19.685 /
%
% Fig POLYLINE object
%
\linethickness= 0.500pt
\setplotsymbol ({\thinlinefont .})
\setdots < 0.0953cm>
\plot 14.287 22.066 14.287 22.062 /
\plot 14.287 22.062 14.285 22.051 /
\plot 14.285 22.051 14.281 22.032 /
\plot 14.281 22.032 14.275 22.005 /
\plot 14.275 22.005 14.268 21.965 /
\plot 14.268 21.965 14.260 21.912 /
\plot 14.260 21.912 14.249 21.848 /
\plot 14.249 21.848 14.237 21.774 /
\plot 14.237 21.774 14.224 21.692 /
\plot 14.224 21.692 14.211 21.601 /
\plot 14.211 21.601 14.199 21.505 /
\plot 14.199 21.505 14.188 21.404 /
\plot 14.188 21.404 14.177 21.304 /
\plot 14.177 21.304 14.171 21.205 /
\plot 14.171 21.205 14.167 21.107 /
\plot 14.167 21.107 14.167 21.014 /
\plot 14.167 21.014 14.173 20.927 /
\plot 14.173 20.927 14.182 20.847 /
\plot 14.182 20.847 14.199 20.773 /
\plot 14.199 20.773 14.220 20.707 /
\plot 14.220 20.707 14.249 20.648 /
\plot 14.249 20.648 14.285 20.595 /
\plot 14.285 20.595 14.330 20.551 /
\plot 14.330 20.551 14.383 20.513 /
\plot 14.383 20.513 14.446 20.479 /
\plot 14.446 20.479 14.495 20.460 /
\plot 14.495 20.460 14.548 20.441 /
\plot 14.548 20.441 14.605 20.426 /
\plot 14.605 20.426 14.668 20.411 /
\plot 14.668 20.411 14.736 20.398 /
\plot 14.736 20.398 14.810 20.388 /
\plot 14.810 20.388 14.889 20.377 /
\plot 14.889 20.377 14.971 20.369 /
\plot 14.971 20.369 15.058 20.360 /
\plot 15.058 20.360 15.151 20.354 /
\plot 15.151 20.354 15.248 20.350 /
\plot 15.248 20.350 15.348 20.343 /
\plot 15.348 20.343 15.452 20.341 /
\plot 15.452 20.341 15.560 20.337 /
\plot 15.560 20.337 15.668 20.335 /
\plot 15.668 20.335 15.780 20.333 /
\plot 15.780 20.333 15.892 20.331 /
\plot 15.892 20.331 16.006 20.331 /
\plot 16.006 20.331 16.121 20.331 /
\plot 16.121 20.331 16.235 20.331 /
\plot 16.235 20.331 16.349 20.331 /
\plot 16.349 20.331 16.461 20.333 /
\plot 16.461 20.333 16.573 20.335 /
\plot 16.573 20.335 16.681 20.337 /
\plot 16.681 20.337 16.787 20.341 /
\plot 16.787 20.341 16.891 20.343 /
\plot 16.891 20.343 16.990 20.350 /
\plot 16.990 20.350 17.088 20.354 /
\plot 17.088 20.354 17.179 20.360 /
\plot 17.179 20.360 17.266 20.369 /
\plot 17.266 20.369 17.348 20.377 /
\plot 17.348 20.377 17.427 20.388 /
\plot 17.427 20.388 17.498 20.398 /
\plot 17.498 20.398 17.564 20.411 /
\plot 17.564 20.411 17.625 20.426 /
\plot 17.625 20.426 17.683 20.441 /
\plot 17.683 20.441 17.733 20.460 /
\plot 17.733 20.460 17.780 20.479 /
\plot 17.780 20.479 17.839 20.513 /
\plot 17.839 20.513 17.888 20.551 /
\plot 17.888 20.551 17.928 20.595 /
\plot 17.928 20.595 17.960 20.648 /
\plot 17.960 20.648 17.983 20.707 /
\plot 17.983 20.707 17.998 20.773 /
\plot 17.998 20.773 18.006 20.847 /
\plot 18.006 20.847 18.009 20.927 /
\plot 18.009 20.927 18.004 21.014 /
\plot 18.004 21.014 17.996 21.107 /
\plot 17.996 21.107 17.983 21.205 /
\plot 17.983 21.205 17.966 21.304 /
\plot 17.966 21.304 17.947 21.404 /
\plot 17.947 21.404 17.926 21.505 /
\plot 17.926 21.505 17.903 21.601 /
\plot 17.903 21.601 17.882 21.692 /
\plot 17.882 21.692 17.860 21.774 /
\plot 17.860 21.774 17.841 21.848 /
\plot 17.841 21.848 17.824 21.912 /
\plot 17.824 21.912 17.810 21.965 /
\plot 17.810 21.965 17.799 22.005 /
\plot 17.799 22.005 17.791 22.032 /
\plot 17.791 22.032 17.784 22.051 /
\plot 17.784 22.051 17.782 22.062 /
\plot 17.782 22.062 17.780 22.066 /
%
% Fig TEXT object
%
\put{\SetFigFont{7}{8.4}{rm}$(I_0^*,I_0^*)$} [lB] at  1.746 23.495
%
% Fig TEXT object
%
\put{\SetFigFont{7}{8.4}{rm}$U_1^4(2D_4)$} [lB] at  8.096 23.178
%
% Fig TEXT object
%
\put{\SetFigFont{7}{8.4}{rm}$U_2^4(3A_1+D_4)$} [lB] at 15.399 23.336
%
% Fig TEXT object
%
\put{\SetFigFont{7}{8.4}{rm}$p$} [lB] at 15.716 20.796
%
% Fig TEXT object
%
\put{\SetFigFont{7}{8.4}{rm}$E_1$} [lB] at  6.985 22.066
%
% Fig TEXT object
%
\put{\SetFigFont{7}{8.4}{rm}$E_2$} [lB] at 10.954 22.066
\linethickness=0pt
\putrectangle corners at  0.768 23.908 and 18.440 19.660
\endpicture}
}
   \vskip 0.5cm
   \begin{center}
    Figure $3$
   \end{center}
  \centerline{\font\thinlinefont=cmr5
\begingroup\makeatletter\ifx\SetFigFont\undefined%
\gdef\SetFigFont#1#2#3#4#5{%
  \reset@font\fontsize{#1}{#2pt}%
  \fontfamily{#3}\fontseries{#4}\fontshape{#5}%
  \selectfont}%
\fi\endgroup%
\mbox{\beginpicture
\setcoordinatesystem units <0.56000cm,0.56000cm>
\unitlength=0.56000cm
\linethickness=1pt
\setplotsymbol ({\makebox(0,0)[l]{\tencirc\symbol{'160}}})
\setshadesymbol ({\thinlinefont .})
\setlinear
%
% Fig POLYLINE object
%
\linethickness= 0.500pt
\setplotsymbol ({\thinlinefont .})
\plot  2.857 23.019  5.239 22.066 /
%
% Fig POLYLINE object
%
\linethickness= 0.500pt
\setplotsymbol ({\thinlinefont .})
\plot  3.969 23.178  1.429 21.907 /
%
% Fig POLYLINE object
%
\linethickness= 0.500pt
\setplotsymbol ({\thinlinefont .})
\putrule from  2.540 23.336 to  2.540 22.066
%
% Fig POLYLINE object
%
\linethickness= 0.500pt
\setplotsymbol ({\thinlinefont .})
\putrule from  1.905 24.448 to  1.905 21.749
%
% Fig POLYLINE object
%
\linethickness= 0.500pt
\setplotsymbol ({\thinlinefont .})
\plot  4.445 23.495  4.128 21.907 /
%
% Fig POLYLINE object
%
\linethickness= 0.500pt
\setplotsymbol ({\thinlinefont .})
\plot  5.080 23.654  4.604 21.590 /
%
% Fig POLYLINE object
%
\linethickness= 0.500pt
\setplotsymbol ({\thinlinefont .})
\plot  6.985 24.606  6.350 22.225 /
%
% Fig POLYLINE object
%
\linethickness= 0.500pt
\setplotsymbol ({\thinlinefont .})
\plot  6.191 22.860  7.779 22.066 /
%
% Fig POLYLINE object
%
\linethickness= 0.500pt
\setplotsymbol ({\thinlinefont .})
\plot  6.350 23.971  8.096 23.178 /
%
% Fig POLYLINE object
%
\linethickness= 0.500pt
\setplotsymbol ({\thinlinefont .})
\plot  7.779 23.812  7.144 21.749 /
%
% Fig POLYLINE object
%
\linethickness= 0.500pt
\setplotsymbol ({\thinlinefont .})
\setdashes < 0.1270cm>
\plot  1.429 24.130  7.461 24.130 /
%
% Fig POLYLINE object
%
\linethickness= 0.500pt
\setplotsymbol ({\thinlinefont .})
\setsolid
\putrule from  9.049 22.860 to  9.842 22.860
%
% arrow head
%
\plot  9.588 22.796  9.842 22.860  9.588 22.924 /
%
%
% Fig POLYLINE object
%
\linethickness= 0.500pt
\setplotsymbol ({\thinlinefont .})
\plot 11.906 23.654 14.446 22.225 /
%
% Fig POLYLINE object
%
\linethickness= 0.500pt
\setplotsymbol ({\thinlinefont .})
\plot 13.811 23.812 12.859 22.384 /
%
% Fig POLYLINE object
%
\linethickness= 0.500pt
\setplotsymbol ({\thinlinefont .})
\plot 14.446 23.654 13.494 22.225 /
%
% Fig POLYLINE object
%
\linethickness= 0.500pt
\setplotsymbol ({\thinlinefont .})
\plot 11.271 23.654 11.906 22.543 /
%
% Fig POLYLINE object
%
\linethickness= 0.500pt
\setplotsymbol ({\thinlinefont .})
\plot 15.558 23.971 17.621 23.178 /
%
% Fig POLYLINE object
%
\linethickness= 0.500pt
\setplotsymbol ({\thinlinefont .})
\plot 15.558 22.543 17.304 22.384 /
%
% Fig POLYLINE object
%
\linethickness= 0.500pt
\setplotsymbol ({\thinlinefont .})
\plot 17.145 23.971 16.828 22.066 /
%
% Fig POLYLINE object
%
\linethickness= 0.500pt
\setplotsymbol ({\thinlinefont .})
\putrule from 18.098 22.701 to 18.891 22.701
%
% arrow head
%
\plot 18.637 22.638 18.891 22.701 18.637 22.765 /
%
%
% Fig POLYLINE object
%
\linethickness= 0.500pt
\setplotsymbol ({\thinlinefont .})
\plot 23.495 23.971 25.400 23.654 /
%
% Fig POLYLINE object
%
\linethickness= 0.500pt
\setplotsymbol ({\thinlinefont .})
\putrule from 23.495 22.701 to 25.400 22.701
%
% Fig POLYLINE object
%
\linethickness= 0.500pt
\setplotsymbol ({\thinlinefont .})
\plot 25.082 24.289 24.924 22.225 /
%
% Fig POLYLINE object
%
\linethickness= 0.500pt
\setplotsymbol ({\thinlinefont .})
\plot 20.955 23.495 22.860 22.384 /
%
% Fig POLYLINE object
%
\linethickness= 0.500pt
\setplotsymbol ({\thinlinefont .})
\plot 22.384 23.336 21.907 22.543 /
%
% Fig POLYLINE object
%
\linethickness= 0.500pt
\setplotsymbol ({\thinlinefont .})
\plot 22.860 23.178 22.384 22.384 /
%
% Fig POLYLINE object
%
\linethickness= 0.500pt
\setplotsymbol ({\thinlinefont .})
\plot 20.479 23.019 20.955 22.225 /
%
% Fig POLYLINE object
%
\linethickness= 0.500pt
\setplotsymbol ({\thinlinefont .})
\setdashes < 0.1270cm>
\plot 21.907 23.812 19.844 22.225 /
%
% Fig POLYLINE object
%
\linethickness= 0.500pt
\setplotsymbol ({\thinlinefont .})
\setsolid
\plot 13.018 24.765 11.430 22.384 /
%
% Fig POLYLINE object
%
\linethickness= 0.500pt
\setplotsymbol ({\thinlinefont .})
\setdashes < 0.1270cm>
\plot 12.065 24.606 16.669 24.130 /
%
% Fig POLYLINE object
%
\linethickness= 0.500pt
\setplotsymbol ({\thinlinefont .})
\plot 15.875 24.606 15.875 22.066 /
%
% Fig POLYLINE object
%
\linethickness= 0.500pt
\setplotsymbol ({\thinlinefont .})
\setdots < 0.0953cm>
\plot 23.971 24.130 23.971 24.126 /
\plot 23.971 24.126 23.971 24.119 /
\plot 23.971 24.119 23.971 24.105 /
\plot 23.971 24.105 23.973 24.083 /
\plot 23.973 24.083 23.973 24.052 /
\plot 23.973 24.052 23.975 24.009 /
\plot 23.975 24.009 23.978 23.956 /
\plot 23.978 23.956 23.980 23.893 /
\plot 23.980 23.893 23.982 23.817 /
\plot 23.982 23.817 23.984 23.732 /
\plot 23.984 23.732 23.986 23.637 /
\plot 23.986 23.637 23.988 23.533 /
\plot 23.988 23.533 23.990 23.423 /
\plot 23.990 23.423 23.990 23.307 /
\plot 23.990 23.307 23.990 23.186 /
\plot 23.990 23.186 23.990 23.063 /
\plot 23.990 23.063 23.988 22.940 /
\plot 23.988 22.940 23.986 22.818 /
\plot 23.986 22.818 23.980 22.695 /
\plot 23.980 22.695 23.975 22.576 /
\plot 23.975 22.576 23.967 22.462 /
\plot 23.967 22.462 23.956 22.354 /
\plot 23.956 22.354 23.946 22.250 /
\plot 23.946 22.250 23.931 22.155 /
\plot 23.931 22.155 23.914 22.064 /
\plot 23.914 22.064 23.893 21.982 /
\plot 23.893 21.982 23.870 21.905 /
\plot 23.870 21.905 23.844 21.836 /
\plot 23.844 21.836 23.815 21.774 /
\plot 23.815 21.774 23.781 21.719 /
\plot 23.781 21.719 23.743 21.670 /
\plot 23.743 21.670 23.700 21.628 /
\plot 23.700 21.628 23.654 21.590 /
\plot 23.654 21.590 23.605 21.558 /
\plot 23.605 21.558 23.552 21.533 /
\plot 23.552 21.533 23.495 21.510 /
\plot 23.495 21.510 23.434 21.488 /
\plot 23.434 21.488 23.366 21.474 /
\plot 23.366 21.474 23.294 21.461 /
\plot 23.294 21.461 23.218 21.450 /
\plot 23.218 21.450 23.137 21.444 /
\plot 23.137 21.444 23.053 21.440 /
\plot 23.053 21.440 22.964 21.440 /
\plot 22.964 21.440 22.873 21.440 /
\plot 22.873 21.440 22.777 21.442 /
\plot 22.777 21.442 22.678 21.446 /
\plot 22.678 21.446 22.578 21.452 /
\plot 22.578 21.452 22.475 21.461 /
\plot 22.475 21.461 22.371 21.469 /
\plot 22.371 21.469 22.265 21.480 /
\plot 22.265 21.480 22.162 21.491 /
\plot 22.162 21.491 22.056 21.503 /
\plot 22.056 21.503 21.950 21.514 /
\plot 21.950 21.514 21.844 21.527 /
\plot 21.844 21.527 21.742 21.541 /
\plot 21.742 21.541 21.641 21.554 /
\plot 21.641 21.554 21.541 21.567 /
\plot 21.541 21.567 21.446 21.582 /
\plot 21.446 21.582 21.353 21.594 /
\plot 21.353 21.594 21.262 21.609 /
\plot 21.262 21.609 21.175 21.622 /
\plot 21.175 21.622 21.095 21.637 /
\plot 21.095 21.637 21.016 21.649 /
\plot 21.016 21.649 20.942 21.664 /
\plot 20.942 21.664 20.872 21.679 /
\plot 20.872 21.679 20.807 21.696 /
\plot 20.807 21.696 20.748 21.713 /
\plot 20.748 21.713 20.690 21.730 /
\plot 20.690 21.730 20.637 21.749 /
\plot 20.637 21.749 20.566 21.780 /
\plot 20.566 21.780 20.502 21.816 /
\plot 20.502 21.816 20.443 21.859 /
\plot 20.443 21.859 20.392 21.907 /
\plot 20.392 21.907 20.345 21.963 /
\plot 20.345 21.963 20.303 22.024 /
\plot 20.303 22.024 20.265 22.092 /
\plot 20.265 22.092 20.229 22.166 /
\plot 20.229 22.166 20.197 22.244 /
\plot 20.197 22.244 20.168 22.327 /
\plot 20.168 22.327 20.140 22.411 /
\plot 20.140 22.411 20.117 22.498 /
\plot 20.117 22.498 20.094 22.583 /
\plot 20.094 22.583 20.074 22.665 /
\plot 20.074 22.665 20.058 22.741 /
\plot 20.058 22.741 20.043 22.809 /
\plot 20.043 22.809 20.030 22.868 /
\plot 20.030 22.868 20.022 22.919 /
\plot 20.022 22.919 20.013 22.957 /
\plot 20.013 22.957 20.009 22.985 /
\plot 20.009 22.985 20.005 23.004 /
\plot 20.005 23.004 20.003 23.015 /
\plot 20.003 23.015 20.003 23.019 /
%
% Fig TEXT object
%
\put{\SetFigFont{7}{8.4}{rm}$(I_1^*,I_4,I_1)$} [lB] at  2.857 25.400
%
% Fig TEXT object
%
\put{\SetFigFont{7}{8.4}{rm}$U_1^9(D_5+A_3)$} [lB] at 15.240 25.241
%
% Fig TEXT object
%
\put{\SetFigFont{7}{8.4}{rm}$U_2^9(A_1+2A_3)$} [lB] at 23.178 25.082
\linethickness=0pt
\putrectangle corners at  1.403 25.813 and 25.425 21.414
\endpicture}
}
  \vskip 0.5cm
  \centerline{\font\thinlinefont=cmr5
\begingroup\makeatletter\ifx\SetFigFont\undefined%
\gdef\SetFigFont#1#2#3#4#5{%
  \reset@font\fontsize{#1}{#2pt}%
  \fontfamily{#3}\fontseries{#4}\fontshape{#5}%
  \selectfont}%
\fi\endgroup%
\mbox{\beginpicture
\setcoordinatesystem units <0.56000cm,0.56000cm>
\unitlength=0.56000cm
\linethickness=1pt
\setplotsymbol ({\makebox(0,0)[l]{\tencirc\symbol{'160}}})
\setshadesymbol ({\thinlinefont .})
\setlinear
%
% Fig POLYLINE object
%
\linethickness= 0.500pt
\setplotsymbol ({\thinlinefont .})
\plot  2.064 23.971  3.493 20.479 /
%
% Fig POLYLINE object
%
\linethickness= 0.500pt
\setplotsymbol ({\thinlinefont .})
\plot  4.445 22.701  5.239 21.907 /
%
% Fig POLYLINE object
%
\linethickness= 0.500pt
\setplotsymbol ({\thinlinefont .})
\plot  4.604 22.066  5.397 21.273 /
%
% Fig POLYLINE object
%
\linethickness= 0.500pt
\setplotsymbol ({\thinlinefont .})
\plot  5.556 21.907  4.604 20.479 /
%
% Fig POLYLINE object
%
\linethickness= 0.500pt
\setplotsymbol ({\thinlinefont .})
\plot  4.604 21.114  5.715 20.320 /
%
% Fig POLYLINE object
%
\linethickness= 0.500pt
\setplotsymbol ({\thinlinefont .})
\plot  7.144 22.860  5.080 20.161 /
%
% Fig POLYLINE object
%
\linethickness= 0.500pt
\setplotsymbol ({\thinlinefont .})
\plot  5.080 23.654  6.985 21.749 /
%
% Fig POLYLINE object
%
\linethickness= 0.500pt
\setplotsymbol ({\thinlinefont .})
\putrule from  8.096 22.543 to  9.049 22.543
%
% arrow head
%
\plot  8.795 22.479  9.049 22.543  8.795 22.606 /
%
%
% Fig POLYLINE object
%
\linethickness= 0.500pt
\setplotsymbol ({\thinlinefont .})
\plot 12.859 22.066 13.970 21.273 /
%
% Fig POLYLINE object
%
\linethickness= 0.500pt
\setplotsymbol ({\thinlinefont .})
\plot 13.970 21.907 13.176 20.320 /
%
% Fig POLYLINE object
%
\linethickness= 0.500pt
\setplotsymbol ({\thinlinefont .})
\plot 13.176 21.114 15.240 20.161 /
%
% Fig POLYLINE object
%
\linethickness= 0.500pt
\setplotsymbol ({\thinlinefont .})
\plot 15.558 22.701 14.287 20.161 /
%
% Fig POLYLINE object
%
\linethickness= 0.500pt
\setplotsymbol ({\thinlinefont .})
\plot 14.287 20.320 14.287 20.320 /
%
% Fig POLYLINE object
%
\linethickness= 0.500pt
\setplotsymbol ({\thinlinefont .})
\plot 12.859 24.130 15.716 21.907 /
%
% Fig POLYLINE object
%
\linethickness= 0.500pt
\setplotsymbol ({\thinlinefont .})
\putrule from 16.510 22.543 to 17.462 22.543
%
% arrow head
%
\plot 17.209 22.479 17.462 22.543 17.209 22.606 /
%
%
% Fig POLYLINE object
%
\linethickness= 0.500pt
\setplotsymbol ({\thinlinefont .})
\plot 20.955 22.543 20.479 21.273 /
%
% Fig POLYLINE object
%
\linethickness= 0.500pt
\setplotsymbol ({\thinlinefont .})
\plot 20.479 21.907 21.431 20.955 /
%
% Fig POLYLINE object
%
\linethickness= 0.500pt
\setplotsymbol ({\thinlinefont .})
\plot 20.479 20.637 23.019 20.003 /
%
% Fig POLYLINE object
%
\linethickness= 0.500pt
\setplotsymbol ({\thinlinefont .})
\plot 23.019 22.860 22.066 19.685 /
%
% Fig POLYLINE object
%
\linethickness= 0.500pt
\setplotsymbol ({\thinlinefont .})
\plot 20.637 23.654 23.336 22.225 /
%
% Fig POLYLINE object
%
\linethickness= 0.500pt
\setplotsymbol ({\thinlinefont .})
\plot  5.080 22.384  4.445 21.590 /
%
% Fig POLYLINE object
%
\linethickness= 0.500pt
\setplotsymbol ({\thinlinefont .})
\setdashes < 0.1270cm>
\plot  1.746 23.654  6.191 24.130 /
%
% Fig POLYLINE object
%
\linethickness= 0.500pt
\setplotsymbol ({\thinlinefont .})
\plot  2.857 22.225  4.763 21.590 /
%
% Fig POLYLINE object
%
\linethickness= 0.500pt
\setplotsymbol ({\thinlinefont .})
\plot  2.857 21.431  5.397 21.114 /
%
% Fig POLYLINE object
%
\linethickness= 0.500pt
\setplotsymbol ({\thinlinefont .})
\setsolid
\plot 12.700 23.336 13.652 21.907 /
%
% Fig POLYLINE object
%
\linethickness= 0.500pt
\setplotsymbol ({\thinlinefont .})
\plot 13.652 22.543 12.700 21.431 /
%
% Fig POLYLINE object
%
\linethickness= 0.500pt
\setplotsymbol ({\thinlinefont .})
\setdashes < 0.1270cm>
\plot 11.906 24.289 10.001 20.796 /
%
% Fig POLYLINE object
%
\linethickness= 0.500pt
\setplotsymbol ({\thinlinefont .})
\plot 10.795 22.701 14.605 23.654 /
%
% Fig POLYLINE object
%
\linethickness= 0.500pt
\setplotsymbol ({\thinlinefont .})
\plot 11.113 22.225 13.176 21.431 /
%
% Fig POLYLINE object
%
\linethickness= 0.500pt
\setplotsymbol ({\thinlinefont .})
\plot 13.335 21.431 13.335 21.431 /
%
% Fig POLYLINE object
%
\linethickness= 0.500pt
\setplotsymbol ({\thinlinefont .})
\plot 10.001 21.749 13.652 20.479 /
%
% Fig POLYLINE object
%
\linethickness= 0.500pt
\setplotsymbol ({\thinlinefont .})
\setsolid
\plot 20.161 23.178 21.114 21.907 /
%
% Fig POLYLINE object
%
\linethickness= 0.500pt
\setplotsymbol ({\thinlinefont .})
\setdashes < 0.1270cm>
\plot 17.780 22.860 22.384 23.178 /
%
% Fig POLYLINE object
%
\linethickness= 0.500pt
\setplotsymbol ({\thinlinefont .})
\plot 22.543 23.178 22.543 23.178 /
%
% Fig POLYLINE object
%
\linethickness= 0.500pt
\setplotsymbol ({\thinlinefont .})
\setsolid
\plot  2.857 23.019  2.853 23.017 /
\plot  2.853 23.017  2.845 23.012 /
\plot  2.845 23.012  2.830 23.004 /
\plot  2.830 23.004  2.807 22.991 /
\plot  2.807 22.991  2.775 22.972 /
\plot  2.775 22.972  2.733 22.949 /
\plot  2.733 22.949  2.680 22.921 /
\plot  2.680 22.921  2.620 22.888 /
\plot  2.620 22.888  2.553 22.852 /
\plot  2.553 22.852  2.479 22.809 /
\plot  2.479 22.809  2.400 22.765 /
\plot  2.400 22.765  2.320 22.720 /
\plot  2.320 22.720  2.239 22.674 /
\plot  2.239 22.674  2.159 22.627 /
\plot  2.159 22.627  2.081 22.581 /
\plot  2.081 22.581  2.007 22.536 /
\plot  2.007 22.536  1.937 22.494 /
\plot  1.937 22.494  1.871 22.451 /
\plot  1.871 22.451  1.812 22.413 /
\plot  1.812 22.413  1.761 22.377 /
\plot  1.761 22.377  1.715 22.344 /
\plot  1.715 22.344  1.674 22.312 /
\plot  1.674 22.312  1.640 22.282 /
\plot  1.640 22.282  1.611 22.253 /
\plot  1.611 22.253  1.587 22.225 /
\plot  1.587 22.225  1.564 22.193 /
\plot  1.564 22.193  1.547 22.162 /
\plot  1.547 22.162  1.535 22.128 /
\plot  1.535 22.128  1.526 22.094 /
\plot  1.526 22.094  1.522 22.060 /
\plot  1.522 22.060  1.520 22.026 /
\putrule from  1.520 22.026 to  1.520 21.988
\plot  1.520 21.988  1.522 21.952 /
\plot  1.522 21.952  1.526 21.916 /
\plot  1.526 21.916  1.532 21.878 /
\plot  1.532 21.878  1.541 21.840 /
\plot  1.541 21.840  1.549 21.804 /
\plot  1.549 21.804  1.560 21.770 /
\plot  1.560 21.770  1.571 21.736 /
\plot  1.571 21.736  1.585 21.704 /
\plot  1.585 21.704  1.600 21.677 /
\plot  1.600 21.677  1.617 21.651 /
\plot  1.617 21.651  1.636 21.630 /
\plot  1.636 21.630  1.657 21.613 /
\plot  1.657 21.613  1.683 21.601 /
\plot  1.683 21.601  1.712 21.592 /
\plot  1.712 21.592  1.746 21.590 /
\plot  1.746 21.590  1.778 21.592 /
\plot  1.778 21.592  1.816 21.598 /
\plot  1.816 21.598  1.858 21.609 /
\plot  1.858 21.609  1.907 21.624 /
\plot  1.907 21.624  1.962 21.643 /
\plot  1.962 21.643  2.021 21.666 /
\plot  2.021 21.666  2.089 21.694 /
\plot  2.089 21.694  2.161 21.725 /
\plot  2.161 21.725  2.239 21.761 /
\plot  2.239 21.761  2.322 21.800 /
\plot  2.322 21.800  2.409 21.842 /
\plot  2.409 21.842  2.498 21.884 /
\plot  2.498 21.884  2.589 21.929 /
\plot  2.589 21.929  2.675 21.973 /
\plot  2.675 21.973  2.762 22.015 /
\plot  2.762 22.015  2.843 22.056 /
\plot  2.843 22.056  2.917 22.094 /
\plot  2.917 22.094  2.982 22.126 /
\plot  2.982 22.126  3.037 22.155 /
\plot  3.037 22.155  3.084 22.178 /
\plot  3.084 22.178  3.120 22.197 /
\plot  3.120 22.197  3.145 22.210 /
\plot  3.145 22.210  3.162 22.219 /
\plot  3.162 22.219  3.171 22.223 /
\plot  3.171 22.223  3.175 22.225 /
%
% Fig POLYLINE object
%
\linethickness= 0.500pt
\setplotsymbol ({\thinlinefont .})
\plot  6.032 24.606  6.035 24.608 /
\putrule from  6.035 24.608 to  6.035 24.610
\plot  6.035 24.610  6.037 24.613 /
\plot  6.037 24.613  6.039 24.617 /
\plot  6.039 24.617  6.043 24.621 /
\plot  6.043 24.621  6.047 24.627 /
\plot  6.047 24.627  6.052 24.634 /
\plot  6.052 24.634  6.056 24.642 /
\plot  6.056 24.642  6.062 24.651 /
\plot  6.062 24.651  6.068 24.659 /
\plot  6.068 24.659  6.075 24.670 /
\plot  6.075 24.670  6.083 24.680 /
\plot  6.083 24.680  6.090 24.691 /
\plot  6.090 24.691  6.098 24.701 /
\plot  6.098 24.701  6.104 24.712 /
\plot  6.104 24.712  6.113 24.723 /
\plot  6.113 24.723  6.119 24.733 /
\plot  6.119 24.733  6.126 24.744 /
\plot  6.126 24.744  6.132 24.752 /
\plot  6.132 24.752  6.138 24.761 /
\plot  6.138 24.761  6.143 24.767 /
\plot  6.143 24.767  6.147 24.773 /
\plot  6.147 24.773  6.151 24.778 /
\putrule from  6.151 24.778 to  6.151 24.780
\putrule from  6.151 24.780 to  6.153 24.780
\putrule from  6.153 24.780 to  6.151 24.780
\plot  6.151 24.780  6.149 24.778 /
\plot  6.149 24.778  6.147 24.771 /
\plot  6.147 24.771  6.140 24.763 /
\plot  6.140 24.763  6.134 24.754 /
\plot  6.134 24.754  6.126 24.742 /
\plot  6.126 24.742  6.115 24.727 /
\plot  6.115 24.727  6.102 24.708 /
\plot  6.102 24.708  6.088 24.687 /
\plot  6.088 24.687  6.073 24.663 /
\plot  6.073 24.663  6.054 24.636 /
\plot  6.054 24.636  6.032 24.606 /
\plot  6.032 24.606  6.009 24.570 /
\plot  6.009 24.570  5.982 24.532 /
\plot  5.982 24.532  5.952 24.488 /
\plot  5.952 24.488  5.918 24.439 /
\plot  5.918 24.439  5.880 24.386 /
\plot  5.880 24.386  5.840 24.327 /
\plot  5.840 24.327  5.795 24.263 /
\plot  5.795 24.263  5.749 24.194 /
\plot  5.749 24.194  5.696 24.117 /
\plot  5.696 24.117  5.641 24.037 /
\plot  5.641 24.037  5.582 23.950 /
\plot  5.582 23.950  5.518 23.859 /
\plot  5.518 23.859  5.453 23.764 /
\plot  5.453 23.764  5.385 23.662 /
\plot  5.385 23.662  5.313 23.561 /
\plot  5.313 23.561  5.239 23.453 /
\plot  5.239 23.453  5.165 23.345 /
\plot  5.165 23.345  5.091 23.235 /
\plot  5.091 23.235  5.014 23.125 /
\plot  5.014 23.125  4.938 23.015 /
\plot  4.938 23.015  4.864 22.909 /
\plot  4.864 22.909  4.792 22.803 /
\plot  4.792 22.803  4.724 22.703 /
\plot  4.724 22.703  4.659 22.608 /
\plot  4.659 22.608  4.597 22.519 /
\plot  4.597 22.519  4.542 22.439 /
\plot  4.542 22.439  4.492 22.365 /
\plot  4.492 22.365  4.445 22.299 /
\plot  4.445 22.299  4.407 22.242 /
\plot  4.407 22.242  4.373 22.193 /
\plot  4.373 22.193  4.348 22.155 /
\plot  4.348 22.155  4.326 22.123 /
\plot  4.326 22.123  4.310 22.100 /
\plot  4.310 22.100  4.299 22.085 /
\plot  4.299 22.085  4.293 22.075 /
\plot  4.293 22.075  4.288 22.068 /
\plot  4.288 22.068  4.286 22.066 /
%
% Fig POLYLINE object
%
\linethickness= 0.500pt
\setplotsymbol ({\thinlinefont .})
\putrule from 11.748 23.336 to 11.743 23.336
\plot 11.743 23.336 11.735 23.338 /
\putrule from 11.735 23.338 to 11.718 23.338
\plot 11.718 23.338 11.692 23.340 /
\plot 11.692 23.340 11.654 23.345 /
\plot 11.654 23.345 11.608 23.349 /
\plot 11.608 23.349 11.549 23.353 /
\plot 11.549 23.353 11.481 23.360 /
\plot 11.481 23.360 11.402 23.366 /
\plot 11.402 23.366 11.314 23.372 /
\plot 11.314 23.372 11.220 23.379 /
\plot 11.220 23.379 11.121 23.385 /
\plot 11.121 23.385 11.019 23.391 /
\plot 11.019 23.391 10.918 23.398 /
\plot 10.918 23.398 10.816 23.404 /
\plot 10.816 23.404 10.715 23.408 /
\plot 10.715 23.408 10.619 23.410 /
\plot 10.619 23.410 10.528 23.412 /
\putrule from 10.528 23.412 to 10.444 23.412
\plot 10.444 23.412 10.365 23.410 /
\plot 10.365 23.410 10.293 23.408 /
\plot 10.293 23.408 10.228 23.402 /
\plot 10.228 23.402 10.171 23.393 /
\plot 10.171 23.393 10.118 23.383 /
\plot 10.118 23.383 10.073 23.370 /
\plot 10.073 23.370 10.035 23.355 /
\plot 10.035 23.355 10.001 23.336 /
\plot 10.001 23.336  9.967 23.311 /
\plot  9.967 23.311  9.940 23.281 /
\plot  9.940 23.281  9.917 23.247 /
\plot  9.917 23.247  9.898 23.209 /
\plot  9.898 23.209  9.881 23.167 /
\plot  9.881 23.167  9.868 23.120 /
\plot  9.868 23.120  9.857 23.070 /
\plot  9.857 23.070  9.851 23.017 /
\plot  9.851 23.017  9.845 22.959 /
\plot  9.845 22.959  9.840 22.900 /
\plot  9.840 22.900  9.838 22.841 /
\putrule from  9.838 22.841 to  9.838 22.782
\putrule from  9.838 22.782 to  9.838 22.720
\plot  9.838 22.720  9.840 22.661 /
\plot  9.840 22.661  9.845 22.602 /
\plot  9.845 22.602  9.851 22.545 /
\plot  9.851 22.545  9.857 22.492 /
\plot  9.857 22.492  9.868 22.441 /
\plot  9.868 22.441  9.881 22.394 /
\plot  9.881 22.394  9.898 22.352 /
\plot  9.898 22.352  9.917 22.314 /
\plot  9.917 22.314  9.940 22.280 /
\plot  9.940 22.280  9.967 22.250 /
\plot  9.967 22.250 10.001 22.225 /
\plot 10.001 22.225 10.035 22.206 /
\plot 10.035 22.206 10.073 22.191 /
\plot 10.073 22.191 10.118 22.178 /
\plot 10.118 22.178 10.171 22.168 /
\plot 10.171 22.168 10.228 22.159 /
\plot 10.228 22.159 10.293 22.153 /
\plot 10.293 22.153 10.365 22.151 /
\plot 10.365 22.151 10.444 22.149 /
\putrule from 10.444 22.149 to 10.528 22.149
\plot 10.528 22.149 10.619 22.151 /
\plot 10.619 22.151 10.715 22.153 /
\plot 10.715 22.153 10.816 22.157 /
\plot 10.816 22.157 10.918 22.164 /
\plot 10.918 22.164 11.019 22.170 /
\plot 11.019 22.170 11.121 22.176 /
\plot 11.121 22.176 11.220 22.183 /
\plot 11.220 22.183 11.314 22.189 /
\plot 11.314 22.189 11.402 22.195 /
\plot 11.402 22.195 11.481 22.202 /
\plot 11.481 22.202 11.549 22.208 /
\plot 11.549 22.208 11.608 22.212 /
\plot 11.608 22.212 11.654 22.217 /
\plot 11.654 22.217 11.692 22.221 /
\plot 11.692 22.221 11.718 22.223 /
\putrule from 11.718 22.223 to 11.735 22.223
\plot 11.735 22.223 11.743 22.225 /
\putrule from 11.743 22.225 to 11.748 22.225
%
% Fig POLYLINE object
%
\linethickness= 0.500pt
\setplotsymbol ({\thinlinefont .})
\putrule from 19.209 23.336 to 19.205 23.336
\plot 19.205 23.336 19.192 23.338 /
\plot 19.192 23.338 19.171 23.340 /
\plot 19.171 23.340 19.139 23.343 /
\plot 19.139 23.343 19.094 23.347 /
\plot 19.094 23.347 19.037 23.351 /
\plot 19.037 23.351 18.970 23.357 /
\plot 18.970 23.357 18.891 23.364 /
\plot 18.891 23.364 18.804 23.370 /
\plot 18.804 23.370 18.711 23.376 /
\plot 18.711 23.376 18.614 23.383 /
\plot 18.614 23.383 18.514 23.389 /
\plot 18.514 23.389 18.419 23.396 /
\plot 18.419 23.396 18.324 23.400 /
\plot 18.324 23.400 18.235 23.402 /
\putrule from 18.235 23.402 to 18.153 23.402
\plot 18.153 23.402 18.078 23.400 /
\plot 18.078 23.400 18.009 23.396 /
\plot 18.009 23.396 17.949 23.389 /
\plot 17.949 23.389 17.896 23.381 /
\plot 17.896 23.381 17.852 23.370 /
\plot 17.852 23.370 17.812 23.355 /
\plot 17.812 23.355 17.780 23.336 /
\plot 17.780 23.336 17.750 23.313 /
\plot 17.750 23.313 17.725 23.285 /
\plot 17.725 23.285 17.706 23.252 /
\plot 17.706 23.252 17.689 23.216 /
\plot 17.689 23.216 17.676 23.173 /
\plot 17.676 23.173 17.666 23.129 /
\plot 17.666 23.129 17.659 23.078 /
\plot 17.659 23.078 17.653 23.027 /
\plot 17.653 23.027 17.649 22.972 /
\plot 17.649 22.972 17.647 22.917 /
\putrule from 17.647 22.917 to 17.647 22.860
\putrule from 17.647 22.860 to 17.647 22.803
\plot 17.647 22.803 17.649 22.748 /
\plot 17.649 22.748 17.653 22.693 /
\plot 17.653 22.693 17.659 22.642 /
\plot 17.659 22.642 17.666 22.591 /
\plot 17.666 22.591 17.676 22.547 /
\plot 17.676 22.547 17.689 22.504 /
\plot 17.689 22.504 17.706 22.468 /
\plot 17.706 22.468 17.725 22.435 /
\plot 17.725 22.435 17.750 22.407 /
\plot 17.750 22.407 17.780 22.384 /
\plot 17.780 22.384 17.812 22.365 /
\plot 17.812 22.365 17.852 22.350 /
\plot 17.852 22.350 17.896 22.339 /
\plot 17.896 22.339 17.949 22.331 /
\plot 17.949 22.331 18.009 22.324 /
\plot 18.009 22.324 18.078 22.320 /
\plot 18.078 22.320 18.153 22.318 /
\putrule from 18.153 22.318 to 18.235 22.318
\plot 18.235 22.318 18.324 22.320 /
\plot 18.324 22.320 18.419 22.324 /
\plot 18.419 22.324 18.514 22.331 /
\plot 18.514 22.331 18.614 22.335 /
\plot 18.614 22.335 18.711 22.344 /
\plot 18.711 22.344 18.804 22.350 /
\plot 18.804 22.350 18.891 22.356 /
\plot 18.891 22.356 18.970 22.363 /
\plot 18.970 22.363 19.037 22.369 /
\plot 19.037 22.369 19.094 22.373 /
\plot 19.094 22.373 19.139 22.377 /
\plot 19.139 22.377 19.171 22.380 /
\plot 19.171 22.380 19.192 22.382 /
\plot 19.192 22.382 19.205 22.384 /
\putrule from 19.205 22.384 to 19.209 22.384
%
% Fig POLYLINE object
%
\linethickness= 0.500pt
\setplotsymbol ({\thinlinefont .})
\setdashes < 0.1270cm>
\plot 18.733 22.543 20.796 21.273 /
%
% Fig POLYLINE object
%
\linethickness= 0.500pt
\setplotsymbol ({\thinlinefont .})
\setdots < 0.0953cm>
\plot 21.590 23.971 21.586 23.971 /
\plot 21.586 23.971 21.579 23.973 /
\plot 21.579 23.973 21.567 23.975 /
\plot 21.567 23.975 21.543 23.978 /
\plot 21.543 23.978 21.514 23.982 /
\plot 21.514 23.982 21.471 23.988 /
\plot 21.471 23.988 21.421 23.995 /
\plot 21.421 23.995 21.357 24.005 /
\plot 21.357 24.005 21.281 24.014 /
\plot 21.281 24.014 21.194 24.026 /
\plot 21.194 24.026 21.097 24.039 /
\plot 21.097 24.039 20.987 24.054 /
\plot 20.987 24.054 20.870 24.069 /
\plot 20.870 24.069 20.743 24.083 /
\plot 20.743 24.083 20.610 24.098 /
\plot 20.610 24.098 20.472 24.115 /
\plot 20.472 24.115 20.331 24.130 /
\plot 20.331 24.130 20.185 24.145 /
\plot 20.185 24.145 20.038 24.160 /
\plot 20.038 24.160 19.895 24.172 /
\plot 19.895 24.172 19.751 24.185 /
\plot 19.751 24.185 19.609 24.194 /
\plot 19.609 24.194 19.473 24.202 /
\plot 19.473 24.202 19.342 24.208 /
\plot 19.342 24.208 19.215 24.210 /
\plot 19.215 24.210 19.097 24.210 /
\plot 19.097 24.210 18.982 24.208 /
\plot 18.982 24.208 18.876 24.202 /
\plot 18.876 24.202 18.779 24.194 /
\plot 18.779 24.194 18.688 24.181 /
\plot 18.688 24.181 18.605 24.164 /
\plot 18.605 24.164 18.529 24.143 /
\plot 18.529 24.143 18.462 24.117 /
\plot 18.462 24.117 18.400 24.088 /
\plot 18.400 24.088 18.347 24.054 /
\plot 18.347 24.054 18.299 24.016 /
\plot 18.299 24.016 18.256 23.971 /
\plot 18.256 23.971 18.220 23.925 /
\plot 18.220 23.925 18.191 23.874 /
\plot 18.191 23.874 18.163 23.819 /
\plot 18.163 23.819 18.140 23.760 /
\plot 18.140 23.760 18.121 23.694 /
\plot 18.121 23.694 18.104 23.624 /
\plot 18.104 23.624 18.091 23.548 /
\plot 18.091 23.548 18.083 23.470 /
\plot 18.083 23.470 18.076 23.385 /
\plot 18.076 23.385 18.072 23.296 /
\plot 18.072 23.296 18.070 23.203 /
\plot 18.070 23.203 18.072 23.106 /
\plot 18.072 23.106 18.074 23.006 /
\plot 18.074 23.006 18.078 22.902 /
\plot 18.078 22.902 18.087 22.797 /
\plot 18.087 22.797 18.093 22.689 /
\plot 18.093 22.689 18.104 22.578 /
\plot 18.104 22.578 18.114 22.466 /
\plot 18.114 22.466 18.127 22.354 /
\plot 18.127 22.354 18.140 22.240 /
\plot 18.140 22.240 18.155 22.128 /
\plot 18.155 22.128 18.169 22.013 /
\plot 18.169 22.013 18.184 21.901 /
\plot 18.184 21.901 18.201 21.791 /
\plot 18.201 21.791 18.216 21.683 /
\plot 18.216 21.683 18.235 21.575 /
\plot 18.235 21.575 18.252 21.471 /
\plot 18.252 21.471 18.271 21.370 /
\plot 18.271 21.370 18.290 21.273 /
\plot 18.290 21.273 18.311 21.179 /
\plot 18.311 21.179 18.330 21.088 /
\plot 18.330 21.088 18.354 21.002 /
\plot 18.354 21.002 18.375 20.921 /
\plot 18.375 20.921 18.398 20.843 /
\plot 18.398 20.843 18.423 20.771 /
\plot 18.423 20.771 18.449 20.703 /
\plot 18.449 20.703 18.479 20.642 /
\plot 18.479 20.642 18.508 20.582 /
\plot 18.508 20.582 18.540 20.527 /
\plot 18.540 20.527 18.574 20.479 /
\plot 18.574 20.479 18.618 20.426 /
\plot 18.618 20.426 18.667 20.377 /
\plot 18.667 20.377 18.720 20.335 /
\plot 18.720 20.335 18.779 20.297 /
\plot 18.779 20.297 18.845 20.263 /
\plot 18.845 20.263 18.917 20.233 /
\plot 18.917 20.233 18.995 20.210 /
\plot 18.995 20.210 19.080 20.187 /
\plot 19.080 20.187 19.171 20.170 /
\plot 19.171 20.170 19.268 20.155 /
\plot 19.268 20.155 19.370 20.142 /
\plot 19.370 20.142 19.480 20.134 /
\plot 19.480 20.134 19.592 20.127 /
\plot 19.592 20.127 19.708 20.123 /
\plot 19.708 20.123 19.829 20.119 /
\plot 19.829 20.119 19.950 20.119 /
\plot 19.950 20.119 20.072 20.119 /
\plot 20.072 20.119 20.191 20.121 /
\plot 20.191 20.121 20.309 20.123 /
\plot 20.309 20.123 20.424 20.127 /
\plot 20.424 20.127 20.532 20.130 /
\plot 20.532 20.130 20.633 20.134 /
\plot 20.633 20.134 20.726 20.138 /
\plot 20.726 20.138 20.809 20.142 /
\plot 20.809 20.142 20.881 20.146 /
\plot 20.881 20.146 20.944 20.151 /
\plot 20.944 20.151 20.995 20.153 /
\plot 20.995 20.153 21.038 20.157 /
\plot 21.038 20.157 21.067 20.159 /
\plot 21.067 20.159 21.088 20.159 /
\plot 21.088 20.159 21.103 20.161 /
\plot 21.103 20.161 21.110 20.161 /
\plot 21.110 20.161 21.114 20.161 /
%
% Fig POLYLINE object
%
\linethickness= 0.500pt
\setplotsymbol ({\thinlinefont .})
\setdashes < 0.1270cm>
\plot 21.273 21.590 20.479 19.685 /
%
% Fig TEXT object
%
\put{\SetFigFont{7}{8.4}{rm}$(I_1,I_1,I_2,I_8)$} [lB] at  2.381 25.400
%
% Fig TEXT object
%
\put{\SetFigFont{7}{8.4}{rm}$U_1^{12}(A_1+A_7)$} [lB] at 10.795 25.241
%
% Fig TEXT object
%
\put{\SetFigFont{7}{8.4}{rm}$U_2^{12}(A_1+2A_3)$} [lB] at 18.574 25.400
\linethickness=0pt
\putrectangle corners at  1.494 25.813 and 23.362 19.660
\endpicture}
}
  \vskip 0.5cm
  \begin{center}
   Figure $4$
  \end{center}
{\bf Configurations of the Singularity types}
\vskip 0.3cm
\hskip -0.6cm $(a)$ The Picard number one Gorenstein log del Pezzo surfaces.
\vskip 0.3cm
\hskip -0.6cm
The following configurations are the figures of the negative curves on $U$ 
which is the minimal resolution of the Picard number one Gorenstein log del Pezzo surface and $2\leq K_U^2 \leq 7$.

\vskip 0.5cm 
\centerline{\font\thinlinefont=cmr5
\begingroup\makeatletter\ifx\SetFigFont\undefined%
\gdef\SetFigFont#1#2#3#4#5{%
  \reset@font\fontsize{#1}{#2pt}%
  \fontfamily{#3}\fontseries{#4}\fontshape{#5}%
  \selectfont}%
\fi\endgroup%
\mbox{\beginpicture
\setcoordinatesystem units <0.78000cm,0.78000cm>
\unitlength=0.78000cm
\linethickness=1pt
\setplotsymbol ({\makebox(0,0)[l]{\tencirc\symbol{'160}}})
\setshadesymbol ({\thinlinefont .})
\setlinear
%
% Fig POLYLINE object
%
\linethickness= 0.500pt
\setplotsymbol ({\thinlinefont .})
\putrule from  2.381 20.003 to  6.032 20.003
%
% Fig POLYLINE object
%
\linethickness= 0.500pt
\setplotsymbol ({\thinlinefont .})
\putrule from  3.016 21.114 to  3.016 19.526
%
% Fig POLYLINE object
%
\linethickness= 0.500pt
\setplotsymbol ({\thinlinefont .})
\putrule from  3.969 21.114 to  3.969 19.526
%
% Fig POLYLINE object
%
\linethickness= 0.500pt
\setplotsymbol ({\thinlinefont .})
\putrule from  3.651 20.796 to  4.604 20.796
%
% Fig POLYLINE object
%
\linethickness= 0.500pt
\setplotsymbol ({\thinlinefont .})
\putrule from  5.397 21.114 to  5.397 19.526
%
% Fig POLYLINE object
%
\linethickness= 0.500pt
\setplotsymbol ({\thinlinefont .})
\putrule from  5.080 20.796 to  6.509 20.796
%
% Fig POLYLINE object
%
\linethickness= 0.500pt
\setplotsymbol ({\thinlinefont .})
\putrule from  6.032 21.907 to  6.032 20.637
%
% Fig POLYLINE object
%
\linethickness= 0.500pt
\setplotsymbol ({\thinlinefont .})
\setdashes < 0.1270cm>
\plot  5.874 21.749  6.826 21.749 /
%
% Fig POLYLINE object
%
\linethickness= 0.500pt
\setplotsymbol ({\thinlinefont .})
\setsolid
\putrule from  8.890 21.114 to  8.890 19.685
%
% Fig POLYLINE object
%
\linethickness= 0.500pt
\setplotsymbol ({\thinlinefont .})
\putrule from  9.684 21.114 to  9.684 19.685
%
% Fig POLYLINE object
%
\linethickness= 0.500pt
\setplotsymbol ({\thinlinefont .})
\putrule from  9.366 20.796 to 10.954 20.796
%
% Fig POLYLINE object
%
\linethickness= 0.500pt
\setplotsymbol ({\thinlinefont .})
\putrule from 10.636 21.114 to 10.636 19.526
%
% Fig POLYLINE object
%
\linethickness= 0.500pt
\setplotsymbol ({\thinlinefont .})
\putrule from 10.319 20.003 to 12.859 20.003
%
% Fig POLYLINE object
%
\linethickness= 0.500pt
\setplotsymbol ({\thinlinefont .})
\putrule from 11.589 21.114 to 11.589 19.526
%
% Fig POLYLINE object
%
\linethickness= 0.500pt
\setplotsymbol ({\thinlinefont .})
\putrule from 12.383 21.114 to 12.383 19.526
%
% Fig POLYLINE object
%
\linethickness= 0.500pt
\setplotsymbol ({\thinlinefont .})
\setdashes < 0.1270cm>
\plot 12.224 20.796 13.176 20.796 /
%
% Fig POLYLINE object
%
\linethickness= 0.500pt
\setplotsymbol ({\thinlinefont .})
\plot 13.176 20.796 13.176 20.796 /
%
% Fig POLYLINE object
%
\linethickness= 0.500pt
\setplotsymbol ({\thinlinefont .})
\plot  8.414 20.003  9.842 20.003 /
%
% Fig POLYLINE object
%
\linethickness= 0.500pt
\setplotsymbol ({\thinlinefont .})
\setsolid
\putrule from 15.399 20.003 to 17.462 20.003
%
% Fig POLYLINE object
%
\linethickness= 0.500pt
\setplotsymbol ({\thinlinefont .})
\putrule from 18.256 20.003 to 20.320 20.003
%
% Fig POLYLINE object
%
\linethickness= 0.500pt
\setplotsymbol ({\thinlinefont .})
\putrule from 16.192 21.114 to 16.192 19.526
%
% Fig POLYLINE object
%
\linethickness= 0.500pt
\setplotsymbol ({\thinlinefont .})
\putrule from 16.828 21.114 to 16.828 19.526
%
% Fig POLYLINE object
%
\linethickness= 0.500pt
\setplotsymbol ({\thinlinefont .})
\putrule from 18.415 21.273 to 18.415 19.526
%
% Fig POLYLINE object
%
\linethickness= 0.500pt
\setplotsymbol ({\thinlinefont .})
\putrule from 19.367 21.273 to 19.367 19.526
%
% Fig POLYLINE object
%
\linethickness= 0.500pt
\setplotsymbol ({\thinlinefont .})
\putrule from 16.669 20.955 to 18.891 20.955
%
% Fig POLYLINE object
%
\linethickness= 0.500pt
\setplotsymbol ({\thinlinefont .})
\setdashes < 0.1270cm>
\plot 15.716 20.955 15.716 19.526 /
%
% Fig POLYLINE object
%
\linethickness= 0.500pt
\setplotsymbol ({\thinlinefont .})
\plot 19.844 21.114 19.844 19.526 /
%
% Fig TEXT object
%
\put{\SetFigFont{9}{10.8}{rm}$E_7$} [lB] at  3.810 18.415
%
% Fig TEXT object
%
\put{\SetFigFont{9}{10.8}{rm}$A_1+D_6$} [lB] at  9.684 18.415
%
% Fig TEXT object
%
\put{\SetFigFont{9}{10.8}{rm}$A_7$} [lB] at 16.828 18.415
\linethickness=0pt
\putrectangle corners at  2.356 21.933 and 20.345 18.415
\endpicture}
}
\vskip 0.5cm 
\centerline{\font\thinlinefont=cmr5
\begingroup\makeatletter\ifx\SetFigFont\undefined%
\gdef\SetFigFont#1#2#3#4#5{%
  \reset@font\fontsize{#1}{#2pt}%
  \fontfamily{#3}\fontseries{#4}\fontshape{#5}%
  \selectfont}%
\fi\endgroup%
\mbox{\beginpicture
\setcoordinatesystem units <0.78000cm,0.78000cm>
\unitlength=0.78000cm
\linethickness=1pt
\setplotsymbol ({\makebox(0,0)[l]{\tencirc\symbol{'160}}})
\setshadesymbol ({\thinlinefont .})
\setlinear
%
% Fig POLYLINE object
%
\linethickness= 0.500pt
\setplotsymbol ({\thinlinefont .})
\putrule from  3.016 21.431 to  3.016 18.415
%
% Fig POLYLINE object
%
\linethickness= 0.500pt
\setplotsymbol ({\thinlinefont .})
\plot  3.969 20.955  2.540 19.050 /
%
% Fig POLYLINE object
%
\linethickness= 0.500pt
\setplotsymbol ({\thinlinefont .})
\putrule from  4.921 19.050 to  7.461 19.050
%
% Fig POLYLINE object
%
\linethickness= 0.500pt
\setplotsymbol ({\thinlinefont .})
\putrule from  5.397 20.003 to  5.397 18.415
%
% Fig POLYLINE object
%
\linethickness= 0.500pt
\setplotsymbol ({\thinlinefont .})
\plot  4.604 20.955  5.715 19.367 /
%
% Fig POLYLINE object
%
\linethickness= 0.500pt
\setplotsymbol ({\thinlinefont .})
\putrule from  6.350 20.637 to  6.350 18.415
%
% Fig POLYLINE object
%
\linethickness= 0.500pt
\setplotsymbol ({\thinlinefont .})
\plot  5.239 21.749  6.668 19.844 /
%
% Fig POLYLINE object
%
\linethickness= 0.500pt
\setplotsymbol ({\thinlinefont .})
\setdashes < 0.1270cm>
\plot  7.144 20.479  7.144 18.415 /
%
% Fig POLYLINE object
%
\linethickness= 0.500pt
\setplotsymbol ({\thinlinefont .})
\plot  3.493 20.637  5.080 20.637 /
%
% Fig POLYLINE object
%
\linethickness= 0.500pt
\setplotsymbol ({\thinlinefont .})
\plot  2.699 21.273  6.032 21.273 /
%
% Fig POLYLINE object
%
\linethickness= 0.500pt
\setplotsymbol ({\thinlinefont .})
\setsolid
\putrule from 12.065 19.526 to 14.446 19.526
%
% Fig POLYLINE object
%
\linethickness= 0.500pt
\setplotsymbol ({\thinlinefont .})
\putrule from 11.113 20.479 to 11.113 19.209
%
% Fig POLYLINE object
%
\linethickness= 0.500pt
\setplotsymbol ({\thinlinefont .})
\putrule from 12.224 20.479 to 12.224 19.209
%
% Fig POLYLINE object
%
\linethickness= 0.500pt
\setplotsymbol ({\thinlinefont .})
\putrule from 10.319 21.114 to 10.319 19.209
%
% Fig POLYLINE object
%
\linethickness= 0.500pt
\setplotsymbol ({\thinlinefont .})
\putrule from 13.018 21.431 to 13.018 19.209
%
% Fig POLYLINE object
%
\linethickness= 0.500pt
\setplotsymbol ({\thinlinefont .})
\putrule from  9.684 20.637 to  9.684 18.415
%
% Fig POLYLINE object
%
\linethickness= 0.500pt
\setplotsymbol ({\thinlinefont .})
\putrule from 13.811 20.796 to 13.811 18.574
%
% Fig POLYLINE object
%
\linethickness= 0.500pt
\setplotsymbol ({\thinlinefont .})
\setdashes < 0.1270cm>
\plot  9.366 19.526 11.589 19.526 /
%
% Fig POLYLINE object
%
\linethickness= 0.500pt
\setplotsymbol ({\thinlinefont .})
\plot 10.795 20.161 12.541 20.161 /
%
% Fig POLYLINE object
%
\linethickness= 0.500pt
\setplotsymbol ({\thinlinefont .})
\plot 10.001 20.955 13.335 20.955 /
%
% Fig POLYLINE object
%
\linethickness= 0.500pt
\setplotsymbol ({\thinlinefont .})
\plot  9.366 18.891 14.287 18.891 /
%
% Fig POLYLINE object
%
\linethickness= 0.500pt
\setplotsymbol ({\thinlinefont .})
\setsolid
\plot 17.780 21.590 16.828 20.320 /
%
% Fig POLYLINE object
%
\linethickness= 0.500pt
\setplotsymbol ({\thinlinefont .})
\putrule from 16.986 20.796 to 16.986 19.050
%
% Fig POLYLINE object
%
\linethickness= 0.500pt
\setplotsymbol ({\thinlinefont .})
\plot 16.828 19.367 18.098 18.256 /
%
% Fig POLYLINE object
%
\linethickness= 0.500pt
\setplotsymbol ({\thinlinefont .})
\plot 19.209 21.749 20.479 20.479 /
%
% Fig POLYLINE object
%
\linethickness= 0.500pt
\setplotsymbol ({\thinlinefont .})
\putrule from 20.320 20.955 to 20.320 19.209
%
% Fig POLYLINE object
%
\linethickness= 0.500pt
\setplotsymbol ({\thinlinefont .})
\plot 20.479 19.685 19.526 18.415 /
%
% Fig POLYLINE object
%
\linethickness= 0.500pt
\setplotsymbol ({\thinlinefont .})
\putrule from 18.574 20.796 to 18.574 19.050
%
% Fig POLYLINE object
%
\linethickness= 0.500pt
\setplotsymbol ({\thinlinefont .})
\setdashes < 0.1270cm>
\plot 17.304 21.431 20.003 21.431 /
%
% Fig POLYLINE object
%
\linethickness= 0.500pt
\setplotsymbol ({\thinlinefont .})
\plot 17.462 18.574 20.161 18.574 /
%
% Fig POLYLINE object
%
\linethickness= 0.500pt
\setplotsymbol ({\thinlinefont .})
\plot 16.669 20.161 18.891 20.161 /
%
% Fig POLYLINE object
%
\linethickness= 0.500pt
\setplotsymbol ({\thinlinefont .})
\plot 18.415 19.844 20.637 19.844 /
%
% Fig TEXT object
%
\put{\SetFigFont{9}{10.8}{rm}$A_2+A_5$} [lB] at  3.969 17.145
%
% Fig TEXT object
%
\put{\SetFigFont{9}{10.8}{rm}$3A_1+D_4$} [lB] at 10.954 17.304
%
% Fig TEXT object
%
\put{\SetFigFont{9}{10.8}{rm}$A_1+2A_3$} [lB] at 17.939 17.304
\linethickness=0pt
\putrectangle corners at  2.515 21.774 and 20.663 17.145
\endpicture}
}
\vskip 0.5cm
\centerline{\font\thinlinefont=cmr5
\begingroup\makeatletter\ifx\SetFigFont\undefined%
\gdef\SetFigFont#1#2#3#4#5{%
  \reset@font\fontsize{#1}{#2pt}%
  \fontfamily{#3}\fontseries{#4}\fontshape{#5}%
  \selectfont}%
\fi\endgroup%
\mbox{\beginpicture
\setcoordinatesystem units <0.78000cm,0.78000cm>
\unitlength=0.78000cm
\linethickness=1pt
\setplotsymbol ({\makebox(0,0)[l]{\tencirc\symbol{'160}}})
\setshadesymbol ({\thinlinefont .})
\setlinear
%
% Fig POLYLINE object
%
\linethickness= 0.500pt
\setplotsymbol ({\thinlinefont .})
\plot  3.016 22.066  1.746 20.955 /
%
% Fig POLYLINE object
%
\linethickness= 0.500pt
\setplotsymbol ({\thinlinefont .})
\putrule from  1.905 21.431 to  1.905 18.891
%
% Fig POLYLINE object
%
\linethickness= 0.500pt
\setplotsymbol ({\thinlinefont .})
\plot  3.651 22.066  4.763 20.955 /
%
% Fig POLYLINE object
%
\linethickness= 0.500pt
\setplotsymbol ({\thinlinefont .})
\plot  4.763 21.273  3.493 20.161 /
%
% Fig POLYLINE object
%
\linethickness= 0.500pt
\setplotsymbol ({\thinlinefont .})
\plot  5.715 20.955  4.445 19.685 /
%
% Fig POLYLINE object
%
\linethickness= 0.500pt
\setplotsymbol ({\thinlinefont .})
\putrule from  4.921 20.479 to  4.921 18.891
%
% Fig POLYLINE object
%
\linethickness= 0.500pt
\setplotsymbol ({\thinlinefont .})
\setdashes < 0.1270cm>
\plot  2.381 21.907  4.286 21.907 /
%
% Fig POLYLINE object
%
\linethickness= 0.500pt
\setplotsymbol ({\thinlinefont .})
\plot  1.746 19.209  5.397 19.209 /
%
% Fig POLYLINE object
%
\linethickness= 0.500pt
\setplotsymbol ({\thinlinefont .})
\plot  3.651 20.796  4.763 19.685 /
%
% Fig POLYLINE object
%
\linethickness= 0.500pt
\setplotsymbol ({\thinlinefont .})
\setsolid
\putrule from  8.572 21.431 to  8.572 19.526
%
% Fig POLYLINE object
%
\linethickness= 0.500pt
\setplotsymbol ({\thinlinefont .})
\putrule from  9.684 21.431 to  9.684 19.526
%
% Fig POLYLINE object
%
\linethickness= 0.500pt
\setplotsymbol ({\thinlinefont .})
\putrule from  9.366 21.114 to 11.271 21.114
%
% Fig POLYLINE object
%
\linethickness= 0.500pt
\setplotsymbol ({\thinlinefont .})
\putrule from 10.954 21.431 to 10.954 19.526
%
% Fig POLYLINE object
%
\linethickness= 0.500pt
\setplotsymbol ({\thinlinefont .})
\putrule from 10.636 20.003 to 12.700 20.003
%
% Fig POLYLINE object
%
\linethickness= 0.500pt
\setplotsymbol ({\thinlinefont .})
\putrule from 11.748 21.431 to 11.748 19.526
%
% Fig POLYLINE object
%
\linethickness= 0.500pt
\setplotsymbol ({\thinlinefont .})
\setdashes < 0.1270cm>
\plot  8.096 20.003 10.160 20.003 /
%
% Fig POLYLINE object
%
\linethickness= 0.500pt
\setplotsymbol ({\thinlinefont .})
\plot 12.224 21.431 12.224 19.526 /
%
% Fig POLYLINE object
%
\linethickness= 0.500pt
\setplotsymbol ({\thinlinefont .})
\setsolid
\putrule from 14.764 20.003 to 17.939 20.003
%
% Fig POLYLINE object
%
\linethickness= 0.500pt
\setplotsymbol ({\thinlinefont .})
\putrule from 15.081 21.273 to 15.081 19.526
%
% Fig POLYLINE object
%
\linethickness= 0.500pt
\setplotsymbol ({\thinlinefont .})
\putrule from 15.716 21.273 to 15.716 19.526
%
% Fig POLYLINE object
%
\linethickness= 0.500pt
\setplotsymbol ({\thinlinefont .})
\putrule from 16.828 21.273 to 16.828 19.526
%
% Fig POLYLINE object
%
\linethickness= 0.500pt
\setplotsymbol ({\thinlinefont .})
\putrule from 15.399 20.955 to 16.192 20.955
%
% Fig POLYLINE object
%
\linethickness= 0.500pt
\setplotsymbol ({\thinlinefont .})
\putrule from 16.510 20.955 to 17.780 20.955
%
% Fig POLYLINE object
%
\linethickness= 0.500pt
\setplotsymbol ({\thinlinefont .})
\setdashes < 0.1270cm>
\plot 17.462 21.907 17.462 20.479 /
%
% Fig TEXT object
%
\put{\SetFigFont{9}{10.8}{rm}$3A_2$} [lB] at  2.699 17.621
%
% Fig TEXT object
%
\put{\SetFigFont{9}{10.8}{rm}$A_1+A_5$} [lB] at  9.525 17.939
%
% Fig TEXT object
%
\put{\SetFigFont{9}{10.8}{rm}$E_6$} [lB] at 15.558 18.098
\linethickness=0pt
\putrectangle corners at  1.721 22.092 and 17.964 17.621
\endpicture}
}
\vskip 0.5cm 
\centerline{\font\thinlinefont=cmr5
\begingroup\makeatletter\ifx\SetFigFont\undefined%
\gdef\SetFigFont#1#2#3#4#5{%
  \reset@font\fontsize{#1}{#2pt}%
  \fontfamily{#3}\fontseries{#4}\fontshape{#5}%
  \selectfont}%
\fi\endgroup%
\mbox{\beginpicture
\setcoordinatesystem units <0.78000cm,0.78000cm>
\unitlength=0.78000cm
\linethickness=1pt
\setplotsymbol ({\makebox(0,0)[l]{\tencirc\symbol{'160}}})
\setshadesymbol ({\thinlinefont .})
\setlinear
%
% Fig POLYLINE object
%
\linethickness= 0.500pt
\setplotsymbol ({\thinlinefont .})
\putrule from  1.746 20.955 to  1.746 19.050
%
% Fig POLYLINE object
%
\linethickness= 0.500pt
\setplotsymbol ({\thinlinefont .})
\putrule from  2.540 20.955 to  2.540 19.050
%
% Fig POLYLINE object
%
\linethickness= 0.500pt
\setplotsymbol ({\thinlinefont .})
\putrule from  2.223 19.526 to  3.651 19.526
%
% Fig POLYLINE object
%
\linethickness= 0.500pt
\setplotsymbol ({\thinlinefont .})
\putrule from  3.334 20.955 to  3.334 19.050
%
% Fig POLYLINE object
%
\linethickness= 0.500pt
\setplotsymbol ({\thinlinefont .})
\putrule from  4.445 20.955 to  4.445 19.050
%
% Fig POLYLINE object
%
\linethickness= 0.500pt
\setplotsymbol ({\thinlinefont .})
\setdashes < 0.1270cm>
\plot  1.429 20.479  2.857 20.479 /
%
% Fig POLYLINE object
%
\linethickness= 0.500pt
\setplotsymbol ({\thinlinefont .})
\plot  3.016 20.003  4.763 20.003 /
%
% Fig POLYLINE object
%
\linethickness= 0.500pt
\setplotsymbol ({\thinlinefont .})
\setsolid
\putrule from  6.826 19.526 to  9.207 19.526
%
% Fig POLYLINE object
%
\linethickness= 0.500pt
\setplotsymbol ({\thinlinefont .})
\putrule from  6.985 20.955 to  6.985 19.209
%
% Fig POLYLINE object
%
\linethickness= 0.500pt
\setplotsymbol ({\thinlinefont .})
\putrule from  7.779 21.114 to  7.779 19.209
%
% Fig POLYLINE object
%
\linethickness= 0.500pt
\setplotsymbol ({\thinlinefont .})
\putrule from  7.461 20.796 to  8.255 20.796
%
% Fig POLYLINE object
%
\linethickness= 0.500pt
\setplotsymbol ({\thinlinefont .})
\putrule from  8.731 21.114 to  8.731 19.209
%
% Fig POLYLINE object
%
\linethickness= 0.500pt
\setplotsymbol ({\thinlinefont .})
\setdashes < 0.1270cm>
\plot  8.572 20.796  9.684 20.796 /
%
% Fig POLYLINE object
%
\linethickness= 0.500pt
\setplotsymbol ({\thinlinefont .})
\setsolid
\putrule from 11.430 19.526 to 13.652 19.526
%
% Fig POLYLINE object
%
\linethickness= 0.500pt
\setplotsymbol ({\thinlinefont .})
\putrule from 11.748 20.955 to 11.748 19.050
%
% Fig POLYLINE object
%
\linethickness= 0.500pt
\setplotsymbol ({\thinlinefont .})
\putrule from 12.383 20.955 to 12.383 19.050
%
% Fig POLYLINE object
%
\linethickness= 0.500pt
\setplotsymbol ({\thinlinefont .})
\plot 12.383 19.209 12.383 19.209 /
%
% Fig POLYLINE object
%
\linethickness= 0.500pt
\setplotsymbol ({\thinlinefont .})
\putrule from 12.065 20.637 to 12.859 20.637
%
% Fig POLYLINE object
%
\linethickness= 0.500pt
\setplotsymbol ({\thinlinefont .})
\setdashes < 0.1270cm>
\plot 13.335 20.955 13.335 19.050 /
%
% Fig POLYLINE object
%
\linethickness= 0.500pt
\setplotsymbol ({\thinlinefont .})
\setsolid
\putrule from 15.875 21.114 to 15.875 19.209
%
% Fig POLYLINE object
%
\linethickness= 0.500pt
\setplotsymbol ({\thinlinefont .})
\putrule from 16.986 21.114 to 16.986 19.209
%
% Fig POLYLINE object
%
\linethickness= 0.500pt
\setplotsymbol ({\thinlinefont .})
\putrule from 16.510 20.637 to 17.939 20.637
%
% Fig POLYLINE object
%
\linethickness= 0.500pt
\setplotsymbol ({\thinlinefont .})
\setdashes < 0.1270cm>
\plot 15.399 19.685 17.780 19.685 /
%
% Fig TEXT object
%
\put{\SetFigFont{9}{10.8}{rm}$2A_1+A_3$} [lB] at  2.540 17.939
%
% Fig TEXT object
%
\put{\SetFigFont{9}{10.8}{rm}$D_5$} [lB] at  7.461 18.098
%
% Fig TEXT object
%
\put{\SetFigFont{9}{10.8}{rm}$A_4$} [lB] at 12.065 17.939
%
% Fig TEXT object
%
\put{\SetFigFont{9}{10.8}{rm}$A_2+A_1$} [lB] at 16.034 18.098
\linethickness=0pt
\putrectangle corners at  1.403 21.139 and 17.964 17.939
\endpicture}
}
\begin{center}
Figure $5$
\end{center}
\vskip 0.3cm
$(b)$ The Picard number two Gorenstein log del Pezzo surface.
\vskip 0.3cm
\hskip -0.6cm
The following configurations are the figures of the negative curves on $U$ 
which is the minimal resolution of the Picard number two Gorenstein log del Pezzo surface and $K_U^2 \geq 3$.

\vskip 0.5cm 
\centerline{\font\thinlinefont=cmr5
\begingroup\makeatletter\ifx\SetFigFont\undefined%
\gdef\SetFigFont#1#2#3#4#5{%
  \reset@font\fontsize{#1}{#2pt}%
  \fontfamily{#3}\fontseries{#4}\fontshape{#5}%
  \selectfont}%
\fi\endgroup%
\mbox
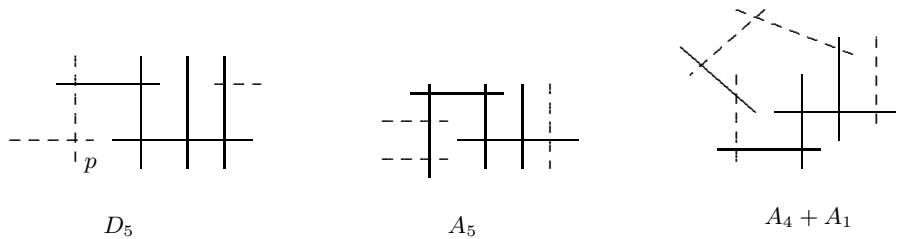
{\beginpicture
\setcoordinatesystem units <0.78000cm,0.78000cm>
\unitlength=0.78000cm
\linethickness=1pt
\setplotsymbol ({\makebox(0,0)[l]{\tencirc\symbol{'160}}})
\setshadesymbol ({\thinlinefont .})
\setlinear
%
% Fig POLYLINE object
%
\linethickness= 0.500pt
\setplotsymbol ({\thinlinefont .})
\putrule from  5.239 21.907 to  7.620 21.907
%
% Fig POLYLINE object
%
\linethickness= 0.500pt
\setplotsymbol ({\thinlinefont .})
\putrule from  5.715 23.336 to  5.715 21.431
%
% Fig POLYLINE object
%
\linethickness= 0.500pt
\setplotsymbol ({\thinlinefont .})
\putrule from  6.509 23.336 to  6.509 21.431
%
% Fig POLYLINE object
%
\linethickness= 0.500pt
\setplotsymbol ({\thinlinefont .})
\putrule from  7.144 23.336 to  7.144 21.431
%
% Fig POLYLINE object
%
\linethickness= 0.500pt
\setplotsymbol ({\thinlinefont .})
\putrule from  4.286 22.860 to  6.032 22.860
%
% Fig POLYLINE object
%
\linethickness= 0.500pt
\setplotsymbol ({\thinlinefont .})
\setdashes < 0.1270cm>
\plot  6.985 22.860  7.779 22.860 /
%
% Fig POLYLINE object
%
\linethickness= 0.500pt
\setplotsymbol ({\thinlinefont .})
\plot  4.604 23.336  4.604 21.431 /
%
% Fig POLYLINE object
%
\linethickness= 0.500pt
\setplotsymbol ({\thinlinefont .})
\plot  3.493 21.907  4.921 21.907 /
%
% Fig POLYLINE object
%
\linethickness= 0.500pt
\setplotsymbol ({\thinlinefont .})
\setsolid
\putrule from 11.113 21.907 to 13.176 21.907
%
% Fig POLYLINE object
%
\linethickness= 0.500pt
\setplotsymbol ({\thinlinefont .})
\putrule from 12.224 22.860 to 12.224 21.431
%
% Fig POLYLINE object
%
\linethickness= 0.500pt
\setplotsymbol ({\thinlinefont .})
\putrule from 11.589 22.860 to 11.589 21.431
%
% Fig POLYLINE object
%
\linethickness= 0.500pt
\setplotsymbol ({\thinlinefont .})
\putrule from 10.319 22.701 to 11.906 22.701
%
% Fig POLYLINE object
%
\linethickness= 0.500pt
\setplotsymbol ({\thinlinefont .})
\putrule from 10.636 22.860 to 10.636 21.273
%
% Fig POLYLINE object
%
\linethickness= 0.500pt
\setplotsymbol ({\thinlinefont .})
\setdashes < 0.1270cm>
\plot 12.700 22.860 12.700 21.431 /
%
% Fig POLYLINE object
%
\linethickness= 0.500pt
\setplotsymbol ({\thinlinefont .})
\plot  9.842 22.225 10.954 22.225 /
%
% Fig POLYLINE object
%
\linethickness= 0.500pt
\setplotsymbol ({\thinlinefont .})
\plot  9.842 21.590 10.954 21.590 /
%
% Fig POLYLINE object
%
\linethickness= 0.500pt
\setplotsymbol ({\thinlinefont .})
\setsolid
\putrule from 16.510 22.384 to 18.574 22.384
%
% Fig POLYLINE object
%
\linethickness= 0.500pt
\setplotsymbol ({\thinlinefont .})
\putrule from 17.621 23.654 to 17.621 21.907
%
% Fig POLYLINE object
%
\linethickness= 0.500pt
\setplotsymbol ({\thinlinefont .})
\putrule from 16.986 23.019 to 16.986 21.431
%
% Fig POLYLINE object
%
\linethickness= 0.500pt
\setplotsymbol ({\thinlinefont .})
\putrule from 15.558 21.749 to 17.304 21.749
%
% Fig POLYLINE object
%
\linethickness= 0.500pt
\setplotsymbol ({\thinlinefont .})
\plot 14.922 23.495 16.192 22.384 /
%
% Fig POLYLINE object
%
\linethickness= 0.500pt
\setplotsymbol ({\thinlinefont .})
\setdashes < 0.1270cm>
\plot 15.875 23.019 15.875 21.431 /
%
% Fig POLYLINE object
%
\linethickness= 0.500pt
\setplotsymbol ({\thinlinefont .})
\plot 18.256 23.654 18.256 22.066 /
%
% Fig POLYLINE object
%
\linethickness= 0.500pt
\setplotsymbol ({\thinlinefont .})
\plot 16.351 24.130 15.081 23.019 /
%
% Fig POLYLINE object
%
\linethickness= 0.500pt
\setplotsymbol ({\thinlinefont .})
\plot 15.875 24.130 17.939 23.336 /
%
% Fig TEXT object
%
\put{\SetFigFont{9}{10.8}{rm}$D_5$} [lB] at  5.080 20.320
%
% Fig TEXT object
%
\put{\SetFigFont{9}{10.8}{rm}$A_5$} [lB] at 10.954 20.320
%
% Fig TEXT object
%
\put{\SetFigFont{9}{10.8}{rm}$A_4+A_1$} [lB] at 16.351 20.479
%
% Fig TEXT object
%
\put{\SetFigFont{9}{10.8}{rm}$p$} [lB] at  4.763 21.431
\linethickness=0pt
\putrectangle corners at  3.467 24.155 and 18.599 20.320
\endpicture}
}
\vskip 0.5cm 
\centerline{\font\thinlinefont=cmr5
\begingroup\makeatletter\ifx\SetFigFont\undefined%
\gdef\SetFigFont#1#2#3#4#5{%
  \reset@font\fontsize{#1}{#2pt}%
  \fontfamily{#3}\fontseries{#4}\fontshape{#5}%
  \selectfont}%
\fi\endgroup%
\mbox{\beginpicture
\setcoordinatesystem units <0.78000cm,0.78000cm>
\unitlength=0.78000cm
\linethickness=1pt
\setplotsymbol ({\makebox(0,0)[l]{\tencirc\symbol{'160}}})
\setshadesymbol ({\thinlinefont .})
\setlinear
%
% Fig POLYLINE object
%
\linethickness= 0.500pt
\setplotsymbol ({\thinlinefont .})
\plot  3.969 23.654  2.699 22.225 /
%
% Fig POLYLINE object
%
\linethickness= 0.500pt
\setplotsymbol ({\thinlinefont .})
\putrule from  3.016 23.019 to  3.016 20.796
%
% Fig POLYLINE object
%
\linethickness= 0.500pt
\setplotsymbol ({\thinlinefont .})
\plot  2.699 21.431  4.286 20.003 /
%
% Fig POLYLINE object
%
\linethickness= 0.500pt
\setplotsymbol ({\thinlinefont .})
\plot  5.556 23.812  7.144 22.701 /
%
% Fig POLYLINE object
%
\linethickness= 0.500pt
\setplotsymbol ({\thinlinefont .})
\plot  7.144 21.590  5.556 20.003 /
%
% Fig POLYLINE object
%
\linethickness= 0.500pt
\setplotsymbol ({\thinlinefont .})
\setdashes < 0.1270cm>
\plot  6.826 23.336  6.826 20.796 /
%
% Fig POLYLINE object
%
\linethickness= 0.500pt
\setplotsymbol ({\thinlinefont .})
\plot  3.493 23.495  6.350 23.654 /
%
% Fig POLYLINE object
%
\linethickness= 0.500pt
\setplotsymbol ({\thinlinefont .})
\plot  3.651 20.161  6.350 20.161 /
%
% Fig POLYLINE object
%
\linethickness= 0.500pt
\setplotsymbol ({\thinlinefont .})
\plot  2.699 21.749  5.080 23.019 /
%
% Fig POLYLINE object
%
\linethickness= 0.500pt
\setplotsymbol ({\thinlinefont .})
\plot  4.286 22.860  7.303 22.066 /
%
% Fig POLYLINE object
%
\linethickness= 0.500pt
\setplotsymbol ({\thinlinefont .})
\setsolid
\plot 11.271 23.812 10.001 22.225 /
%
% Fig POLYLINE object
%
\linethickness= 0.500pt
\setplotsymbol ({\thinlinefont .})
\plot 10.001 22.701 11.271 20.796 /
%
% Fig POLYLINE object
%
\linethickness= 0.500pt
\setplotsymbol ({\thinlinefont .})
\plot 14.129 23.019 12.859 20.796 /
%
% Fig POLYLINE object
%
\linethickness= 0.500pt
\setplotsymbol ({\thinlinefont .})
\plot 12.700 23.812 14.129 22.225 /
%
% Fig POLYLINE object
%
\linethickness= 0.500pt
\setplotsymbol ({\thinlinefont .})
\putrule from 11.906 24.765 to 11.906 23.019
%
% Fig POLYLINE object
%
\linethickness= 0.500pt
\setplotsymbol ({\thinlinefont .})
\setdashes < 0.1270cm>
\plot 10.636 21.114 13.494 21.114 /
%
% Fig POLYLINE object
%
\linethickness= 0.500pt
\setplotsymbol ({\thinlinefont .})
\plot 10.319 23.178 12.065 21.749 /
%
% Fig POLYLINE object
%
\linethickness= 0.500pt
\setplotsymbol ({\thinlinefont .})
\plot 13.652 23.495 11.589 21.749 /
%
% Fig POLYLINE object
%
\linethickness= 0.500pt
\setplotsymbol ({\thinlinefont .})
\plot 10.636 23.336 12.224 23.178 /
%
% Fig POLYLINE object
%
\linethickness= 0.500pt
\setplotsymbol ({\thinlinefont .})
\plot 11.748 23.495 13.494 23.654 /
%
% Fig POLYLINE object
%
\linethickness= 0.500pt
\setplotsymbol ({\thinlinefont .})
\setsolid
\putrule from 17.304 21.273 to 20.637 21.273
%
% Fig POLYLINE object
%
\linethickness= 0.500pt
\setplotsymbol ({\thinlinefont .})
\putrule from 17.780 23.019 to 17.780 20.796
%
% Fig POLYLINE object
%
\linethickness= 0.500pt
\setplotsymbol ({\thinlinefont .})
\putrule from 18.733 23.019 to 18.733 20.796
%
% Fig POLYLINE object
%
\linethickness= 0.500pt
\setplotsymbol ({\thinlinefont .})
\putrule from 19.844 23.178 to 19.844 20.796
%
% Fig POLYLINE object
%
\linethickness= 0.500pt
\setplotsymbol ({\thinlinefont .})
\setdashes < 0.1270cm>
\plot 16.986 22.543 18.098 22.543 /
%
% Fig POLYLINE object
%
\linethickness= 0.500pt
\setplotsymbol ({\thinlinefont .})
\plot 19.526 22.543 20.796 22.543 /
%
% Fig TEXT object
%
\put{\SetFigFont{9}{10.8}{rm}$A_3+2A_1$} [lB] at  4.286 18.891
%
% Fig TEXT object
%
\put{\SetFigFont{9}{10.8}{rm}$2A_2+A_1$} [lB] at 11.271 19.209
%
% Fig TEXT object
%
\put{\SetFigFont{9}{10.8}{rm}$D_4$} [lB] at 18.256 19.367
\linethickness=0pt
\putrectangle corners at  2.673 24.790 and 20.822 18.891
\endpicture}
}
\vskip 0.5cm
\centerline{\font\thinlinefont=cmr5
\begingroup\makeatletter\ifx\SetFigFont\undefined%
\gdef\SetFigFont#1#2#3#4#5{%
  \reset@font\fontsize{#1}{#2pt}%
  \fontfamily{#3}\fontseries{#4}\fontshape{#5}%
  \selectfont}%
\fi\endgroup%
\mbox{\beginpicture
\setcoordinatesystem units <0.66000cm,0.66000cm>
\unitlength=0.66000cm
\linethickness=1pt
\setplotsymbol ({\makebox(0,0)[l]{\tencirc\symbol{'160}}})
\setshadesymbol ({\thinlinefont .})
\setlinear
%
% Fig POLYLINE object
%
\linethickness= 0.500pt
\setplotsymbol ({\thinlinefont .})
\plot 22.066 23.654 21.273 22.384 /
%
% Fig POLYLINE object
%
\linethickness= 0.500pt
\setplotsymbol ({\thinlinefont .})
\plot 21.273 21.749 22.384 20.637 /
%
% Fig POLYLINE object
%
\linethickness= 0.500pt
\setplotsymbol ({\thinlinefont .})
\plot 23.178 23.812 24.130 22.701 /
%
% Fig POLYLINE object
%
\linethickness= 0.500pt
\setplotsymbol ({\thinlinefont .})
\plot 24.130 21.907 23.336 20.637 /
%
% Fig POLYLINE object
%
\linethickness= 0.500pt
\setplotsymbol ({\thinlinefont .})
\setdashes < 0.1270cm>
\plot 21.590 23.495 23.812 23.495 /
%
% Fig POLYLINE object
%
\linethickness= 0.500pt
\setplotsymbol ({\thinlinefont .})
\plot 21.431 23.019 21.431 21.114 /
%
% Fig POLYLINE object
%
\linethickness= 0.500pt
\setplotsymbol ({\thinlinefont .})
\plot 23.812 23.336 23.971 21.273 /
%
% Fig POLYLINE object
%
\linethickness= 0.500pt
\setplotsymbol ({\thinlinefont .})
\plot 21.907 20.796 23.812 20.796 /
%
% Fig POLYLINE object
%
\linethickness= 0.500pt
\setplotsymbol ({\thinlinefont .})
\setsolid
\putrule from  4.604 21.273 to  6.826 21.273
%
% Fig POLYLINE object
%
\linethickness= 0.500pt
\setplotsymbol ({\thinlinefont .})
\putrule from  5.874 22.543 to  5.874 20.796
%
% Fig POLYLINE object
%
\linethickness= 0.500pt
\setplotsymbol ({\thinlinefont .})
\putrule from  4.921 22.543 to  4.921 20.796
%
% Fig POLYLINE object
%
\linethickness= 0.500pt
\setplotsymbol ({\thinlinefont .})
\putrule from  3.651 22.225 to  5.239 22.225
%
% Fig POLYLINE object
%
\linethickness= 0.500pt
\setplotsymbol ({\thinlinefont .})
\setdashes < 0.1270cm>
\plot  3.810 22.543  3.810 20.637 /
%
% Fig POLYLINE object
%
\linethickness= 0.500pt
\setplotsymbol ({\thinlinefont .})
\plot  2.381 21.273  4.128 21.273 /
%
% Fig POLYLINE object
%
\linethickness= 0.500pt
\setplotsymbol ({\thinlinefont .})
\plot  6.509 22.543  6.509 20.796 /
%
% Fig POLYLINE object
%
\linethickness= 0.500pt
\setplotsymbol ({\thinlinefont .})
\setsolid
\putrule from 11.906 21.431 to 14.129 21.431
%
% Fig POLYLINE object
%
\linethickness= 0.500pt
\setplotsymbol ({\thinlinefont .})
\putrule from  9.684 21.431 to 11.271 21.431
%
% Fig POLYLINE object
%
\linethickness= 0.500pt
\setplotsymbol ({\thinlinefont .})
\putrule from 12.224 22.701 to 12.224 21.114
%
% Fig POLYLINE object
%
\linethickness= 0.500pt
\setplotsymbol ({\thinlinefont .})
\putrule from 10.636 22.384 to 12.541 22.384
%
% Fig POLYLINE object
%
\linethickness= 0.500pt
\setplotsymbol ({\thinlinefont .})
\setdashes < 0.1270cm>
\plot 13.018 22.701 13.018 21.114 /
%
% Fig POLYLINE object
%
\linethickness= 0.500pt
\setplotsymbol ({\thinlinefont .})
\plot 13.652 22.701 13.652 21.114 /
%
% Fig POLYLINE object
%
\linethickness= 0.500pt
\setplotsymbol ({\thinlinefont .})
\plot 10.954 22.701 10.954 21.114 /
%
% Fig POLYLINE object
%
\linethickness= 0.500pt
\setplotsymbol ({\thinlinefont .})
\setsolid
\plot 17.939 23.178 18.891 22.384 /
%
% Fig POLYLINE object
%
\linethickness= 0.500pt
\setplotsymbol ({\thinlinefont .})
\plot 19.050 22.066 18.415 20.955 /
%
% Fig POLYLINE object
%
\linethickness= 0.500pt
\setplotsymbol ({\thinlinefont .})
\setdashes < 0.1270cm>
\plot 16.828 23.019 18.415 23.019 /
%
% Fig POLYLINE object
%
\linethickness= 0.500pt
\setplotsymbol ({\thinlinefont .})
\plot 16.986 21.114 18.733 21.114 /
%
% Fig POLYLINE object
%
\linethickness= 0.500pt
\setplotsymbol ({\thinlinefont .})
\plot 18.733 22.860 17.780 21.907 /
%
% Fig POLYLINE object
%
\linethickness= 0.500pt
\setplotsymbol ({\thinlinefont .})
\plot 17.621 22.225 19.050 21.431 /
%
% Fig POLYLINE object
%
\linethickness= 0.500pt
\setplotsymbol ({\thinlinefont .})
\setsolid
\plot 16.294 22.018 17.405 20.748 /
%
% Fig POLYLINE object
%
\linethickness= 0.500pt
\setplotsymbol ({\thinlinefont .})
\plot 17.145 23.178 16.351 21.749 /
%
% Fig TEXT object
%
\put{\SetFigFont{8}{9.6}{rm}$A_4$} [lB] at  4.128 19.685
%
% Fig TEXT object
%
\put{\SetFigFont{8}{9.6}{rm}$A_3+A_1$} [lB] at 11.113 20.003
%
% Fig TEXT object
%
\put{\SetFigFont{8}{9.6}{rm}$A_2+2A_1$} [lB] at 17.304 20.003
%
% Fig TEXT object
%
\put{\SetFigFont{8}{9.6}{rm}$4A_1$} [lB] at 22.066 19.844
\linethickness=0pt
\putrectangle corners at  2.356 23.838 and 24.155 19.685
\endpicture}
}
\vskip 0.5cm 
\centerline{\font\thinlinefont=cmr5
\begingroup\makeatletter\ifx\SetFigFont\undefined%
\gdef\SetFigFont#1#2#3#4#5{%
  \reset@font\fontsize{#1}{#2pt}%
  \fontfamily{#3}\fontseries{#4}\fontshape{#5}%
  \selectfont}%
\fi\endgroup%
\mbox{\beginpicture
\setcoordinatesystem units <0.83000cm,0.83000cm>
\unitlength=0.83000cm
\linethickness=1pt
\setplotsymbol ({\makebox(0,0)[l]{\tencirc\symbol{'160}}})
\setshadesymbol ({\thinlinefont .})
\setlinear
%
% Fig POLYLINE object
%
\linethickness= 0.500pt
\setplotsymbol ({\thinlinefont .})
\putrule from  2.064 20.637 to  4.445 20.637
%
% Fig POLYLINE object
%
\linethickness= 0.500pt
\setplotsymbol ({\thinlinefont .})
\putrule from  3.493 21.907 to  3.493 20.320
%
% Fig POLYLINE object
%
\linethickness= 0.500pt
\setplotsymbol ({\thinlinefont .})
\putrule from  2.381 21.907 to  2.381 20.320
%
% Fig POLYLINE object
%
\linethickness= 0.500pt
\setplotsymbol ({\thinlinefont .})
\setdashes < 0.1270cm>
\plot  4.128 21.907  4.128 20.320 /
%
% Fig POLYLINE object
%
\linethickness= 0.500pt
\setplotsymbol ({\thinlinefont .})
\plot  1.746 21.590  3.016 21.590 /
%
% Fig POLYLINE object
%
\linethickness= 0.500pt
\setplotsymbol ({\thinlinefont .})
\setsolid
\putrule from  7.779 20.796 to  9.525 20.796
%
% Fig POLYLINE object
%
\linethickness= 0.500pt
\setplotsymbol ({\thinlinefont .})
\putrule from  9.049 22.066 to  9.049 20.479
%
% Fig POLYLINE object
%
\linethickness= 0.500pt
\setplotsymbol ({\thinlinefont .})
\putrule from 10.160 22.225 to 10.160 20.320
%
% Fig POLYLINE object
%
\linethickness= 0.500pt
\setplotsymbol ({\thinlinefont .})
\setdashes < 0.1270cm>
\plot  8.731 21.749 10.478 21.749 /
%
% Fig POLYLINE object
%
\linethickness= 0.500pt
\setplotsymbol ({\thinlinefont .})
\plot  7.938 22.066  7.938 20.320 /
%
% Fig POLYLINE object
%
\linethickness= 0.500pt
\setplotsymbol ({\thinlinefont .})
\plot  6.826 21.590  8.255 21.590 /
%
% Fig POLYLINE object
%
\linethickness= 0.500pt
\setplotsymbol ({\thinlinefont .})
\setsolid
\putrule from 12.859 20.796 to 15.240 20.796
%
% Fig POLYLINE object
%
\linethickness= 0.500pt
\setplotsymbol ({\thinlinefont .})
\putrule from 13.335 22.066 to 13.335 20.320
%
% Fig POLYLINE object
%
\linethickness= 0.500pt
\setplotsymbol ({\thinlinefont .})
\setdashes < 0.1270cm>
\plot 13.970 22.066 13.970 20.320 /
%
% Fig POLYLINE object
%
\linethickness= 0.500pt
\setplotsymbol ({\thinlinefont .})
\plot 14.605 22.066 14.605 20.479 /
%
% Fig POLYLINE object
%
\linethickness= 0.500pt
\setplotsymbol ({\thinlinefont .})
\plot  3.175 17.145  5.080 17.145 /
%
% Fig POLYLINE object
%
\linethickness= 0.500pt
\setplotsymbol ({\thinlinefont .})
\plot  2.064 15.875  3.651 15.875 /
%
% Fig POLYLINE object
%
\linethickness= 0.500pt
\setplotsymbol ({\thinlinefont .})
\setsolid
\putrule from  3.493 17.304 to  3.493 15.558
%
% Fig POLYLINE object
%
\linethickness= 0.500pt
\setplotsymbol ({\thinlinefont .})
\putrule from  4.763 17.304 to  4.763 15.558
%
% Fig POLYLINE object
%
\linethickness= 0.500pt
\setplotsymbol ({\thinlinefont .})
\setdashes < 0.1270cm>
\plot  8.572 17.462  8.572 15.558 /
%
% Fig POLYLINE object
%
\linethickness= 0.500pt
\setplotsymbol ({\thinlinefont .})
\setsolid
\putrule from 10.001 17.462 to 10.001 15.558
%
% Fig POLYLINE object
%
\linethickness= 0.500pt
\setplotsymbol ({\thinlinefont .})
\setdashes < 0.1270cm>
\plot  8.255 16.034 10.319 16.034 /
%
% Fig TEXT object
%
\put{\SetFigFont{10}{12.0}{rm}$A_3$} [lB] at  2.699 19.050
%
% Fig TEXT object
%
\put{\SetFigFont{10}{12.0}{rm}$A_2+A_1$} [lB] at  8.096 19.050
%
% Fig TEXT object
%
\put{\SetFigFont{10}{12.0}{rm}$A_2$} [lB] at 13.494 19.050
%
% Fig TEXT object
%
\put{\SetFigFont{10}{12.0}{rm}$2A_1$} [lB] at  2.699 14.287
%
% Fig TEXT object
%
\put{\SetFigFont{10}{12.0}{rm}$A_1$} [lB] at  8.731 14.446
\linethickness=0pt
\putrectangle corners at  1.721 22.250 and 15.265 14.287
\endpicture}
}
\begin{center}
Figure $6$
\end{center}

\vskip 0.5cm
\hskip -0.6cm {\bf  Table $1$. \ Mordell-Weil Rank one Rational Elliptic Surfaces}
\vskip 0.3cm
In the following table, $Y \longrightarrow U_1 (D)$ is the blow-down of the zero section ${\cal O}$ and $U_1(D) \longrightarrow V$ is the minimal resolution of the singularities of $V$ with $D$ the exceptional divisor; this map is nothing but the contraction of all $(-2)$-curves. The symble $\langle m \rangle$ stands for a rank one lattice ${\bf Z}P$ with the generator $P$ such that $\langle P,P \rangle =m$. We also denote $P_n:=nP$ where $n \in {\bf Z}$. The $A_1^*$ indicates the lattice generated by the element $P$ such that $<P,P>=1/2$. In the last column, $E$ is a $(-1)$-curve and $R$ is a connected component of $D$ such that $E+R$ is a linear chain on $U_1(D)$. We can't find such a linear chain for the first $6$ cases.
\begin{center}
\begin{tabular}{| c|c|c|c|c|}\hline
 $\sharp$  &  Fibre type on $Y$ &  MW(Y) &  $U_1 (D)$   &  $E+R$ on $U_1(D)$   \\  \hline  
$1$ &  $I_3^*, II, I_1$    &  $<1/4>$   &  $U^1_1(D_7)$  &   \\ \hline
$2$ &$I_3^*, 3I_1$      &  $<1/4>$   &   $U_1^2( D_7)$    &   \\ \hline
$3$ & $I_4, 4I_2$     &  $<1/4>\oplus ({\bf Z}/2{\bf Z})^2$   &  $U_1^3( A_3+4A_1)$    &     \\ \hline
$4$ & $I_1^*, III, I_2$      &  $<1/4>\oplus ({\bf Z}/2{\bf Z}) $   &  $U_1^4( D_5+2A_1)$   &     \\ \hline
$5$ &  $I_1^*, 2I_2, I_1$     &  $ <1/4>\oplus ({\bf Z}/2{\bf Z})  $   & $U_1^5( D_5+2A_1)$     &   \\ \hline
$6$ & $I_0^*, I_4, 2I_1$     &  $  <1/4>\oplus ({\bf Z}/2{\bf Z})  $   & $U_1^6( D_4+A_3)$     &   \\ \hline 
$7$ & $III^*, II, I_1$     &  $A_1^*$   &  $U_1^7(E_7)$   &  $P_2$    \\ \hline
 $8$ &  $III^*, 3I_1$    &  $A_1^*$   &  $U_1^8(E_7)$     &   $P_2$ \\ \hline
$9$ & $I_8, 4I_1$     &  $A_1^*\oplus ({\bf Z}/2{\bf Z}) $   &   $U_1^9(A_7)$   &   $P_2$  \\ \hline
$10$ & $I_8, II, 2I_1$      &  $<1/8>$   &  $U_1^{10}(A_7)$     &   $P_3+A_7$  \\ \hline
$11$ &  $I_8, 4I_1$     &  $<1/8>$   &  $U_1^{11}(A_7)$    &   $P_3+A_7$  \\ \hline
$12$ & $2I_4, I_2, 2I_1$     &  $ A_1^* \oplus ({\bf Z}/4{\bf Z})  $   &  $U_1^{12}(2A_3+A_1)$    &   $P_2$ \\ \hline
$13$ &$I_2^*, III, I_1$      &  $ A_1^* \oplus {\bf Z}/2{\bf Z}$   &  $U_1^{13}(D_6+A_1)$    &  $P_2$ \\ \hline
$14$ & $I^*_2, I_2, 2I_1$      &  $A_1^* \oplus {\bf Z}/2{\bf Z} $   & $U_1^{14}(D_6+A_1)$       & $P_2$  \\ \hline
$15$ &  $I_6, IV, 2I_1$    &  $ A_1^* \oplus {\bf Z}/3{\bf Z}$   &  $U_1^{15}(A_5+A_2)$     &  $P_2$ \\ \hline
$16$ &  $I_6, I_3, 3I_1$     &  $A_1^* \oplus {\bf Z}/3{\bf Z}$   &  $U_1^{16}(A_5+A_2)$     & $P_2$  \\ \hline
$17$ &   $I_0^*, 3I_2$     &  $ A_1^* \oplus ({\bf Z}/2{\bf Z})^2$   & $U_1^{17}(D_4+3A_1)$ & $P_2$ \\ \hline
$18$ &  $IV^*, III, I_1$    &  $<1/6>$   &  $U_1^{18}(E_6+A_1)$    & $P_3+A_1$ \\ \hline
$19$ & $IV^*, I_2, II$     &  $<1/6>$   &  $U_1^{19}(E_6+A_1)$    &  $P_3+A_1$ \\ \hline
$20$ & $IV^*, I_2, 2I_1$     &  $<1/6>$   &  $U_1^{20}(E_6+A_1)$    & $P_3+A_1$  \\ \hline
$21$ & $I_4, I_3, III, I_2$     &  $<1/12> \oplus {\bf Z}/2{\bf Z}$   &  $U_1^{21}(A_3+A_2+2A_1)$ &  $P_4+A_2$  \\ \hline
$22$ & $I_4, I_3, 2I_2, I_1$     &  $ <1/12> \oplus {\bf Z}/2{\bf Z} $   &  $U_1^{22}(A_3+A_2+2A_1)$    & $P_4+A_2$ \\ \hline
$23$ & $IV, 2I_3, I_2$     &  $ <1/6> \oplus {\bf Z}/3{\bf Z}$   & 
$U_1^{23}(3A_2+A_1)$      & $P_3+A_1$  \\ \hline
$24$ &  $3I_3, I_2, I_1$    &  $ <1/6> \oplus {\bf Z}/3{\bf Z}$   &  $U_1^{24}(3A_2+A_1)$    & $P_3+A_1$ \\ \hline
$25$ &$I_7, III, 2I_1$      &  $ <1/14>$   &  $U_1^{25}(A_6+A_1)$  &  $P_4+A_6$ \\ \hline
$26$ &  $I_7, I_2, II, I_1$    &  $<1/14>$   &  $U_1^{26}(A_6+A_1)$  & $P_4+A_6$  \\ \hline
$27$ &  $I_7, I_2, 3I_1$    &  $<1/14>$   & $U_1^{27}(A_6+A_1)$       &  $P_4+A_6$ \\ \hline
$28$ & $I_1^*, IV, I_1$     &  $<1/12>$   & $U_1^{28}(D_5+A_2)$    &   $P_4+A_2$ \\ \hline
 $29$ & $I_1^*, I_3, II$      &  $<1/12>$   &  $U_1^{29}(D_5+A_2)$      &  $P_4+A_2$    \\ \hline
$30$ & $I_1^*, I_3, 2I_1$      &  $<1/12>$   & $U_1^{30}(D_5+A_2)$   & $P_4+A_2$    \\ \hline
$31$ &  $I_6, III, I_2, I_1$     &  $<1/6> \oplus {\bf Z}/2{\bf Z}$   &  $U_1^{31}( A_5+2A_1)$  & $P_3+A_1$ \\ \hline
$32$ &  $I_6, 2I_2, 2I_1$     &  $ <1/6> \oplus {\bf Z}/2{\bf Z}$   &   $ U_1^{32}(A_5+2A_1)$    & $P_3+A_1$     \\ \hline
$33$ & $I_5, I_4, II, I_1$     &  $ <1/20>$   &  $U_1^{33}(A_4+A_3)$       & $ P_5+A_3$  \\ \hline
$34$ &  $I_5, I_4, 3I_1$     &  $<1/20>$   &  $U_1^{34}(A_4+A_3)$    &  $ P_6+A_4$ \\ \hline
$35$ &   $I_5, IV, I_2, I_1$    &  $<1/30>$   &   $U_1^{35}(A_4+A_2+A_1)$   & $ P_6+A_4$   \\ \hline
$36$ &  $I_5, I_3, III, I_1$     &  $<1/30>$   &  $U_1^{36}(A_4+A_2+A_1)$    &$ P_6+A_4$   \\ \hline
$37$ &  $I_5, I_3, I_2, II$     &  $<1/30>$   &  $U_1^{37}(A_4+A_2+A_1)$      & $ P_6+A_4$  \\ \hline
$38$ &   $I_5, I_3, I_2, 2I_1$    &  $<1/30>$   &  $U_1^{38}(A_4+A_2+A_1)$      & $ P_6+A_4$   \\ \hline 
 
\end{tabular} 
\end{center}
 
\vskip 0.6cm
\hskip -0.6cm {\bf Table $2$. \  Nice Exceptinal Curves on $U_1$}
\vskip 0.5cm
\hskip -0.6cm
In the following table, we employ the following notation and convention:
\vskip 0.3cm
\hskip -0.6cm
We let $G_i$, $H_i$, $J_i$, $L_i$ and $M_i$ be the $i$-th
component in the corresponding fibres $F_1$, $F_2$, $F_3$ $F_4$ and $F_5$
respectively. The numbering of the singular fibre is defined in 
the following diagrams.
For the three types of singular fibres, the $0$-th component in a fibre intersects the zero section. For the other types, the numbering is natural. The components of the singular fibre are numbered as in [Ye]. 
%\vskip 0.3cm
%\centerline {\input mw1.tex} 
%\vskip 0.3cm
%\begin{center}
%Figure $7$
%\end{center}
 For a lattice $L=<m>\oplus {\bf Z}/k{\bf Z}$, we let $P_1$ to be the generator of $<m>$ and $P_n=nP_1$ where $n\in {\bf Z}$. If $Q\in {\bf Z}/k{\bf Z}$,
i.e., the torsion part of $L$, then we denote $P_{n,Q}=P_n+Q$.  $(R_1 R_2)$ is the intersection number of the sections $R_1$ and $R_2$.
 \vskip 0.5cm
\begin{center}
\begin{tabular}{| c|c|c|c|c|}\hline
 $\sharp$ &  Fibre type  & Sections disjoint from & $NECs$ & sections which intersects  \\
  & on $Y$    &  ${\cal O}$ on $Y$  &    on $U_1$    &   NEC  on $U_1$  \\ \hline
$1$ &  $I_3^*, II, I_1$  &  $ P_1=(2)$, $P_{-1}=(3) $  & $P_{\pm1}$, $P_{\pm 2}$ & $(P_{\pm 1} P_{\mp 2})=1$   \\ 

$2$ &$I_3^*, 3I_1$     &   $ P_{\pm 2}=(1)$  & $G_0$ & $(P_2 P_{-2})=1$    \\ \hline

$3$ & $I_4, 4I_2$       &  $Q_1=(2,1,1,0,0)$,& $P_{\pm 2}$&  $(P_2 P_{-2})=1 $ \\
    & &  $Q_2=(2,0,0,1,1)$, &  & $(P_{\pm 2} P_{{\mp 2},Q_1})=2 $     \\
    &      & $Q_3=(0,1,1,1,1)$,& & $ (P_{\pm 2} P_{\mp 1})=1 $ \\
&      &  $P_1=(1,0,1,1,0)$, &  &$ (P_{\pm 2} P_{{\mp 1},Q_2})=1$ \\

   &   & $P_{1,Q_1}=(3,1,0,1,0)$,& & $ (P_{\pm 2} P_{\mp 1 ,Q_1})=1 $  \\
&   &  $P_{1,Q_2}=(3,0,1,0,1)$, & & $(P_{\pm 2} P_{\mp 1 ,Q_2})=1 $  \\

 & &          $P_{1,Q_3}=(1,1,0,0,1)$,& & $(P_{\pm 2}  Q_3)=1$  \\
& &          $P_{-1}=(1,0,1,0,1)$, &   &\\

& &     $P_{{-1},Q_1}=(3,1,0,0,1)$,  &  &\\ 
& &      $P_{{-1},Q_2}=(3,0,1,1,0)$, &  &\\

& &          $P_{-1,Q_3}=(1,1,0,1,0)$,   &   & \\
& &             $P_{\pm 2}=(2,0,0,0,0)$, &  & \\

& &          $P_{\pm 2 ,Q_1}=(0,1,1,0,0)$,&  &  \\ 
& &        $P_{\pm 2 ,Q_2}=(0,0,0,1,1)$. & & \\ \hline

\end{tabular} 
\end{center} 

\begin{center}
\begin{tabular}{| c|c|c|c|c|}\hline

$4$ & $I_1^*, III, I_2$  & $Q=(1,1,1)$, $P_1=(2,1,0)$,& $P_{\pm 2}$ & $(P_{\pm 2} P_{\mp 1})=1$      \\
                           
 $5$ &  $I_1^*, 2I_2, I_1$    &     $P_{1,Q}=(3,0,1)$, $P_{-1}=(2,0,1)$,& $G_0$ &   $(P_{\pm 2} P_{\mp 1 ,Q})=1$    \\
 & &   $P_{-1 ,Q}=(3,1,0)$,  $P_{\pm 2}=(1,0,0)$,& &
 $(P_{\pm 2} P_{\mp 2,Q})=2$ \\
& & $P_{\pm 2,Q}=(0,1,1)$. &     & $(P_2  P_{-2})=1$ \\ 
& & & & $(G_0 P_{\pm 2,Q})=1$ \\
\hline
$6$ & $I_0^*, I_4, 2I_1$ &    $Q=(1,2)$, $P_1=(2,1)$  & $P_{\pm 2}$, &$(P_{\pm 2} P_{\mp 1})=1$   \\
     &   &         $P_{1,Q}=(3,3)$, $P_{-1}=(3,1)$ & $P_{\pm 2,Q}$  &$(P_{\pm 2} P_{\mp 1,Q})=1$    \\
    &  &    $P_{-1,Q}=(2,3)$, $P_{\pm 2}=(1,0)$ & $G_0$ & $(P_2P_{-2})=1$    \\                                   &  &  $P_{\pm 2,Q}=(0,2)$ & & $(P_{\pm 2,Q} P_{\mp 1})=1$   \\
 & &                        & & $(P_{\pm 2} P_{\mp 2,Q})=2$   \\
 & &                        & & $(P_{2,Q} P_{-2,Q})=1$    \\ 
& & & & $(G_0  P_{\pm 2,Q})=1$ \\
\hline

$7$ & $III^*, II, I_1$      &  $P_{\pm 1}=(1)$, $P_{\mp 2}=(0)$ & $P_{\pm 1}$ &    \\  
$8$ &  $III^*, 3I_1$    &              & $G_0$    & \\ \hline
$9$ & $I_8, 4I_1$      &   $Q=(4)$, $P_{\pm 1}=(2)$ & $Q$, $P_{\pm 1}$  & $(P_{\pm 1} P_{\mp 1,Q})=1$ \\    
    & &  $P_{\pm 1,Q}=(6)$, $P_{\pm 2}=(0)$  & $P_{\pm 1,Q}$  & $(P_{\pm 2}Q)=1$   \\ \hline

$10$ & $I_8, II, 2I_1$  &   $P_1=(3)$, $P_2=(6)$ & $P_{\pm 1}$    &$(P_{\pm 1} P_{\mp 3})=1$ \\
  $11$ &  $I_8, 4I_1$ & $P_3=(1)$, $P_{-1}=(5)$&$P_{\pm 2}$  &$(P_2 P_{-2})=1$   \\
       & &  $P_{-2}=(2)$, $P_{-3}=(7)$ & $P_{\pm 3}$    &$(P_{\pm 3} P_{\mp 2})=1$    \\ \hline

$12$ & $2I_4, I_2, 2I_1$  & $Q_1=(1,1,1)$, $Q_2=(2,2,0)$  & $\sharp \{ NEC \}$ & \\
    &    &  $Q_3=(3,3,1)$, $P_{\pm 1}=(2,0,1)$ & $=\phi$  &        \\
      &   &  $P_{\pm 1,Q_1}=(3,1,0)$, $P_{\pm 1,Q_2}=(0,2,1)$ &   &    \\
      &  &  $P_{\pm 1,Q_3}=(1,3,0)$, $P_{\pm 2}=(0,0,0)$  &  &   \\ \hline

$13$ &$I_2^*, III, I_1$ &  $Q=(2,1)$, $P_{\pm 1}=(3,0)$&  $P_{\pm 1}$  &    $(P_{\pm 1} P_{\mp 1,Q})=1$ \\
 $14$ & $I^*_2, I_2, 2I_1$       &  $P_{\pm 1, Q}=(1,1)$, $P_{\pm 2}=(0,0)$    & $G_0$  & \\ \hline

$15$ &  $I_6, IV, 2I_1$ & $Q_1=(2,1)$, $Q_2=(4,2)$ & $P_{\pm 1}$   & $(P_{\pm 1} P_{\mp 1,Q_1})=1$ \\
  $16$ &  $I_6, I_3, 3I_1$   &        $P_{\pm 1}=(3,0)$, $P_{\pm 1,Q_1}=(5,1)$ &  & $(P_{\pm 1} P_{\mp 1,Q_2})=1$ \\
     & &             $P_{\pm 1,Q_2}=(1,2)$, $P_{\pm 2}=(0,0)$   &    & \\ \hline

$17$ &   $I_0^*, 3I_2$ & $Q_1=(1,1,1,0)$, $Q_2=(2,0,1,1)$ &   $G_0$ & $(G_0 P_{\pm1})=1$\\
    &       & $Q_3=(3,1,0,1)$, $P_{\pm 1}=(0,1,1,1)$ & &$(G_0 P_{\pm 2})=1$   \\
     &     &$P_{\pm 1, Q_1}=(1,0,0,1)$, $P_{\pm 1,Q_2}=(2,1,0,0)$ & &  \\
   &      & $P_{\pm 1,Q_3}=(3,0,1,0)$, $P_{\pm 2}=(0,0,0,0)$ & &   \\ \hline

$18$ &  $IV^*, III, I_1$ & $P_1=(1,1)$, $P_2=(2,0)$ & $P_{\pm 2}$ & $(P_{\pm 3} P_{\mp 1})=1$ \\         
  $19$ & $IV^*, I_2, II$  & $P_{\pm 3}=(0,1)$, $P_{-1}=(2,1)$ & $P_{\pm 3}$    &$(P_{2} P_{-2})=1$ \\    
  $20$ & $IV^*, I_2, 2I_1$ &   $P_{-2}=(1,0)$      & $G_0$      & $(P_{\pm 3} P_{\mp 2})=2$ \\    
     & & &   & $(P_{3} P_{-3})=2$     \\ 
& & & & $(G_0 P_{\pm 3})=1$ \\
\hline

$21$ & $I_4, I_3, III, I_2$ & $Q=(2,0,1,1)$, $P_1=(1,1,1,0)$& $P_{\pm 4}$  & $(P_{\pm 4} P_{\mp 1})=1$ \\    
  $22$ & $I_4, I_3, 2I_2, I_1$        & $P_{1,Q}=(3,1,0,1)$, $P_2=(2,2,0,0)$  &    &$(P_{\pm 4} P_{\mp 1,Q})=1$ \\    
 &     & $P_{2,Q}=(0,2,1,1)$, $P_3=(3,0,1,0)$  &  & $(P_{\pm 4} P_{\mp 2})=1$ \\    
 & & $P_{3,Q}=(1,0,0,1)$, $P_4=(0,1,0,0)$  &  & $(P_{\pm 4} P_{\mp 2,Q})=1$ \\    
 & & $P_{-1}=(1,2,0,1)$, $P_{-1,Q}=(3,2,1,0)$  &  & $(P_{\pm 4} P_{\mp 3})=2$ \\    
 &  & $P_{-2}=(2,1,0,0)$, $P_{-2,Q}=(0,1,1,1)$  &  & $(P_{\pm 4} P_{\mp 3,Q})=2$ \\    
 & & $P_{-3}=(3,0,0,1)$, $P_{-3,Q}=(1,0,1,0)$  &  &   $(P_{4} P_{-4})=2$ \\    
  &   & $P_{-4}=(0,2,0,0)$  &  &  \\ \hline

\end{tabular} 
\end{center}

\begin{center}
\begin{tabular}{| c|c|c|c|c|}\hline

 $23$ & $IV, 2I_3, I_2$ &  $Q_1=(1,1,1,0)$, $Q_2=(2,2,2,0)$ &$P_{\pm 3}$  & $(P_{\pm 3} P_{\mp 1})=1$ \\      
 $24$ &  $3I_3, I_2, I_1$    &  $P_1=(1,2,0,1)$, $P_{1,Q_1}=(2,0,1,1)$ &     &  $(P_{\pm 3} P_{\mp 1,Q_1})=1$ \\ 
&  &  $P_{1,Q_2}=(0,1,2,1)$, $P_2=(1,0,2,0)$ &   & $(P_{\pm 3} P_{\mp 1,Q_2})=1$ \\ 
&  &  $P_{2,Q_1}=(2,1,0,0)$, $P_{2,Q_2}=(0,2,1,0)$ &    &$(P_{\pm 3} P_{\mp 2})=2$ \\ 
&  &  $P_{-1}=(1,0,2,1)$, $P_{-1,Q_1}=(2,1,0,1)$ &   & $(P_{\pm 3} P_{\mp 2,Q_1})=2$ \\ 
&  &  $P_{-1,Q_2}=(0,2,1,1)$, $P_{-2}=(1,2,0,0)$ &   &$(P_{\pm 3} P_{\mp 2,Q_2})=2$ \\ 
&  &  $P_{-2,Q_1}=(2,0,1,0)$, $P_{-2,Q_2}=(0,1,2,0)$ &     &$(P_{3} P_{-3})=2$ \\
&    &  $P_{\pm 3}=(0,0,0,1)$ &    & \\ \hline

$25$ &$I_7, III, 2I_1$ &  $P_1=(2,1)$, $P_2=(4,0)$&$P_{\pm 2}$ &$(P_{\pm 3} P_{\mp 2})=1$ \\  
  $26$ &  $I_7, I_2, II, I_1$  &  $P_3=(6,1)$, $P_4=(1,0)$& $P_{\pm 4}$ & $(P_{\pm 4} P_{\mp 1})=1$ \\  
  $27$ &  $I_7, I_2, 3I_1$  &  $P_{-1}=(5,1)$, $P_{-2}=(3,0)$&  & $(P_{\pm 4} P_{\mp 2})=1$ \\  
&  &  $P_{-3}=(1,1)$, $P_{-4}=(6,0)$& & $(P_{\pm 4} P_{\mp 3})=1$ \\
&  &   &   &  $(P_4  P_{-4})=2$ \\ \hline

$28$ & $I_1^*, IV, I_1$ &    $P_1=(2,1)$, $P_2=(1,2)$ &  $P_{\pm 3}$  &$(P_{\pm 3} P_{\mp 2})=1$ \\  
  $29$ & $I_1^*, I_3, II$    &  $P_3=(3,0)$, $P_4=(0,1)$  & $P_{\pm 4}$    &  $(P_{\pm 4} P_{\mp 2})=1$ \\ 
  $30$ & $I_1^*, I_3, 2I_1$    &  $P_{-1}=(3,2)$, $P_{-2}=(1,1)$  & $G_0$ & $(P_{\pm 4} P_{\mp 3})=2$ \\ 
   &  &  $P_{-3}=(2,0)$, $P_{-4}=(0,2)$  &   &$(P_{-4} P_{4})=2$ \\  
   &    &         &   & $(P_{-3} P_{3})=1$    \\ 
& & && $(G_0  P_{\pm 4})=1$ \\
\hline

$31$ &  $I_6, III, I_2, I_1$ & $Q=(3,1,0)$, $P_1=(2,0,1)$ & $P_{\pm 2}$ & $(P_{\pm 3} P_{\mp 1})=1$ \\  
  $32$ &  $I_6, 2I_2, 2I_1$        & $P_{1,Q}=(5,1,1)$, $P_2=(4,0,0)$& $P_{\pm 3}$      &  $(P_{2} P_{-2})=1$ \\  
     &   & $P_{2,Q}=(1,1,0)$, $P_{\pm 3}=(0,0,1)$ &  & $(P_{\pm 2} P_{\mp 2,Q})=1$ \\  
&         & $P_{-1}=(4,0,1)$, $P_{-1}+Q=(1,1,1)$ &  & $(P_{\pm 3} P_{\mp 2})=2$ \\  
&        & $P_{-2}=(2,0,0)$, $P_{-2}+Q=(5,1,0)$ &  & $(P_{\pm 3} P_{\mp 2,Q})=2$ \\  
&        &     &  & $(P_{3} P_{-3})=2$    \\ \hline

$33$ & $I_5, I_4, II, I_1$      & $P_1=(2,1)$, $P_2=(4,2)$ & $P_{\pm 4}$ &$(P_{\pm 4} P_{\mp 3})=1$ \\  
$34$ &  $I_5, I_4, 3I_1$        &  $P_3=(1,3)$, $P_4=(3,0)$ & $P_{\pm 5}$ & $(P_{\pm 5} P_{\mp 3})=1$ \\  
&   & $P_5=(0,1)$, $P_{-1}=(3,3)$ & &  $(P_{4} P_{-4})=1$ \\  
  &   & $P_{-2}=(1,2)$, $P_{-3}=(4,1)$ &  & $(P_{\pm 5} P_{\mp 4})=2$ \\  
 &    & $P_{-4}=(2,0)$, $P_{-5}=(0,3)$ &  &  $(P_{5} P_{-5})=2$ \\ \hline

$35$ &   $I_5, IV, I_2, I_1$  &  $P_1=(1,1,1)$, $P_2=(2,2,0)$ &$P_{\pm 6}$   & $(P_{\pm 6} P_{\mp 1})=1$ \\     
$36$ &  $I_5, I_3, III, I_1$   &  $P_3=(3,0,1)$, $P_4=(4,1,0)$ &   &$(P_{\pm 6} P_{\mp 2})=1$ \\    
  $37$ &  $I_5, I_3, I_2, II$ &  $P_5=(0,2,1)$, $P_6=(1,0,0)$ &  &$(P_{\pm 6} P_{\mp 3})=1$ \\  
  $38$ &   $I_5, I_3, I_2, 2I_1$ &  $P_{-1}=(4,2,1)$, $P_{-2}=(3,1,0)$ &    &  $(P_{\pm 6} P_{\mp 4})=1$ \\  
&  &  $P_{-3}=(2,0,1)$, $P_{-4}=(1,2,0)$ &  &$(P_{\pm 6} P_{\mp 5})=2$ \\  
&  &  $P_{-5}=(0,1,1)$, $P_{-6}=(4,0,0)$ &   &$(P_{6} P_{-6})=2$    \\ \hline

\end{tabular} 
\end{center} 
 
\vskip 1cm
\hskip -0.6cm {\bf Table $3$. \  The Singularity types of the Picard number two Gorenstein log del Pezzo surfaces}
$$
$$
To find out the singularity types of the Picard number two Gorenstein log del Pezzo surfaces,
we start blowing down the zero-section ${\cal O}$ of  the Mordell-Weil rank one rational elliptic surface (cf. Theorem $3.5$) and then blow down an NEC (if exists) on the resulting surface,
and continue this process. The following table indicates the all singularity types found by this method. 
\vskip 0.5cm
\hskip -0.6cm
In this table, the numbers $i$ in the first column indicates that we start from the  No.$i$ Mordell-Weil rank one rational elliptic surface (cf. Table $1$ of  Appendix). Since we pay attention to the singularity types, if there are more than one surfaces with the same Mordell-Weil group which may produce the same singularity type, then we only pick up one of those surfaces. For example, blowing down the zero section of No.1 and No.2 Mordell-Weil rank one rational elliptic surfaces will give rise to the  same singularity type $D_7$ and the same configuration of the negative curves on the surface; thus we may only choose one surface, i.e., the No.$1$ Mordell-Weil rank one rational elliptic surface. In the last column of the table, the number $m$ in `` Exceptional divisor of $U_1 \longrightarrow U_m$'' is the largest number $m$ such that $U_m$ is the minimal resolution of the singularities of the Picard number two Gorenstein log del Pezzo surface. For the notation of $P_i$, $G_j$, $H_k$ and $Q$, we refer to the Table $2$ of  Appendix. The suffix $n$ of $U_n$ singifies that $K_{U_n}^2=n$. 
\vskip 0.5cm
\hskip -0.6cm Let's give an example, say, the type $D_7$, to explain how to read this table. There are five NECs on $U_1(D_7)$: $P_{\pm 1}$, $P_{\pm 2}$ and $G_0$ (cf. Table $2$ of Appendix). Thus there are five possibilities to do operation (blowing down) on $U_1(D_7)$. But by the symmetry, we only need to blow down one of the three NECs: $P_1$, $P_2$ or $G_0$. Thus blowing down $P_1$ (resp. $P_2$, or $G_0$) of $U_1(D_7)$ will give rise to $U_2(A_6)$ (resp. $U_2(D_6)$ and $U_2(A_1+D_5)$) and the first component of ``Exceptional divisor of $U_1$ to $U_m$'' in the last column is $P_1$ (resp. $P_2$, or $G_0$). When we reach $U_2(A_6)$, we may do the similiar operation as on $U_1(D_7)$ and get $U_3(A_5)$ and $U_3(A_4+A_1)$. 
$$
$$
\begin{center}
\begin{tabular}{| c|c|c|c|c|c|c|c|c|}\hline
$No.$ & $U_1$  &  $U_2$ & $U_3$ & $U_4$ & $U_5$ & $U_6$ & $U_7$ & Exceptional divisor   \\
 & & & & & & & & of $U_1 \longrightarrow U_m$
 \\ \hline 
$1$   &  $D_7$ & $A_6$ & $A_5$ & $A_3+A_1$ & $A_2+A_1$ & $2A_1$ & $A_1$ & $P_1+P_2+G_2+$  \\
 $2$ &  &    &     &         &       &   & & $G_1+G_4+G_5$ \\
 \cline{5-9}
 &   &  &  &  $A_4$& $A_2+A_1$ &  $2A_1$  & $A_1$ & $P_1+P_2+G_1$ \\
 &    &   &    &       &   & &  & $G_2+G_4+G_5$  \\ \cline{6-9}
&  &  &  &     & $A_3$ & $A_2$   & $A_1$ &$P_1+P_2+G_1$ \\
&        &       &     &         &       &   & & $G_4+G_5+G_2$ \\ \cline{7-9}
&  & &     &    &   & $2A_1$ & $A_1$ & $P_1+P_2+G_1$ \\

&  &  &  &    &       &   & & $G_4+G_2+G_5$ \\ \cline{4-9} 

   & &   & $A_4+A_1$ & $A_3+A_1$ & $A_2+A_1$ & $2A_1$ & $A_1$ &  
  $P_1+G_0+P_2+$ \\
& & & & & & & & $G_1+G_4+G_5$ 
\\ \cline{5-9}

 & &   &   &$A_4$ & $A_2+A_1$& $2A_1$ & $A_1$ & $P_1+P_0+P_2+$ \\
& & & & & & & & $G_2+G_4+G_5$ \\ \cline{6-9}

& & & &  & $A_3$  & $A_2$ & $A_1$ 
& $P_1+G_2+P_{-1}+$ \\ 
&  & & & & & & & $G_7+G_6+G_5$
\\  \cline{7-9}

& & & & & &$2A_1$ &$A_1$ & $P_1+G_0+P_2+$ \\
& & & & & & & & $G_4+G_2+G_5$  \\ \cline{3-9}

&     & $D_6$ & $D_5$ & $D_4$ & $A_3$ & $2A_1$ & $A_1$ & $P_2+G_1+G_4+$\\ 
& & & & & & & & $G_5+G_6+P_1$ \\ \cline{7-9}

&    & &              & &             & $A_2$ & $A_1$ & $P_2+G_1+G_4+$  \\ 
& & & & & & & & $G_5+P_1+G_6$ \\ \cline{5-9}

&    & &              & $A_4$ & $A_2+A_1$ & $2A_1$ & $A_1$ & $P_2+G_1+P_1+$ \\ 
& & & & & & & & $G_2+G_4+G_5$ \\ \cline{6-9}

&    & &              &       & $A_3$ & $A_2$ & $A_1$ & $P_2+G_1+P_1+$ \\ 
& & & & & & & & $G_4+G_5+G_6$  \\ \cline{7-9}

&     & &             &       &       & $2A_1$ & $A_1$ & $P_2+G_1+P_1+$  \\ 

& & & & & & & & $G_4+G_2+G_6$   \\ \cline{3-9}

&  & $A_1+D_5$ & $D_5$ & $D_4$ & $A_3$ & $A_2$ & $A_1$ & $G_0+P_2+G_4+$ \\

& & & & & & & & $G_5+P_1+G_6$\\ \cline{7-9}       

&  & &                 & &             & $2A_1$ & $A_1$ & $G_0+P_2+G_4+$ \\

& & & & & & & & $G_5+G_6+P_1$\\ \cline{5-9}

&  & &                 & $A_4$ & $A_2+A_1$ & $2A_1$ & $A_1$&  $G_0+P_2+P_1+$ \\

& & & & & & & & $G_2+G_4+G_5$\\ \cline{6-9}

&  & &                 &       & $A_3$ & $A_2$ & $A_1$ & $G_0+P_2+P_1+$ \\ 

& & & & & & & & $G_4+G_5+G_6$\\ \cline{7-9}

&  & &                 & &             & $2A_1$ & $A_1$ & $G_0+P_2+P_1+$ \\

& & & & & & & & $G_4+G_2+G_5$  \\ \cline{1-9}
\end{tabular} 
\end{center} 

\begin{center}
\begin{tabular}{| c|c|c|c|c|c|c|c|c|}\hline

&     &  & $A_5$ & $A_3+A_1$ & $A_2+A_1$ & $2A_1$ & $A_1$ & $P_2+P_1+G_2+$\\ 
& & & & & & & & $G_1+G_4+G_5$ \\ \cline{5-9}

&     &  &   & $A_4$ & $A_2+A_1$ & $2A_1$ & $A_1$ & $P_2+P_1+G_1+$\\ 
& & & & & & & & $G_2+G_4+G_5$ \\ \cline{6-9}

&     &  &  &  & $A_3$ & $A_2$ & $A_1$ & $P_2+P_1+G_1+$\\ 
& & & & & & & & $G_4+G_5+G_2$ \\ \cline{6-9}

&     &  &  &  &  & $2A_1$ & $A_1$ & $P_2+P_1+G_1+$\\ 
& & & & & & & & $G_4+G_2+G_5$ \\ \cline{4-9}

&  &           & $A_1+A_4$ & $A_2+2A_1$ & $A_2+A_1$ & $2A_1$ & $A_1  $ & $G_0+P_1+G_2+$ \\

& & & & & & & & $P_2+G_4+G_5$ \\ \cline{5-9}

&  & &  & $A_3+A_1$ & $A_2+A_1$ & $2A_1$ & $A_1$ &  $G_0+P_1+P_{-1}+$ \\
& & & & & & & & $G_2+G_7+G_6$ 
\\ \cline{5-9}
 &  &           & & $A_4$ & $A_2+A_1$ & $2A_1$ & $A_1  $ & $G_0+P_1+P_2+$ \\
& & & & & & & & $G_2+G_4+G_5$ \\ \cline{6-9}
&  &      &  & & $A_3$ & $A_2$ & $A_1  $ & $G_0+P_1+P_2+$ \\
& & & & & & & & $G_4+G_5+G_2$ \\ \cline{6-9}

&  &    & & &  & $2A_1$ & $A_1  $ & $G_0+P_1+P_2+$ \\

& & & & & & & & $G_4+G_2+G_5$ \\ \cline{1-9}

$3$  & $A_3+4A_1$ & $6A_1$ &  & & & &   & $P_2$ 
\\ \hline

$4$ & $D_5+2A_1$ & $D_4+2A_1$ & $A_3+2A_1$ & $4A_1$ &  & &  & $P_2+G_1+G_4$
\\ \cline{3-9}

$5$   & & $A_3+3A_1$ & $A_3+2A_1$ & $4A_1$ &  & &  &  $G_0+P_2+G_4$
\\ \hline
$6$ & $D_4+A_3$ & $2A_3$ & $A_3+2A_1$ & $4A_1$ & & &  & $P_2+G_1+P_{2, Q}$ \\  \cline{3-9}
    &           & $D_4+2A_1$ & $A_3+2A_1$ & $4A_1$ &  & &  & $P_{2, Q}+P_2+G_1$ \\ \cline{3-9}
     &           & $A_3+3A_1$ & $A_3+2A_1$ & $4A_1$ &  & &  & $G_0+P_2+P_{2, Q}$\\ \hline
$7$  & $E_7$ & $E_6$ & $D_5$  & $A_4$ & $A_2+A_1$& $2A_1$ & $A_1$ & $P_1+G_1+G_7+$ \\
$8$& & & & & & & & $G_6+G_0+G_5$ \\ \cline{6-9}

  & & &                 &     &$A_3$& $2A_1$&$A_1$ & $P_1+G_1+G_7+$  \\
& & & & & & & & $G_0+G_6+G_5$ \\ \cline{7-9}

     & &&                  &&            & $A_2$&$A_1$ & $P_1+G_1+G_7+$  \\ 
& & & & & & & & $G_0+G_5+G_6$ \\ \cline{5-9}

     & &&                  &$D_4$&$A_3$& $2A_1$&$A_1$ & $P_1+G_1+G_0+$  \\

& & & & & & & & $G_5+G_4+G_7$ \\ \cline{7-9}
     & &&                  &&          & $A_2$&$A_1$ & $P_1+G_1+G_0+$  \\

& & & & & & & & $G_5+G_7+G_4$ \\ \cline{3-9}

&   &$D_6$&$D_5$ &$A_4$&$A_2+A_1$& $2A_1$&$A_1$ & $G_0+P_1+G_5+$  \\
 & & & & & & & & $G_4+G_1+G_7$ \\ \cline{6-9}

&   &&              &     &$A_3$& $2A_1$&$A_1$ & $G_0+P_1+G_5+$  \\
& & & & & & & & $G_1+G_4+G_7$ \\ \cline{7-9}

& &&                &&          & $A_2$&$A_1$ & $G_0+P_1+G_5+$  \\
& & & & & & & & $G_1+G_7+G_4$ \\ \cline{5-9}

& &&                &$D_4$    &$A_3$& $2A_1$&$A_1$ & $G_0+P_1+G_1+$  \\

& & & & & & & & $G_5+G_4+G_7$ \\ \cline{7-9}

& &&                &&             & $A_2$&$A_1$ & $G_0+P_1+G_1+$  \\
& & & & & & & & $G_5+G_7+G_4$ \\ \cline{1-9}

\end{tabular} 
\end{center} 

\begin{center}
\begin{tabular}{| c|c|c|c|c|c|c|c|c|}\hline

&    &  &$A_5$  &$A_3+A_1$  &$A_2+A_1$& $2A_1$&$A_1$ & $G_0+G_5+G_4+$  \\

& & & & & & & & $P_1+G_1+G_7$ \\ \cline{5-9}

&    &  &   &$A_4$  &$A_2+A_1$& $2A_1$&$A_1$ & $G_0+G_5+G_1+$  \\

& & & & & & & & $G_4+G_1+G_7$ \\ \cline{6-9}

& &&             &      &$A_3$& $2A_1$&$A_1$ & $G_0+G_5+G_1+$  \\

 & & & & & & & & $G_1+G_4+G_7$ \\ \cline{7-9}

& &&             &&            & $A_2$&$A_1$ & $G_0+G_5+G_1+$  \\

& & & & & & & & $G_1+G_7+G_4$ \\ \cline{1-9}

$9$ & $A_7$& $A_5+A_1$ &$A_3+2A_1$&$4A_1$ & & & & $P_1+Q+P_{1,Q}$  \\ \cline{3-9}

&          & $2A_3   $ &$A_3+2A_1$&$4A_1$ & & & & $Q+P_1+P_{1,Q}$  \\
\hline

$10$ & $A_7$ & $A_5+A_1$ & $2A_2+A_1$ & $A_2+2A_1$ & $A_2+A_1$ & $2A_1$ & $A_1$ & $P_2+P_1+P_3+$  \\
$11$ & & & & & & & & $G_1+G_3+G_4$ \\ \cline{4-9}

 &&   & $A_1+A_4$ & $A_2+2A_1$ & $A_2+A_1$ & $2A_1$ & $A_1$ & $P_2+P_3+P_1+$  \\
& & & & & & & & $G_1+G_3+G_4$ \\ \cline{5-9}

&& &  & $A_3+A_1$ & $A_2+A_1$ & $2A_1$ & $A_1$  & $P_2+P_3+G_1+$  \\
& & & & & & & & $G_2+G_3+G_4$ \\ \cline{5-9}

&& &  & $A_4$ & $A_2+A_1$ & $2A_1$ & $A_1$  & $P_2+P_3+G_0+$  \\
& & & & & & & & $P_1+G_5+G_6$ \\ \cline{6-9}

&& &  &  & $A_3$ & $A_2$ & $A_1$  & $P_2+P_3+G_0+$  \\
& & & & & & & & $G_6+G_5+P_1$ \\ \cline{7-9}

&& &  &   &   & $2A_1$ & $A_1$  & $P_2+P_3+G_0+$  \\
& & & & & & & & $G_6+P_1+G_5$ \\ \cline{3-9}

 &  & $A_2+A_4$ & $2A_2+A_1$ & $A_2+2A_1$ & $A_2+A_1$ & $2A_1$ & $A_1$  & $P_1+P_2+P_3+$  \\

& & & & & & & & $G_1+G_3+G_4$ \\ \cline{4-9}

&&             & $A_4+A_1$ & $A_2+2A_1$ & $A_2+A_1$ & $2A_1$ & $A_1$ & $P_1+P_3+P_2+$  \\ 

& & & & & & & & $G_1+G_3+G_4$ \\ \cline{5-9}

&&             &   & $A_3+A_1$ & $A_2+A_1$ & $2A_1$ & $A_1$ & $P_1+P_3+G_0+$  \\ 

& & & & & & & & $P_2+G_6+G_5$ \\ \cline{5-9}

&&             &           & $A_4$ & $A_2+A_1$ & $2A_1$ & $A_1$ & $P_1+P_3+G_1+$  \\ 

& & & & & & & & $P_2+G_3+G_4$ \\ \cline{6-9}
&&              &          &       & $A_3$ & $2A_1$ & $A_1$ & $P_1+P_3+G_1+$  \\

 & & & & & & & & $G_3+G_2+G_4$ \\ \cline{7-9}
&&            &            &       &       & $ A_2$ & $A_1$  & $P_1+P_3+G_1+$  \\

 & & & & & & & & $G_3+G_4+G_5$ \\ \cline{3-9}

& & $A_6$ & $A_4+A_1$ & $A_2+2A_1$ & $A_2+A_1$ & $2A_1$ & $A_1$  & $P_3+P_1+P_2+$  \\
& & & & & & & & $G_2+G_3+G_4$ \\ \cline{5-9}
&&  & & $A_3+A_1$ & $A_2+A_1$ & $2A_1$ & $A_1$  & $P_3+P_1+G_1+$  \\
& & & & & & & & $G_1+G_3+G_4$ \\ \cline{5-9}

&&             &           & $A_4$ & $A_2+A_1$ & $2A_1$ & $A_1$  & $P_3+P_1+G_1+$  \\
& & & & & & & & $P_2+G_3+G_4$ \\ \cline{6-9}

&&              &          &       & $A_3$ & $2A_1$ & $A_1$  & $P_3+P_1+G_1+$  \\
& & & & & & & & $G_3+P_2+G_4$ \\ \cline{7-9}

&&            &            &       &       & $ A_2$ & $A_1$  & $P_3+P_1+G_1+$  \\

& & & & & & & & $G_3+G_4+G_5$ \\ \cline{1-9}

\end{tabular} 
\end{center} 

\begin{center}
\begin{tabular}{| c|c|c|c|c|c|c|c|c|}\hline

&&        & $A_5$ &   $A_3+A_1$ & $A_2+A_1$ & $2A_1$ & $A_1$  & $P_3+G_1+P_2+$  \\
& & & & & & & & $G_2+G_3+G_4$ \\ \cline{5-9}

&&        &  & $A_4$ & $A_2+A_1$ & $2A_1$ & $A_1$  & $P_3+G_1+G_2+$  \\

& & & & & & & & $P_2+G_3+G_4$ \\ \cline{6-9}

&&        & &             & $A_3$     & $2A_1$ & $A_1$  & $P_3+G_1+G_2+$  \\

 & & & & & & & & $G_3+G_2+G_4$ \\ \cline{7-9}

&&        & &             &           & $A_2$ & $A_1$ & $P_3+G_1+G_2+$  \\
 
& & & & & & & & $G_3+G_4+G_5$ \\ \cline{1-9}

$12$   &  $2A_3+A_1$ &  &   &  &  &   &  & \\ \hline
$13$   &  $D_6+A_1$ & $A_5+A_1$ & $A_3+2A_1$ & $4A_1$ &  &  &  & $P_1+G_3+G_0$ \\ \cline{3-9}
 $14$  &    & $D_4+2A_1$ & $A_3+2A_1$ & $4A_1$ &   &   &  & $G_0+P_1+G_3$  \\ \hline

$15$  &  $A_5+A_2$ & $3A_2$ &   &   &  &  &  & $P_1$ \\

$16$ & & & & & & & &  \\ \hline

$17$   &  $D_4+3A_1$ & $6A_1$  &   &  &  &  &   & $G_0$ \\ \hline

$18$   &  $E_6+A_1$ & $E_6$ & $D_5$ & $A_4$ & $A_2+A_1$ & $2A_1$ & $A_1$ & $P_3+P_1+G_1+$  \\
$19$ & & & & & & & & $G_4+P_2+G_2$ \\ \cline{6-9}
$20$&&&     &   & $A_3$ & $2A_1$  & $A_1$ & $P_3+P_1+G_1+$  \\
 & & & & & & & & $P_2+G_4+G_2$ \\ \cline{7-9}
&&&    & & & $A_2$ & $A_1$ &$P_3+P_1+G_1+$  \\
 & & & & & & & & $P_2+G_2+G_4$ \\ \cline{5-9}
&&&    & $D_4$ & $A_3$ & $2A_1$ & $A_1$& $P_3+P_1+P_2+$  \\
 & & & & & & & & $G_2+G_5+G_1$ \\ \cline{7-9}
&&&     & &  & $A_2$ & $A_1$& $P_3+P_1+P_2+$  \\
 & & & & & & & & $G_2+G_1+G_5$ \\ \cline{3-9}
&    & $D_5+A_1$ & $D_5$ & $A_4$ & $A_2+A_1$ & $2A_1$ & $A_1$ & $P_2+P_3+G_2+$  \\
& & & & & & & & $G_5+p_1+G_1$ \\ \cline{6-9}
&&&                      &                     & $A_3$ & $2A_1$ & $A_1$ & $P_2+P_3+G_2+$  \\
& & & & & & & & $P_1+G_5+G_1$
\\ \cline{7-9}
&&&    &  &       & $A_2$  &  $A_1$ & $P_2+P_3+G_2+$  \\
& & & & & & & & $P_1+G_1+G_5$ \\ \cline{5-9}

&&&  & $D_4$ & $A_3$ & $2A_1$ & $A_1$& $P_2+P_3+P_1+$  \\ 

& & & & & & & & $G_2+G_5+G_1$ \\ \cline{7-9}

&&&         & &           & $A_2$ & $A_1$  & $P_2+P_3+P_1+$  \\

& & & & & & & & $G_2+G_1+G_5$ \\ \cline{4-9}

&    &   & $A_4+A_1$ & $A_2+2A_1$ & $A_2+A_1$ & $2A_1$ & $A_1$ & $P_2+G_2+G_5+$  \\

& & & & & & & & $P_3+P_1+G_1$ \\ \cline{5-9}

&    &    && $A_3+A_1$ & $A_2+A_1$ & $2A_1$ & $A_1$ & $P_2+G_2+G_0+$  \\

& & & & & & & & $G_3+G_6+G_4$ \\ \cline{5-9}

&    &   &           & $A_4$ & $A_2+A_1$ & $2A_1$ & $A_1$ & $P_2+G_2+P_3+$  \\

& & & & & & & & $G_5+P_1+G_1$ \\ \cline{6-9}

&    &       &       &     & $A_3$ & $2A_1$ & $A_1$ & $P_2++G_2+P_3+$  \\
& & & & & & & & $P_1+G_5+G_1$ \\ \cline{7-9}
&    & &             &      &   & $A_2$ & $A_1$ & $P_2++G_2+P_3+$  \\
& & & & & & & & $P_1+G_1+G_5$ \\ \cline{1-9}

\end{tabular} 
\end{center} 

\begin{center}
\begin{tabular}{| c|c|c|c|c|c|c|c|c|}\hline

& & $A_5+A_1$ & $2A_2+A_1$ & $A_2+2A_1$  &$A_2+A_1$ & $2A_1$  &  $A_1$ & $G_0+G_3+P_2+$  \\
& & & & & & & & $G_2+G_6+G_4$ \\ \cline{4-9}

& & & $A_4+A_1$ & $A_2+2A_1$  &$A_2+A_1$ & $2A_1$  &  $A_1$ & $G_0+P_2+G_3+$  \\
& & & & & & & & $G_2+G_6+G_4$ \\ \cline{5-9}

& &  & & $A_3+A_1$  &$A_2+A_1$ & $2A_1$  &  $A_1$ & $G_0+P_2+G_2+$  \\
& & & & & & & & $G_3+G_6+G_4$ \\ \cline{5-9}

& &  & & $A_4$  &$A_2+A_1$ & $2A_1$  &  $A_1$ & $G_0+P_2+P_{-1}+$  \\
& & & & & & & & $G_3+P_1+G_1$ \\ \cline{6-9}

& &  & & &$A_3$ & $A_2$&$A_1$ & $G_0+P_2+P_{-1}+$  \\
& & & & & & & & $P_1+G_1+G_3$ \\ \cline{7-9}

& &  & &  && $2A_1$  &  $A_1$ & $G_0+P_2+P_{-1}+$  \\
& & & & & & & & $P_1+G_3+G_1$ \\ \cline{1-9}

$21$   &  $A_3+A_2+$ & $A_3+3A_1$ & $A_3+2A_1$ & $4A_1$ &  &  & & $P_4+H_1+P_2$ \\
$22$& $2A_1$ & & & & & & &  \\ \hline

$23$   &  $3A_2+A_1$ & $3A_2$ &  & &  &  &  & $P_3$ \\
$24$ & & & & & & & &  
\\ \hline

$25$&$A_6+A_1$ & $A_3+A_2+$ & $2A_2+A_1$ & $A_2+2A_1$ & $A_2+A_1$ & $2A_1$  & $A_1$ & $P_2+P_4+G_1+$  \\
$26$ & &  $A_1$   & & & & & & $G_2+G_4+G_5$\\ \cline{3-9}

$27$ & & $A_5+A_1$ & $2A_2+A_1$ & $A_2+2A_1$ & $A_2+A_1$ & $2A_1$  & $A_1$& $P_4+P_2+G_1+$  \\ 

& & & & & & & & $G_2+G_4+G_5$
\\ \cline{4-9}

&&   & $A_4+A_1$ & $A_2+2A_1$ & $A_2+A_1$ & $2A_1$  & $A_1$& $P_4+G_2+P_2+$  \\ 

& & & & & & & & $G_2+G_4+G_5$\\ \cline{5-9}

&& &   & $A_3+A_1$ & $A_2+A_1$ & $2A_1$  & $A_1$& $P_4+G_1+G_2+$  \\ 
& & & & & & & & $G_3+G_4+G_5$\\ \cline{5-9}

&& &   & $A_4$ & $A_2+A_1$ & $2A_1$  & $A_1$& $P_4+G_1+P_1+$  \\ 
& & & & & & & & $P_2+P_3+G_6$\\ \cline{6-9}

&& &   &  & $A_3$ & $A_2$  & $A_1$& $P_4+G_1+P_1+$  \\ 
& & & & & & & & $P_2+G_6+P_2$\\ \cline{7-9}

&& && &  & $2A_1$  & $A_1$& $P_4+G_1+P_1+$  \\ 
& & & & & & & & $P_3+G_6+P_2$\\ \cline{1-9}

$28 $ &$D_5+A_2$ & $A_3+A_2+$& $2A_2+A_1$ & $A_2+2A_1$ & $A_2+A_1$ & $2A_1$ & $A_1$ & $G_0+P_3+G_3+$  \\ 

$29$ & & $A_1$ & & & & & & $G_5+P_1+H_1$ \\ \cline{3-9}

$30 $ & & $D_5+A_1$& $D_5$ & $A_4$ & $A_2+A_1$ & $2A_1$ & $A_1$ &$P_4+H_1+P_3+$  \\ 

& & & & & & & & $G_3+P_2+G_1$
\\ \cline{6-9}
 & && & & $A_3$ & $2A_1$ & $A_1$ &$P_4+H_1+P_3+$  \\ 

& & & & & & & & $P_2+G_3+G_1$
\\ \cline{7-9}
&    &     &    &     &  & $A_2$ & $A_1$ &$P_4+H_1+P_3+$  \\ 

& & & & & & & & $P_2+G_1+G_3$ \\ \cline{4-9}

& & &  & $D_4$ & $A_3$ & $A_2$ & $A_1$ &$P_4+H_1+P_2+$  \\ 

& & & & & & & & $G_1+P_3+G_4$
\\ \cline{7-9}
& & & & &  & $2A_1$ & $A_1$ &$P_4+H_1+P_2+$  \\ 
& & & & & & & & $G_1+P_3+G_4$
\\ \cline{1-9}

\end{tabular} 
\end{center}

\begin{center}
\begin{tabular}{| c|c|c|c|c|c|c|c|c|}\hline

&&         & $A_4+A_1$ & $A_2+2A_1$ & $A_2+A_1$ & $2A_1$ & $A_1$&$P_4+P_3+G_3$  \\ 
& & & & & & & & $H_1+P_2+G_1$ \\ \cline{5-9}

&&  &  & $A_3+2A_1$ & $A_2+A_1$ & $2A_1$ & $A_1$&$P_4+P_3+P_1$  \\ 
& & & & & & & & $G_2+G_5+G_4$ \\ \cline{5-9}

&&  &  & $A_4$ & $A_2+A_1$ & $2A_1$ & $A_1$&$P_4+P_3+H_1+$  \\ 
& & & & & & & & $G_3+P_2+G_1$

\\ \cline{6-9}

 &&  & &  & $A_3$ & $2A_1$ & $A_1$&$P_4+P_3+H_1+$  \\ 

& & & & & & & & $P_2+G_3+G_1$
\\ \cline{7-9}

 &&         &   &              &   & $ A_2$ & $A_1$ &$P_4+P_3+H_1+$  \\ 
& & & & & & & & $P_2+G_1+G_3$
\\ \cline{3-9}

& & $A_4+A_2$& $2A_2+A_1$ & $A_2+2A_1$ & $A_2+A_1$ & $2A_1$ & $A_1$&$P_3+G_3+P_4+$  \\ 
& & & & & & & & $H_1+P_2+G_1$
\\ \cline{4-9}

 && & $A_4+A_1$ & $A_2+2A_1$ & $A_2+A_1$ & $2A_1$ & $A_1$&$P_3+P_4+G_3+$  \\ 
& & & & & & & & $H_1+P_2+G_1$
\\ \cline{5-9}

&& &  & $A_3+A_1$ & $A_2+A_1$ & $2A_1$ & $A_1$&$P_4+P_3+P_1+$  \\ 
& & & & & & & & $G_2+G_5+G_4$
\\ \cline{5-9}

&&         &  & $A_4$ & $A_2+A_1$ & $2A_1$ & $A_1$&$P_3+P_4+H_1+$  \\ 
& & & & & & & & $G_3+P_2+G_1$

\\ \cline{6-9}
 &&         &    &  & $A_3$ & $2A_1$ & $A_1$&$P_3+P_4+H_1+$  \\ 
& & & & & & & & $P_2+G_3+G_1$
\\ \cline{7-9}
 &&         &   & &   & $ A_2$ & $A_1$&$P_3+P_4+H_1+$  \\ 
& & & & & & & & $P_2+G_1+G_3$
\\ \cline{1-9}

$31$   &  $A_5+2A_1$ & $A_5+A_1$ & $A_3+2A_1$ & $4A_1$ &  &  &  & $P_3+P_2+P_1$\\ \cline{3-9}
  $32$  &    & $A_3+3A_1$ & $A_3+2A_1$ & $4A_1$ &   &  & &$P_2+P_3+P_1$ \\ \hline

$33$   &  $A_4+A_3$ & $A_3+A_2+$ & $2A_2+A_1$ & $A_2+2A_1$ & $A_2+A_1$ & $2A_1$ & $A_1$&$P_4+P_5+H_1+$  \\ 
$34$& & $A_1$ & & & & & & $H_2+P_3+G_1$ \\ \cline{3-9}

   &    & $A_4+A_2$ & $2A_2+A_1$ & $A_2+2A_1$ & $A_2+A_1$ & $2A_1$ & $A_1$&$P_5+P_1+P_4+$  \\

& & & & & & & & $G_2+P_3+H_3$ \\ \cline{4-9}
    &    &          & $A_4+A_1$ & $A_3+A_1$ & $A_2+A_1$ & $2A_1$ & $A_1$&$P_5+H_1+P_2+$  \\
& & & & & & & & $G_4+g_3+G_2$ \\ \cline{5-9}
  &    &            &      & $A_2+2A_1$ & $A_2+A_1$ & $2A_1$& $A_1$&$P_5+H_1+P_4+$  \\
& & & & & & & & $P_2+G_3+G_2$ \\ \cline{5-9}
    &   &          &   & $A_4$ & $A_2+A_1$ & $2A_1$ & $A_1$&$P_5+H_1+H_2+$  \\
& & & & & & & & $P_4+P_3+G_1$\\ \cline{6-9}
    &    &         &   &       & $A_3$ & $2A_1$ & $A_1$&$P_5+H_1+H_2$  \\
 & & & & & & & & $P_3+P_4+G_1$ \\ \cline{7-9}
   &    &           &   &      &       & $ A_2$ & $A_1$ &$P_5+H_1+H_2$  \\

& & & & & & & & $P_3+G_1+G_2$ \\ \cline{1-9}

$35$   &  $A_4+A_2+$ & $A_3+A_2+A_1$ & $2A_2+A_1$ & $A_2+2A_1$ & $A_2+A_1$ & $2A_1$ 
 &$A_1$ & $P_6+G_1+G_2$ \\

$36$ & $A_1$& & & & & & & $G_3+P_4+H_1$  \\
$37$& & & & & & & &   \\
$38$ & & & & & & & &   \\ \hline

\end{tabular} 
\end{center} 
$$
$$
\section*{References}

\hskip 0.6cm [Artin1] \ M.Artin: {\em Some numerical criteria for contractibility of curves on algebraic surfaces}, Amer.J.Math {\bf 84} (1962), 485-496.

[Artin2]\ M.Artin: {\em On isolated rational singularities of surfaces}, Amer.J.Math {\bf 88} (1966), 129-136.

[ATZ]\ E.Artal Bartolo, H.Tokunaga and D.-Q.Zhang: 
{\em Miranda-Persson's Problem on Extremal Elliptic $K3$ Surfaces},
Pacific J.Math., to appear.

[Brenton] \ L.Brenton: {\em On Singular complex surfaces with negative canonical bundle, with applications to singular  compactifications of  ${C}^2$ and to $3$-dimensional rational singularities}, Math.Ann. {\bf  248} (1980), 117-124.

[D] \ M.Demazure: {\em Surfaces de del Pezzo}, in: Lecture Notes in Mathematics, Vol. {\bf  777}, Springer, 1980.

[DZ]\ I.Dolgachev and D.-Q. Zhang: {\em Coble rational surfaces}, Amer.J.Math.{\bf 123}(2001), 79-114.

[Fujiki] \ A.Fujiki, R.Kobayashi, and S.Lu: {\em  On the Fundamental group of Certain Open normal Surfaces}, Saitama Math.J. {\bf  11} (1993), 15-20.

[Fu] \ M.Furushima: {\em Singular del Pezzo surfaces and analytic compactifications of $3$-dimensional complex affine space ${\bf  C}^3$}, Nagoya Math.J.{\bf  104} (1986), 1-28.

[GPZ]\ R.V.Gurjar, C.R.Pradeep and D.-Q. Zhang: {\em  On Gorenstein Surfaces 
isomorphic to ${\bf  P}^2/G$ }, preprint, 2000.

[GZ]\ R.V.Gurjar and D.-Q. Zhang:  {\em $\pi_1$ of smooth points of a log del Pezzo surface is finite: I }, J. Math. Sci. Univ. Tokyo {\bf  1} (1994), 137-180.

[GZ1]\ R.V.Gurjar and D.-Q. Zhang:  {\em $\pi_1$ of smooth points of a log del Pezzo surface is finite: II }, J. Math. Sci. Univ. Tokyo  {\bf  2} (1995), 165-196.

[Her]\ S.Herfurtner : {\em  Elliptic surfaces with four singular fibres}, Math.Ann. {\bf  291}(1991), 319-342.

[Hir]\ U.Schmickler-Hirzebruch : {\em  Elliptische Fl${\rm \ddot{a}}$chen ${\rm \ddot{u}}$ber ${\bf  P}_1{\bf  C}$ mit drei Ausnahmefasern und die hypergeometrische Differentialgleichung}. Schriftreihe des Mathematischen Instituts der Universit${\rm \ddot{a}}$t  M${\rm \ddot{u}}$nster, 2. Ser. {\bf  33}(1985).

[Kas] \ A.Kas : {\em Weierstrass normal forms and invariants of elliptic surfaces}, Trans.A.M.S. {\bf  225}(1977), 259-266.

[K1]\ K.Kodaira: {\em On compact complex 
analytic surfaces II}, Ann.Math. {\bf  77}(1963), 563-626.

[K2] \ K.Kodaira: {\em Pluricanonical systems on algebraic surfaces of general type}, J.Math.Soc.Japan Vol.{\bf  20}(1968), 170-192.

[KM] \ S.Keel and J.McKernan: {\em  Rational curves on quasi-projective surfaces}, Memoirs of the A.M.S. {\bf  669}(1999).

[KMM] \ Y.Kawamata, K.Matsuda, K.Matsuki:  {\em Introduction to the Minimal Model Problem}, Advanced Studies in Pure Mathematics {\bf  10}(1985), Algebraic Geometry, Sendai.

[M]\ R.Miranda:  {\em  Persson's list of singular fibres for a rational elliptic surface}, Math.Z.
{\bf  190}(1990), 191-211.

[MP]\ R.Miranda and U.Persson:  {\em On extremal rational elliptic surfaces}, Math.Z.
{\bf  193}(1986), 537-558.

[MP1]\ R.Miranda and U.Persson: {\em  Configurations of $I_n$ fibers
on elliptic $K3$ surfaces}, Math.Z. {\bf  201}(1989), 339-361.

[MP2] \ R.Miranda and U.Persson: {\em Mordell-Weil Groups of extremal
elliptic $K3$ surfaces}, Problems in the theory of
surfaces and their classification (Cortona, 1988), Symposia 
Mathematica, XXXII, Academic Press, London,1991, 167-192.

[M1]\ M.Miyanishi: {\em Open Algebraic Surfaces}, CRM monograph series, Amer.Math.Soc., 2000.

[MZ1]\ M.Miyanishi and D.-Q. Zhang : {\em  Gorenstein log del Pezzo surfaces of 
rank one}, J. of Algebra, {\bf  118} (1988), 63-84.

[MZ2]\ M.Miyanishi and D.-Q. Zhang : {\em Gorenstein log del Pezzo surfaces,II}, J. of Algebra, {\bf  156} (1993), 183-193.

[OS]\ K.Oguiso and T.Shioda: {\em The Mordell-Weil Lattice of a Rational Elliptic Surface}, Commentarii Mathematici, Universitatis Sancti Pauli, Vol.{\bf  40}(1991), 83-99.

[P]\ U.Persson : {\em Configurations of Kordaira fibers on rational elliptic surfaces}, Math.Z.{\bf  205} (1990), 1-47.

[Reid]\ M.Reid : {\em Chapters on algebraic surfaces}, in "Complex Algebraic Geometry", ed. J. Koll\'{a}r, IAS/Park City Mathematics Series {\bf  3} (1997), 5-159. 

[S] \ T.Shioda: {\em Elliptic modular surfaces},J.Math.Soc.Japan
{\bf  24}(1972),20-59.

[Sh] \ T.Shioda: {\em  On the Mordell-Weil Lattices}, Commentarii
Mathematici, Universitatis Sancti Pauli, Vol.{\bf  39}(1990),
211-240.

[SZ]\ I.Shimada and D.-Q.Zhang : {\em Classification of extremal
elliptic $K3$ surfaces and fundamental groups of open $K3$ surfaces}, Nagoya Math.J.{\bf 161}(2001), 23-54.

[T1]\ S.Takayama : {\em Simple connectedness of weak Fano varities}, 
J.Alg.Geom., Vol.{\bf 9} (2000), 403-407.

[Ye]\ Q.Ye: {\em On extremal elliptic $K3$  
surfaces}, J.Korean Math.Soc. {\bf  36}(1999), 1091-1113.

[Z1]\ D.-Q. Zhang : {\em On Iitaka surfaces}, Osaka J. Math. 
{\bf  24}(1987), 417-460.

[Z2]\ D.-Q. Zhang : {\em Logarithmic del Pezzo surfaces of rank one with contractible boundaries}, Osaka J. Math.  {\bf  25}(1988), no.2, 461-497.

[Z3]\ D.-Q. Zhang : {\em  Logarithmic del Pezzo  surfaces with rational double and triple singular points}, T${\hat o}$hoku Math. J. $(2)$ {\bf  41}(1989), no.3, 399-452.

[Z4] \ D.-Q. Zhang : {\em Algebraic surfaces with nef and big anti-canonical divisor    }, Math.Proc.Cambridge Philos.Soc. 
{\bf  117}(1995),no.1,  161-163.

[Z5]\ D.-Q. Zhang : {\em The fundamental group of the smooth part of a log Fano varity}, Osaka J. Math. 
{\bf 32}(1995), no.3, 637-644.

[Z6]\ D.-Q. Zhang : {\em  Algebraic surfaces with log canonical singularities and the fundamental groups of their smooth parts}, Trans.Amer.Math.Soc. 
{\bf 348}(1996), no.10, 4175-4184.
$$
$$

{\rm Departmant of Mathematics

National University of Singapore

2 Science Drive 2

Singapore 117543}

{\it e-mail}: mathyeq@hotmail.com

\end{document}